\documentclass[twoside,11pt]{article}
\usepackage{amsthm}
\usepackage[preprint]{formatting}
\usepackage{caption}

\usepackage{mathtools,bm}
\usepackage{tikz}
\usetikzlibrary{
  arrows.meta,
  positioning
}

\usepackage[T1]{fontenc}
\usepackage{booktabs}
\usepackage{microtype}
\usepackage{tabularx}
\usepackage{multirow}
\usepackage{cleveref}
\usepackage{enumitem}
\usepackage{geometry}
\usepackage[dvipsnames,table]{xcolor}
\definecolor{grayshade}{RGB}{240,247,255}

\hypersetup{
  colorlinks=true,
  linkcolor=blue,
  citecolor=blue,
  urlcolor=blue,
}

\newcommand{\wgray}[1]{\textcolor{.}{#1}}
\newcommand{\bvar}[1]{\textcolor{blue}{#1}}
\newcommand{\Lean}{(\ensuremath{\checkmark}\ Formalized and verified in Lean)}
\newcommand{\cures}{c_{u,\mathrm{res}}}

\definecolor{certiforange}{RGB}{180,68,2}
\definecolor{certifpurple}{RGB}{15,25,100}
\definecolor{certpipelinetextmuted}{RGB}{82,82,88}
\definecolor{certpipelinerulegray}{RGB}{140,140,148}
\definecolor{certpipelinerulemid}{RGB}{176,176,183}

\tikzset{
  certpipeline flow/.style={
    -{Stealth[length=1.15mm,width=0.80mm]},
    draw=certpipelinerulegray,
    line width=0.34pt
  },
  certpipeline transition/.style={
    midway,
    right=2.9pt,
    fill=white,
    inner xsep=0.75pt,
    inner ysep=0.35pt,
    text=certpipelinetextmuted,
    font=\fontsize{5.85}{6.2}\selectfont\itshape
  },
  certpipeline stage/.style={
    draw=certpipelinerulemid,
    fill=white,
    line width=0.31pt,
    minimum width=8.15cm,
    text width=7.60cm,
    minimum height=11.45mm,
    align=center,
    inner xsep=6.2pt,
    inner ysep=3.2pt,
    outer sep=0pt
  }
}

\newcommand{\certpipelinetitle}[2]{{\fontsize{8.1}{8.55}\selectfont
   \bfseries\color{#1}\textsc{#2}\par}}

\newcommand{\certpipelinesubtitle}[1]{{\vspace{0.45pt}\fontsize{6.2}{6.55}\selectfont
   \color{certpipelinetextmuted}#1\par}}

\newcommand{\certpipelinephasetag}[2]{{\fontsize{4.8}{5.15}\selectfont
   \bfseries\scshape\color{#1}#2}}

\newcommand{\certpipelinetagpanel}[6]{\node[certpipeline stage,#6] (#1) {\certpipelinetitle{#2}{#4}\certpipelinesubtitle{#5}};
  \draw[#2!78,line width=0.54pt]
    ([xshift=0.2pt]#1.north west) --
    ([xshift=-0.2pt]#1.north east);
  \node[
    anchor=west,
    fill=white,
    inner xsep=2.15pt,
    inner ysep=0.55pt,
    text=#2
  ] at ([xshift=6.4mm]#1.north west) {\certpipelinephasetag{#2}{#3}};
}

\definecolor{lightproofgray}{gray}{0.55}
\definecolor{darksalmon}{RGB}{240,80,70}
\newcommand{\dsalmon}[1]{\textcolor{darksalmon}{#1}}

\newcommand{\DLamL}{\ensuremath{\dsalmon{\Lambda_{\mathcal L}}}}
\newcommand{\DLamLlow}{\ensuremath{\dsalmon{\Lambda_{\mathcal L}^{\rm low}}}}
\newcommand{\DLamLhigh}{\ensuremath{\dsalmon{\Lambda_{\mathcal L}^{\rm high}}}}
\newcommand{\DLamDLhigh}{\ensuremath{\dsalmon{\Lambda_{D^\alpha\mathcal L}^{\rm high}}}}
\newcommand{\DLamDLlow}{\ensuremath{\dsalmon{\Lambda_{D^\alpha\mathcal L}^{\rm low}}}}
\newcommand{\DLamNthree}{\ensuremath{\dsalmon{\Lambda_{\mathcal N,3}}}}
\newcommand{\DLamDNthree}{\ensuremath{\dsalmon{\Lambda_{D^\alpha\mathcal N,3}}}}
\newcommand{\DLamThree}{\ensuremath{\dsalmon{\Lambda_{3}}}}
\newcommand{\DLamRPhi}{\ensuremath{\dsalmon{\Lambda_{\mathcal R,\Phi}}}}
\newcommand{\DLamDR}{\ensuremath{\dsalmon{\Lambda_{D^\alpha\mathcal R}}}}
\newcommand{\DLamR}{\ensuremath{\dsalmon{\Lambda_{\mathcal R}}}}
\newcommand{\DLamStab}{\ensuremath{\dsalmon{\Lambda_{\rm stab}}}}

\newcounter{certalg}
\newenvironment{certalgorithm}[1]{\refstepcounter{certalg}\par\medskip
\noindent\rule{\linewidth}{0.4pt}\par\smallskip
\noindent\textbf{Algorithm \thecertalg. #1}\par\smallskip
\begin{enumerate}[label=\textbf{\arabic*.},leftmargin=2.2em,itemsep=0.35em,topsep=0.25em]
}{\end{enumerate}
\smallskip\noindent\rule{\linewidth}{0.4pt}
\par\medskip
}

\newenvironment{grayproof}
  {\hfill 
  \begin{proof} \hspace{0.5mm} \color{lightproofgray}}
  {\end{proof}}

\newtheorem*{remark*}{Remark}
\newtheorem{proposition}{Proposition}

\title{Stability Framework for the \\ Singularity of the Euler Equations on $\mathbb{R}^3$}

\author{\name Valentin Duruisseaux* \email vduruiss@caltech.edu   \AND
\name Adarsh Ganeshram* \email aganeshram@berkeley.edu 
   \AND
\name Robert J. George \email rgeorge@caltech.edu 
\AND 
   \name Anima Anandkumar \email anima@caltech.edu  \AND
   \\
    \addr Department of Mathematics,     University of California Berkeley, Berkeley, CA, USA \AND
    \addr Department of Computing and Mathematical Sciences, California Institute of Technology, Pasadena, CA, USA
}

\begin{document}

\maketitle

\hfill 

\begin{abstract}
In a recent numerical study, we found a high-precision singular profile for the Euler equations on the unbounded domain $\mathbb{R}^3$. The present manuscript complements that study by establishing in detail a preliminary framework for proving (nonlinear) stability of the approximate self-similar profile, reducing the analysis to a large but finite collection of explicit estimates and computable constants. Conditional on rigorous certification of the estimates and constants appearing in the argument, and on the candidate profile  satisfying the required nonlinear stability conditions, the framework closes the stability proof and, crucially, allows the resulting stable rescaled profile to be reconstructed as an admissible solution in the original variables that becomes singular in finite physical time. With the overall stability and reconstruction mechanisms  formulated, the remaining work within this approach is largely quantitative:  determining whether the explicit constants and margins can be rigorously certified with sufficient positive margin and, where necessary, sharpening selected analytic  estimates.
\end{abstract}

\hfill \\

\tableofcontents

\clearpage

\section{Introduction}

A central problem in the analysis of partial differential equations (PDE) is whether all solutions arising from smooth initial data remain regular for all time or whether they may lose regularity in finite time. Such a loss of regularity is known as \emph{singularity formation}, or \emph{blowup}. In particular, establishing whether the Navier--Stokes equations have a finite-time blowup is among the most prominent open problems in fluid dynamics, and is one of the six unsolved Millennium problems~\citep{FeffermanNavierStokesCMI} and the $15^{{\textit th}}$ Smale problem~\citep{smale1998mathematical}. Closely related is the question whether smooth initial conditions can lead to finite-time singularities in the Euler equations.\\

\noindent Blowups have been established in related simpler problems or under modified assumptions. For instance, finite-time blowups have been established for the Euler equations in the presence of a boundary given smooth initial conditions~\citep{HouJiaje,Chen2026Analysis}, or under non-smooth initial conditions on the unbounded domain $\mathbb{R}^3$~\citep{nonSmoothProof,shkoller2026incompressible,JiajeIncompI,JiajeIncompII}.  
However, it is not clear how these additional assumptions can be removed. For instance, boundary-driven blowup mechanisms exploit the confinement imposed by solid walls to generate rapid vorticity growth~\citep{HouLuoBoundary,HouLuo1,HouJiaje}, without which discovering an approximate singular profile for the Euler equations on $\mathbb{R}^3$  has so far been elusive. Adding viscosity to the Euler equations results in the Navier--Stokes problem, which makes finite-time blowup substantially more difficult to realize due to dissipation of small-scale structures associated with any potential singularity~\citep{tao2007NSBlogpot,Tao2015}. \\

\noindent A common strategy for establishing finite-time blowup is to construct an approximate solution numerically, and then to use stability theory to prove the existence of a singular solution to the governing PDEs that remains sufficiently close to this approximation. A self-similar ansatz offers a tractable path to both the numerical search for such an approximate profile and its stability analysis. It allows us to continuously zoom in or out at the specified rate,  as the solution evolves, so that its overall ``shape'' is unchanged over time. More precisely,  a blowup solution is \emph{nearly self-similar} if, under a suitable dynamic rescaling, it remains close to an approximate blowup profile all the way up to the blowup time. The dynamic rescaling equations thus serve as a starting point for stability analysis around a given approximate profile.\\

\noindent In a recent numerical study~\citep{EulerBlowupMain}, we used Physics-Informed Neural Networks (PINNs)~\citep{PINN_OG} to obtain a highly accurate approximate self-similar profile for the Euler equations on $\mathbb{R}^3$. The PINN profile was then converted into piecewise polynomial splines suitable for rigorous certification with tight bounds. The spline representation facilitates exact interval-wise differentiation and evaluation of profile-dependent quantities. We also used \texttt{Arb}~\citep{Arb}, an arbitrary-precision interval arithmetic library, to rigorously control numerical and rounding errors throughout the certification procedure. This allowed us in particular to bound the residual errors of the spline profiles by $\mathcal{O}(10^{-5})$ in $L^2$ and $\mathcal{O}(10^{-3})$ in $L^\infty$. \\

\noindent In addition to finding an approximate profile, establishing stability is crucial and complex. This becomes simpler when a global outgoing property is satisfied~\citep{DGBlowup}, which pushes the self-similar flow away from the singular core. This expansive effect of global outgoing flow is also closely connected to the favorable higher-order damping that underlies the stability analysis. However, for the Euler system with a self-similar ansatz, \citet{constantin2026putative} provides strong theoretical evidence that nontrivial fixed points of the meridional flow are likely to exist and ``well known to cause tremendous headaches''. This motivates replacing the stringent global outgoing condition with a weaker local one that pushes the self-similar flow away from these fixed points, preventing perturbations from concentrating there and providing a favorable structure for the stability analysis. Establishing stability with this weaker outgoing property however requires a substantially more delicate analysis than in prior works. \\

\paragraph{In this paper,} we establish a preliminary proof framework for the nonlinear stability of the approximate blowup profiles, reducing the argument to a (large) finite collection of explicit, rigorously verifiable estimates.\\

\noindent Starting from the spline representation of the approximate profiles, we establish a stability framework: following prior works~\citep{DGBlowup,HouJiaje,Chen2026Analysis} and related works on stability of blowups, we linearize the appropriate dynamically rescaled equations around the spline approximate profiles and formulate a stability theorem in suitable weighted $L^2$ and Sobolev norms. We then derive complementary low-order and high-order energy estimates for both the linear and nonlinear terms, with explicit computable bounds. These estimates are incorporated into a single energy functional that  provides the quantitative  structure needed to  determine whether the full stability  argument can be closed. The  stability question is thereby reduced to  determining whether a large finite collection of explicit  estimates can be rigorously certified with sufficient margin using interval arithmetic or \texttt{Arb} computations.\\

\hfill 

\noindent \emph{Remaining Work to Complete the Proof.} The remaining work within the proposed stability framework is primarily  computational: rigorously certifying the profile-dependent constants and margins  using spline representations, interval arithmetic, certified matrix bounds, and finite-dimensional  optimization, and determining whether they provide sufficient positive margin to close the stability argument. Some estimates may require further refinement to obtain sufficient margin, but the modular structure of the argument allows such adjustments without changing the underlying stability mechanism.  If these  certifications yield the  required positive margins and thereby establish nonlinear stability of the rescaled profile  under the required class of perturbations,  this stability can then be transferred through the dynamic rescaling to finite-time singularity formation in the original variables. At the critical scaling $\lambda=1/2$, the corresponding control of the modulation parameters near their self-similar values would  then allow the  stable rescaled profile  to be reconstructed as an admissible solution that becomes singular in finite time.\\

\noindent In parallel, the \texttt{LeanPDE} paper develops a formalization in Lean~\citep{Lean4} of key components of the argument, with particular emphasis on the symbolic derivations and proof steps underlying the analysis. Using \texttt{Mathlib}~\citep{mathlib} and \texttt{TorchLean}~\citep{TorchLean}, this formalization connects symbolic verification in Lean with certified numerical computations performed separately in \texttt{Julia} using \texttt{Arb} interval arithmetic~\citep{Julia,Arb}. The resulting pipeline assigns complementary roles to formal symbolic reasoning and rigorous numerical certification, and forms part of the broader \texttt{LeanPDE} effort to develop an integrated framework for machine-checked PDE analysis. Further details of the formalization and verification workflow are provided in the \texttt{LeanPDE} paper.

\clearpage

\section{Approximate Self-Similar Profiles of the Euler Equations} \label{sec: Euler}

\subsection{The 3D Axisymmetric Equations} \label{sec: The 3D Axisymmetric Equations}

Starting from the incompressible 3D Euler equations written by \citet{EulerOriginalPaper},
\begin{equation}
\label{eq:Original_Euler}
\partial_t \tilde u + (\tilde u\cdot\nabla) \tilde u = -\nabla \tilde p,
\qquad
\nabla\cdot \tilde u = 0,
\end{equation}
we solve for the velocity field $\tilde u(x,t)\in\mathbb{R}^3$ and the scalar pressure $\tilde p(x,t)\in\mathbb{R}$, for $x\in\mathbb{R}^3$ and $t\ge 0$. Taking the curl of the above equation yields the vorticity formulation
\begin{equation}
\label{eq:vorticity}
\partial_t \tilde \omega + \tilde u\cdot\nabla \tilde \omega = \tilde \omega\cdot\nabla \tilde u,
\qquad
\tilde \omega := \nabla\times \tilde u .
\end{equation}
In $\mathbb{R}^3$, one can recover $\tilde u$ from $\tilde \omega$ via the Biot--Savart law. Equivalently, one can introduce a vector streamfunction $\tilde \psi(x,t)\in\mathbb{R}^3$ satisfying
\begin{equation}
\label{eq:vector_potential}
-\Delta \tilde \psi = \tilde \omega,
\qquad
\tilde u = \nabla\times \tilde \psi,
\end{equation}
together with a gauge condition (e.g.\ $\nabla\cdot \tilde \psi=0$) to fix $\tilde \psi$. \\ 

\noindent For flows that are symmetric about a fixed spatial axis, such as the $\mathring z$-axis, it is convenient to rewrite equation~\eqref{eq:Original_Euler} in cylindrical coordinates
\begin{equation}
x=(\mathring{r}\cos\theta,\mathring{r}\sin\theta,\mathring{z}),\qquad \mathring{r} \ge 0,\ \  \theta\in[0,2\pi),\ \  \mathring{z}\in\mathbb{R}.
\end{equation}
In these coordinates, the velocity and vorticity decompose as
\begin{align}
\tilde u = u^r(\mathring{r},\mathring{z},t)\,e_r + u^\theta(\mathring{r},\mathring{z},t)\,e_\theta + u^z(\mathring{r},\mathring{z},t)\,e_z, \\
\tilde \omega = \omega^r(\mathring{r},\mathring{z},t)\,e_r + \omega^\theta(\mathring{r},\mathring{z},t)\,e_\theta + \omega^z(\mathring{r},\mathring{z},t)\,e_z,
\end{align}
where $e_r, e_\theta, e_z$ are the unit vectors in cylindrical coordinates. In addition, by definition of axisymmetry, we have that $\partial_\theta(\cdot)=0$. \\ 

\noindent Following \citet{HouLuoBoundary}, it is convenient to introduce the rescaled swirl and streamfunction variables
\begin{equation}
\label{eq:scaled_vars}
u_{\bullet} := \frac{u^\theta}{\mathring{r}},\qquad
\omega_{\bullet} := \frac{\omega^\theta}{\mathring{r}},\qquad
\psi_{\bullet} := \frac{\psi^\theta}{\mathring{r}},
\end{equation}
where $\psi^\theta$ denotes the $\theta$-component of the vector potential $\psi$ in~\eqref{eq:vector_potential}. This reformulation has the advantage of removing the $1/\mathring{r}$ singularity from the cylindrical coordinates. \\

\noindent We also define the radial component $u^r$ and axial component $u^z$ of the velocity in cylindrical coordinates. These meridian velocity components $(u^r,u^z)$ can be recovered from $\psi_{\bullet}$ via the (axisymmetric) streamfunction relations
\begin{equation}
\label{eq:stream}
u^{r} = -\,\mathring{r}\,\partial_{\mathring{z}} \psi_{\bullet},
\qquad
u^{z} = 2\psi_{\bullet} + \mathring{r}\,\partial_{\mathring{r}} \psi_{\bullet}.
\end{equation}

\hfill 

\paragraph{3D Axisymmetric Euler Equations.} The axisymmetric 3D Euler equations can be written as
\begin{align}
\label{eq:Euler1}
\partial_t u_{\bullet} + u^r\,\partial_{\mathring{r}} u_{\bullet} + u^z\,\partial_{\mathring{z}} u_{\bullet}
&= 2u_{\bullet}\,\partial_{\mathring{z}} \psi_{\bullet}, \\
\label{eq:Euler2}
\partial_t \omega_{\bullet} + u^r\,\partial_{\mathring{r}} \omega_{\bullet} + u^z\,\partial_{\mathring{z}} \omega_{\bullet}
&= 2u_{\bullet}\,\partial_{\mathring{z}}u_{\bullet}, \\
\label{eq:Euler3}
-\Bigl[\partial_{\mathring{r}}^2 + \frac{3}{\mathring{r}}\partial_{\mathring{r}} + \partial_{\mathring{z}}^2\Bigr]\psi_{\bullet}
&= \omega_{\bullet}.
\end{align}

\noindent Equations \eqref{eq:stream}, \eqref{eq:Euler1}--\eqref{eq:Euler3} form the reduced axisymmetric-with-swirl system in $(\mathring{r},\mathring{z})$ for $(u_{\bullet},\omega_{\bullet},\psi_{\bullet})$.

\clearpage

\noindent \textbf{Axis interpretation.}
Every occurrence of \(\partial_{\mathring r}^2+3\mathring r^{-1}\partial_{\mathring r} +\partial_{\mathring z}^2\) at \(\mathring r=0\) denotes its axis-regular continuation, not literal division by zero.  For a \(C^2\) profile even in
\(\mathring r\),
\begin{equation}
\left(\partial_{\mathring r}^2+\frac{3}{\mathring r}
\partial_{\mathring r}+\partial_{\mathring z}^2\right)f(0,\mathring z)
=4\,\partial_{\mathring r\mathring r}f(0,\mathring z)
+\partial_{\mathring z\mathring z}f(0,\mathring z).
\end{equation}

\hfill 

\subsection{Traveling-Wave Ansatz} \label{sec: Ansatz}

\noindent
We consider a \emph{traveling self-similar} ansatz. 
\begin{itemize}
    \item \emph{Self-similar} refers to the fact that the geometry of the solution, after a proper rescaling of space, keeps an approximately fixed shape. 
    \item \emph{Traveling} means that the center of this singularity is not fixed. Instead, it moves along the symmetry ($\mathring{r} = 0$) axis while the solution converges to the self-similar geometry.
\end{itemize} 

\noindent Since the ansatz is written in terms of powers of $\tau = T-t$, it is also a \emph{backward self-similar ansatz}: the possible singular profile is described by rescaling backward from the blowup time $T$.\\

\noindent
The translation is introduced only in the axial ($\mathring{z}$) variable. This is consistent with the axisymmetric setting, where $\mathring r=0$ is the fixed symmetry axis and $\mathring r\geq 0$ measures distance to that axis. A radial ($\mathring{r}$) translation would move the singularity away from the symmetry axis and would generally destroy the parity and regularity structure imposed at $\mathring r=0$. By contrast, an axial translation in $\mathring z$ preserves the cylindrical symmetry and provides a natural mechanism for a coherent structure to drift while it concentrates.\\

\noindent
Given the original axisymmetric variables as $(u_\bullet,\omega_\bullet,\psi_\bullet)$, we introduce the \emph{rescaled profile variables} $(u,\omega,\psi)$ via the \emph{traveling self-similar} ans
\begin{align}
u_\bullet(\mathring r,\mathring z,t)
&= (T-t)^{c_u}\,
 u\!\left(\frac{\mathring r}{(T-t)^\lambda},\frac{\mathring z-\mathring z_c(t)}{(T-t)^\lambda}\right),
\label{ansantz1}\\
\omega_\bullet(\mathring r,\mathring z,t)
&= (T-t)^{c_\omega}\,
 \omega\!\left(\frac{\mathring r}{(T-t)^\lambda},\frac{\mathring z-\mathring z_c(t)}{(T-t)^\lambda}\right),
\label{ansantz2}\\
\psi_\bullet(\mathring r,\mathring z,t)
&= (T-t)^{c_\psi}\,
 \psi\!\left(\frac{\mathring r}{(T-t)^\lambda},\frac{\mathring z-\mathring z_c(t)}{(T-t)^\lambda}\right).
\label{ansantz3}
\end{align}
\noindent
Here $T$ is the \emph{blowup time}, $\lambda$ is the \emph{spatial blowup rate}, $\mathring z_c(t)$ is the center of the traveling profile, and $c_u,c_\omega,c_\psi$ are \emph{amplitude exponents}. \\

\noindent The rescaled coordinates are
\begin{equation}
    r=\frac{\mathring r}{(T-t)^\lambda},
    \qquad
    z=\frac{\mathring z-\mathring z_c(t)}{(T-t)^\lambda},
\label{eq:rescaled_coordinates_ansatz}
\end{equation}
\noindent
and the rescaled meridian velocity components $u^r$ and $u^z$, can be recovered from $\psi$ via
\begin{equation}
    u^r=-r\partial_z\psi,
    \qquad
    u^z=2\psi+r\partial_r\psi.
    \label{eq:rescaled_stream_relations}
\end{equation}

\noindent The traveling speed is encoded by the constant $C$ through
\begin{equation}
    \frac{d}{dt}\mathring z_c(t)=-C(T-t)^{\lambda-1}.
    \label{eq:center_speed}
\end{equation}
\noindent
With this convention, the moving coordinate $z$ satisfies
\begin{equation}
    \left(\partial_t z\right)_{\mathring r,\mathring z}
    =\frac{\lambda z+C}{T-t}.
    \label{eq:z_coordinate_speed_ansatz}
\end{equation}
\noindent
Thus $C$ is the \emph{traveling speed} of the profile in the rescaled axial coordinate. If $C=0$, the singularity is centered at a fixed axial point in the rescaled frame. If $C\neq0$, the singularity drifts in the original variables while remaining stationary in the moving rescaled frame.\\

\paragraph{Scaling exponents.}
\noindent
The exponents are fixed by balancing the powers of $T-t$ in the axisymmetric system. This balance enforces that time differentiation, transport, stretching, and the elliptic streamfunction relation all contribute at the same leading order. It gives
\begin{equation}
    c_u=-1,
    \qquad
    c_\omega=-1-\lambda,
    \qquad
    c_\psi=-1+\lambda.
    \label{eq:scaling_exponents}
\end{equation}
\noindent
These exponents have a simple interpretation. The variable $u_\bullet$ grows like $(T-t)^{-1}$. The vorticity variable $\omega_\bullet$ contains one additional spatial derivative at scale $(T-t)^\lambda$, and therefore grows like $(T-t)^{-1-\lambda}$. The streamfunction is two spatial derivatives smoother than $\omega_\bullet$, which leads to the exponent $-1+\lambda$.\\

\hfill

\subsection{Steady-State Profile Equations}
\label{sec: Steady-State Profile Equations}

We next present the profile equations generated by the traveling self-similar ansatz. This reduction transforms the original singular evolution into a time-independent nonlinear system for the rescaled variables $(U,\Omega,\Psi)$, together with the scaling exponent $\lambda$ and the drift speed $C$. \\ 

\noindent For notation simplicity, we define the elliptic operator
\begin{equation}
    \mathcal{E} \coloneqq \partial_r^2 + \frac{3}{r}\partial_r + \partial_z^2.
\end{equation}

\hfill

\begin{proposition}[Traveling-wave profile equations]
\label{thm:traveling_wave_profile_eps1}
Assume that $(u_\bullet,\omega_\bullet,\psi_\bullet)$ is given by the ansatz \eqref{ansantz1}--\eqref{ansantz3}, that the center satisfies \eqref{eq:center_speed}, and that the exponents are chosen as in \eqref{eq:scaling_exponents}.  Set $\lambda=\frac12.$ Then $(U,\Omega,\Psi)$ satisfies
\begin{align}
U + \bigl(\lambda r+U^r\bigr)\partial_rU + \bigl(C+\lambda z+U^z\bigr)\partial_zU
&=2U\partial_z\Psi,
\label{eq:ProfileEq1_eps1}\\
(1+\lambda)\Omega + \bigl(\lambda r+U^r\bigr)\partial_r\Omega
+ \bigl(C+\lambda z+U^z\bigr)\partial_z\Omega
&=2U\partial_zU,
\label{eq:ProfileEq2_eps1}\\
-\, \mathcal E \,\Psi&=\Omega,
\label{eq:ProfileEq3_eps1}
\end{align}
where
\begin{equation}
U^r=-r\partial_z\Psi,
\qquad
U^z=2\Psi+r\partial_r\Psi.
\label{eq:profile_velocity_recovery_eps1}
\end{equation}
Equivalently,
\begin{align}
U+r\bigl(\lambda-\partial_z\Psi\bigr)\partial_rU
+\bigl(C+\lambda z+2\Psi+r\partial_r\Psi\bigr)\partial_zU
&=2U\partial_z\Psi,
\label{eq:ProfileEq1_expanded_eps1}\\
(1+\lambda)\Omega+r\bigl(\lambda-\partial_z\Psi\bigr)\partial_r\Omega
+\bigl(C+\lambda z+2\Psi+r\partial_r\Psi\bigr)\partial_z\Omega
&=2U\partial_zU,
\label{eq:ProfileEq2_expanded_eps1}\\
-\, \mathcal E \,\Psi&=\Omega.
\label{eq:ProfileEq3_expanded_eps1}
\end{align}
\end{proposition}

\begin{grayproof}  The proof is provided in the numerical paper~\citep{EulerBlowupMain}.
\end{grayproof}

\hfill  

\noindent Equations \eqref{eq:ProfileEq1_eps1}--\eqref{eq:ProfileEq3_eps1} are the steady-state equations for the traveling self-similar profile. Thus, the traveling-wave ansatz reduces the study of a possible finite-time singularity to an autonomous nonlinear profile problem for $(U,\Omega,\Psi)$ and the parameters $\lambda$ and $C$. The parameter $\lambda$ determines the rate at which the spatial scale collapses, while $C$ determines the speed of the axial drift in the rescaled frame. \\

\clearpage 

\subsection{Numerical Profiles} \label{sec: numerical profiles}

\noindent In a recent numerical study~\citep{EulerBlowupMain}, we parametrized and optimized a PINN solution for the Euler self-similar profile system. This profile was obtained to high accuracy, as can be seen from the profile equations residuals reported in \Cref{tab:euler_ns_errors}. The PINN computation serves as a numerical discovery tool, producing an accurate candidate self-similar profile, providing a practical and well-conditioned means of obtaining an accurate profile when direct optimization over analytic ansatzes or alternative basis expansions would be considerably more difficult. \\

\begin{table*}[htbp]
\centering

\caption{Residual errors for the $U$, $\partial_r U$, $\partial_z U$, $\Omega$, and $\Psi$ profile equations for the Euler system, comparing PINN and spline representations. Note that the residuals are estimated empirically on an independent set of $100{,}000{,}000$ randomly sampled collocation points for the PINN representations, whereas they are evaluated analytically and exactly for the spline representations. in particular, cRMSE refers to the continuous root-mean squared error $ \sqrt{
\frac{1}{\left|\mathbb{D}\right|}
\int_{\mathbb{D}}
\mathsf{residual}\!\left(r_{\!\mathfrak{comp}}, z_{\mathfrak{comp}}\right)^2
\,\mathrm{d}r_{\mathfrak{comp}}\,\mathrm{d}z_{\mathfrak{comp}}
} $. Here, $\mathbb{D}$ denotes the full computational domain (corresponding to the half plane), $B_{0.1}(0^{\star})$, $B_1(0^{\star})$, and $B_{10}(0^{\star})$ denote balls of radii $0.1$, $1$, and $10$, respectively, centered at the origin $0^{\star}$ of the stability analysis (i.e. the on-axis meridional fixed point), $\mathsf{axis}$ denotes the $r=0$ axis, and $\mathsf{FF}$ denotes the far-field region $\left(r_{\!\mathfrak{comp}}, z_{\mathfrak{comp}}\right) \in [10, 30] \times [-30, -10] \cup [10, 30] \times [10, 30] $.}
\label{tab:euler_ns_errors}

\vspace{0.7mm}
\renewcommand{\arraystretch}{1.04}

\begin{tabularx}{0.9\linewidth}{
@{}
c@{\hspace{2mm}}
l
*{4}{>{\centering\arraybackslash}X}
@{}
}

&
&
\multicolumn{2}{c}{\textbf{PINNs}}
&
\multicolumn{2}{c}{\textbf{Splines}}
\\[-1.1mm]

\cmidrule(lr){3-4}
\cmidrule(lr){5-6}

&
\textbf{Domain}
& RMSE
& $\ell^\infty$
& cRMSE
& $L^\infty$
\\

\midrule
\addlinespace[0.7mm]


\multirow{6}{*}{
\rotatebox[origin=c]{90}{\textbf{\small $U$-equation}}
}
&
$B_{0.1}(0^{\star})$
& $2.8\!\cdot\!10^{-7}$
& $3.4\!\cdot\!10^{-7}$
& $2.8\!\cdot\!10^{-7}$
& $3.5\!\cdot\!10^{-7}$
\\

&
$B_1(0^\star)$
& $7.0\!\cdot\!10^{-6}$
& $3.3\!\cdot\!10^{-5}$
& $6.7\!\cdot\!10^{-6}$
& $3.4\!\cdot\!10^{-5}$
\\

&
$B_{10}(0^{\star})$
& $1.2\!\cdot\!10^{-5}$
& $5.7\!\cdot\!10^{-5}$
& $1.7\!\cdot\!10^{-5}$
& $5.8\!\cdot\!10^{-5}$
\\

&
$\mathsf{axis}$
& $9.0\!\cdot\!10^{-5}$
& $1.4\!\cdot\!10^{-3}$
& $9.1\!\cdot\!10^{-5}$
& $1.4\!\cdot\!10^{-3}$
\\

&
$\mathsf{FF}$
& $6.2\!\cdot\!10^{-6}$
& $2.2\!\cdot\!10^{-5}$
& $6.3\!\cdot\!10^{-6}$
& $2.3\!\cdot\!10^{-5}$
\\

&
\cellcolor{grayshade}$\mathbb{D}$
& \cellcolor{grayshade}$4.2\!\cdot\!10^{-5}$
& \cellcolor{grayshade}$3.6\!\cdot\!10^{-3}$
& \cellcolor{grayshade}$5.7\!\cdot\!10^{-5}$
& \cellcolor{grayshade}$3.7\!\cdot\!10^{-3}$
\\

\addlinespace[1.0mm]
\midrule
\addlinespace[1.0mm]


\multirow{6}{*}{
\rotatebox[origin=c]{90}{\textbf{\small $\partial_r U$-equation}}
}
&
\rule[-0.55ex]{0pt}{3.15ex}$B_{0.1}(0^{\star})$
& $1.5\!\cdot\!10^{-6}$
& $3.3\!\cdot\!10^{-6}$
& $1.6\!\cdot\!10^{-6}$
& $3.4\!\cdot\!10^{-6}$
\\

&
\rule[-0.55ex]{0pt}{3.15ex}$B_1(0^\star)$
& $5.2\!\cdot\!10^{-5}$
& $1.9\!\cdot\!10^{-4}$
& $5.0\!\cdot\!10^{-5}$
& $1.9\!\cdot\!10^{-4}$
\\

&
\rule[-0.55ex]{0pt}{3.15ex}$B_{10}(0^{\star})$
& $4.2\!\cdot\!10^{-5}$
& $3.1\!\cdot\!10^{-4}$
& $6.9\!\cdot\!10^{-5}$
& $3.2\!\cdot\!10^{-4}$
\\

&
\rule[-0.55ex]{0pt}{3.15ex}$\mathsf{axis}$
& $0$
& $0$
& $0$
& $0$
\\

&
\rule[-0.55ex]{0pt}{3.15ex}$\mathsf{FF}$
& $1.3\!\cdot\!10^{-10}$
& $2.4\!\cdot\!10^{-9}$
& $1.3\!\cdot\!10^{-10}$
& $2.5\!\cdot\!10^{-9}$
\\

&
\cellcolor{grayshade}\rule[-0.55ex]{0pt}{3.15ex}$\mathbb{D}$
& \cellcolor{grayshade}$1.2\!\cdot\!10^{-4}$
& \cellcolor{grayshade}$8.3\!\cdot\!10^{-3}$
& \cellcolor{grayshade}$1.5\!\cdot\!10^{-4}$
& \cellcolor{grayshade}$8.5\!\cdot\!10^{-3}$
\\

\addlinespace[1.0mm]
\midrule
\addlinespace[1.0mm]


\multirow{6}{*}{
\rotatebox[origin=c]{90}{\textbf{\small $\partial_z U$-equation}}
}
&
\rule[-0.55ex]{0pt}{3.15ex}$B_{0.1}(0^{\star})$
& $1.1\!\cdot\!10^{-6}$
& $1.9\!\cdot\!10^{-6}$
& $1.1\!\cdot\!10^{-6}$
& $2.3\!\cdot\!10^{-6}$
\\

&
\rule[-0.55ex]{0pt}{3.15ex}$B_1(0^\star)$
& $2.3\!\cdot\!10^{-5}$
& $1.0\!\cdot\!10^{-4}$
& $2.3\!\cdot\!10^{-5}$
& $1.1\!\cdot\!10^{-4}$
\\

&
\rule[-0.55ex]{0pt}{3.15ex}$B_{10}(0^{\star})$
& $1.3\!\cdot\!10^{-5}$
& $1.3\!\cdot\!10^{-4}$
& $2.8\!\cdot\!10^{-5}$
& $1.3\!\cdot\!10^{-4}$
\\

&
\rule[-0.55ex]{0pt}{3.15ex}$\mathsf{axis}$
& $1.1\!\cdot\!10^{-6}$
& $1.7\!\cdot\!10^{-5}$
& $1.2\!\cdot\!10^{-6}$
& $1.8\!\cdot\!10^{-5}$
\\

&
\rule[-0.55ex]{0pt}{3.15ex}$\mathsf{FF}$
& $1.5\!\cdot\!10^{-11}$
& $1.8\!\cdot\!10^{-10}$
& $1.6\!\cdot\!10^{-11}$
& $2.7\!\cdot\!10^{-10}$
\\

&
\cellcolor{grayshade}\rule[-0.55ex]{0pt}{3.15ex}$\mathbb{D}$
& \cellcolor{grayshade}$2.9\!\cdot\!10^{-6}$
& \cellcolor{grayshade}$1.3\!\cdot\!10^{-4}$
& \cellcolor{grayshade}$4.7\!\cdot\!10^{-6}$
& \cellcolor{grayshade}$1.3\!\cdot\!10^{-4}$
\\

\addlinespace[1.0mm]
\midrule
\addlinespace[1.0mm]


\multirow{6}{*}{
\rotatebox[origin=c]{90}{\textbf{\small $\Omega$-equation}}
}
&
$B_{0.1}(0^{\star})$
& $5.0\!\cdot\!10^{-7}$
& $5.8\!\cdot\!10^{-7}$
& $5.1\!\cdot\!10^{-7}$
& $5.9\!\cdot\!10^{-7}$
\\

&
$B_1(0^\star)$
& $3.4\!\cdot\!10^{-6}$
& $1.1\!\cdot\!10^{-5}$
& $3.3\!\cdot\!10^{-6}$
& $1.1\!\cdot\!10^{-5}$
\\

&
$B_{10}(0^{\star})$
& $4.4\!\cdot\!10^{-6}$
& $2.7\!\cdot\!10^{-5}$
& $5.2\!\cdot\!10^{-6}$
& $2.8\!\cdot\!10^{-5}$
\\

&
$\mathsf{axis}$
& $4.8\!\cdot\!10^{-5}$
& $3.7\!\cdot\!10^{-4}$
& $4.9\!\cdot\!10^{-5}$
& $3.8\!\cdot\!10^{-4}$
\\

&
$\mathsf{FF}$
& $3.1\!\cdot\!10^{-6}$
& $1.4\!\cdot\!10^{-5}$
& $3.2\!\cdot\!10^{-6}$
& $1.4\!\cdot\!10^{-5}$
\\

&
\cellcolor{grayshade}$\mathbb{D}$
& \cellcolor{grayshade}$1.7\!\cdot\!10^{-5}$
& \cellcolor{grayshade}$9.5\!\cdot\!10^{-4}$
& \cellcolor{grayshade}$1.9\!\cdot\!10^{-5}$
& \cellcolor{grayshade}$9.6\!\cdot\!10^{-4}$
\\

\addlinespace[1.0mm]
\midrule
\addlinespace[1.0mm]


\multirow{6}{*}{
\rotatebox[origin=c]{90}{\textbf{\small $\Psi$-equation}}
}
&
$B_{0.1}(0^{\star})$
& $3.6\!\cdot\!10^{-7}$
& $1.2\!\cdot\!10^{-6}$
& $3.7\!\cdot\!10^{-7}$
& $4.2\!\cdot\!10^{-7}$
\\

&
$B_1(0^\star)$
& $2.8\!\cdot\!10^{-6}$
& $8.9\!\cdot\!10^{-6}$
& $2.7\!\cdot\!10^{-6}$
& $9.0\!\cdot\!10^{-6}$
\\

&
$B_{10}(0^{\star})$
& $6.9\!\cdot\!10^{-6}$
& $2.6\!\cdot\!10^{-5}$
& $8.0\!\cdot\!10^{-6}$
& $2.7\!\cdot\!10^{-5}$
\\

&
$\mathsf{axis}$
& $8.5\!\cdot\!10^{-5}$
& $1.4\!\cdot\!10^{-3}$
& $8.6\!\cdot\!10^{-5}$
& $1.5\!\cdot\!10^{-3}$
\\

&
$\mathsf{FF}$
& $2.1\!\cdot\!10^{-6}$
& $8.3\!\cdot\!10^{-6}$
& $2.1\!\cdot\!10^{-6}$
& $8.3\!\cdot\!10^{-6}$
\\

&
\cellcolor{grayshade}$\mathbb{D}$
& \cellcolor{grayshade}$3.8\!\cdot\!10^{-5}$
& \cellcolor{grayshade}$1.4\!\cdot\!10^{-3}$
& \cellcolor{grayshade}$5.3\!\cdot\!10^{-5}$
& \cellcolor{grayshade}$1.5\!\cdot\!10^{-3}$
\\

\addlinespace[0.7mm]
\bottomrule

\end{tabularx}

\vspace{1mm}
\end{table*}

\noindent For the stability analysis, we replace the PINN representation by a simpler explicit analytic profile $(\bar u, \, \bar \psi, \, \bar \omega)$, together with certified intervals for the fitted parameters. All residuals, energy estimates, and nonlinear constants are then evaluated using this fitted profile rather than the original neural network. This different representation allows for certified residuals and norms, which is essential in our stability argument as it involves the repeatedly differentiating, calculating weighted perturbations, and finding $L^\infty$ bounds of the numerical profile. These operations can be certified for explicit analytic representations (via interval arithmetic), whereas certifying them directly for the PINN would require a separate verification of the network evaluation, automatic differentiation procedure, and $L^\infty$ bounds.  

\hfill 

\noindent We must now choose an analytic representation that is sufficiently expressive to reproduce the computed profile with high accuracy, while remaining sufficiently explicit for rigorous certification. For this purpose, we use \emph{piecewise polynomial splines}. This choice is particularly well suited to verification because, on each subinterval, the approximation is given by an explicit low-degree polynomial, so that its derivatives, extrema, and residuals can be evaluated and bounded using standard rigorous tools, including interval arithmetic. Moreover, the local nature of the representation allows the mesh to be refined only where the profile exhibits rapid variation without unnecessarily increasing the complexity of the approximation elsewhere. In contrast, global polynomial or spectral representations have to find a tradeoff between global and local behavior which may require a relatively high degree of polynomials, leading to larger expressions and less convenient rigorous bounds. Piecewise splines avoid this, achieving high local accuracy with moderate polynomial degree, preserve smoothness across subintervals, and yield a representation whose errors and derivatives can be certified in an more interpretable and computationally efficient manner. \\ 

\hfill 

\noindent We define $\mathcal R_u^{\mathsf{raw}}$ and $\mathcal R_\omega^{\mathsf{raw}}$ to be the PDE residuals accounting for the fact that the numerical profiles do not satisfy exactly the Euler equations:
\begin{align}
\mathcal R_u^{\mathsf{raw}}
&:=-(\lambda r+\bar u^r)\partial_r\bar u-(\lambda z+\bar C+\bar u^z)\partial_z\bar u+2\bar u\partial_z\bar\psi+\bar c_u\bar u ,\\
\mathcal R_\omega^{\mathsf{raw}}
&:=-(\lambda r+\bar u^r)\partial_r\bar\omega-(\lambda z+\bar C+\bar u^z)\partial_z\bar\omega+2\bar u\partial_z\bar u+(\bar c_u-\lambda)\bar\omega .
\end{align}

\noindent The spline fit preserves the small profile-equation residuals observed in the numerical solution, as reported in \Cref{tab:euler_ns_errors}, so the passage to an analytic representation incurs only a minor loss of accuracy. For the stability analysis, we actually construct the surrogate through a two-head fit: \((\bar u,\bar\psi)\) are fitted separately, while \(\bar\omega\) is defined by applying the elliptic operator exactly to \(\bar\psi\). Consequently, the residual in the elliptic equation is exactly~0, eliminating a source of error and significantly simplifying the subsequent analysis.

\newpage

\noindent Overall, the certification pipeline has four separate layers, as described in the diagram below, where only the first arrow is a numerical approximation step: 

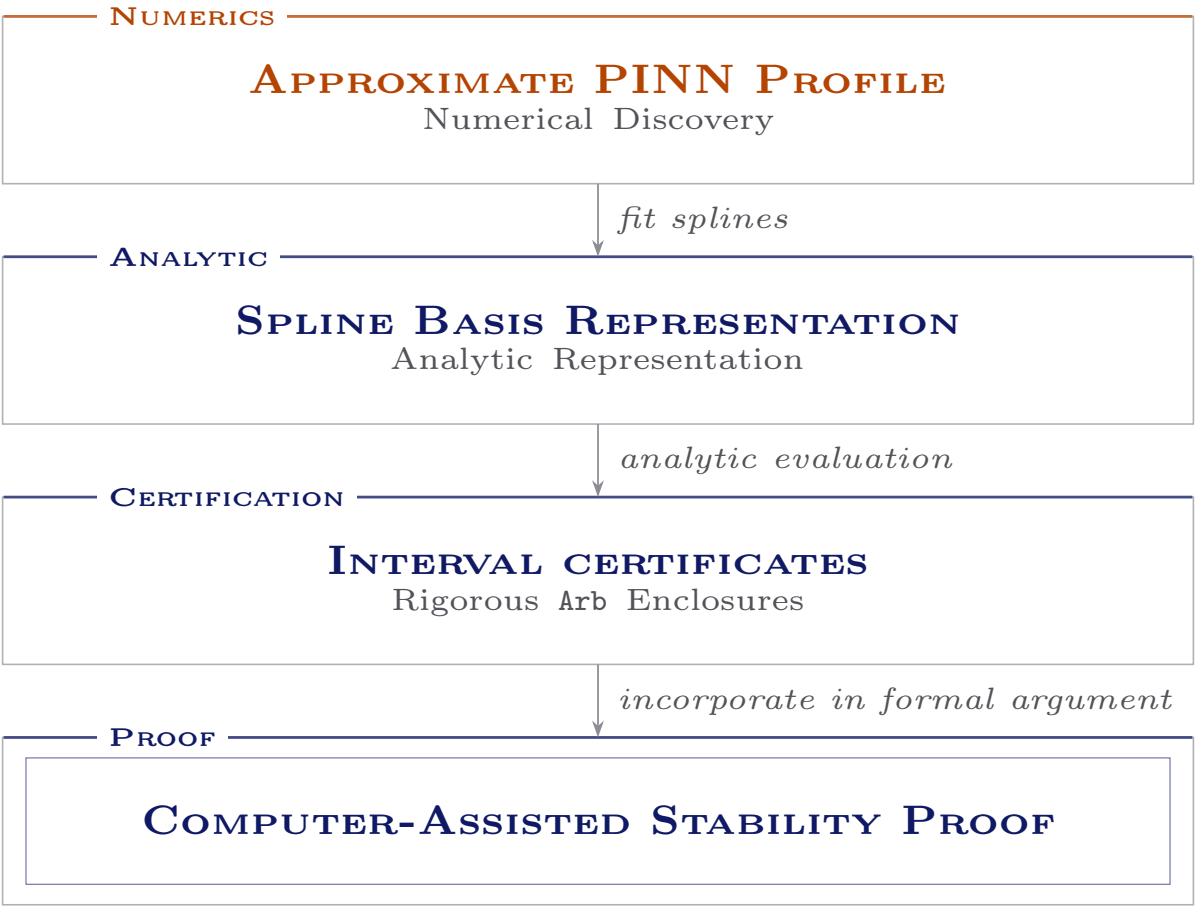
\vspace{4mm}

\begin{figure}[h]
  \centering

  \resizebox{0.9\linewidth}{!}{%
  \begin{tikzpicture}[node distance=5.0mm,font=\rmfamily]

    \certpipelinetagpanel
      {certprofile}
      {certiforange}
      {Numerics}
      {Approximate PINN Profile}
      {Numerical Discovery}
      {}

    \certpipelinetagpanel
      {certspline}
      {certifpurple}
      {Analytic}
      {Spline Basis Representation}
      {Analytic Representation}
      {below=of certprofile}

    \certpipelinetagpanel
      {certinterval}
      {certifpurple}
      {Certification}
      {Interval certificates}
      {Rigorous \texttt{Arb} Enclosures}
      {below=of certspline}

    \certpipelinetagpanel
      {certlean}
      {certifpurple}
      {Proof}
      {Computer-Assisted Stability Proof}
      {}
      {below=of certinterval}

    \draw[certifpurple!62,line width=0.22pt]
      ([xshift=1.6mm,yshift=-1.4mm]certlean.north west)
      rectangle
      ([xshift=-1.6mm,yshift=1.4mm]certlean.south east);

    \draw[certpipeline flow]
      (certprofile.south) --
      node[certpipeline transition] {fit splines}
      (certspline.north);

    \draw[certpipeline flow]
      (certspline.south) --
      node[certpipeline transition] {analytic evaluation}
      (certinterval.north);

    \draw[certpipeline flow]
      (certinterval.south) --
      node[certpipeline transition] {incorporate in formal argument}
      (certlean.north);

  \end{tikzpicture}%
  }

\end{figure}

\hfill \\

\clearpage

\section{Stability}

A singularity is \emph{stable} if the associated blowup mechanism persists under small perturbations of the initial data, in the sense that the same mechanism occurs for an open neighborhood of nearby initial conditions. Conversely, an \emph{unstable} singularity requires infinitely precise initial conditions: small perturbations deflect the evolution away from the singularity mechanism.

\hfill 

\noindent The dynamic rescaling formulation presented in \Cref{sec: Dynamic Rescaling} is the natural starting point for such a stability analysis. In the original physical variables, the solution may grow and concentrate as $t\to T^-$, so there is no fixed object around which to linearize. After rescaling, the approximate blowup profile becomes a steady state of the rescaled equations, while the blowup rate, translation rate, and amplitude normalization are recorded by modulation parameters. Linear stability can therefore be studied by perturbing this steady rescaled profile and analyzing the resulting linearized rescaled dynamics. \\

\noindent We follow the same overall structure as \citep{DGBlowup,HouJiaje,Chen2026Analysis} and related works on stability of blowups. We first compute the blowup profile numerically, then  develop a framework for establishing linear stability  of the dynamically rescaled equation in suitable weighted norms and  for proving nonlinear stability using Sobolev embeddings and weighted estimates. The main distinction is that the weight functions and other tunable parameters in the stability proof are optimized numerically to  seek sufficient margin for the argument  to close, rather than being chosen and refined through expert intuition and a potentially long and technically challenging process of trial and error. The overall proof architecture is summarized in the diagram on page~\pageref{diagram: Stability Proof}, while the optimization and certification strategy used to seek sufficient margin to close the proof is described in detail in \Cref{sec: Computer-Assisted Strategy for Proving Stability}.

\hfill

\begin{remark*}
An alternative sufficient criterion is to prove invertibility of the full linearized profile operator~\citep{elgindi2025dynamics}. In this case, the self-similar profile equation is viewed perturbatively around an approximate profile, and the goal is to invert the linearized operator after removing the neutral scaling direction. This is stronger than the criterion used here. In contrast, the numerical eigenvalue analysis in \citep[Section~4.4]{PengfeiThesis} gives a more limited, finite-dimensional test of the same stability mechanism: it consists of discretizing the dynamically rescaled linearized equation around the computed profile and inspecting the spectrum of the resulting Jacobian. This provides useful evidence for stability, but is not by itself a rigorous substitute for the required linear estimate. It only concerns a finite-dimensional discretization and therefore controls, at best, resolved low-frequency perturbations. As noted in \citep[Section~1.3.1]{chen2022singularity}, the proof must control the infinite-dimensional operator, including high-frequency modes, nonlocal velocity terms, possible neutral directions, and transient growth caused by non-normality. Thus, negative real parts of the computed eigenvalues do not yield the uniform estimate required for nonlinear stability.
\end{remark*}

\subsection{Dynamic Rescaling Equations}
\label{sec: Dynamic Rescaling}

\subsubsection{Introduction}

\noindent
The steady profile equations above describe an exact traveling self-similar solution. In computations and stability arguments, however, one usually does not know the correct profile, or traveling speed in advance. The \emph{dynamic rescaling formulation} addresses this by letting the rescaled solution evolve while continuously adjusting its location, scale, and amplitude. The key idea is to replace the finite-time singular behavior in physical variables by a long-time evolution in normalized variables. The singularity is kept at order-one size and near a fixed location in the rescaled coordinates. In this frame, blowup no longer appears primarily as unbounded growth or shrinking length scale. Instead, a stable blowup scenario should appear as convergence toward a steady state. That steady-state is precisely the traveling-wave profile described in \Cref{thm:traveling_wave_profile_eps1}. \\

\paragraph{Modulated coordinates and amplitudes.}
\noindent
In the dynamic formulation, the fixed self-similar parameters in the ansatz are replaced by time-dependent \emph{modulation parameters}. Rather than prescribing a blowup time, a spatial scale, and a traveling center in advance, the rescaled frame is adjusted as the solution evolves. More precisely, let $(\mathring r,\mathring z,\mathring t)$ denote the physical variables and let $(r,z,t)$ denote the dynamically rescaled variables.  We introduce a common spatial scale $\mathrm{s}_r(t)>0$, amplitude scales $\mathrm{s}_u(t), \, \mathrm{s}_\omega(t)>0$, and an axial center $\mathring z_c(t)$ by
\begin{equation}
\mathring r=\mathrm{s}_r(t) \, r, \qquad\quad  \mathring z-\mathring z_c(t)=\mathrm{s}_r(t) \, z,
\label{eq:dynamic_spatial_rescaling_setup}
\end{equation}
and
\begin{equation}
u_\bullet(\mathring r,\mathring z,\mathring t)
=\mathrm{s}_u(t) \, u(r,z,t), \qquad \quad 
\omega_\bullet(\mathring r,\mathring z,\mathring t) =\mathrm{s}_\omega(t) \, \omega(r,z,t).
\label{eq:dynamic_u_omega_rescaling_setup}
\end{equation}
Then,
\begin{equation}
\psi_\bullet(\mathring r,\mathring z,\mathring t)
=\mathrm{s}_\omega(t) \, \mathrm{s}_r(t)^2 \, \psi(r,z,t).
\label{eq:dynamic_psi_rescaling_setup}
\end{equation}
The normalization of the transport, stretching, and time-derivative terms is fixed by
\begin{equation}
\mathrm{s}_u(t)=\mathrm{s}_\omega(t) \, \mathrm{s}_r(t),
\qquad
\frac{dt}{d\mathring t}=\mathrm{s}_\omega(t) \, \mathrm{s}_r(t).
\label{eq:dynamic_amplitude_time_rescaling_setup}
\end{equation}
We then define the rescaling rates by
\begin{equation}
\partial_t\mathrm{s}_r(t)=-\lambda \,\mathrm{s}_r(t),
\qquad
\partial_t\mathrm{s}_u(t)=-c_u(t) \, \mathrm{s}_u(t),
\qquad
\partial_t\mathrm{s}_\omega(t)=-c_\omega(t) \, \mathrm{s}_\omega(t).
\label{eq:dynamic_rate_definitions_setup}
\end{equation}
and define the axial drift via
\begin{equation}
\partial_t\mathring z_c(t)=-C(t) \, \mathrm{s}_r(t).
\label{eq:dynamic_axial_drift_setup}
\end{equation}

\noindent Overall, this adjustment is encoded by the coordinate drift
\begin{equation}
    \left(\partial_t r\right)_{\mathring r,\mathring z}=\lambda r,
    \qquad
    \left(\partial_t z\right)_{\mathring r,\mathring z}=\lambda z+C(t).
    \label{eq:dynamic_coordinate_drift}
\end{equation}
where
\begin{itemize}
    \item $\lambda$ is the \emph{spatial rescaling rate}. It measures, in the current rescaled time variable, how quickly the computational frame expands
    \item the function $C(t)$ is the \emph{instantaneous traveling speed} of the axial center in the same rescaled frame.
    \item  the term $\lambda z$ is the axial part of the dilation, while $C(t)$ is the additional axial translation. 
\end{itemize}

\noindent The amplitude of the solution must also be normalized. This is achieved via two scalar functions $c_u(t)$ and $c_\omega(t)$, which are the logarithmic amplitude rates for $u$ and $\omega$. The $c_u(t)u$ and $c_\omega(t)\omega$ terms appearing in the rescaled evolution equations do not represent new forcing, they only record the rate at which the normalization of the rescaled variables is changed in order to keep the profile at order-one amplitude. \\

\paragraph{Use and connection with stability.}
\noindent
The profile formulation and the dynamic rescaling formulation play complementary roles. The profile equations describe the limiting singularity once the blowup rate, traveling speed, and amplitude rates have settled to constants, while the dynamic equations describe the evolution toward that limiting object. Numerically, this is essential: instead of resolving a structure whose amplitude diverges and whose length scale collapses in physical variables, one follows a normalized solution whose main profile remains order one in a fixed computational frame. The modulation parameters provide direct diagnostics for the blowup scenario: if $C(t)$, $c_u(t)$, and $c_\omega(t)$ converge, their limits give the axial traveling speed, and the amplitude growth rates. If the rescaled profiles converge at the same time, the computation gives more than evidence of singular growth. It gives a candidate self-similar profile, the values of the scaling parameters, and the profile equations that the limiting object should satisfy.\\

\noindent
The dynamic rescaling formulation is the natural setting for stability. In the original physical variables, a perturbation of a blowup solution can change the blowup time, shift the center location, or alter the amplitude normalization. Such changes may appear as growth or drift even when the underlying profile is stable. The normalization conditions remove these neutral directions by continuously recentering, rescaling, and renormalizing the solution. After this gauge freedom has been fixed, stability becomes a question about the evolution of genuine perturbations in the dynamically rescaled variables. A stable traveling self-similar singularity is represented by a stable steady state of the dynamic rescaling equations: solutions starting close to the profile should remain controlled, the modulation parameters should approach limiting constants, and the rescaled solution should converge to the steady profile. In this formulation, finite-time blowup in the original variables is linked to convergence toward a stable fixed point of the evolution.\\

\subsubsection{Dynamic Rescaling Equations for the Euler Equations}

\noindent 
As before, define the elliptic operator
\begin{equation}
    \mathcal{E} \coloneqq \partial_r^2 + \frac{3}{r}\partial_r + \partial_z^2.
\end{equation}

\hfill  

\noindent The proposition below provides the axisymmetric Euler equations in the dynamically rescaled variables.  \\

\begin{proposition}[Dynamic rescaling equations for the Euler equations]
\label{thm:dynamic_rescaling_equations eps1}
In the dynamically rescaled variables, the axisymmetric system takes the form
\begin{align}
\partial_t u
+\bigl(\lambda r+ u^r\bigr)\partial_r u
+\bigl(\lambda z+C(t)+ u^z\bigr)\partial_z u
&=2u\partial_z\psi+c_u(t)u 
\label{eq:dynamic_rescaled_u eps1}\\
\partial_t\omega
+\bigl(\lambda r+ u^r\bigr)\partial_r\omega
+\bigl(\lambda z+C(t)+ u^z\bigr)\partial_z\omega
&=2u\partial_z u+c_\omega(t)\omega
\label{eq:dynamic_rescaled_omega_eps1}\\
- \mathcal{E} \psi
&=\omega,
\label{eq:dynamic_rescaled_psi eps1}
\end{align}
with
\begin{equation}
    u^r=-r\partial_z\psi,
    \qquad
    u^z=2\psi+r\partial_r\psi,
    \label{eq:dynamic_rescaled_stream_intro eps1}
\end{equation}
and $\lambda = 1/2$. 
\end{proposition}
\begin{grayproof}
The derivation is provided in the numerical paper~\citep{EulerBlowupMain}.
\end{grayproof}

\hfill

\paragraph{Normalization conditions.}

\noindent
When $\lambda=1/2$ is fixed, the independent modulation parameters are
$C(t)$ and $c_u(t)$. The vorticity amplitude rate is determined by
\begin{equation}
    c_\omega(t)=c_u(t)-\frac{1}{2}.
\end{equation}
Before the stability analysis, write the radial and axial transport coefficients of the numerically obtained profile as
\begin{equation}
G_r(r,z)=\lambda r+u^r_{\mathsf{num}}(r,z),
\qquad
G_z(r,z)=\lambda z+C_{\mathsf{num}}+u^z_{\mathsf{num}}(r,z).
\label{eq:numerical-transport-intro-add}
\end{equation}
We choose an axial shift $z_{\mathsf{shift}}$ so that the point $(0,z_{\mathsf{shift}})$ of the numerical profile becomes the analysis origin. For this point to be a fixed point of the meridional flow, we require
\begin{equation}
G_r(0,z_{\mathsf{shift}})=0,
\qquad
G_z(0,z_{\mathsf{shift}})=0.
\label{eq:numerical-transport-stagnation-intro-add}
\end{equation}
The radial condition is automatic because $u^r_{\mathsf{num}}(0,z)=0$ on the axis. Using $u^z_{\mathsf{num}}(0,z)=2\psi_{\mathsf{num}}(0,z)$, the axial condition $G_z(0,z_{\mathsf{shift}})=0$ becomes the scalar equation
\begin{equation}
C_{\mathsf{num}}+\lambda z_{\mathsf{shift}}+2\psi_{\mathsf{num}}(0,z_{\mathsf{shift}})=0.
\label{eq:shift-root-intro-add}
\end{equation}
We impose this condition and choose $z_{\mathsf{shift}}$ as a numerically certified root. Having fixed $z_{\mathsf{shift}}$, we recenter each numerical profile field by
\begin{equation}
\bar f(r,z):=f_{\mathsf{num}}(r,z+z_{\mathsf{shift}}).
\label{eq:profile-recentering-intro-add}
\end{equation}
Denoting the constant term in the recentered axial transport by $\bar C$, we continue to write the reference transport coefficients as
\begin{equation}
G_r(r,z)=\lambda r+\bar u^r(r,z),
\qquad
G_z(r,z)=\lambda z+\bar C+\bar u^z(r,z).
\label{eq:recentered-transport-intro-add}
\end{equation}
With this notation, the imposed root condition gives
\begin{equation}
G_r(0,0)=0,
\qquad
G_z(0,0)=\bar C+2\bar\psi_0=0.
\label{eq:recentered-stagnation-intro-add}
\end{equation}

\noindent After applying this shift, set $\bar u_0:=\bar u(0,0)$. We preserve the meridional fixed point and the amplitude at the origin by imposing
\begin{equation}
C(t)+u^z(0,0,t)=0,
\qquad
u(0,0,t)=\bar u_0.
\label{eq:normalization_conditions_intro_eps1}
\end{equation}
The first condition keeps the reference meridional fixed point at the origin. The second keeps the amplitude fixed there and therefore gives $\partial_tu(0,0,t)=0$.  \\

\paragraph{Limiting behaviour.}
\noindent
If the dynamic rescaling captures a stable traveling self-similar blowup, the rescaled solution should converge as rescaled time advances to a time-independent profile,
\begin{equation}
    u(r,z,t)\to u_\infty(r,z),
    \qquad
    \omega(r,z,t)\to \omega_\infty(r,z),
    \qquad
    \psi(r,z,t)\to \psi_\infty(r,z),\label{eq:profile_convergence_intro}
\end{equation}
\noindent
and the modulation parameters should converge to constants,
\begin{equation}
    C(t)\to C_\infty,
    \qquad
    c_u(t)\to c_{u,\infty},
    \qquad
    c_\omega(t)\to c_{\omega,\infty}.
\label{eq:modulation_convergence_intro}
\end{equation}
\noindent
In that limit, the time derivatives in the dynamic rescaling equations disappear. The steady dynamic equations are
\begin{align}
\Bigl(\lambda \, r+ u^r_\infty\Bigr)\partial_r u_\infty
+\Bigl(\lambda \, z+ C_\infty + u^z_\infty \Bigr)\partial_z u_\infty
&=2 u_\infty \, \partial_z \psi_\infty+ c_{u,\infty} \, u_\infty 
\label{eq:steady_dynamic_u}\\
\Bigl(\lambda \, r+ u^r_\infty\Bigr)\partial_r\omega_\infty
+\Bigl(\lambda \, z+C_\infty+ u^z_\infty\Bigr)\partial_z\omega_\infty
&=2 u_\infty \, \partial_z u_\infty+ c_{\omega,\infty} \, \omega_\infty  
\label{eq:steady_dynamic_w}
\end{align}
\noindent
with $u^r_\infty =-r\partial_z \psi_\infty$ and $u^z_\infty=2 \psi_\infty+r\partial_r\psi_\infty$. \\ 

\noindent
For the traveling-wave ansatz in \eqref{ansantz1}--\eqref{ansantz3}, the limiting amplitude rates are
\begin{equation}
    c_{u,\infty}=-1,
    \qquad
    c_{\omega,\infty} = -(1+\lambda)=-3/2.
    \label{eq:limiting_amplitude_rates}
\end{equation}

\noindent By writing
\begin{equation}
\bigl(u_\infty, \, , \, \omega_\infty ,\psi_\infty  ,\, C_\infty  \bigr) = \bigl(U, \, \Omega, \, \Psi, \, C\bigr),
\end{equation} 
and substituting \eqref{eq:limiting_amplitude_rates} into \eqref{eq:steady_dynamic_u}--\eqref{eq:steady_dynamic_w}, we see that the limiting dynamic equations agree with the corresponding profile equations. In this sense, the traveling-wave profiles are fixed points of the appropriate dynamically rescaled evolution. \\

\subsection{Linearization}

\noindent We now introduce the perturbation equations around a steady state of the dynamically rescaled system. We then separate the terms that are linear in the perturbation from the higher-order remainders. The linearized equations identify the operator whose stability determines whether the candidate blowup profile is attracting to first order in the rescaled variables.

\subsubsection{Perturbations}

\noindent We denote perturbations using \textcolor{blue}{$\delta$}'s in \textcolor{blue}{blue}. We use the subscript $_0$ to denote evaluation at the origin $(r,z)=(0,0)$.  We perturb a given steady state
\begin{equation}
    \bar u,\bar\omega,\bar\psi,\bar C,\bar c_u,\bar c_\omega
\end{equation}
of the dynamic rescaling equations, by writing
\begin{equation}
    u = \bar u + \textcolor{blue}{\delta u},
    \qquad
    \omega=\bar\omega+\textcolor{blue}{\delta\omega},
    \qquad
    \psi = \bar\psi + \textcolor{blue}{\delta\psi},
\label{eq:pert_ansatz_1}
\end{equation}
\begin{equation}
    C=\bar C+\textcolor{blue}{\delta C},
    \qquad
    c_u = \bar c_u + \textcolor{blue}{\delta c_u},
    \qquad
    c_\omega=\bar c_\omega+\textcolor{blue}{\delta c_\omega}. \label{eq:pert_ansatz_2}
\end{equation}

\hfill 

\noindent The perturbation variables measure the difference between the evolving rescaled solution and the steady blowup profile. The goal of linear stability analysis is to determine whether these perturbations remain bounded, or preferably decay, under the dynamics obtained by linearizing around the steady state. \\

\noindent The perturbations must inherit the symmetries of the underlying steady state and therefore belong to the same functional class. In particular, the parity and regularity conditions imposed on the rescaled variables are inherited by $\textcolor{blue}{\delta u}$, $\textcolor{blue}{\delta\omega}$, and $\textcolor{blue}{\delta\psi}$, with the latter two coupled by the elliptic equation relating vorticity and streamfunction. The perturbations must also decay at the far-field. The admissible perturbation class includes the weighted regularity, integrability, trace, and boundary-flux assumptions used in the energy estimates. In particular, the cutoff boundary-flux terms vanish as the cutoffs are removed. This ensures that the perturbation equations remain in the same functional class as the steady profile. \\

\noindent For notational convenience, we will use $\textcolor{blue}{\delta}$ as shorthand for the full list of perturbation variables,
\begin{equation}
\label{eq:delta_shorthand_intro}
    \textcolor{blue}{\delta}
    :=
    \{ \textcolor{blue}{\delta u},\textcolor{blue}{\delta\omega},\textcolor{blue}{\delta\psi},\textcolor{blue}{\delta C},
    \textcolor{blue}{\delta c_u} \}.
\end{equation}

\noindent We also write the meridional velocity as
\begin{equation} \label{eq: meridional}
    u^r=\bar u^r+\textcolor{blue}{\delta u^r},
    \qquad
    u^z=\bar u^z+\textcolor{blue}{\delta u^z},
\end{equation}
where
\begin{equation} \label{eq:velocity_perturbation_intro}
\textcolor{blue}{\delta u^r} = -r\,\partial_z\textcolor{blue}{\delta\psi},
\qquad
\textcolor{blue}{\delta u^z} = 2\,\textcolor{blue}{\delta\psi} + r\,\partial_r\textcolor{blue}{\delta\psi}, 
\end{equation}
and also have
\begin{equation}
\bar c_\omega=\bar c_u-\lambda,
\qquad
\textcolor{blue}{\delta c_\omega}=\textcolor{blue}{\delta c_u}.
\end{equation}

\hfill \\

\noindent \textbf{Modulations and Normalization Conditions.} \, The modulation parameters are fixed by normalization conditions at the origin. We impose
\begin{equation}
    C(t)+u^z_0(t)=0,
    \qquad
    u_0(t)=\bar u_0.
\end{equation}

\hfill 

\noindent For the perturbations, these conditions give
\begin{equation} \label{eq: perturbation normalization}
    \textcolor{blue}{\delta C}+2\textcolor{blue}{\delta\psi}_0=0,
    \qquad
    \textcolor{blue}{\delta u}_0 = 0,
    \qquad
    (\partial_t\textcolor{blue}{\delta u})_0=0.
\end{equation}

\hfill

\noindent These normalization conditions fix the axial translation and amplitude normalization associated with $C(t)$ and $c_u(t)$. The spatial rate $\lambda=1/2$ has already been fixed. They are essential for formulating linear stability in the modulated variables, since otherwise the perturbation could drift along symmetry directions rather than measuring a genuine change in the profile.\\

\hfill 

\subsubsection{Choice of Variables for the Linear Stability Analysis}

The reduced variables $u$ and $\omega$ do not transform in the same way under rescaling~\citep{EulerBlowupMain}. Indeed, the scaling of the full three-dimensional fields induces the following scaling of the reduced variables:
\begin{equation}
    u_\eta(r,z,t)=\eta^{\alpha+1}u(\eta r,\eta z,\eta^{\alpha+1}t),
    \qquad
    \omega_\eta(r,z,t)=\eta^{\alpha+2}\omega(\eta r,\eta z,\eta^{\alpha+1}t).
    \label{eq:euler-main-scaling}
\end{equation}

\noindent The relative scaling of $u$, $\omega$, and $\nabla u$ is unchanged. Thus $\omega$ carries one extra factor of $\eta$ relative to $u$, reflecting the one-derivative relation between the underlying full velocity and vorticity fields. By contrast,
\begin{equation}
    \omega_\eta(r,z,t)=\eta^{\alpha+2}\omega(\eta r,\eta z,\eta^{\alpha+1}t),
    \qquad
    \nabla u_\eta(r,z,t)=\eta^{\alpha+2}(\nabla u)(\eta r,\eta z,\eta^{\alpha+1}t).
\end{equation}

\hfill

\noindent Therefore $\omega$ and $\nabla u$ are the natural pair of quantities in a scale-consistent linear stability analysis. Quantities with different scaling correspond to different differential orders, so comparing them directly mixes objects of different analytic strength. For this reason, instead of using the equations for $(\partial_t\textcolor{blue}{\delta \omega},\partial_t\textcolor{blue}{\delta u})$, the stability analysis is formulated using the equations for $    \partial_t\textcolor{blue}{\delta\omega},
 \,
    \partial_t\partial_z\textcolor{blue}{\delta u},
 \,
    \partial_t\partial_r\textcolor{blue}{\delta u}$. We collect the associated perturbations into the scale-consistent perturbation vector
\begin{equation}
\label{eq:ns-deltav-class}
\textcolor{blue}{\delta v}
:=
(\textcolor{blue}{\delta\omega},\nabla\textcolor{blue}{\delta u})
=
(\textcolor{blue}{\delta\omega},\partial_r\textcolor{blue}{\delta u},\partial_z\textcolor{blue}{\delta u}),
\end{equation}
and for compact notation, we also set
\begin{equation}
\textcolor{blue}{\delta v_\omega}:=\textcolor{blue}{\delta\omega},
\qquad
\textcolor{blue}{\delta v_r}:=\partial_r\textcolor{blue}{\delta u},
\qquad
\textcolor{blue}{\delta v_z}:=\partial_z\textcolor{blue}{\delta u}.
\end{equation}

\hfill  

\subsubsection{Perturbation Equations}

Linearization consists of expanding the perturbation equation around the candidate blowup profile and retaining only the terms that are first order in the perturbation. The resulting linearized equation gives the leading-order dynamics of small disturbances near the profile. The discarded terms are collected in nonlinear remainders, which record the higher-order interactions among the perturbations. The full derivations of the linearized and modulation equations are given in the numerical paper~\citep{EulerBlowupMain}. \\

\noindent We start from the dynamic rescaling equations~\eqref{eq:dynamic_rescaled_u eps1}--\eqref{eq:dynamic_rescaled_psi eps1}:
\begin{align}
\partial_t u
+\bigl(\lambda r+ u^r\bigr)\partial_r u
+\bigl(\lambda z+C(t)+ u^z\bigr)\partial_z u
&=2u\partial_z\psi+c_u(t)u 
\label{eq:dynamic_rescaled_u stab}\\
\partial_t\omega
+\bigl(\lambda r+ u^r\bigr)\partial_r\omega
+\bigl(\lambda z+C(t)+ u^z\bigr)\partial_z\omega
&=2u\partial_z u+c_\omega(t)\omega 
\label{eq:dynamic_rescaled_omega stab}\\
- \mathcal{E} \psi
&=\omega.
\label{eq:dynamic_rescaled_psi stab}
\end{align}
Here,
\begin{equation}
    u^r=-r\partial_z\psi,
    \qquad
    u^z=2\psi+r\partial_r\psi,
    \label{eq:dynamic_rescaled_stream_intro stab}
\end{equation}
and $\mathcal{E}$ is the elliptic operator
\begin{equation}
    \mathcal{E} \coloneqq \partial_r^2 + \frac{3}{r}\partial_r + \partial_z^2.
\end{equation}

\hfill \\

\noindent Substituting \eqref{eq:pert_ansatz_1}--\eqref{eq:pert_ansatz_2}, subtracting the steady-state equations, re-centering the modulation coefficients as specified below, and separating linear terms from nonlinear terms in the perturbations gives
\begin{equation}
\label{eq:linearized_system_intro}
    \partial_t\textcolor{blue}{\delta u}
    =
    \mathcal L_u(\textcolor{blue}{\delta}) + \mathcal N_u(\textcolor{blue}{\delta}) \ + \ \mathcal R_u,
    \qquad
    \partial_t\textcolor{blue}{\delta\omega}
    =
    \mathcal L_\omega(\textcolor{blue}{\delta}) + \mathcal N_\omega(\textcolor{blue}{\delta}) \ + \ \mathcal R_\omega.
\end{equation}

\hfill 

\noindent The raw profile residuals $\mathcal R_u^{\mathsf{raw}}$ and $\mathcal R_\omega^{\mathsf{raw}}$ account for the fact that the numerical profiles are not exact:
\begin{align}
\mathcal R_u^{\mathsf{raw}}
&=-(\lambda r+\bar u^r)\partial_r\bar u
-(\lambda z+\bar C+\bar u^z)\partial_z\bar u
+2\bar u\partial_z\bar\psi+\bar c_u\bar u ,\\
\mathcal R_\omega^{\mathsf{raw}}
&=-(\lambda r+\bar u^r)\partial_r\bar\omega
-(\lambda z+\bar C+\bar u^z)\partial_z\bar\omega
+2\bar u\partial_z\bar u+(\bar c_u-\lambda)\bar\omega .
\end{align}

\hfill

\noindent We also define the residual amplitude offset
\begin{equation}
c_{u,\mathrm{res}}
:=-\frac{\mathcal R_u^{\mathsf{raw}}(0)}{\bar u_0},
\label{eq:residual_phase_offsets_intro}
\end{equation}
 under the nondegeneracy condition $\bar u_0\ne0$. \\
 
 \noindent We absorb this fixed offset into the reference amplitude rate. From this point onward, $\bvar{\delta C}$ and $\bvar{\delta c_u}$ denote the variable modulation coefficients, so that
\begin{equation}
C=\bar C+\bvar{\delta C},
\qquad
c_u=\bar c_u+c_{u,\mathrm{res}}+\bvar{\delta c_u},
\qquad
c_\omega=\bar c_\omega+c_{u,\mathrm{res}}+\bvar{\delta c_u}.
\label{eq:recentered_modulation_variables_intro}
\end{equation}

\hfill \\

\noindent The residuals used in the perturbation equations are
\begin{align}
\mathcal R_u
&:=\mathcal R_u^{\mathsf{raw}}+c_{u,\mathrm{res}}\,\bar u,\\
\mathcal R_\omega
&:=\mathcal R_\omega^{\mathsf{raw}}+c_{u,\mathrm{res}}\,\bar\omega.
\label{eq:recentered_profile_residuals_intro}
\end{align}

\hfill \\

\noindent The linearized $u$-operator and $\omega$-operator are
\vspace{2mm}

\begin{align}
\mathcal L_u(\textcolor{blue}{\delta})
:={}&
\bigl(2\partial_z\bar\psi+\bar c_u+c_{u,\mathrm{res}}\bigr)\,\textcolor{blue}{\delta u}
-\bigl(\lambda r+\bar u^r\bigr)\,\partial_r \textcolor{blue}{\delta u}
-\bigl(\lambda z + \bar C + \bar u^z\bigr)\,\partial_z \textcolor{blue}{\delta u}
+\bigl(2\bar u+ r\,\partial_r\bar u\bigr)\partial_z\textcolor{blue}{\delta\psi} \nonumber\\
&\quad
-\bigl( r\,\partial_z\bar u\bigr)\partial_r\textcolor{blue}{\delta\psi}
-\bigl(2\partial_z\bar u\bigr)\textcolor{blue}{\delta\psi}
-(\partial_z\bar u)\,\bvar{\delta C}
+\bar u\,\bvar{\delta c_u},
\label{eq:Lu_Euler_definition}
\end{align}
\begin{align}
\mathcal L_\omega(\textcolor{blue}{\delta})
:={}&
(\bar c_u+c_{u,\mathrm{res}}-\lambda)\,\textcolor{blue}{\delta\omega}
-\bigl(\lambda r + \bar u^r \bigr)\partial_r\textcolor{blue}{\delta\omega}
-\left(\lambda z+\bar C+ \bar u^z \right)\partial_z\textcolor{blue}{\delta\omega}
+2(\partial_z\bar u)\,\textcolor{blue}{\delta u}
+2\bar u\,\partial_z\textcolor{blue}{\delta u}
\nonumber\\
&\quad
-2(\partial_z\bar\omega)\,\textcolor{blue}{\delta\psi}
+ r\,(\partial_r\bar\omega)\,\partial_z\textcolor{blue}{\delta\psi}
- r\,(\partial_z\bar\omega)\,\partial_r\textcolor{blue}{\delta\psi}
-(\partial_z\bar\omega)\,\bvar{\delta C}
+\bar\omega\,\bvar{\delta c_u}.
\label{eq:Lw_Euler_definition}
\end{align}

\hfill \\

\noindent The nonlinear remainders are
\vspace{1mm}

\begin{align}
\mathcal N_u(\textcolor{blue}{\delta})
={}&
 r\,\partial_z\textcolor{blue}{\delta\psi}\,\partial_r\textcolor{blue}{\delta u}
-\left(\bvar{\delta C}+\left(2\textcolor{blue}{\delta\psi}+r\partial_r\textcolor{blue}{\delta\psi}\right)\right)\partial_z\textcolor{blue}{\delta u}
+2\,\textcolor{blue}{\delta u}\,\partial_z\textcolor{blue}{\delta\psi}
+\bvar{\delta c_u}\,\textcolor{blue}{\delta u},
\label{eq:Nu_Euler_definition} \\ 
\mathcal N_\omega(\textcolor{blue}{\delta})
={}&
- \textcolor{blue}{\delta u^r} \partial_r\textcolor{blue}{\delta\omega}
-\left(\bvar{\delta C}+\textcolor{blue}{\delta u^z}\right)\partial_z\textcolor{blue}{\delta\omega}
+2\,\textcolor{blue}{\delta u}\,\partial_z\textcolor{blue}{\delta u}
+\bvar{\delta c_u}\,\textcolor{blue}{\delta\omega},
\label{eq:Nw_Euler_definition}
\end{align}

\hfill

\noindent The elliptic equation is already linear, so
\begin{equation}
\label{eq:linearized_elliptic_intro}
 -\mathcal E  \textcolor{blue}{\delta\psi}  \coloneqq   -\left(\partial_r^2 + \frac{3}{r}\partial_r + \partial_z^2\right)\textcolor{blue}{\delta\psi}
    = \textcolor{blue}{\delta\omega}.
\end{equation}

\hfill 

\noindent Consequently, $\textcolor{blue}{\delta\psi}$ is not an independent dynamical unknown: it is recovered from $\textcolor{blue}{\delta\omega}$ through this elliptic equation, and then $\textcolor{blue}{\delta u^r}$ and $\textcolor{blue}{\delta u^z}$ are recovered from \eqref{eq:velocity_perturbation_intro}. 

\hfill \\

\subsubsection{Differentiated Equations for the Scale-Consistent Variables}

Since the scale-consistent variables are $\omega$ and $\nabla u$, we differentiate the $\delta u$-equation:
\begin{align}
\label{eq: differentiated linearized}
\partial_t\partial_r\textcolor{blue}{\delta u}
&=\mathcal L_r(\textcolor{blue}{\delta})+
\mathcal N_{r}(\textcolor{blue}{\delta}) \ +\ \mathcal R_r,
\qquad \quad 
\partial_t\partial_z\textcolor{blue}{\delta u}
=\mathcal L_z(\textcolor{blue}{\delta})+
\mathcal N_{z}(\textcolor{blue}{\delta})  \ +\ \mathcal R_z,
\end{align}
where we denote
\begin{equation}
\mathcal L_r(\textcolor{blue}{\delta}):=\partial_r\mathcal L_u(\textcolor{blue}{\delta}),
\qquad
\mathcal N_r(\textcolor{blue}{\delta}):=\partial_r\mathcal N_u(\textcolor{blue}{\delta}), \qquad \mathcal L_z(\textcolor{blue}{\delta}):=\partial_z\mathcal L_u(\textcolor{blue}{\delta}),
\qquad
\mathcal N_z(\textcolor{blue}{\delta}):=\partial_z\mathcal N_u(\textcolor{blue}{\delta}),
\end{equation}
and 
\begin{equation}
    \mathcal R_z:=\partial_z\mathcal R_u,
\qquad\quad 
\,\mathcal R_r:=\,\partial_r\mathcal R_u.
\end{equation}

\hfill 

\noindent Differentiating the $\mathcal{L}_u$ expressions gives
\vspace{1mm}

\begin{align}
\mathcal L_r(\textcolor{blue}{\delta})
={}&
-(\lambda r+\bar u^r)\,\partial_{rr}\textcolor{blue}{\delta u}
-(\lambda z+\bar C+\bar u^z)\,\partial_{zr}\textcolor{blue}{\delta u}
+\bigl(2\partial_z\bar\psi+\bar c_u+c_{u,\mathrm{res}}-\lambda-\partial_r\bar u^r\bigr)
\partial_r\textcolor{blue}{\delta u}
\nonumber\\
&
-(\partial_r\bar u^z)\,\partial_z\textcolor{blue}{\delta u}  +2\partial_{rz}\bar\psi\,\textcolor{blue}{\delta u} +(2\bar u+ r\partial_r\bar u)\,\partial_{zr}\textcolor{blue}{\delta\psi}
+\bigl(3\partial_r\bar u+ r\partial_{rr}\bar u\bigr)
\partial_z\textcolor{blue}{\delta\psi}
- r(\partial_z\bar u)\,\partial_{rr}\textcolor{blue}{\delta\psi}
\nonumber\\
& -\bigl(3\partial_z\bar u+ r\partial_{rz}\bar u\bigr)
\partial_r\textcolor{blue}{\delta\psi}
-2\partial_{rz}\bar u\,\textcolor{blue}{\delta\psi}
-(\partial_{rz}\bar u)\,\bvar{\delta C}
+(\partial_r\bar u)\,\bvar{\delta c_u},
\label{eq:Lr_Euler}
\end{align}

\begin{align}
\mathcal L_z(\textcolor{blue}{\delta})
={}&
-(\lambda r+\bar u^r)\,\partial_{rz}\textcolor{blue}{\delta u}
-(\lambda z+\bar C+\bar u^z)\,\partial_{zz}\textcolor{blue}{\delta u}
+\bigl(2\partial_z\bar\psi+\bar c_u+c_{u,\mathrm{res}}-\lambda-\partial_z\bar u^z\bigr)
\partial_z\textcolor{blue}{\delta u}
\nonumber\\
& -(\partial_z\bar u^r)\,\partial_r\textcolor{blue}{\delta u} +2\partial_{zz}\bar\psi\,\textcolor{blue}{\delta u}
+(2\bar u+ r\partial_r\bar u)\,\partial_{zz}\textcolor{blue}{\delta\psi}
+\bigl(r\partial_{rz}\bar u\bigr)
\partial_z\textcolor{blue}{\delta\psi}
- r(\partial_z\bar u)\,\partial_{rz}\textcolor{blue}{\delta\psi}
\nonumber\\
& - r(\partial_{zz}\bar u)\,\partial_r\textcolor{blue}{\delta\psi}
-2\partial_{zz}\bar u\,\textcolor{blue}{\delta\psi}
-(\partial_{zz}\bar u)\,\bvar{\delta C}
+(\partial_z\bar u)\,\bvar{\delta c_u}.
\label{eq:Lz_Euler_explicit}
\end{align}

\hfill 

\noindent The differentiated nonlinear remainders are
\vspace{1mm}

\begin{align}
\mathcal N_r(\textcolor{blue}{\delta})
={}&
\left(\partial_z\textcolor{blue}{\delta\psi}+r\partial_{rz}\textcolor{blue}{\delta\psi}\right)
\partial_r\textcolor{blue}{\delta u}
+ r\partial_z\textcolor{blue}{\delta\psi}\,\partial_{rr}\textcolor{blue}{\delta u}
-\left(3\partial_r\textcolor{blue}{\delta\psi}+r\partial_{rr}\textcolor{blue}{\delta\psi}\right)
\partial_z\textcolor{blue}{\delta u}
\label{eq:Nr_Euler}\\
& -\left(\bvar{\delta C}+\left(2\textcolor{blue}{\delta\psi}+r\partial_r\textcolor{blue}{\delta\psi}\right)\right)
\partial_{zr}\textcolor{blue}{\delta u}
+2(\partial_r\textcolor{blue}{\delta u})(\partial_z\textcolor{blue}{\delta\psi})
+2\textcolor{blue}{\delta u}\,\partial_{zr}\textcolor{blue}{\delta\psi}
+\bvar{\delta c_u}\,\partial_r\textcolor{blue}{\delta u},
\nonumber \end{align}
\begin{align} 
\mathcal N_z(\textcolor{blue}{\delta})
={}&
 r(\partial_{zz}\textcolor{blue}{\delta\psi})\,\partial_r\textcolor{blue}{\delta u}
+ r(\partial_z\textcolor{blue}{\delta\psi})\,\partial_{rz}\textcolor{blue}{\delta u}
-\left(2\partial_z\textcolor{blue}{\delta\psi}
+r\partial_{rz}\textcolor{blue}{\delta\psi}\right)\partial_z\textcolor{blue}{\delta u}
\label{eq:Nz_Euler} \\
& -\left(\bvar{\delta C}+\left(2\textcolor{blue}{\delta\psi}+r\partial_r\textcolor{blue}{\delta\psi}\right)\right)
\partial_{zz}\textcolor{blue}{\delta u}
+2(\partial_z\textcolor{blue}{\delta u})(\partial_z\textcolor{blue}{\delta\psi})
+2\textcolor{blue}{\delta u}\,\partial_{zz}\textcolor{blue}{\delta\psi}
+\bvar{\delta c_u}\,\partial_z\textcolor{blue}{\delta u}, \nonumber
\end{align}

\hfill 

\subsubsection{Modulation Equations}

The modulation equations are derived in the numerical paper~\citep{EulerBlowupMain} by imposing the normalization conditions \eqref{eq: perturbation normalization}. They express the perturbations of the rescaling parameters in terms of the profile and its perturbations evaluated at the origin. The translation formula is the perturbation of the stagnation condition $C+u_0^z=0$. For the amplitude formula, we differentiate the normalization $u_0=\bar u_0$ in time, evaluate the $u$-equation at the meridional fixed point, and solve the resulting identity for $\bvar{\delta c_u}$.  \\

 \noindent With the fixed residual amplitude offset already absorbed into the reference coefficient, the modulation coefficients are
\begin{align}
\bvar{\delta C}
&=-2\bvar{\delta\psi}_0,
\label{eq:ns_deltaC_cu}, \qquad 
\bvar{\delta c_u}
=-2(\partial_z\bvar{\delta\psi})_0
\end{align}
Both $\bvar{\delta C}$ and $\bvar{\delta c_u}$ are linear in the perturbation. Products such as $\bvar{\delta C}\,\partial_z\bvar{\delta u}$ and $\bvar{\delta c_u}\,\bvar{\delta u}$ remain ordinary quadratic nonlinear terms.\\

\noindent Equation \eqref{eq:ns_deltaC_cu} determines the two modulation variables from the perturbations, provided the nondegeneracy condition $\bar u_0\ne0$ holds true.

\hfill

\subsubsection{Closed Linearized Variables}

 In terms of $\textcolor{blue}{\delta v}$, the closed linearized system may be written schematically as
\begin{equation}
\label{eq:closed-linearized-deltav}
\partial_t\textcolor{blue}{\delta v}
=
\mathcal L_{\nabla}\textcolor{blue}{\delta v},
\end{equation}
or equivalently in components:
\begin{equation}
\label{eq:closed-linearized-components}
\partial_t
\begin{pmatrix}
\textcolor{blue}{\delta\omega}\\
\partial_r\textcolor{blue}{\delta u}\\
\partial_z\textcolor{blue}{\delta u}
\end{pmatrix}
=
\mathcal L_{\nabla}
\begin{pmatrix}
\textcolor{blue}{\delta\omega}\\
\partial_r\textcolor{blue}{\delta u}\\
\partial_z\textcolor{blue}{\delta u}
\end{pmatrix}.
\end{equation}
The operator $\mathcal L_{\nabla}$ includes transport by the steady rescaled flow, coupling between $\textcolor{blue}{\delta u}$ and $\textcolor{blue}{\delta\omega}$, elliptic recovery of $\textcolor{blue}{\delta\psi}$, velocity reconstruction, and the linearized modulation terms. 

\hfill \\

\subsection{Weighted Norms and Energy}

\subsubsection{Weighted Norms}

\noindent The stability estimates are measured in weighted norms adapted to the candidate blowup profile and to the damping structure of the linearized operator. For this purpose, define
\begin{equation}
\wgray{\Phi_\omega},\ \wgray{\Phi_r},\ \wgray{\Phi_z}:\mathbb D\to(0,\infty)
\label{eq:low-order-weight-functions-main}
\end{equation}
to be the positive weight functions assigned respectively to the three components of $\textcolor{blue}{\delta v}$, on the half-plane
\begin{equation}
\mathbb D=\{(r,z):r\ge0,\ z\in\mathbb R\}.
\end{equation}

\hfill 

\noindent More explicitly, $\wgray{\Phi_\omega}$ weights the vorticity perturbation $\textcolor{blue}{\delta\omega}$, $\wgray{\Phi_r}$ weights $\partial_r\textcolor{blue}{\delta u}$, and $\wgray{\Phi_z}$ weights $\partial_z\textcolor{blue}{\delta u}$. We denote the collection of weights by
\begin{equation}
\label{eq:Phi-vector}
\wgray{\Phi}
:=
(\wgray{\Phi_\omega},\wgray{\Phi_r},\wgray{\Phi_z}).
\end{equation}

\hfill

\noindent For a positive scalar weight $\wgray{\phi}$, define the associated inner product and norm by
\begin{equation}
\label{eq:ns-weighted-inner}
\|f\|_{\wgray{\phi}}^2 := \int_\mathbb{D} |f|^2 \wgray{\phi}\,dr\,dz,
\qquad
\langle f,g\rangle_{\wgray{\phi}}:= \int_\mathbb{D} fg\,\wgray{\phi}\,dr\,dz.
\end{equation}

\hfill 

\noindent Given vector quantities $q := (q_\omega,q_r,q_z)$, and   $\tilde q := (\tilde q_\omega,\tilde q_r,\tilde q_z)$, we define the corresponding $\wgray{\Phi}$-weighted norm and $\wgray{\Phi}$-weighted inner product by 
\begin{equation}
\label{eq:Phi-vector-norm}
\|q\|_{\wgray{\Phi}}^2 := \|q_\omega\|_{\wgray{\Phi_\omega}}^2 + \|q_r\|_{\wgray{\Phi_r}}^2 + \|q_z\|_{\wgray{\Phi_z}}^2. \end{equation}
\begin{equation}
\label{eq:Phi-vector-inner}
\langle q,\tilde q\rangle_{\wgray{\Phi}} := \langle q_\omega,\tilde q_\omega\rangle_{\wgray{\Phi_\omega}} + \langle q_r,\tilde q_r\rangle_{\wgray{\Phi_r}} + \langle q_z,\tilde q_z\rangle_{\wgray{\Phi_z}}.
\end{equation}

\hfill

\subsubsection{Low-Order Weighted Energy.} 

\noindent The low-order energy uses the three singular weights $\wgray{\Phi_\omega}$, $\wgray{\Phi_r}$, $\wgray{\Phi_z}$, and is defined as the $\wgray{\Phi}$-norm of $\textcolor{blue}{\delta v}$:
\begin{equation}
\label{eq:ns-E0}
\mathfrak{E}_0(t) := \|\textcolor{blue}{\delta v}(t)\|_{\wgray{\Phi}}, \qquad \text{i.e.} \quad  \mathfrak{E}_0(t)^2 = \|\textcolor{blue}{\delta\omega}(t)\|_{\wgray{\Phi_\omega}}^2 + \|\partial_r\textcolor{blue}{\delta u}(t)\|_{\wgray{\Phi_r}}^2 + \|\partial_z\textcolor{blue}{\delta u}(t)\|_{\wgray{\Phi_z}}^2.
\end{equation}

\noindent The three weights $\wgray{\Phi_\omega}$, $\wgray{\Phi_r}$, $\wgray{\Phi_z}$ are kept separate because the damping calculation acts on the three components of $\textcolor{blue}{\delta v}$ in different ways. The purpose of the weighted energy is to make the coercive part of the linear part produce a positive damping constant $\DLamL>0$. A discussion of the constraints that the weight functions $\wgray{\Phi_\omega}$, $\wgray{\Phi_r}$, and $\wgray{\Phi_z}$ must satisfy in our setting is provided in \Cref{sec: Constraints on phi}. \\

\noindent While the low-order energy remains the main source of damping, it is not sufficient by itself to close the estimates since the linear and nonlinear terms contain terms which cannot be controlled using a weighted $L^2$ norm, such as terms that must be controlled in $L^\infty$, as well as modulation factors and pointwise quantities evaluated at the origin, such as $(\partial_z\bvar{\delta u})_0$ and $(\partial_{z}\bvar{\delta\omega})_0$.

\hfill

\subsubsection{High-Order Weighted Energy.} 

To control the linear and nonlinear terms which cannot be controlled in the low-order weighted energy $\mathfrak E_0$, we supplement $\mathfrak E_0$ with a high-order weighted energy $\mathfrak H_k$, chosen so that the top derivatives provide the required pointwise, interpolation, elliptic, and modulation bounds. \\

\noindent The high-order part uses the \emph{higher-order weight functions}
\begin{equation} 
  \label{eq:ns-compensated-weights-intro}
 \wgray{\Phi_{i,k}}:=1+\rho^{\,k}\wgray{\Phi_i},
\qquad \rho(r,z):=r^2+z^2,
\qquad i\in\{\omega,r,z\}. 
\end{equation}

\noindent The factor $\rho^{\,k}=|(r,z)|^{2k}$ compensates for the loss of $k$ powers under differentiation near the origin. If $f$ has local order $|(r,z)|^m$, then $D^k f$ has local order $|(r,z)|^{m-k}$; consequently, the factors $|D^k f|^2\rho^{,k}\wgray{\Phi_i}$ and $|f|^2\wgray{\Phi_i}$ have the same local power of $|(r,z)|$. The added $1$ prevents the high-order weight from degenerating at the origin and gives the direct bound
\begin{equation}
\wgray{\Phi_{i,k}}=1+\rho^k\wgray{\Phi_i}\ge 1,
\qquad
\|D^k f\|_{L^2(\mathbb D)}\le \|D^k f\|_{\wgray{\Phi_{i,k}}}.
\label{eq:high-weight-unit-floor-main}
\end{equation}

\noindent We then define the \emph{high-order weighted energy}
\begin{equation}
\label{eq:ns-Hk}
\mathfrak H_k^2
:={}\sum_{|\alpha|=k}\|D^\alpha\bvar{\delta\omega}\|_{\wgray{\Phi_{\omega,k}}}^2
+\sum_{|\alpha|=k}\|D^\alpha\partial_r\bvar{\delta u}\|_{\wgray{\Phi_{r,k}}}^2 +\sum_{|\alpha|=k}\|D^\alpha\partial_z\bvar{\delta u}\|_{\wgray{\Phi_{z,k}}}^2. 
\end{equation}

\hfill  

\subsubsection{Full Energy.} The full energy used for the stability argument is obtained as
\begin{equation} 
\label{eq:ns-Ek}
\mathfrak E_k^2:=\mathfrak E_0^2+\mu_k\mathfrak H_k^2,
\qquad 0<\mu_k\le 1.
\end{equation}
The parameter $\mu_k$ is chosen only small enough so the low-order damping absorbs the low-order loss coming from the differentiated linear estimate, but also large enough so that the differentiated damping absorbs the top-order loss left by the low-order linear estimate. \\

\paragraph{Choosing $k$.} \noindent The order of the top derivative $k$ should be chosen large enough so that the energy controls all point values, products, commutators, and elliptic or streamfunction terms used in the estimates. The modulation formulas require the origin values
\begin{equation}
\begin{aligned}
&\bvar{\delta\psi}_0,
\qquad
(\partial_z\bvar{\delta\psi})_0.
\end{aligned}
\label{eq:modulation-origin-values-k-choice}
\end{equation}
\noindent In two dimensions, a point value of a derivative of order \(m\) is controlled by Sobolev norms through order \(m+2\).  The largest origin derivative in \eqref{eq:modulation-origin-values-k-choice} has order one, so the modulation formulas require \(k\ge3\).  Independently, the reverse allocation in the top-order transport commutators places a factor of derivative order at most two in \(L^\infty\), and the corresponding \(H^2\)-to-\(L^\infty\) estimate again uses derivatives through order four. \textbf{The lowest admissible choice is thus} $\boldsymbol{k=4}$. Larger choices of $k$ are possible, but they increase combinatorially the number of profile derivatives and analytic constants that must be certified. \\

\paragraph{Remark on intermediate derivative orders.} \noindent Although commutator and product estimates generate intermediate derivative orders $0<j<k$, these orders do not need to be included explicitly in the definition of the full energy. Indeed, Appendix~\ref{app:high-order-energy} provides, for each such $j$, the certified interpolation bound \begin{equation} \sum_{|\alpha|=j}\|D^\alpha f\|_{\wgray{\Phi_{i,j}}} \le \eta_{i,j}\sum_{|\alpha|=k}\|D^\alpha f\|_{\wgray{\Phi_{i,k}}} +I_{i,j,k}(\eta_{i,j})\|f\|_{\wgray{\Phi_i}}. \label{eq:intermediate-order-roadmap-main} \end{equation} \noindent Thus, every intermediate-order term is controlled by a combination of the top-order and low-order energy levels. In each product estimate, these two contributions are kept separate until Young's inequality is applied: the top-order term is absorbed by the available $\mathfrak H_k^2$ damping, while the low-order term is controlled by $\mathfrak E_0^2$. The finite family of parameters $\eta_{i,j}$ is chosen so that the total top-order interpolation contribution remains strictly below the available differentiated damping.

\hfill \\

\subsubsection{Constraints on the Weight Functions} \label{sec: Constraints on phi}

\noindent Each weight function \(\wgray{\Phi_i}\), \(i\in\{\omega,r,z\}\), is required to satisfy the following conditions.
\begin{enumerate}[leftmargin=2em]
\vspace{3mm}

\item \textbf{Positivity and coercivity.}
\noindent There is a constant \(C_{\rm floor}>0\) such that
\begin{equation}
\wgray{\Phi_i}(r,z)\ge C_{\rm floor}
\qquad\text{for almost every }(r,z)\in\mathbb D.
\label{eq:weight_coercivity}
\end{equation}
\noindent This prevents the weighted norm from becoming weak in a region where the perturbation must still be controlled.

\vspace{5mm}

\item \textbf{Near-origin admissibility.}
\noindent Each weight function $\wgray{\Phi_i}$ may have its own singularity rate at the origin. We require only local integrability, without imposing a common exponent or a comparison between different components. For certification, one convenient sufficient condition is
\begin{equation}
\Phi_i(r,z)\le C_i(r^2+z^2)^{-p_i},
\qquad
0<r^2+z^2\le \rho_i,
\qquad
0\le p_i<1,
\label{eq:origin-weight-power-model-revision}
\end{equation}
\noindent where $C_i$, $\rho_i$, and $p_i$ are chosen separately for each component. The condition $p_i<1$ guarantees local integrability in two dimensions. This power-law model is only a sufficient criterion, and a direct local-integrability certificate may be used instead.

\vspace{5mm}

\item \textbf{Regularity.}
\noindent Away from the origin, $C^1$ regularity is enough for the basic weighted integrations by parts. When a later estimate differentiates $\wgray{\Phi_i}$ or $\Phi_{i,k}$ further, the required derivatives are assumed to exist on that region and to satisfy the corresponding certified bounds. In the computer-assisted implementation, we use smooth parametrizations with $C^\infty$ activation functions, which meet these requirements and are convenient for gradient-based optimization.

\vspace{5mm}

\item \textbf{Far-field growth needed for velocity control.}
\noindent There are constants \(\rho_{\mathsf{growth}},\mathsf c_{\mathsf{growth}}>0\) such that
\begin{equation}
\min\{\Phi_r,\Phi_z\}\ge \mathsf c_{\mathsf{growth}} \, \rho
\qquad\text{whenever }\rho\ge\rho_{\mathsf{growth}},
\qquad
\rho=r^2+z^2.
\label{eq:delta-u-tail-growth-main}
\end{equation}
\noindent This lower growth condition is used to recover the unweighted \(L^2\) control of \(\bvar{\delta u}\) needed in the global pointwise estimate. If a particular transport, matrix, or elliptic estimate needs a stronger far-field bound, we state and certify that additional bound with the estimate where it is used rather than imposing it as a global relation among the weights.

\vspace{5mm}

\item (Optional) \textbf{Normalization.}
\noindent When the weight functions are optimized via gradient-based optimization, it can be useful to fix one normalization for each weight to remove the irrelevant scaling freedom.

\vspace{5mm}

\item (Optional) \textbf{Evenness in \(r\).}
\noindent Since the pointwise estimates use reflection across the axis, it is convenient to choose the weights with the same even symmetry in $r$:
\begin{equation}
\wgray{\Phi_i}(r,z)=\wgray{\Phi_i}(-r,z).
\end{equation}
\noindent This symmetry can be built directly into the weight parametrization. This requirement concerns only the weights. Each perturbation component retains its own prescribed even or odd parity independently of the fact that its weight is even in $r$.
\end{enumerate}

\hfill

\subsection{Linear Stability Estimates} \label{sec: meaning linear stability}

\subsubsection{Low-Order Energy Linear Stability and the Need for Higher-Order Energy} 

Linear stability asks whether perturbations of solutions of the closed linearized system remain bounded or decay in a suitable norm. In terms of $\textcolor{blue}{\delta v}$, the desired estimate is
\begin{equation}
\label{eq:semigroup_decay_intro_deltav}
\|\textcolor{blue}{\delta v}(t)\|_{\wgray{\Phi}} \leq C e^{-\gamma t} \|\textcolor{blue}{\delta v}(0)\|_{\wgray{\Phi}}, \qquad \gamma>0.
\end{equation}

\hfill 

\noindent A standard rigorous way to prove \eqref{eq:semigroup_decay_intro_deltav} is to establish the differential energy estimate
\begin{equation}
\label{eq:energy_decay_intro_E0}
\frac{1}{2}\frac{d}{dt}\mathfrak{E}_0(t)^2
=
\left\langle
\mathcal L_{\nabla}\textcolor{blue}{\delta v},
\textcolor{blue}{\delta v}
\right\rangle_{\wgray{\Phi}}
\leq
-\DLamL \mathfrak{E}_0(t)^2.
\end{equation}

\noindent This estimate gives a quantitative estimate for every solution of the closed linearized perturbation equation in the specified weighted space. Integrating \eqref{eq:energy_decay_intro_E0} gives exponential decay of $\mathfrak{E}_0(t)$, and hence of the scale-consistent perturbation $\textcolor{blue}{\delta v}$, rather than merely a numerical or visual indication of convergence. \\ 

\noindent For a blowup profile, this is the relevant notion of stability: the original physical solution may still blow up, but the normalized spatial profile is stable in the rescaled variables. Linear stability means that $\textcolor{blue}{\delta u}(t) \rightarrow0$ and $\textcolor{blue}{\delta\omega}(t) \rightarrow 0$ under the linearized rescaled flow. Note that the blowup time, blowup location, and scaling parameters may change under perturbation, but these changes are absorbed by the modulation parameters. \\

\noindent In spectral language, the same idea is that the relevant spectrum of the closed linearized operator, denoted above by $\mathcal L_{\nabla}$, should lie in the stable half-plane. In a numerical stability calculation, one discretizes $\mathcal L_{\nabla}$ and studies the eigenvalues of the resulting matrix. Eigenvalues with negative real part correspond to decaying linear modes, while eigenvalues with positive real part correspond to growing directions. Neutral modes may arise from translation or amplitude normalization. The modulation conditions are designed to fix these directions, so that the stability question concerns perturbations transverse to the artificial symmetry directions. Thus, the role of the linearization is to produce the operator $\mathcal L_{\nabla}$ governing the first-order evolution of the scale-consistent perturbation $\textcolor{blue}{\delta v}$. The role of the linear stability analysis is to show that this operator has no growing directions in the weighted energy $\mathfrak{E}_0$, once the elliptic constraint, the velocity reconstruction, and the modulation conditions are imposed. If such an estimate is proved, then the proposed rescaled profile is linearly stable as a blowup profile. \\

\noindent However, as noted earlier, the low-order energy is not sufficient by itself to close the linear estimates since the linear operators contain terms which cannot be controlled using a weighted $L^2$ norm. It is thus supplemented with the high-order weighted energy $\mathfrak H_k$, to get an estimate of the form
\begin{equation}  
\frac{1}{2}\frac{d}{dt}\mathfrak{E}_0(t)^2
=
\left\langle
\mathcal L_{\nabla}\textcolor{blue}{\delta v},
\textcolor{blue}{\delta v}
\right\rangle_{\wgray{\Phi}}
\leq -\DLamLlow\mathfrak E_0^2+\DLamLhigh\mathfrak H_k^2.
\end{equation}

\hfill \\

\subsubsection{Low-Order Linear Estimates}

\noindent The low-order linear estimate to certify is
\begin{equation}  \ \ 
\label{eq:ns-CL}
\langle \mathcal L_\omega,\bvar{\delta\omega}\rangle_{\wgray{\Phi_\omega}}
+\langle \partial_r\mathcal L_u,\partial_r\bvar{\delta u}\rangle_{\wgray{\Phi_r}}
+\langle \partial_z\mathcal L_u,\partial_z\bvar{\delta u}\rangle_{\wgray{\Phi_z}}
\le -\DLamLlow \, \mathfrak E_0^2+\DLamLhigh \, \mathfrak H_k^2. \ \ 
\end{equation}
Here the constant \(\DLamLlow>0\) is the net low-order linear damping and \(\DLamLhigh\ge0\) is the top-order loss from direct estimates that require \(\mathfrak H_k\).

\hfill \\

\paragraph{High-level steps to obtain the estimate:}
\begin{enumerate}[leftmargin=.35in, itemsep=0.8em]
\item We split the low-order linearized equations into the transport terms, the zeroth-order couplings kept in the signed matrix, and the remaining terms that are estimated directly. After weighted integration by parts, the retained signed matrix terms give
\begin{equation}
\int_{\mathbb D}\textcolor{black}{d\xi}^\top\mathbb M_L\textcolor{black}{d\xi}\,dr\,dz,
\qquad
\textcolor{black}{d\xi}:=\bigl(\Phi_\omega^{1/2}\bvar{\delta\omega},\Phi_r^{1/2}\partial_r\bvar{\delta u},\Phi_z^{1/2}\partial_z\bvar{\delta u}\bigr)^\top.
\end{equation}

\item We certify the matrix damping on the computational box and close the estimate analytically on the tail. If a finite collection of bounded regions does not satisfy the uniform pointwise bound, we use the localized certificate in Appendix~\ref{app:localized-matrix-certificate} to combine the damping away from those regions with certified concentration bounds inside them.

\item We estimate the terms outside the matrix according to their structure. The coupling $2\bar u\,\partial_z\bvar{\delta u}$ uses the pointwise bound for $\partial_z\bvar{\delta u}$. The fixed-profile terms containing undifferentiated $\bvar{\delta u}$ use the ray estimate. Streamfunction terms use the weighted elliptic estimates, and modulation terms use the certified origin-value bounds and the explicit modulation formulas.

\item We collect the direct losses in $A_L^0$ and $A_L^k$, and then apply the Cauchy--Schwarz inequality and Young's inequality. For the global certificate,
\begin{equation}
\DLamLlow=\Lambda_L^{\rm mat}-A_L^0,
\qquad
\DLamLhigh=A_L^k.
\end{equation}
If the localized matrix certificate is used, we replace the global matrix margin by its effective low-order damping and include the associated top-order concentration cost. The estimate closes when the resulting certified constant satisfies $\DLamLlow>0$.

\end{enumerate}

\hfill

\noindent Appendix~\ref{app:low-order-linear} gives the detailed low-order linear estimate and the numerical construction of the required constants.

\hfill \\

\subsubsection{High-Order Linear Estimates}

\noindent We need to certify a top-order damping constant $\DLamDLhigh>0$ and a low-order constant $\DLamDLlow\ge0$ such that
\begin{equation} \ \ 
\label{eq:ns-high-lin}
\begin{aligned}
&\sum_{|\alpha|=k}\langle D^\alpha\mathcal L_\omega,D^\alpha\bvar{\delta\omega}\rangle_{\wgray{\Phi_{\omega,k}}}
+\sum_{|\alpha|=k}\langle D^\alpha\partial_r\mathcal L_u,D^\alpha\partial_r\bvar{\delta u}\rangle_{\wgray{\Phi_{r,k}}}
+\sum_{|\alpha|=k}\langle D^\alpha\partial_z\mathcal L_u,D^\alpha\partial_z\bvar{\delta u}\rangle_{\wgray{\Phi_{z,k}}}
 \\
&\qquad\le \DLamDLlow \, \mathfrak E_0^2 -\DLamDLhigh \, \mathfrak H_k^2.
\end{aligned} \ \ 
\end{equation}

\hfill \\

{
\paragraph{High-level steps to obtain the estimate:}
\begin{enumerate}[leftmargin=.35in, itemsep=0.8em]
\item We apply $D^\alpha$ for all $|\alpha|=k$ to the linearized equations and test against the top-order derivatives of the corresponding perturbations. For the transport terms, we integrate by parts the term in which all $k$ derivatives remain on the perturbation. The commutators in which exactly one derivative hits $G_r$ or $G_z$ still contain $k$ derivatives of the perturbation. Together with the zeroth-order cross-component terms collected in the next step, these contributions form a signed higher-order damping matrix analogous to the low-order damping matrix. Its global damping margin controls the top-order energy, and the localized fallback uses this higher-order coercivity in regions where the low-order matrix does not have a uniform damping margin, in particular near stagnation regions.

\item We add the zeroth-order cross-component terms that also keep all $k$ derivatives on the perturbation. We collect the $3(k+1)$ weighted top derivatives into the symmetric matrix, keeping every cross-component coefficient together with its complete weight conversion. The finite-box enclosure and analytic tail estimate give the principal top-order damping margin.

\item We estimate all remaining differentiated terms by derivative count. If $|\beta|\ge2$ derivatives hit a transport coefficient, the perturbation factor has order $k+1-|\beta|$. We measure that factor in its derivative-level weight and use the certified interpolation estimate when the order lies strictly between $0$ and $k$. We treat the other differentiated coefficients in the same way, always keeping the complete source-to-target weight conversion with the coefficient. Streamfunction terms use the elliptic estimates, and modulation terms use the origin-value bounds. We keep low-order and top-order losses separate until Young's inequality is applied.

\item We subtract the certified top-order losses from the matrix damping margin and collect the lower-order losses. This yields $\DLamDLhigh>0$ and $\DLamDLlow\ge0$.
\end{enumerate}
}

\hfill

\noindent The detailed estimate is proved in Appendix~\ref{app:high-order-linear}.

\hfill \\ 

\subsection{Nonlinear Stability Estimates}

The linear stability analysis  identifies the estimates required to establish stability of the rescaled profile  under the first-order perturbation dynamics, measured in the weighted low-order norm associated with $\wgray{\Phi_\omega}$, $\wgray{\Phi_r}$, and $\wgray{\Phi_z}$. The nonlinear stability argument  seeks to upgrade this  control to the full perturbation equation. Its goal is to prove that, if the initial perturbation is sufficiently small in a suitable energy, then the perturbation remains small for as long as the estimates apply. In the case of an approximate numerical profile, the solution is allowed to remain in a small neighborhood whose size also depends on the certified residual error. \\ 

\noindent The nonlinear stability  argument is reduced to an energy estimate for the full energy $\mathfrak E_k$. Differentiating $\mathfrak E_k^2$ in time and inserting the perturbation equations separates the argument into five distinct components: the residual error, the low-order linear part, the high-order linear part, the low-order nonlinear part, and the high-order nonlinear part. The residual records the defect of the approximate steady profile. The linear terms contain the stabilizing mechanism: the low-order estimate gives the principal damping, while the high-order estimate controls the differentiated equations, allowing only losses that can be absorbed by the low-order energy after choosing $\mu_k$. The nonlinear terms are then estimated perturbatively, with the low-order and high-order contributions bounded by cubic and quartic powers of the full energy.

\hfill \\

\subsubsection{Low-Order Nonlinear Estimates}

We start from the expressions \eqref{eq:Nw_Euler_definition}, \eqref{eq:Nr_Euler}, \eqref{eq:Nz_Euler}, for the nonlinear parts $\mathcal{N}_i$ for $i\in \{\omega,r,z\}$ . \\

\noindent We substitute the exact modulation formulas before separating linear and nonlinear terms:
\begin{equation}
\bvar{\delta C}=-2\bvar{\delta\psi}_0,
\qquad\quad
\bvar{\delta c_u}
=-2(\partial_z\bvar{\delta\psi})_0.
\end{equation}

\hfill

\noindent The precise low-order nonlinear estimate to prove is
\begin{equation}  \ \ 
\label{eq:ns-low-nonlinear-final}
\left|\langle\mathcal N_\omega,\bvar{\delta\omega}\rangle_{\wgray{\Phi_\omega}}
+\langle\mathcal N_{r},\partial_r\bvar{\delta u}\rangle_{\wgray{\Phi_r}}
+\langle\mathcal N_{z},\partial_z\bvar{\delta u}\rangle_{\wgray{\Phi_z}}
\right|
\le \DLamNthree \, \mathfrak E_k^3. 
\end{equation}

\hfill

{
\paragraph{High-level steps to obtain the estimate:}
\begin{enumerate}[leftmargin=.35in, itemsep=0.8em]
\item We separate transport pairings from the remaining nonlinear products. For transport terms, we use the weighted transport identity and certify the complete multiplier, including the weight derivative. In the axial direction, we keep $\bvar{\delta C}+\bvar{\delta u^z}$ as a single centered coefficient. The normalization makes this coefficient vanish at the origin, so we use that cancellation before estimating its product with a singular weight derivative.

\item For each non-transport product, we keep the two factors from the energy pairing in their weighted $L^2$ norms and place the remaining coefficient or lower-order perturbative factor in $L^\infty$. If the two weighted factors carry different component weights, we keep the complete source-to-target weight conversion with the pointwise multiplier. We first use the velocity or streamfunction estimates to control the $L^\infty$ factor by the energy and then take the final product bound.

\item We substitute the exact modulation formulas before counting powers of the energy.The modulation coefficients are linear, so the resulting nonlinear energy pairings are cubic. 

\item We sum the finitely many certified termwise constants to obtain $\DLamNthree$.
\end{enumerate}
}

\hfill

\noindent The termwise proof and the definition of $\DLamNthree$ are given in Appendix~\ref{app:low-order-nonlinear}.

\hfill \\

\subsubsection{High-Order Nonlinear Estimates}

\hfill 

\noindent The precise high-order nonlinear estimate to prove is

\begin{equation}
\label{eq:ns-high-nonlinear}
\mu_k\left|
\sum_{i\in\{\omega,r,z\}}\sum_{|\alpha|=k}
\left\langle D^\alpha\mathcal N_i,D^\alpha\bvar{\delta v_i}\right\rangle_{\wgray{\Phi_{i,k}}}
\right|
\le \DLamDNthree \, \mathfrak E_k^3.
\end{equation}

\hfill \\

{
\paragraph{High-level steps to obtain the estimate:}
\begin{enumerate}[leftmargin=.35in, itemsep=0.8em]
\item We distribute $D^\alpha$ for all $|\alpha|=k$, over each nonlinear transport term. The principal terms in which all $k$ derivatives remain on the transported perturbation are integrated by parts with the top-order weight. In every commutator, we compare the derivative orders of the two perturbative factors, place the lower-order factor in $L^\infty$, and keep the higher-order factor in weighted $L^2$. If the weighted factor has order $j<k$, we measure it with its order-$j$ weight and keep the full conversion to $\wgray{\Phi_{i,k}}$  in the same pointwise multiplier. For axial transport, we keep $\bvar{\delta C}+\bvar{\delta u^z}$ centered so that its vanishing at the origin is used before a singular top-order weight derivative is estimated.

\item We treat each differentiated source product by the same rule. We keep a perturbative factor of maximal derivative order in weighted $L^2$ and place the lower-order perturbative factor in $L^\infty$. The pointwise multiplier retains the complete source-to-target weight conversion. We first reduce streamfunction factors with the elliptic estimates, then use the certified interpolation inequalities for any intermediate derivative order, keeping the resulting low-order and top-order contributions separate.

\item We substitute the modulation formulas before assigning nonlinear order. The modulation coefficients are linear and all resulting nonlinear energy pairings are cubic.

\item We sum the transport commutators, source products, and modulation terms to obtain the constants in the high-order nonlinear estimate. 
\end{enumerate}
}

\hfill

\noindent The detailed proof is given in Appendix~\ref{app:high-order-nonlinear}. 

\hfill \\

\subsection{PDE Residuals}

If the $(\bar u,\bar\omega,\bar\psi,\bar\lambda,\bar C,\bar c_u)$ profile is exact, the raw residuals below vanish. If the profile is numerical, these are the defects obtained by substituting the profile into the steady equations:
\begin{align}
\label{eq:ns-Ru-raw}
\mathcal R_u^{\mathsf{raw}}
&:=-(\lambda r+\bar u^r)\partial_r\bar u
-(\lambda z+\bar C+\bar u^z)\partial_z\bar u
+2\bar u\partial_z\bar\psi+\bar c_u\bar u ,\\
\label{eq:ns-Romega-raw}
\mathcal R_\omega^{\mathsf{raw}}
&:=-(\lambda r+\bar u^r)\partial_r\bar\omega
-(\lambda z+\bar C+\bar u^z)\partial_z\bar\omega
+2\bar u\partial_z\bar u+(\bar c_u-\lambda)\bar\omega .
\end{align}

\hfill

\noindent The residuals entering the perturbation equations are
\begin{align}
\label{eq:ns-Ru}
\mathcal R_u
&:=\mathcal R_u^{\mathsf{raw}}
+\cures\,\bar u,\\
\label{eq:ns-Romega}
\mathcal R_\omega
&:=\mathcal R_\omega^{\mathsf{raw}}
+\cures\,\bar\omega.
\end{align}

\hfill

\noindent The elliptic residual is not listed here because the spline representation of 
\(\bar\omega\) is obtained by applying the elliptic operator to \(\bar\psi\) exactly, as discussed in \Cref{sec: numerical profiles}. As a result, the steady-state profile satisfies the elliptic equation exactly, and the elliptic residual is exactly 0.

\hfill \\

\noindent Define
\begin{equation}
\label{eq:ns-Rrz}
\mathcal R_r:=\partial_r\mathcal R_u,
\qquad
\mathcal R_z:=\partial_z\mathcal R_u.
\end{equation}

\hfill 

\noindent Then the low-order perturbation equations are
\begin{align}
\label{eq:ns-component-eqns}
\partial_t\bvar{\delta\omega}
&=\mathcal L_\omega(\bvar{\delta})+
\mathcal N_\omega(\bvar{\delta})+\mathcal R_\omega,\\
\partial_t\partial_r\bvar{\delta u}
&=\partial_r\mathcal L_u(\bvar{\delta})+
\mathcal N_{r}(\bvar{\delta})+\mathcal R_r,\\
\partial_t\partial_z\bvar{\delta u}
&=\partial_z\mathcal L_u(\bvar{\delta})+
\mathcal N_{z}(\bvar{\delta})+\mathcal R_z.
\end{align}

\hfill 

\noindent Define the residual constants by separating the low-order and top-order pieces:
\begin{align}
\label{eq:ns-eps0}
\DLamRPhi^2
&:=\|\mathcal R_\omega\|_{\wgray{\Phi_\omega}}^2
+\|\mathcal R_r\|_{\wgray{\Phi_r}}^2
+\|\mathcal R_z\|_{\wgray{\Phi_z}}^2,\\
\label{eq:ns-epsk}
\DLamDR^2
&:=\sum_{|\alpha|=k}
\Bigl(
\|D^\alpha\mathcal R_\omega\|_{\wgray{\Phi_{\omega,k}}}^2
+\|D^\alpha\mathcal R_r\|_{\wgray{\Phi_{r,k}}}^2
+\|D^\alpha\mathcal R_z\|_{\wgray{\Phi_{z,k}}}^2
\Bigr),\\
\label{eq:ns-eps-total}
\DLamR^2
&:=\DLamRPhi^2+\mu_k\DLamDR^2.
\end{align}

\hfill 

\noindent Applying the Cauchy--Schwarz inequality gives
\begin{equation} 
\label{eq:ns-res-pair}
\left| \,
\sum_{i\in\{\omega,r,z\}}
\langle\mathcal R_i,\bvar{\delta v_i}\rangle_{\wgray{\Phi_i}}
+\mu_k\sum_{i\in\{\omega,r,z\}}\sum_{|\alpha|=k}
\langle D^\alpha\mathcal R_i,D^\alpha\bvar{\delta v_i}\rangle_{\wgray{\Phi_{i,k}}} \,
\right|
\le \DLamR \, \mathfrak E_k. 
\end{equation}

\hfill \\

\subsection{Stability Theorems}

\hfill

\noindent We are now ready to state the stability results for the Euler equations. Importantly, every named certification constant attached to an analytic estimate is understood to be a rigorously validated finite numerical value. Conceptually, the analytic argument identifies the constants that are needed, while the numerical step certifies them once for the chosen profiles, weights, and parameters. The stability proofs are not complete until appropriate constants have been certified numerically. \\

\subsubsection{Single-Radius Stability Theorems}

\noindent We now state the full stability theorem.\\

\begin{theorem}[Single-radius stability]\label{thm:main-stability}
Assume that the following perturbation equations hold for admissible perturbations:
\begin{align}
\partial_t\bvar{\delta\omega}
&=\mathcal L_\omega(\bvar{\delta})+
\mathcal N_\omega(\bvar{\delta})+\mathcal R_\omega, \label{eq: Stability theorem perturbation 1}\\
\partial_t\partial_r\bvar{\delta u}
&=\partial_r\mathcal L_u(\bvar{\delta})+
\mathcal N_{r}(\bvar{\delta})+\mathcal R_r,\label{eq: Stability theorem perturbation 2}\\
\partial_t\partial_z\bvar{\delta u}
&=\partial_z\mathcal L_u(\bvar{\delta})+
\mathcal N_{z}(\bvar{\delta})+\mathcal R_z. \label{eq: Stability theorem perturbation 3}
\end{align}
Assume that the following estimates have been certified with finite constants:

\begin{equation}
\begin{aligned}
&\langle \mathcal L_\omega,\bvar{\delta\omega}\rangle_{\wgray{\Phi_\omega}}
+\langle \partial_r\mathcal L_u,\partial_r\bvar{\delta u}\rangle_{\wgray{\Phi_r}}
+\langle \partial_z\mathcal L_u,\partial_z\bvar{\delta u}\rangle_{\wgray{\Phi_z}}
\le -\DLamLlow \, \mathfrak E_0^2+\DLamLhigh \, \mathfrak H_k^2.
\end{aligned}
\label{eq:ns-thm-low-linear}
\end{equation}

\vspace{-1mm}

\begin{equation}
\begin{aligned}
&\sum_{i\in\{\omega,r,z\}}\sum_{|\alpha|=k}
\langle D^\alpha\mathcal L_i,D^\alpha\bvar{\delta v_i}\rangle_{\wgray{\Phi_{i,k}}}
\le \DLamDLlow \, \mathfrak E_0^2-\DLamDLhigh \, \mathfrak H_k^2.
\end{aligned}
\label{eq:ns-thm-high-linear}
\end{equation}

\vspace{-1mm}

\begin{equation}
\begin{aligned}
&\left|\langle\mathcal N_\omega,\bvar{\delta\omega}\rangle_{\wgray{\Phi_\omega}}
+\langle\mathcal N_{r},\partial_r\bvar{\delta u}\rangle_{\wgray{\Phi_r}}
+\langle\mathcal N_{z},\partial_z\bvar{\delta u}\rangle_{\wgray{\Phi_z}}
\right|
\le \DLamNthree \, \mathfrak E_k^3.
\end{aligned}
\label{eq:ns-thm-low-nonlinear}
\end{equation}

\vspace{-1mm}

\begin{equation}
\mu_k \, \Bigg|
\sum_{i\in\{\omega,r,z\}}\sum_{|\alpha|=k}
\langle D^\alpha\mathcal N_i,D^\alpha\bvar{\delta v_i}\rangle_{\wgray{\Phi_{i,k}}}
\Bigg|
\le \DLamDNthree \, \mathfrak E_k^3.
\label{eq:ns-thm-high-nonlinear}
\end{equation}

\vspace{-1mm}

\begin{equation}
\begin{aligned}
&\Bigg|
\sum_{i\in\{\omega,r,z\}}\langle \mathcal R_i,\bvar{\delta v_i}\rangle_{\wgray{\Phi_i}}
+\mu_k\sum_{i\in\{\omega,r,z\}}\sum_{|\alpha|=k}
\langle D^\alpha\mathcal R_i,D^\alpha\bvar{\delta v_i}\rangle_{\wgray{\Phi_{i,k}}}
\Bigg|
\le \DLamR \, \mathfrak E_k.
\end{aligned}
\label{eq:ns-thm-residual}
\end{equation}

\hfill 

\vspace{3mm}

\noindent Choose $\mu_k$ so that $0<\mu_k\le 1$ and
\begin{equation}
\label{eq:ns-mu-choice}
\DLamStab(\mu_k):=
\min\left\{
\DLamLlow-\mu_k\DLamDLlow,
\ \DLamDLhigh-\frac{\DLamLhigh}{\mu_k}
\right\}>0.
\end{equation}
After this choice of $\mu_k$ is fixed, write $\DLamStab$ for the positive number $\DLamStab(\mu_k)$.  Define the combined nonlinear constants by
\begin{equation}
\label{eq:ns-N2N3}
\DLamThree:=\DLamNthree+\DLamDNthree.
\end{equation}

\noindent Assume that $\delta_*>0$ satisfies

\begin{equation}
\label{eq:ns-smallness}
\DLamThree\delta_*+\frac{\DLamR}{\delta_*}<\DLamStab.
\end{equation}

\noindent Then every solution with $\mathfrak E_k(0)<\delta_*$ remains in the ball $\mathfrak E_k(t)<\delta_*$ for as long as the solution exists and the certified estimates apply.
\end{theorem}

\begin{grayproof}
The proof is presented in Appendix~\ref{appx: nonlinear stability proof}. 
\end{grayproof}

\clearpage 

\noindent The diagram below summarizes the strategy used for proving nonlinear stability in the Euler equations.

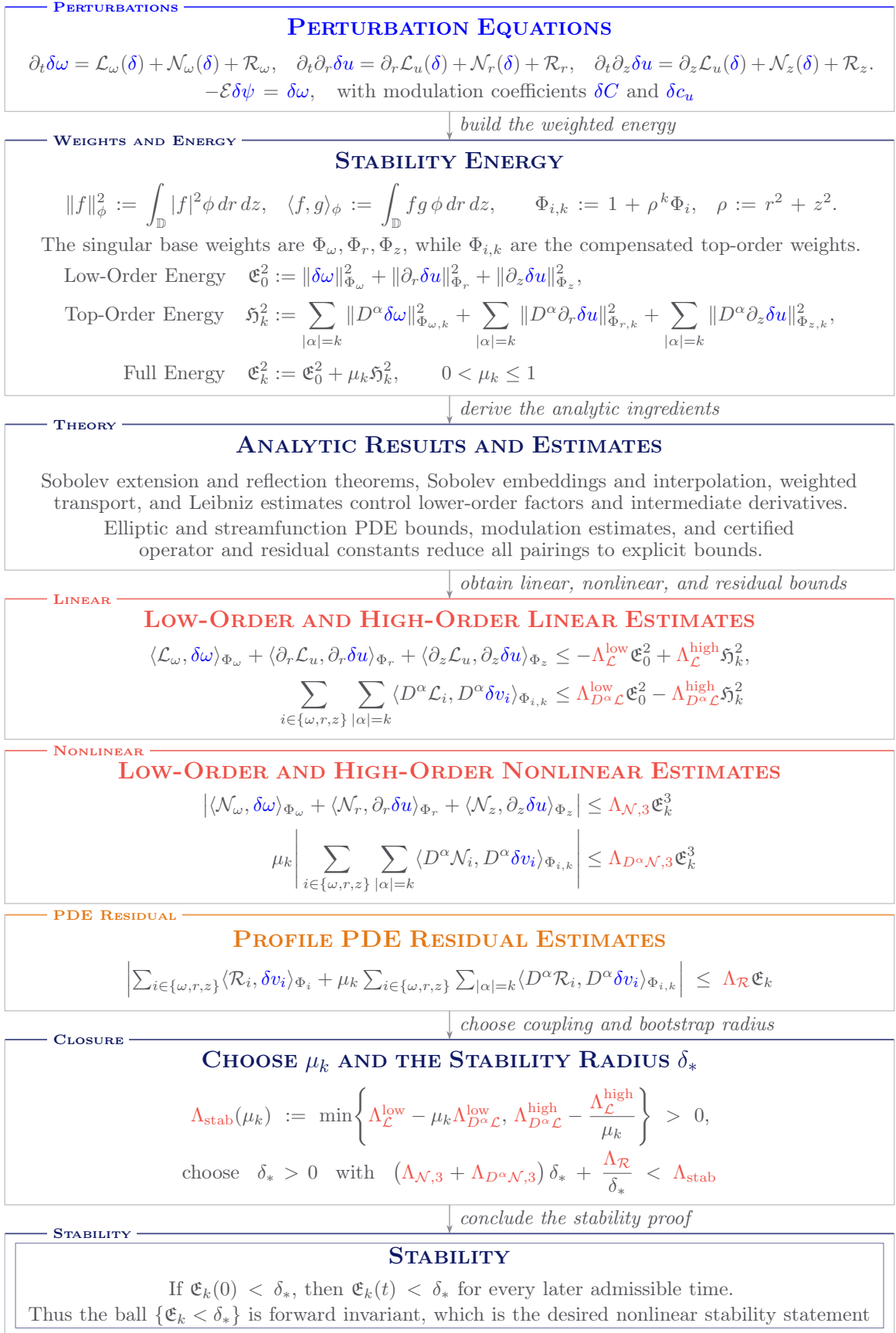
\begin{figure}[h]
\centering
\begingroup

\definecolor{eulerstaborange}{RGB}{230,120,20}
\definecolor{eulerstabpurple}{RGB}{15,25,100}
\definecolor{eulerstabsalmon}{RGB}{240,80,70}
\definecolor{eulerstabtextmuted}{RGB}{82,82,88}
\definecolor{eulerstabrulegray}{RGB}{140,140,148}
\definecolor{eulerstabrulemid}{RGB}{176,176,183}

\newcommand{\eulerstabLamLlow}{%
  \ensuremath{\textcolor{eulerstabsalmon}
  {\Lambda_{\mathcal L}^{\rm low}}}%
}

\newcommand{\eulerstabLamLhigh}{%
  \ensuremath{\textcolor{eulerstabsalmon}
  {\Lambda_{\mathcal L}^{\rm high}}}%
}

\newcommand{\eulerstabLamDLhigh}{%
  \ensuremath{\textcolor{eulerstabsalmon}
  {\Lambda_{D^\alpha\mathcal L}^{\rm high}}}%
}

\newcommand{\eulerstabLamDLlow}{%
  \ensuremath{\textcolor{eulerstabsalmon}
  {\Lambda_{D^\alpha\mathcal L}^{\rm low}}}%
}

\newcommand{\eulerstabLamNthree}{%
  \ensuremath{\textcolor{eulerstabsalmon}
  {\Lambda_{\mathcal N,3}}}%
}

\newcommand{\eulerstabLamNfour}{%
  \ensuremath{\textcolor{eulerstabsalmon}
  {\Lambda_{\mathcal N,4}}}%
}

\newcommand{\eulerstabLamDNthree}{%
  \ensuremath{\textcolor{eulerstabsalmon}
  {\Lambda_{D^\alpha\mathcal N,3}}}%
}

\newcommand{\eulerstabLamDNfour}{%
  \ensuremath{\textcolor{eulerstabsalmon}
  {\Lambda_{D^\alpha\mathcal N,4}}}%
}

\newcommand{\eulerstabLamR}{%
  \ensuremath{\textcolor{eulerstabsalmon}
  {\Lambda_{\mathcal R}}}%
}

\newcommand{\eulerstabLamStab}{%
  \ensuremath{\textcolor{eulerstabsalmon}
  {\Lambda_{\rm stab}}}%
}

\tikzset{
  eulerstab flow/.style={
    -{Stealth[length=1.55mm,width=1.08mm]},
    draw=eulerstabrulegray,
    line width=0.46pt
  },
  eulerstab transition/.style={
    midway,
    right=3.5pt,
    fill=white,
    inner xsep=1.10pt,
    inner ysep=0.65pt,
    text=eulerstabtextmuted,
    font=\fontsize{9.6}{10.6}\selectfont\itshape
  },
  eulerstab stage/.style={
    draw=eulerstabrulemid,
    fill=white,
    line width=0.40pt,
    minimum width=15.75cm,
    text width=14.90cm,
    minimum height=17.0mm,
    align=center,
    inner xsep=8.0pt,
    inner ysep=5.8pt,
    outer sep=0pt
  }
}

\newcommand{\eulerstabtitle}[2]{%
  {\fontsize{12.4}{13.4}\selectfont
   \bfseries\color{#1}\textsc{#2}\par
   \vspace{6.6pt}}%
}

\newcommand{\eulerstabtitlecompact}[2]{%
  {\fontsize{12.4}{13.4}\selectfont
   \bfseries\color{#1}\textsc{#2}\par
   \vspace{5.0pt}}%
}

\newcommand{\eulerstabsubtitle}[1]{%
  {\fontsize{9.8}{11.0}\selectfont
   \color{eulerstabtextmuted}#1\par
   \vspace{-0.8pt}}%
}

\newcommand{\eulerstabphasetag}[2]{%
  {\fontsize{7.4}{8.2}\selectfont
   \bfseries\scshape\color{#1}#2}%
}

\newcommand{\eulerstabtagpanel}[6]{%
  \node[eulerstab stage,#6] (#1) {%
    \eulerstabtitle{#2}{#4}%
    \eulerstabsubtitle{#5}%
  };
  \draw[#2!78,line width=0.70pt]
    ([xshift=0.2pt]#1.north west) --
    ([xshift=-0.2pt]#1.north east);
  \node[
    anchor=west,
    fill=white,
    inner xsep=2.7pt,
    inner ysep=1.0pt,
    text=#2
  ] at ([xshift=7.5mm]#1.north west) {%
    \eulerstabphasetag{#2}{#3}%
  };
}

\newcommand{\eulerstabtagpanelcompact}[6]{%
  \node[eulerstab stage,#6] (#1) {%
    \eulerstabtitlecompact{#2}{#4}%
    \eulerstabsubtitle{#5}%
  };
  \draw[#2!78,line width=0.70pt]
    ([xshift=0.2pt]#1.north west) --
    ([xshift=-0.2pt]#1.north east);
  \node[
    anchor=west,
    fill=white,
    inner xsep=2.7pt,
    inner ysep=1.0pt,
    text=#2
  ] at ([xshift=7.5mm]#1.north west) {%
    \eulerstabphasetag{#2}{#3}%
  };
}

\hyphenpenalty=10000
\exhyphenpenalty=10000
\emergencystretch=1em

\resizebox{0.87\linewidth}{!}{%
\begin{tikzpicture}[node distance=5.1mm,font=\rmfamily]

  \eulerstabtagpanel
    {eulersetup}
    {blue}
    {Perturbations}
    {Perturbation Equations}
    {$\partial_t\bvar{\delta\omega}
      =\mathcal L_\omega(\bvar{\delta})
      +\mathcal N_\omega(\bvar{\delta})
      +\mathcal R_\omega$,\quad
     $\partial_t\partial_r\bvar{\delta u}
      =\partial_r\mathcal L_u(\bvar{\delta})
      +\mathcal N_r(\bvar{\delta})
      +\mathcal R_r$,\quad
     $\partial_t\partial_z\bvar{\delta u}
      =\partial_z\mathcal L_u(\bvar{\delta})
      +\mathcal N_z(\bvar{\delta})
      +\mathcal R_z$.\\[3pt]
     $-\mathcal E\bvar{\delta\psi} \!\!\!\!\!\! = \!\!\!\!\!\! \bvar{\delta\omega}$, \ \ 
     with modulation coefficients
     $\bvar{\delta C}$ and $\bvar{\delta c_u}$}
    {}

  \eulerstabtagpanel
    {weightedenergy}
    {eulerstabpurple}
    {Weights and Energy}
    {Stability Energy}
    {$\|f\|_{\phi}^{2} \!
      := \! \displaystyle\int_{\mathbb{D}}|f|^{2}\phi\,dr\,dz$,\quad  
     $\langle f,g\rangle_{\phi} \! 
      := \! \displaystyle\int_{\mathbb{D}}fg\,\phi\,dr\,dz$,\qquad 
     $\wgray{\Phi_{i,k}} \! 
      := \! 1+\rho^{\,k}\wgray{\Phi_i}$,\quad 
     $\rho \! := \! r^2+z^2$.\\[4pt]
     The singular base weights are
     $\wgray{\Phi_\omega},\wgray{\Phi_r},\wgray{\Phi_z}$,
     while $\wgray{\Phi_{i,k}}$ are the compensated top-order weights.\\[3pt]
     $\begin{aligned} \text{Low-Order Energy } \ \ 
      \mathfrak E_0^2
      &:=\|\bvar{\delta\omega}\|_{\wgray{\Phi_\omega}}^2
      +\|\partial_r\bvar{\delta u}\|_{\wgray{\Phi_r}}^2
      +\|\partial_z\bvar{\delta u}\|_{\wgray{\Phi_z}}^2,
      \\[1.8pt]
     \text{Top-Order Energy } \ \  \mathfrak H_k^2
      &:=\sum_{|\alpha|=k}
      \|D^\alpha\bvar{\delta\omega}\|_{\wgray{\Phi_{\omega,k}}}^2
      +\sum_{|\alpha|=k}
      \|D^\alpha\partial_r\bvar{\delta u}\|_{\wgray{\Phi_{r,k}}}^2
      +\sum_{|\alpha|=k}
      \|D^\alpha\partial_z\bvar{\delta u}\|_{\wgray{\Phi_{z,k}}}^2,
      \\[1.8pt]
      \text{Full Energy } \ \  \mathfrak E_k^2
      &:=\mathfrak E_0^2+\mu_k\mathfrak H_k^2,
      \qquad 0<\mu_k\le1
     \end{aligned}$}
    {below=of eulersetup,
     text width=15.35cm,
     inner xsep=4.5pt}

  \eulerstabtagpanel
    {eulertoolbox}
    {eulerstabpurple}
    {Theory}
    {Analytic Results and Estimates}
    {Sobolev extension and reflection theorems,
     Sobolev embeddings and interpolation,
     weighted transport, and Leibniz estimates control
     lower-order factors and intermediate derivatives.\\[1.8pt]
     Elliptic and streamfunction PDE bounds,
     modulation estimates, and certified operator and residual
     constants reduce all pairings to explicit bounds. }
    {below=of weightedenergy}

  \eulerstabtagpanel
    {eulerlinear}
    {eulerstabsalmon}
    {Linear}
    {Low-Order and High-Order Linear Estimates}
    {$\begin{aligned}
      \langle\mathcal L_\omega,\bvar{\delta\omega}\rangle_
        {\wgray{\Phi_\omega}}
      +\langle\partial_r\mathcal L_u,
        \partial_r\bvar{\delta u}\rangle_{\wgray{\Phi_r}}
      +\langle\partial_z\mathcal L_u,
        \partial_z\bvar{\delta u}\rangle_{\wgray{\Phi_z}}
      &\le
      -\eulerstabLamLlow\mathfrak E_0^2
      +\eulerstabLamLhigh\mathfrak H_k^2,
      \\[1.8pt]
      \sum_{i\in\{\omega,r,z\}}\sum_{|\alpha|=k}
      \langle D^\alpha\mathcal L_i,
        D^\alpha\bvar{\delta v_i}\rangle_{\wgray{\Phi_{i,k}}}
      &\le
      \eulerstabLamDLlow\mathfrak E_0^2
      -\eulerstabLamDLhigh\mathfrak H_k^2
    \end{aligned}$}
    {below=of eulertoolbox}

  \eulerstabtagpanel
    {eulernonlinear}
    {eulerstabsalmon}
    {Nonlinear}
    {Low-Order and High-Order Nonlinear Estimates}
    {$\begin{aligned}
      \bigl|
      \langle\mathcal N_\omega,\bvar{\delta\omega}\rangle_
        {\wgray{\Phi_\omega}}
      +\langle\mathcal N_r,\partial_r\bvar{\delta u}\rangle_
        {\wgray{\Phi_r}}
      +\langle\mathcal N_z,\partial_z\bvar{\delta u}\rangle_
        {\wgray{\Phi_z}}
      \bigr|
      &\le
      \eulerstabLamNthree\mathfrak E_k^3
      \\[1.8pt]
      \mu_k\Bigg|
      \sum_{i\in\{\omega,r,z\}}\sum_{|\alpha|=k}
      \langle D^\alpha\mathcal N_i,
        D^\alpha\bvar{\delta v_i}\rangle_{\wgray{\Phi_{i,k}}}
      \Bigg|
      &\le
      \eulerstabLamDNthree\mathfrak E_k^3
    \end{aligned}$}
    {below=2.0mm of eulerlinear}

  \eulerstabtagpanel
    {eulerresidual}
    {eulerstaborange}
    {PDE Residual}
    {Profile PDE Residual Estimates}
    {$\left|
      \sum_{i\in\{\omega,r,z\}}
      \langle\mathcal R_i,\bvar{\delta v_i}\rangle_{\wgray{\Phi_i}}
      +\mu_k
      \sum_{i\in\{\omega,r,z\}}\sum_{|\alpha|=k}
      \langle D^\alpha\mathcal R_i,
        D^\alpha\bvar{\delta v_i}\rangle_{\wgray{\Phi_{i,k}}}
      \right|
      \!\! \le \!\! \eulerstabLamR\mathfrak E_k$}
    {below=2.0mm of eulernonlinear}

  \eulerstabtagpanel
    {eulerclosure}
    {eulerstabpurple}
    {Closure}
    {Choose $\mu_k$ and the Stability Radius $\delta_*$}
    {$\eulerstabLamStab(\mu_k)
      \!\!\!\!:=\!\!\!\!\min\!\left\{
      \eulerstabLamLlow-\mu_k\eulerstabLamDLlow,\,
      \eulerstabLamDLhigh
      -\dfrac{\eulerstabLamLhigh}{\mu_k}
      \right\} \!\!\!\!>\!\!\!\!0$,\\[1.8pt]
     choose \, $\delta_* \!\!\!\!\!>\!\!\!\!\!0$ \,  with \, 
     $\bigl(\eulerstabLamNthree \!\! + \!\! \eulerstabLamDNthree\bigr) \, \delta_*
      \!+\!\dfrac{\eulerstabLamR }{\delta_*}
     \!\! \!\!< \!\!\!\! \eulerstabLamStab$}
    {below=of eulerresidual}

  \eulerstabtagpanelcompact
    {eulerproof}
    {eulerstabpurple}
    {Stability}
    {Stability}
    {If $\mathfrak E_k(0) \!\!\!\!<\!\!\!\! \delta_*$, then
     $\mathfrak E_k(t) \!\!\!\! < \!\!\!\! \delta_*$ for every later admissible time.\\[1.8pt]
     Thus the ball $\{\mathfrak E_k<\delta_*\}$ is forward invariant,
     which is the desired nonlinear stability statement}
    {below=of eulerclosure,
     inner ysep=6.8pt}

  \draw[eulerstabpurple!62,line width=0.28pt]
    ([xshift=1.9mm,yshift=-1.7mm]eulerproof.north west)
    rectangle
    ([xshift=-1.9mm,yshift=1.7mm]eulerproof.south east);

  \draw[eulerstab flow]
    (eulersetup.south) --
    node[eulerstab transition] {build the weighted energy}
    (weightedenergy.north);

  \draw[eulerstab flow]
    (weightedenergy.south) --
    node[eulerstab transition] {derive the analytic ingredients}
    (eulertoolbox.north);

  \draw[eulerstab flow]
    (eulertoolbox.south) --
    node[eulerstab transition] {obtain linear, nonlinear, and residual bounds}
    (eulerlinear.north);

  \draw[eulerstab flow]
    (eulerresidual.south) --
    node[eulerstab transition]
      {choose coupling and bootstrap radius}
    (eulerclosure.north);

  \draw[eulerstab flow]
    (eulerclosure.south) --
    node[eulerstab transition] {conclude the stability proof}
    (eulerproof.north);

\end{tikzpicture}%
}

\endgroup
\vspace{-18mm}
\label{diagram: Stability Proof}
\end{figure}

\clearpage

\subsubsection{Two-Radii Stability Theorem}

\noindent The single-radius certificate can be relaxed by prescribing separate target bounds for the low-order and top-order energies.  The same weighted energy argument still applies, provided the weighted ball is chosen small enough to lie inside the desired two-radii region,
\begin{equation}
\label{eq:two-radii-target-radii}
\mathfrak E_0(t)<\delta_0,
\qquad
\mathfrak H_k(t)<\delta_k,
\end{equation}
with $\delta_0>0$ and $\delta_k>0$.  \\

\noindent These bounds are obtained from the same full energy $
\mathfrak E_k^2=
\mathfrak E_0^2+\mu_k\mathfrak H_k^2$. \\

\noindent For a chosen weight $\mu_k>0$, the invariant set is the weighted ball $\mathfrak E_k(t)<\delta.$ This ball is contained in the target set \eqref{eq:two-radii-target-radii} provided that
\begin{equation}
\label{eq:two-radii-containment}
0<\delta\le \min\{\delta_0,\sqrt{\mu_k}\,\delta_k\},
\end{equation}
since this implies
\begin{equation}
\mathfrak E_0(t)\le \mathfrak E_k(t)<\delta\le\delta_0,
\qquad
\mathfrak H_k(t)\le \frac{\mathfrak E_k(t)}{\sqrt{\mu_k}}
<\frac{\delta}{\sqrt{\mu_k}}\le\delta_k.
\end{equation}

\hfill \\ 

\begin{theorem}[Two-radii nonlinear stability]
\label{thm:two-radii-stability}
Assume that the perturbation equations \eqref{eq: Stability theorem perturbation 1}--\eqref{eq: Stability theorem perturbation 3} hold for admissible perturbations. Fix $\delta_0>0$ and $\delta_k>0$.  Choose $\mu_k>0$ and $\delta>0$ satisfying
\begin{equation}
0<\delta\le \min\{\delta_0,\sqrt{\mu_k}\,\delta_k\}.
\end{equation}
Assume that the same certified estimates of \Cref{thm:main-stability} hold, namely \eqref{eq:ns-thm-low-linear}--\eqref{eq:ns-thm-residual}, with the same weighted energy
\begin{equation}
\mathfrak E_k^2=\mathfrak E_0^2+\mu_k\mathfrak H_k^2.
\end{equation}

\noindent Assume that $\mu_k$ satisfies the stability condition
\begin{equation}
\label{eq:two-radii-mu-choice}
\DLamStab(\mu_k):=
\min\left\{
\DLamLlow-\mu_k\DLamDLlow,
\ \DLamDLhigh-\frac{\DLamLhigh}{\mu_k}
\right\}>0.
\end{equation}
After this choice of $\mu_k$ is fixed, write $\DLamStab$ for the positive number $\DLamStab(\mu_k)$.  Define $\DLamThree:=\DLamNthree+\DLamDNthree$ and assume that
\begin{equation}
\label{eq:two-radii-smallness}
\DLamThree\delta+\frac{\DLamR}{\delta}<\DLamStab.
\end{equation}
Then every solution with
\begin{equation}
\label{eq:two-radii-initial-data}
\mathfrak E_0(0)^2+\mu_k\mathfrak H_k(0)^2<\delta^2
\end{equation}
satisfies
\begin{equation}
\label{eq:two-radii-conclusion}
\mathfrak E_0(t)<\delta_0,
\qquad
\mathfrak H_k(t)<\delta_k
\end{equation}
for as long as the solution exists and the certified estimates apply.
\end{theorem}

\begin{grayproof}
\noindent The weighted-energy part of the argument is identical to the single-radius argument in \Cref{thm:main-stability}. Thus, $\mathfrak E_k(t)<\delta$ for all times under consideration. It remains only to translate this invariant weighted ball into the two separate target bounds.  By the containment condition
\eqref{eq:two-radii-containment},
\begin{equation}
    \mathfrak E_0(t)\le \mathfrak E_k(t)<\delta\le\delta_0,
    \qquad
    \mathfrak H_k(t)\le\frac{\mathfrak E_k(t)}{\sqrt{\mu_k}}
    <\frac{\delta}{\sqrt{\mu_k}}\le\delta_k.
\end{equation}
\end{grayproof}

\clearpage

\subsection{Computer-Assisted Strategy for Proving Stability}
\label{sec: Computer-Assisted Strategy for Proving Stability}

\noindent Theorems~\ref{thm:main-stability} and~\ref{thm:two-radii-stability} reduce the nonlinear stability problem to a finite-dimensional optimization problem.  The aim is to choose the weights so that the linear damping margin is larger than the optimized nonlinear and residual error.  For each candidate choice of weights, we compute the certified constants, optimize the radius, and then update the weights until the closing inequality becomes strict.

\hfill

\noindent The key simplification is that, once the weights are fixed, the remaining optimization is only over the scalar radius \(\delta\).  Decreasing \(\delta\) reduces the cubic and quartic contributions, while increasing \(\delta\) reduces the residual contribution \(\DLamR/\delta\).  Thus the best radius balances these two effects, and the outer search changes the weights to improve this balance.

\hfill

\noindent The optimization loop proceeds follows.

\hfill \\ 

\noindent\textbf{Step 1. Choose the class of admissible weight functions.} \\ 

\noindent Start with a finite-dimensional parametrization of a class of admissible weight functions, chosen so that the weights have the required positivity, regularity, and symmetries. \\

\noindent For each parameter value~\(\wgray{\vartheta}\), write the low-order weights as
\begin{equation}
\wgray{\Phi_{\vartheta}}
=
\bigl(\wgray{\Phi_{\omega,\vartheta}},
\wgray{\Phi_{r,\vartheta}},
\wgray{\Phi_{z,\vartheta}}\bigr),
\end{equation}
and define the order-\(k\) weights by
\begin{equation}
\wgray{\Phi_{i,k,\vartheta}}
=
1+\rho^k\wgray{\Phi_{i,\vartheta}},
\qquad
\rho(r,z)=r^2+z^2,
\qquad
 i\in\{\omega,r,z\}.
\end{equation}
The optimization variables are the weight parameter \(\wgray{\vartheta}\) and the relative top-order weight \(\mu_k\).  In the single-radius case of \Cref{thm:main-stability}, one imposes \(0<\mu_k\le1\).  In the two-radii case of \Cref{thm:two-radii-stability}, one only imposes \(\mu_k>0\).  In both cases the weighted energy is
\begin{equation}
\mathfrak E_k^2
=
\mathfrak E_0^2+
\mu_k\mathfrak H_k^2.
\end{equation}
Choose an initial admissible candidate \((\wgray{\vartheta},\mu_k)\).

\hfill \\ 

\noindent\textbf{Step 2. Compute the linear damping for the current candidate.} \\

\noindent For the current values of \((\wgray{\vartheta},\mu_k)\), compute
\begin{equation}
\DLamLlow,
\qquad
\DLamLhigh,
\qquad
\DLamDLlow,
\qquad
\DLamDLhigh.
\end{equation}
The damping margin available for the combined energy is
\begin{equation}
\DLamStab
:=
\min\left\{
\DLamLlow-\mu_k\DLamDLlow,
\ \DLamDLhigh-\frac{\DLamLhigh}{\mu_k}
\right\}.
\label{eq:strategy-stab-margin}
\end{equation}
If \(\DLamStab\le0\), this candidate cannot close the theorem, so we penalize this term in the outer optimization and move to a new candidate.

\hfill \\ 

\noindent\textbf{Step 3. Compute the nonlinear and residual constants.}\\

\noindent If \(\DLamStab>0\), compute
\begin{equation}
\DLamThree
=
\DLamNthree+\DLamDNthree.
\end{equation}

\hfill 

\noindent For this candidate, the error that must be dominated by the damping is
\begin{equation}
\DLamThree\delta
+
\frac{\DLamR}{\delta}.
\label{eq:strategy-error-curve}
\end{equation}

\hfill \\ 

\noindent\textbf{Step 4. Optimize the radius and evaluate the gap.}\\

\noindent For the current candidate, let \(\mathcal A_\delta(\wgray{\vartheta},\mu_k)\subset(0,\infty)\) denote the nonempty set of radii for which every radius-dependent hypothesis used in the proof is valid. \\

\noindent In the single-radius case, the admissible condition is
\begin{equation}
\delta\in\mathcal A_\delta(\wgray{\vartheta},\mu_k).
\label{eq:strategy-single-radius-range}
\end{equation}

\hfill 

\noindent In the two-radii case with prescribed target bounds \(\delta_0>0\) and \(\delta_k>0\), the admissible conditions are
\begin{equation}
\delta\in\mathcal A_\delta(\wgray{\vartheta},\mu_k),
\qquad
0<\delta\le \min\{\delta_0,\sqrt{\mu_k}\,\delta_k\}.
\label{eq:strategy-two-radii-range}
\end{equation}

\hfill

\noindent Define the optimized gap
\begin{equation}
\mathcal J(\wgray{\vartheta},\mu_k)
:=
\inf_{\delta\ \mathrm{admissible}}
\left(
\DLamThree\delta
+
\frac{\DLamR}{\delta}
\right)
-
\DLamStab.
\label{eq:strategy-weight-objective}
\end{equation}

\hfill \\

\noindent The candidate succeeds exactly when \(\mathcal J(\wgray{\vartheta},\mu_k)<0\).  In this case, choose an admissible radius \(\delta_*\) such that
\begin{equation}
\DLamThree\delta_*
+
\frac{\DLamR}{\delta_*}
<
\DLamStab.
\label{eq:strategy-optimized-closing}
\end{equation}

\hfill \\ 

\hfill 

\noindent\textbf{Step 5. Iterate the outer optimization.} \\ 

\noindent The outer objective is the gap \(\mathcal J(\wgray{\vartheta},\mu_k)\). \\

\noindent One evaluation of this objective consists of Steps 2--4: compute the damping margin, compute the error constants, optimize \(\delta\), and subtract the margin. \\

\noindent If \(\mathcal J(\wgray{\vartheta},\mu_k)<0\), the  candidate satisfies the numerical closing condition and is retained for rigorous certification.  \\

\noindent If \(\mathcal J(\wgray{\vartheta},\mu_k)\ge0\), update \((\wgray{\vartheta},\mu_k)\) to decrease \(\mathcal J\), and repeat Steps 2--4.

\hfill \\ 

\hfill

\noindent\textbf{Step 6. Record the tentative certificate.}\\

\noindent Once \(\mathcal J<0\), record the weights and a selected admissible radius \(\delta_*\) satisfying \eqref{eq:strategy-optimized-closing}. The single-radius certificate gives
\begin{equation}
\mathfrak E_k(t)<\delta_*.
\end{equation}

\noindent Provided that $    0 <  \delta_*\le \min\{\delta_0,\sqrt{\mu_k}\delta_k\}$, the two-radii certificate gives
\begin{equation}
\mathfrak E_0(t)<\delta_0,
\qquad
\mathfrak H_k(t)<\delta_k.
\end{equation}

\hfill  

\noindent\textbf{Step 7. Validate rigorously the tentative certificate.}\\

\noindent Validate rigorously and analytically that the  candidate weights , radius, and associated constants satisfy the closing condition and all required estimates. The final certificate also verifies every radius-dependent admissibility guard at the selected radius \(\delta_*\).

\hfill  

\noindent \dotfill

\hfill 

\noindent We summarize below the numerical strategy used to  identify and rigorously certify the conditions required by the stability theorems.

\begin{certalgorithm}{Loop for optimizing the nonlinear stability certificate}

\item Initialize an admissible candidate \((\wgray{\vartheta},\mu_k)\). 

\vspace{1mm}

\item Repeat the following loop.
\begin{enumerate}[label=\textbf{\alph*.},leftmargin=2.5em,itemsep=1.3em,topsep=0.5em]
\item Compute
\begin{equation}
    \DLamStab
=
\min\left\{
\DLamLlow-\mu_k\DLamDLlow,
\ \DLamDLhigh-\frac{\DLamLhigh}{\mu_k}
\right\}.
\end{equation}
If \(\DLamStab\le0\), penalize the candidate and update \((\wgray{\vartheta},\mu_k)\).
\item If \(\DLamStab>0\), compute \(\DLamR\) and  $
    \DLamThree=\DLamNthree+\DLamDNthree$.
\item Determine the nonempty admissible set of radii for the current candidate, and compute
\begin{equation}
\mathcal J(\wgray{\vartheta},\mu_k)
=
\inf_{\delta\ \mathrm{admissible}}
\left(
\DLamThree\delta+
\frac{\DLamR}{\delta}
\right)
-
\DLamStab.
\end{equation}
\item If \(\mathcal J<0\), select an admissible radius \(\delta_*\) satisfying \eqref{eq:strategy-optimized-closing}, record the resulting tentative certificate, and  proceed to rigorous validation.
\item If \(\mathcal J\ge0\), update \((\wgray{\vartheta},\mu_k)\) to decrease \(\mathcal J\), and repeat the loop.
\end{enumerate}

\vspace{1mm}

\item Record the final weights, the selected admissible radius \(\delta_*\), and the resulting bounds.

\vspace{1mm}

\item Validate rigorously the tentative proof using analytic certificates, including every radius-dependent admissibility guard at \(\delta_*\).

\end{certalgorithm}

\clearpage

\section*{Conclusion}
\addcontentsline{toc}{section}{Conclusion}

\noindent In a recent numerical study~\citep{EulerBlowupMain}, we provided the first numerical evidence that the Euler equations admit a self-similar singularity in the unbounded setting using PINNs. The resulting PINN approximate profile can serve as the starting point for a computer-assisted proof, after being converted to a simpler analytic representation more suitable for rigorous certification. \\

\noindent This paper establishes a preliminary framework for proving nonlinear stability of the numerically discovered self-similar profile, reducing the remaining gap to rigorous certification and, if necessary, targeted analysis refinement. Starting from the exact perturbation and modulation equations, we construct a coupled low- and high-order energy and establish the linear and nonlinear estimates needed to control the perturbation dynamics, with every constant either explicit or computable by a prescribed rigorous procedure. Together with the residual and auxiliary bounds, these estimates form a unified stability framework. Nonlinear stability is therefore no longer left as a qualitative perturbative scenario, but reduced to a finite collection of concrete quantitative conditions whose rigorous verification would determine whether the candidate profile is stable and, if so, complete the nonlinear stability proof. \\

\noindent With the proof architecture in place, the remaining work is largely quantitative: to rigorously certify the profile-dependent constants and margins appearing in the established estimates and determine whether they are sufficient to close the stability argument. The spline representation, interval arithmetic, certified matrix bounds, and finite-dimensional optimization provide a direct framework for doing so with rigorous error control. If some margins prove insufficient, the modular structure of the argument allows individual estimates, weights, spline resolutions, or parameter choices to be sharpened without altering the underlying stability mechanism or proof strategy. If the required quantitative conditions can ultimately be certified with positive margin, then the candidate profile is nonlinearly stable, and the finite-time reconstruction yields an admissible solution that remains asymptotic to the profile and develops a singularity in finite physical time. Thus, within the framework developed here, the principal  immediate task is to determine quantitatively whether the candidate profile satisfies the conditions required for stability. If sufficient positive margin cannot be obtained, further analytic refinement may be required, and the possibility that the candidate profile is not stable under the required class of perturbations would remain open. \\ 

\noindent In parallel, the \texttt{LeanPDE} paper develops a formalization in Lean~\citep{Lean4} of key components of the argument, with particular emphasis on the symbolic derivations and proof steps underlying the stability analysis. Using \texttt{TorchLean}~\citep{TorchLean} together with \texttt{Mathlib}~\citep{mathlib}, the formalization connects machine-learning-based computations with Lean verification workflows. Numerical certificates are generated separately in \texttt{Julia} using \texttt{Arb} interval arithmetic~\citep{Julia,Arb}, while Lean verifies the algebraic identities and symbolic relations entering the certified estimates. The resulting pipeline assigns complementary roles to its different components: physics-informed optimization for discovery, splines and interval arithmetic for rigorous numerical certification, and Lean for formal symbolic verification and, ultimately, machine-checked analysis. In this way, \texttt{LeanPDE} provides a systematic route toward integrating the computational certificates developed here into a formally verified PDE proof.

\clearpage

\bibliographystyle{unsrtnat}
\bibliography{bib}

\clearpage

\appendix

\addcontentsline{toc}{part}{Appendix}

\section{Interpolation Between Low-Order and High-Order Energies}\label{app:high-order-energy}

The high-order energy contains only the top derivative level because all intermediate derivative levels are recovered by a weighted interpolation argument. \\ 

\noindent We first introduce the weights used at each derivative level.  For $1\le j\le k$, set
\begin{equation}
\label{eq:ns-Phiij}
\wgray{\Phi_{i,j}}:=1+\rho^{\,j}\wgray{\Phi_i},
\qquad
\wgray{\Phi_{i,0}}:=\wgray{\Phi_i},
\qquad i\in\{\omega,r,z\}.
\end{equation}
Thus the level $j$ weight has the same scaling as the top-order weight at derivative order $j$, while the level $0$ weight is the original low-order weight. \\

\hfill 

\noindent The one-step interpolation estimate compares the three adjacent weights $\wgray{\Phi_{i,j-1}}$, $\wgray{\Phi_{i,j}}$, and $\wgray{\Phi_{i,j+1}}$ for a fixed component $i$, so we next record the corresponding comparison constants:

\begin{equation}
\label{eq:ns-Aij}
A^{(0)}_{i,j}:=   \operatorname*{ess\,sup}_{\mathbb D}
 \frac{\wgray{\Phi_{i,j}}}{(\wgray{\Phi_{i,j-1}}\wgray{\Phi_{i,j+1}})^{1/2}},
\qquad
A^{(1)}_{i,j}:=   \operatorname*{ess\,sup}_{\mathbb D}
 \frac{|\nabla\wgray{\Phi_{i,j}}|}{(\wgray{\Phi_{i,j-1}}\wgray{\Phi_{i,j}})^{1/2}}.
\end{equation}

The first constant compares adjacent weights directly.  The second constant controls the weight derivative that appears after integration by parts. \\
\\

\noindent We note that the number of distinct derivatives of total order \(j\) in the two spatial variables \(r\) and \(z\) is $j+1$. Given parameters $\eta>0$ and $0<\theta<1$, define the one-step coefficient
\begin{equation}
\label{eq:ns-Bij}
B^{\rm step}_{i,j}(\eta,\theta):=
\frac{(j+1)^2\bigl(A^{(0)}_{i,j}\bigr)^2}{4(1-\theta)^2\eta}
+\frac{(j+1)^2\bigl(A^{(1)}_{i,j}\bigr)^2}{4\theta(1-\theta)}.
\end{equation}

\hfill \\

\noindent Finiteness of these constants is the concrete certification needed for the interpolation step.   

\hfill  \\

\begin{proposition}[One-step weighted interpolation]\label{prop:weighted-interpolation}
Let $\wgray{\Phi_{i,j}}$ be defined by \eqref{eq:ns-Phiij}, and define
\begin{equation}
X_{i,j}(f):=\sum_{|\alpha|=j}\|D^\alpha f\|_{\wgray{\Phi_{i,j}}}.
\end{equation}
Fix $1\le j\le k-1$. Assume that the three adjacent weights $\wgray{\Phi_{i,j-1}}$, $\wgray{\Phi_{i,j}}$, and $\wgray{\Phi_{i,j+1}}$ are locally $W^{1,\infty}$ away from the origin and satisfy positive lower bounds there, as follows for the manuscript weights from \eqref{eq:weight_coercivity} and \eqref{eq:ns-Phiij}. Assume also that the constants in \eqref{eq:ns-Aij} are finite. \\

\noindent Let $f$ have weighted weak derivatives through order $j+1$ with $X_{i,j-1}(f)$, $X_{i,j}(f)$, and $X_{i,j+1}(f)$ finite. Assume that integration by parts is justified after inserting cutoffs at the origin, the axis, and the far field, and that the resulting boundary-flux terms tend to zero as those cutoffs are removed. \\

\noindent Then, for every $\eta>0$ and every $0<\theta<1$,
\begin{equation}
\label{eq:app-step-interp}
X_{i,j}(f)^2
\le \eta X_{i,j+1}(f)^2+
B^{\rm step}_{i,j}(\eta,\theta)X_{i,j-1}(f)^2.
\end{equation}
\end{proposition}

\begin{grayproof}
\textbf{Reduction to adjacent derivative orders.} Fix $|\alpha|=j$ and write $D^\alpha=\partial_\ell D^\beta$ with $|\beta|=j-1$. We insert smooth cutoffs that exclude the origin and truncate the axis and far field. On each cutoff region the relevant weights are ordinary $W^{1,\infty}$ weights bounded above and below, so standard local mollification justifies the integration by parts. Passing first through the mollification limit and then removing the cutoffs using the stated boundary-flux assumption gives
\begin{equation}
\int_{\mathbb D}|D^\alpha f|^2\Phi_{i,j}
=-\int_{\mathbb D}D^\beta f\,\partial_\ell^2D^\beta f\,\Phi_{i,j}
-\int_{\mathbb D}D^\beta f\,\partial_\ell D^\beta f\,\partial_\ell\Phi_{i,j}.
\end{equation}

\hfill \\

\noindent\textbf{Estimating the adjacent-order terms.}  For the second-derivative term, we insert the adjacent weights:
\begin{align}
\left|\int D^\beta f\,\partial_\ell^2D^\beta f\,\Phi_{i,j}\right|
&\le A^{(0)}_{i,j}
\|D^\beta f\|_{\Phi_{i,j-1}}
\|\partial_\ell^2D^\beta f\|_{\Phi_{i,j+1}}.
\end{align}
For the weight-derivative term,
\begin{align}
\left|\int D^\beta f\,\partial_\ell D^\beta f\,\partial_\ell\Phi_{i,j}\right|
&\le A^{(1)}_{i,j}
\|D^\beta f\|_{\Phi_{i,j-1}}
\|\partial_\ell D^\beta f\|_{\Phi_{i,j}}.
\end{align}
Summing over all derivatives of order $j$ introduces the finite factor of $j+1$, and therefore
\begin{equation}
X_{i,j}^2
\le (j+1)A^{(0)}_{i,j}X_{i,j-1}X_{i,j+1}
+(j+1)A^{(1)}_{i,j}X_{i,j-1}X_{i,j}.
\end{equation}

\hfill \\ 

\noindent\textbf{Conclude using Young's inequalities.}  We apply Young's inequality to the first product with parameter $(1-\theta)\eta$.
\begin{equation}
 (j+1)A^{(0)}_{i,j}X_{i,j-1}X_{i,j+1}
 \le (1-\theta)\eta X_{i,j+1}^2
 +\frac{(j+1)^2(A^{(0)}_{i,j})^2}{4(1-\theta)\eta}X_{i,j-1}^2.
\end{equation}
We apply Young's inequality to the second product with parameter $\theta$:
\begin{equation}
 (j+1)A^{(1)}_{i,j}X_{i,j-1}X_{i,j}
 \le \theta X_{i,j}^2+
 \frac{(j+1)^2(A^{(1)}_{i,j})^2}{4\theta}X_{i,j-1}^2.
\end{equation}
Move $\theta X_{i,j}^2$ to the left and divide by $1-\theta$.  This gives \eqref{eq:app-step-interp}.
\end{grayproof}

\hfill \\

\noindent Proposition~\ref{prop:weighted-interpolation} is used for weighted intermediate derivative norms.  Such norms arise when a derivative falls on a coefficient or on a multiplier and the remaining perturbation derivative has order $j$ with $0<j<k$.  The estimate converts that intermediate weighted norm into two pieces: a small multiple of the next higher weighted norm and a controlled multiple of the previous weighted norm.  Iterating this one-step inequality gives the finite interpolation certificate in Proposition~\ref{prop:finite-interpolation-certificate}, which is the form used in the high-order linear and nonlinear estimates.  When the intermediate norm is unweighted, Proposition~\ref{prop:whole-plane-standard-constants} is used instead. 

\hfill

\begin{proposition}[Finite interpolation recursion certificate with exact parameters]\label{prop:finite-interpolation-certificate}
Fix $i\in\{\omega,r,z\}$ and $k\ge2$.  Assume that the one-step estimate \eqref{eq:app-step-interp} has been certified for every $1\le j\le k-1$.  Suppose that a concrete choice of the recursion parameters in the proof has been certified and produces constants \(\eta_{i,j}\) and    \(I_{i,j,k}(\eta_{i,j})\) . Then, for every \(1\le j\le k-1\),
\begin{equation}
\label{eq:interp-additive-certified}
X_{i,j}(f)\le \eta_{i,j} X_{i,k}(f)+I_{i,j,k}(\eta_{i,j})X_{i,0}(f).
\end{equation} 
\end{proposition}

\hfill

\begin{grayproof}
\textbf{Finite system of one-step estimates.}  Choose parameters $\eta_j>0$ and $\theta_j\in(0,1)$ and set
\begin{equation}\label{eq:interp-beta-def}
\beta_j:=B^{\rm step}_{i,j}(\eta_j,\theta_j),\qquad 1\le j\le k-1.
\end{equation}
Define the $(k-1)\times(k-1)$ tridiagonal matrix $Q_i=Q_i(\eta,\theta)$ by
\begin{equation}\label{eq:interp-Q-def}
(Q_i)_{ab}:=
\begin{cases}
1, & a=b,\\
-\eta_a, & b=a+1,\\
-\beta_a, & b=a-1,\\
0, & |a-b|\ge2,
\end{cases}
\qquad 1\le a,b\le k-1,
\end{equation}
where endpoint terms involving $y_0$ and $y_k$ are kept on the right-hand side. \\

\noindent Let
\begin{equation}\label{eq:interp-endpoint-vectors}
b^{(0)}:=\beta_1 e_1,
\qquad
b^{(k)}:=\eta_{k-1}e_{k-1},
\end{equation}
where $e_1,\ldots,e_{k-1}$ are the coordinate vectors in $\mathbb R^{k-1}$. \\

\noindent Assume that $Q_i$ is invertible, that $Q_i^{-1}$ has nonnegative entries, and define
\begin{equation}\label{eq:interp-algebraic-constants}
A^{\rm alg}_{i,j,k}:=\bigl(Q_i^{-1}b^{(k)}\bigr)_j,
\qquad
B^{\rm alg}_{i,j,k}:=\bigl(Q_i^{-1}b^{(0)}\bigr)_j.
\end{equation}
The certified parameters are admissible when
\begin{equation}
\sqrt{A^{\rm alg}_{i,j,k}}\le \eta_{i,j},
\qquad
I_{i,j,k}(\eta_{i,j}):=\sqrt{B^{\rm alg}_{i,j,k}}.
\label{eq:ns-Iijk-alg}
\end{equation} 

\noindent Let $y_j:=X_{i,j}(f)^2$ for $0\le j\le k$.  The one-step inequality gives
\begin{equation}
y_j-\eta_jy_{j+1}-\beta_jy_{j-1}\le0,
\qquad 1\le j\le k-1.
\end{equation}
For $j=1$, the term $\beta_1y_0$ contains the low-order endpoint.  For $j=k-1$, the term $\eta_{k-1}y_k$ contains the top-order endpoint.  Moving these endpoint terms to the right and keeping the internal variables $y=(y_1,\ldots,y_{k-1})$ on the left gives
\begin{equation}
\label{eq:interp-matrix}
Q_i(\eta,\theta)y\le b^{(0)}y_0+b^{(k)}y_k .
\end{equation}

\hfill

\noindent\textbf{Applying the finite interpolation certificate.}  Since $Q_i^{-1}$ has nonnegative entries, multiplying \eqref{eq:interp-matrix} by $Q_i^{-1}$ preserves the componentwise inequality.  The $j$-th component gives
\begin{equation}
y_j\le \bigl(Q_i^{-1}b^{(k)}\bigr)_j y_k+\bigl(Q_i^{-1}b^{(0)}\bigr)_j y_0
=A^{\rm alg}_{i,j,k}y_k+B^{\rm alg}_{i,j,k}y_0.
\end{equation}
Taking square roots and using $\sqrt{a+b}\le\sqrt a+\sqrt b$ gives
\begin{equation}
X_{i,j}(f)
\le \sqrt{A^{\rm alg}_{i,j,k}}X_{i,k}(f)
+\sqrt{B^{\rm alg}_{i,j,k}}X_{i,0}(f).
\end{equation}
If the certified parameters satisfy $\sqrt{A^{\rm alg}_{i,j,k}}\le\eta_{i,j}$, then the preceding inequality becomes \eqref{eq:interp-additive-certified}. 
\end{grayproof}

\hfill

\noindent The conclusion of Proposition~\ref{prop:finite-interpolation-certificate} is available only after the recursion data have been certified for the weight family under consideration.  A certificate must provide precise values or rigorous interval enclosures for the entries of $Q_i$, verify that $Q_i$ is invertible and that $Q_i^{-1}$ has nonnegative entries, for example through an explicit nonsingular M-matrix criterion, and check the inequalities in    \eqref{eq:ns-Iijk-alg}  for every required value of $j$.

\hfill \\ 

\noindent Proposition~\ref{prop:finite-interpolation-certificate} is used when the estimate requires a weighted interpolation inequality with an explicitly certified coefficient in front of the top-order norm.  For an unweighted reflected norm, Proposition~\ref{prop:whole-plane-standard-constants} supplies the corresponding interpolation constant directly. \\

\hfill 

\paragraph{Interpolation statement.}
Proposition~\ref{prop:finite-interpolation-certificate} gives directly the additive estimate

\begin{equation}
\label{eq:ns-lower-from-top}
\sum_{|\alpha|=j}\|D^\alpha f\|_{\wgray{\Phi_{i,j}}}
\le
\eta_{i,j}\sum_{|\alpha|=k}\|D^\alpha f\|_{\wgray{\Phi_{i,k}}}
+I_{i,j,k}(\eta_{i,j})\|f\|_{\wgray{\Phi_i}},
\qquad 0<j<k.
\end{equation}
Here $I_{i,j,k}(\eta_{i,j})$ is the certified constant defined in \eqref{eq:ns-Iijk-alg}.

\hfill  \\

\paragraph{Consequence for the stability energy.}
Applying \eqref{eq:ns-lower-from-top} to $\bvar{\delta\omega}$, $\partial_r\bvar{\delta u}$, or $\partial_z\bvar{\delta u}$ gives
\begin{equation}
\label{eq:ns-all-lower-controlled}
\sum_{|\alpha|=j}\|D^\alpha f\|_{\wgray{\Phi_{i,j}}}
\le
\bigl(\eta_{i,j}\sqrt{k+1}\,\mu_k^{-1/2}+I_{i,j,k}(\eta_{i,j})\bigr)\mathfrak E_k,
\qquad 0<j<k.
\end{equation}
The factor $\sqrt{k+1}$ is part of the certified coefficient.  It appears because \eqref{eq:ns-lower-from-top} contains the sum of the $k+1$ derivatives of total order $k$ in two variables, while $\mathfrak H_k$ contains their square-sum.  Hence
\begin{equation}
\sum_{|\alpha|=k}\|D^\alpha f\|_{\wgray{\Phi_{i,k}}}
\le \sqrt{k+1}\left(\sum_{|\alpha|=k}\|D^\alpha f\|_{\wgray{\Phi_{i,k}}}^2\right)^{1/2}
\le \sqrt{k+1}\,\mathfrak H_k
\le \sqrt{k+1}\,\mu_k^{-1/2}\mathfrak E_k.
\end{equation}

\hfill

\noindent In the stability estimates, \eqref{eq:ns-lower-from-top} is used when the low- and top-order contributions are kept separate until Young's inequality is applied. The parameter $\eta_{i,j}$ is chosen so that the top-order piece fits the damping reserved for that product, while the term with $I_{i,j,k}(\eta_{i,j})$ is recorded in the low-order contribution. Equation~\eqref{eq:ns-all-lower-controlled} is used when a single $\mathfrak E_k$ bound is sufficient.

\clearpage

\section{Analytical Theorems for the Stability Proof}\label{app:analytical-theorems}

\hfill

\subsection{Preliminaries} \label{sec: analytic preliminaries}

\hfill 

 \noindent This subsection fixes the analytic language used throughout the estimates.  The argument first specifies the weak-derivative and multi-index conventions, then records the whole-plane Fourier and Sobolev normalizations, and finally explains how half-plane quantities are reflected to the whole plane and converted back to the weighted norms in the energy.  

\hfill

\noindent In whole-plane estimates the spatial variable is written as $x=(r,z)\in\mathbb R^2$.  In half-plane estimates the same coordinates are used with $r\ge0$.  When a half-plane function is reflected, the reflected whole-plane function is denoted by an uppercase letter, such as $F$, and the original half-plane function is denoted by the corresponding lowercase letter, such as $f$.

\hfill \\

\noindent\textbf{Multi-index notation.}  A multi-index is a vector $\alpha$ of integers, written
\begin{equation}
\alpha=(\alpha_1,\ldots,\alpha_n).
\end{equation}

\hfill 

\noindent Its order is
\begin{equation}
|\alpha|:=\alpha_1+\cdots+\alpha_n,
\end{equation}
and the corresponding weak or classical derivative is
\begin{equation}
D^\alpha u:=\partial_{x_1}^{\alpha_1}\cdots\partial_{x_n}^{\alpha_n}u.
\end{equation}

\hfill 

\noindent If $\beta=(\beta_1,\ldots,\beta_n)$, then $\beta\le\alpha$ means $\beta_j\le\alpha_j$ for every $j$.  The multinomial coefficient is
\begin{equation}
\binom{\alpha}{\beta}:=\prod_{j=1}^n\binom{\alpha_j}{\beta_j}.
\end{equation}

\hfill 

\noindent With this notation, the ordinary Leibniz expansion takes the form
\begin{equation}
D^\alpha(AB)=\sum_{\beta\le\alpha}\binom{\alpha}{\beta}(D^\beta A)(D^{\alpha-\beta}B).
\end{equation}

\clearpage

\noindent\textbf{Weak-derivative framework.}  \\ 

\noindent The next theorem fixes the definition of weak derivatives, the elementary calculus rules for weak derivatives, and the Sobolev norm convention.  These facts justify the differentiated equations, the integration-by-parts identities, and the Sobolev norms used in the energy estimates.

\hfill

\begin{theorem}[Weak derivatives and Sobolev spaces, {\citep[Sections~5.2.1--5.2.3, Theorems~1--2]{evans2010partial}}]\label{thm:evans-weak-sobolev}
Let $U\subset\mathbb R^n$ be open.
\begin{enumerate}[label=(E\arabic*),leftmargin=.38in]
\item A locally integrable function $v$ is the weak derivative $D^\alpha u$ if
\begin{equation}
\int_U uD^\alpha\varphi\,dx=(-1)^{|\alpha|}\int_U v\varphi\,dx
\qquad \text{for all }\varphi\in C_c^\infty(U).
\end{equation}
\item If a weak derivative exists, it is unique up to equality almost everywhere.
\item If $u$ has classical derivatives through order $|\alpha|$, then the classical derivative $D^\alpha u$ is also the weak derivative in the preceding sense.
\item Linearity: for constants $a,b$, whenever the displayed weak derivatives exist,
\begin{equation}
D^\alpha(au+bw)=aD^\alpha u+bD^\alpha w.
\end{equation}
\item Higher weak derivatives are independent of their order: if the relevant derivatives exist locally integrably, then $D^\alpha D^\beta u=D^{\alpha+\beta}u$ in the weak sense.
\item For $1\le p\le\infty$ and an integer $k\ge0$,
\begin{equation}
W^{k,p}(U):=\{u\in L^p(U):D^\alpha u\in L^p(U)\text{ for every }|\alpha|\le k\},
\end{equation}
with norm
\begin{equation}
\|u\|_{W^{k,p}(U)}:=\left(\sum_{|\alpha|\le k}\|D^\alpha u\|_{L^p(U)}^p\right)^{1/p}
\end{equation}
for $1\le p<\infty$, with the usual essential-supremum modification when $p=\infty$.  For $p=2$ one writes $H^k(U)=W^{k,2}(U)$.
\end{enumerate}
\end{theorem}

\hfill \\

\noindent We use the standard $L^\infty$ norm. For a measurable function on $U$,
\begin{equation}
\|f\|_{L^\infty(U)}
:=
\operatorname*{ess\,sup}_{x\in U}|f(x)|
=
\inf\left\{M\ge0:\ |f(x)|\le M\text{ for almost every }x\in U\right\}.
\label{eq:Linfty-essential-supremum}
\end{equation}
\noindent The essential supremum ignores values on sets of measure zero. Unless a pointwise supremum is written explicitly, $\|\cdot\|_{L^\infty}$ below always means the essential supremum. Point evaluations in the modulation conditions are justified separately by the regularity estimates that make them well defined.

\hfill \\ 

\noindent In the stability proof, \Cref{thm:evans-weak-sobolev} justifies using the same derivative notation for smooth approximations and Sobolev perturbations.  It also ensures that the differentiated identities obtained formally from the equations remain valid in the weak sense once the required Sobolev norms are finite. 

\clearpage

\noindent We also write $K \Subset \mathbb{D}$ to mean that \(K\) is compactly contained in \(\mathbb{D}\), i.e.,
\begin{equation}
K \text{ is compact}, \qquad K \subset \mathbb{D}, \qquad \operatorname{dist}(K,\partial\mathbb{D})>0.
\end{equation}

\hfill 

\noindent The local \(L^p\) space is defined by
\begin{equation}
L^p_{\mathrm{loc}}(\mathbb{D}) := \left\{u : u|_K \in L^p(K) \text{ for every } K \Subset \mathbb{D}\right\}.
\end{equation}
Thus, \(u \in L^p_{\mathrm{loc}}(\mathbb{D})\) means that \(u\) belongs to \(L^p\) on every compact subset lying strictly inside \(\mathbb{D}\). In contrast, \(u \in L^p(\mathbb{D})\) requires \(u\) to belong to \(L^p\) on the entire domain. Hence \(L^p_{\mathrm{loc}}(\mathbb{D})\) imposes no integrability condition near the boundary of \(\partial\mathbb{D}\), or at infinity since \(\mathbb{D}\) is unbounded. \\ 

\noindent Similarly, the local Sobolev space is defined by
\begin{equation}
W^{k,p}_{\mathrm{loc}}(\mathbb{D}) := \left\{u : D^\alpha u \in L^p_{\mathrm{loc}}(\mathbb{D}) \text{ for every multi-index } \alpha \text{ with } |\alpha| \leq k\right\}.
\end{equation}
Thus, \(u \in W^{k,p}_{\mathrm{loc}}(\mathbb{D})\) means that \(u\) and all of its weak derivatives up to order \(k\) belong to \(L^p_{\mathrm{loc}}(\mathbb{D})\). Equivalently,
\begin{equation}
u \in W^{k,p}_{\mathrm{loc}}(\mathbb{D}) \qquad \text{if and only if}\qquad u|_U \in W^{k,p}(U) \text{ for every open } U \Subset \mathbb{D}.
\end{equation}

\hfill \\

\noindent We now state the weak product rule, which is used because the differentiated linear and nonlinear estimates repeatedly expand derivatives of products of coefficients, perturbations, and streamfunction multipliers. \\

\begin{proposition}[Weak Leibniz rule for cutoffs and commutation]\label{thm:weak-calculus-used}
Let $\zeta\in C^\infty$ have bounded derivatives up to order~$k$, and let $u\in H^k_{\rm loc}(U)$.  Then for every $|\alpha|\le k$, we have in the weak sense that
\begin{equation}
\label{eq:app-leibniz-cutoff}
D^\alpha(\zeta u)=\sum_{\beta\le\alpha}\binom{\alpha}{\beta}D^\beta\zeta\,D^{\alpha-\beta}u.
\end{equation} 
\end{proposition}

\begin{grayproof} 
\noindent The proof reduces the weak statement to the classical product rule.  We choose a compact set $K\Subset U$ containing the supports of the test functions used in the weak formulation.  By the standard approximation theorem for Sobolev functions on compact subsets, there are functions $u_n\in C^\infty(K)$ such that
\begin{equation}
u_n\to u
\qquad\text{in }H^k(K).
\end{equation}

\vspace{1mm}

\noindent\textbf{Classical product rule on the approximants.}
For each smooth $u_n$, the classical Leibniz rule gives
\begin{equation}
D^\alpha(\zeta u_n)=\sum_{\beta\le\alpha}\binom{\alpha}{\beta}D^\beta\zeta\,D^{\alpha-\beta}u_n.
\end{equation}
Since every coefficient $D^\beta\zeta$ is bounded on $K$, multiplication by $D^\beta\zeta$ is a bounded map on $L^2(K)$.  Therefore
\begin{equation}
D^\beta\zeta\,D^{\alpha-\beta}u_n
\to
D^\beta\zeta\,D^{\alpha-\beta}u
\qquad\text{in }L^2(K)
\end{equation}
for each $\beta\le\alpha$.  Hence the right-hand side of the smooth product rule converges in $L^2(K)$ to
\begin{equation}
\sum_{\beta\le\alpha}\binom{\alpha}{\beta}D^\beta\zeta\,D^{\alpha-\beta}u.
\end{equation}

\vspace{1mm}

\noindent\textbf{Passing to the Sobolev limit.}
At the same time, $\zeta u_n\to\zeta u$ in $H^k(K)$, so
\begin{equation}
D^\alpha(\zeta u_n)\to D^\alpha(\zeta u)
\qquad\text{in }L^2(K).
\end{equation}
The two limits must agree, and this gives \eqref{eq:app-leibniz-cutoff} in $L^2(K)$ and hence in the weak sense on $U$.
\end{grayproof}

\clearpage

\hfill 

\noindent\textbf{Continuous embeddings.}  If $X$ and $Y$ are normed function spaces, the notation
\begin{equation}
X\hookrightarrow Y
\end{equation}
means that $X\subset Y$ and that there is a finite constant $C_{X,Y}$ such that
\begin{equation}
\|F\|_Y\le C_{X,Y}\|F\|_X
\qquad\text{for every }F\in X.
\end{equation}
Thus every embedding used below is an estimate with a constant.  For example, $H^s(\mathbb R^2)\hookrightarrow L^\infty(\mathbb R^2)$ means that a whole-plane Sobolev norm controls the \(L^\infty\) norm.  After reflection, this is the basic mechanism that gives $L^\infty$ bounds for the half-plane perturbation variables. 

\hfill \\ 

\hfill 

\noindent\textbf{Fourier convention.}  The whole-plane Fourier transform is normalized by
\begin{equation}\label{eq:fourier-convention}
\widehat F(\xi):=\int_{\mathbb R^2}e^{-ix\cdot\xi}F(x)\,dx,
\qquad
F(x)=\frac{1}{(2\pi)^2}\int_{\mathbb R^2}e^{ix\cdot\xi}\widehat F(\xi)\,d\xi.
\end{equation}

\hfill 

\noindent With this convention, Plancherel's identity is
\begin{equation}\label{eq:plancherel-convention}
\|F\|_{L^2(\mathbb R^2)}^2
=\frac{1}{(2\pi)^2}\|\widehat F\|_{L^2(\mathbb R^2)}^2.
\end{equation}

\hfill 

\noindent All Bessel-potential norms in this paper use the same normalization.  

\hfill \\

\hfill 

\noindent\textbf{Bessel-potential Sobolev norm.}  For $s\ge0$, define
\begin{equation}\label{eq:Bessel-norm-prelim}
\|F\|_{H^s_B(\mathbb R^2)}
:=\left(\frac{1}{(2\pi)^2}\int_{\mathbb R^2}(1+|\xi|^2)^s|\widehat F(\xi)|^2\,d\xi\right)^{1/2}.
\end{equation}

\hfill 

\noindent This is the norm associated with the Bessel potential operator \begin{equation}
    J^s=(1-\Delta)^{s/2}.
\end{equation}  
Indeed, since the Fourier multiplier of $J^s$ is $(1+|\xi|^2)^{s/2}$, the Plancherel identity gives
\begin{equation}
\|F\|_{H^s_B(\mathbb R^2)}=\|J^sF\|_{L^2(\mathbb R^2)}.
\end{equation}

\hfill 

\noindent The sharp whole-plane Sobolev constants quoted below are stated in this Bessel-potential normalization.   \\

\noindent The following comparison gives the explicit conversion between the two descriptions.

\clearpage

\noindent\textbf{Derivative Sobolev norm and comparison constants.}

\hfill

\begin{proposition}[Sobolev--Bessel norm comparison on \(\mathbb R^2\)]\label{thm:Bessel-derivative-comparison}
For every integer \(k\ge1\),
\begin{equation}\label{eq:Bessel-derivative-equivalence}
\frac{1}{\sqrt{k+1}}\|F\|_{H^k}
\le \|F\|_{H^k_B}
\le 3^{k/2}\|F\|_{H^k}.
\end{equation}
\end{proposition}

\vspace{4mm}

\begin{grayproof} 
\noindent\textbf{Fourier representation of the derivative norm.}
By Plancherel in the convention \eqref{eq:plancherel-convention},
\begin{equation}
\|F\|_{H^k}^2
=\frac{1}{(2\pi)^2}\int_{\mathbb R^2}
\bigg( \,\sum_{|\alpha|\le k}\xi^{2\alpha} \,\bigg)|\widehat F(\xi)|^2\,d\xi.
\end{equation}
Set
\begin{equation}
N_k(\xi):=\sum_{|\alpha|\le k}\xi^{2\alpha},
\qquad
 t:=|\xi|^2=\xi_1^2+\xi_2^2.
\end{equation}

\hfill 

\noindent The Bessel norm is
\begin{equation}
\|F\|_{H^k_B}^2
=\frac{1}{(2\pi)^2}\int_{\mathbb R^2}(1+t)^k|\widehat F(\xi)|^2\,d\xi.
\end{equation}

\hfill 

\noindent The two norms differ only through the factors multiplying \(|\widehat F(\xi)|^2\) under the same Fourier integral.  Therefore pointwise estimates of the form
\begin{equation}
\label{eq:Bessel-polynomial-comparison-target}
\frac{1}{3^k}(1+t)^k\le N_k(\xi)\le (k+1)(1+t)^k
\end{equation}
valid for every \(\xi\), can be multiplied by \(|\widehat F(\xi)|^2\) and integrated over \(\mathbb R^2\).  This converts the comparison of the two norms into the elementary comparison of the two polynomial weights \(N_k(\xi)\) and \((1+t)^k\).

\hfill \\

\noindent\textbf{Lower bound for the Bessel norm.}
For each \(0\le j\le k\), every monomial \(\xi^{2\alpha}\) with \(|\alpha|=j\) is one term in the expansion of \((\xi_1^2+\xi_2^2)^j=t^j\), with a positive coefficient.  Hence
\begin{equation}
\sum_{|\alpha|=j}\xi^{2\alpha}\le t^j\le (1+t)^k.
\end{equation}

\hfill 

\noindent Summing this inequality over \(j=0,\ldots,k\) gives
\begin{equation}
N_k(\xi)
=\sum_{j=0}^k\sum_{|\alpha|=j}\xi^{2\alpha}
\le (k+1)(1+t)^k.
\end{equation}

\hfill 

\noindent After multiplication by \(|\widehat F(\xi)|^2\), integration, and Plancherel's theorem, this gives
\begin{equation}
\|F\|_{H^k}^2\le (k+1)\|F\|_{H^k_B}^2,
\end{equation}
which is the left inequality in \eqref{eq:Bessel-derivative-equivalence}.

\hfill \\

\noindent\textbf{Upper bound for the Bessel norm.}
The multinomial formula gives
\begin{equation}
t^j=(\xi_1^2+\xi_2^2)^j
=\sum_{a=0}^j\binom{j}{a}\xi_1^{2a}\xi_2^{2(j-a)}.
\end{equation}
Since \(\binom{j}{a}\le 2^j\),
\begin{equation}
t^j\le 2^j\sum_{|\alpha|=j}\xi^{2\alpha}.
\end{equation}
Therefore
\begin{equation}
\begin{aligned}
(1+t)^k
&=\sum_{j=0}^k\binom{k}{j}t^j\le\sum_{j=0}^k\binom{k}{j}2^j\sum_{|\alpha|=j}\xi^{2\alpha}\le\left(\sum_{j=0}^k\binom{k}{j}2^j\right)N_k(\xi)
=3^kN_k(\xi).
\end{aligned}
\end{equation}

\hfill 

\noindent Multiplication by \(|\widehat F(\xi)|^2\), integration, and Plancherel's theorem yield
\begin{equation}
\|F\|_{H^k_B}^2\le 3^k\|F\|_{H^k}^2,
\end{equation}
which proves the right inequality in \eqref{eq:Bessel-derivative-equivalence}.

\end{grayproof}

\clearpage

\clearpage

\subsection{Extension, Sobolev, and Morrey estimates}

\noindent The next estimates explain how Sobolev control enters the energy argument.  We first recall the standard extension, Morrey, and Sobolev results that justify pointwise evaluation and full-space comparison.  The parity assumptions then give an explicit reflection from the half-plane $r\ge0$ to $\mathbb R^2$, and the whole-plane estimates provide the $L^\infty$ and interpolation bounds used in the energy estimates.  Streamfunction quantities are handled separately by lifting the elliptic equation to a five-dimensional Poisson problem and reducing the Newtonian kernel back to the variables $(r,z)$.

\hfill \\ 

\noindent The extension theorem is the standard result that lets Sobolev functions on a regular domain be compared with Sobolev functions on the full ambient space. \\

\begin{theorem}[Extension theorem, {\citep[Section 5.4, Theorem 1]{evans2010partial}}]
\label{thm:evans-extension-input}
Let $1<p<\infty$.  If $U\subset\mathbb R^n$ is bounded with $C^1$ boundary and $U\Subset V$, then there is a bounded linear extension operator
\begin{equation}
E:W^{1,p}(U)\to W^{1,p}(\mathbb R^n)
\end{equation}
such that $Eu=u$ almost everywhere on $U$, $Eu$ is supported in $V$, and
\begin{equation}
\|Eu\|_{W^{1,p}(\mathbb R^n)}\le C\|u\|_{W^{1,p}(U)},
\end{equation}
where $C$ depends only on $p$, $U$, and $V$.
\end{theorem}

\hfill 

\noindent This result explains why it is legitimate to move from a half-plane Sobolev problem to a full-space Sobolev problem once the boundary behavior has been fixed.  

\hfill  \\ 

\hfill

\noindent Morrey's estimate is the standard theorem that turns sufficiently strong Sobolev control into pointwise control.  It is used here to identify the regularity mechanism that makes the point values in the modulation equations meaningful once enough derivatives are controlled.\\

\begin{theorem}[Morrey estimate, {\citep[Section 5.6.2, Theorem 5]{evans2010partial}}]
\label{thm:evans-morrey-input}
Let $U\subset\mathbb R^n$ be bounded and open with $C^1$ boundary, let $n<p<\infty$, and let $u\in W^{1,p}(U)$.  Then $u$ has a version $u^*\in C^{0,\gamma}(\overline U)$, where
\begin{equation}
\gamma=1-\frac{n}p,
\end{equation}
and
\begin{equation}
\|u^*\|_{C^{0,\gamma}(\overline U)}\le C\|u\|_{W^{1,p}(U)},
\end{equation}
where $C$ depends only on $p$, $n$, and $U$.
\end{theorem}

\clearpage 

\noindent The Sobolev inequality is the embedding result behind the conversion from derivative control to integrability and Hölder control.  The theorem is first stated for the bounded domain, and later in the full-plane version with explicit constants, which is the form used for the stability estimates.

\hfill

\begin{theorem}[Sobolev inequalities, {\citep[Section 5.6.3, Theorem 6]{evans2010partial}}]
\label{thm:evans-sobolev-input}
Let $U\subset\mathbb R^n$ be bounded with $C^1$ boundary and let $u\in W^{k,p}(U)$. \\

\noindent If $k<n/p$, then $u\in L^q(U)$ with
\begin{equation}
\frac{1}q=\frac{1}p-\frac{k}n,
\end{equation}
and $\|u\|_{L^q(U)}\le C\|u\|_{W^{k,p}(U)}$.  If $k>n/p$, then
\begin{equation}
u\in C^{k-[n/p]-1,\gamma}(\overline U),
\end{equation}
with $\gamma=[n/p]+1-n/p$ when $n/p$ is not an integer, and with any $0<\gamma<1$ when $n/p$ is an integer.  The corresponding Hölder norm is bounded by $C\|u\|_{W^{k,p}(U)}$.  In particular, when $n=2$, $p=2$, and $k=2$,
\begin{equation}
H^2(U)=W^{2,2}(U)\hookrightarrow C^{0,\gamma}(\overline U)\hookrightarrow L^\infty(U),
\qquad 0<\gamma<1.
\end{equation}
\end{theorem}

\hfill

\noindent The special case $H^2\hookrightarrow L^\infty$ is of particular interest for the stability proof.  The Euler modulation formulas require the streamfunction origin values \(\bvar{\delta\psi}_0\) and \((\partial_z\bvar{\delta\psi})_0\), which are controlled by the origin-kernel elliptic estimates below.

\clearpage

\subsection{Whole-plane endpoint Sobolev estimate}\label{subsec:full-space-estimates}

\noindent The next result gives the sharp whole-plane Sobolev constant used after parity reflection to place lower-order factors in $L^\infty$.

\hfill

\begin{proposition}[Morosi--Pizzocchero endpoint Sobolev constant, {\citep[Propositions 2.1 and 3.1]{MorosiPizzoccheroSobolev}}]\label{thm:mp-upper-bound}\label{thm:mp-sharp-endpoint}
Let \(S_{\infty,n,d}\) be the sharp constant in
\begin{equation}
\|f\|_{L^\infty(\mathbb R^d)}\le S_{\infty,n,d}\|f\|_{H^n(\mathbb R^d)},
\qquad n>\frac{d}{2}.
\end{equation} 
Then, the upper bound
\begin{equation}
S^+_{\infty,n,d}
=
\frac{1}{(4\pi)^{d/4}}
\left(\frac{\Gamma(n-d/2)}{\Gamma(n)}\right)^{1/2}
\end{equation}
is sharp at the endpoint \(r=\infty\).  Hence
\begin{equation}
\|f\|_{L^\infty(\mathbb R^d)}
\le
\frac{1}{(4\pi)^{d/4}}
\left(\frac{\Gamma(n-d/2)}{\Gamma(n)}\right)^{1/2}
\|f\|_{H^n(\mathbb R^d)}.
\end{equation} 
\end{proposition}

\hfill 

\hfill \\ 

\noindent In the stability proof, this is used only with \(d=2\) and \(n=s>1\).  Substituting \(d=2\) gives
\begin{equation}
\frac{1}{(4\pi)^{1/2}}
\left(\frac{\Gamma(s-1)}{\Gamma(s)}\right)^{1/2}
=\frac{1}{2\sqrt\pi}\left(\frac{\Gamma(s-1)}{\Gamma(s)}\right)^{1/2}
=\frac{1}{2\sqrt{\pi(s-1)}}.
\end{equation}

\hfill \\

\noindent With the Bessel norm defined in \eqref{eq:Bessel-norm-prelim}, the precise full-plane estimate used is
\begin{equation}\label{eq:MP-sharp-constant-current-convention}
\|F\|_{L^\infty(\mathbb R^2)}
\le \frac{1}{2\sqrt{\pi(s-1)}}\|F\|_{H^s_B(\mathbb R^2)},
\qquad s>1.
\end{equation}
The Morosi--Pizzocchero formula is quoted for the sharp constant in their Bessel-potential convention.  The conversion to the Fourier convention fixed in \eqref{eq:fourier-convention}, including the two-dimensional specialization used later, is carried out more carefully in Proposition~\ref{prop:whole-plane-standard-constants}.

\clearpage

\clearpage

\hfill 

\subsection{Whole-plane estimates used in the energy argument}\label{subsec:whole-plane-estimates}

\noindent The next results are the pointwise and interpolation inequalities used repeatedly in the energy estimates.  Proposition~\ref{prop:whole-plane-reflection} converts half-plane Sobolev norms into whole-plane Sobolev norms.  Proposition~\ref{prop:whole-plane-standard-constants} supplies the explicit $H^s\to L^\infty$ and interpolation constants. \\

\noindent In the energy estimates, each half-plane factor is first reflected to \(\mathbb R^2\).  The whole-plane Sobolev bound then places lower-order factors in \(L^\infty\), and interpolation reduces intermediate derivative levels to top-order and low-order pieces.  

\hfill 

\begin{proposition}[Parity reflection to the whole plane]\label{prop:whole-plane-reflection}
\noindent Let \(m\ge0\) be an integer and let
\begin{equation}
\mathbb{D}:=\{(r,z):r>0,\ z\in\mathbb R\}.
\label{eq:reflection-half-plane}
\end{equation}
\noindent For \(p\in\{0,1\}\), define the parity extension by
\begin{equation}
(\mathcal E_p f)(r,z):=(\operatorname{sgn}r)^p f(|r|,z)
\qquad (r\ne0),
\label{eq:parity-extension-definition}
\end{equation}
\noindent where \(p=0\) and \(p=1\) give the even and odd extensions.  Assume \(f\in H^m(\mathbb{D})\) and
\begin{equation}
\partial_r^q\partial_z^\ell f(0,z)=0
\quad\text{in the trace sense whenever }q+p\text{ is odd and }q+\ell\le m-1.
\label{eq:parity-extension-trace-condition}
\end{equation}
\noindent Then \(\mathcal E_p f\in H^m(\mathbb R^2)\).  With the half-plane norm
\begin{equation}
\|f\|_{H^m(\mathbb{D})}^2
:=\sum_{q+\ell\le m}
\|\partial_r^q\partial_z^\ell f\|_{L^2(\mathbb{D})}^2,
\label{eq:half-plane-intrinsic-Hm-norm}
\end{equation}
\noindent one has the exact identity
\begin{equation}
\|\mathcal E_p f\|_{H^m(\mathbb R^2)}^2
=2\|f\|_{H^m(\mathbb{D})}^2.
\label{eq:parity-extension-Hm-norm}
\end{equation}

\hfill 

\noindent If \(\Phi\) is a measurable nonnegative weight and \(\Phi^{\rm even}(r,z):=\Phi(|r|,z)\), then for every \(q+\ell\le m\),
\begin{equation}
\|\partial_r^q\partial_z^\ell\mathcal E_p f\|_{L^2(\Phi^{\rm even},\mathbb R^2)}^2
=2\|\partial_r^q\partial_z^\ell f\|_{L^2(\Phi,\mathbb{D})}^2,
\label{eq:parity-extension-weighted-norm}
\end{equation}
\noindent whenever the right side is finite. 
\end{proposition}

\begin{grayproof}
\noindent\textbf{Smooth reflection and pointwise derivatives.}  We first assume that \(f\in C^\infty(\overline{\mathbb{D}})\) and that the trace conditions in \eqref{eq:parity-extension-trace-condition} hold classically.  For \(r>0\), the extension equals \(f(r,z)\).  For \(r<0\), it equals \((-1)^p f(-r,z)\).  Differentiating away from the axis gives
\begin{equation}
\partial_r^q\partial_z^\ell(\mathcal E_p f)(r,z)
=
\begin{cases}
\partial_r^q\partial_z^\ell f(r,z),&r>0,\\
(-1)^{q+p}\partial_r^q\partial_z^\ell f(-r,z),&r<0.
\end{cases}
\label{eq:parity-extension-derivative-formula}
\end{equation}

\hfill \\

\noindent\textbf{Matching traces at the axis.}  For \(q+\ell\le m-1\), the jump of the derivative in \eqref{eq:parity-extension-derivative-formula} is
\begin{equation}
\big(\partial_r^q\partial_z^\ell\mathcal E_p f\big)(0^+,z)
-
\big(\partial_r^q\partial_z^\ell\mathcal E_p f\big)(0^-,z)
=
\left(1-(-1)^{q+p}\right)
\partial_r^q\partial_z^\ell f(0,z)
=0.
\label{eq:parity-extension-trace-jump}
\end{equation}
\noindent If \(q+p\) is even, the factor in parentheses vanishes.  If \(q+p\) is odd, the trace vanishes by \eqref{eq:parity-extension-trace-condition}.  Thus every derivative of total order at most \(m-1\) has matching traces across the axis.  This prevents the next normal derivative from creating a distribution supported on \(r=0\).

\hfill \\

\noindent\textbf{Weak gluing across the axis.}  We denote the piecewise function on the right side of \eqref{eq:parity-extension-derivative-formula} by \(G_{q,\ell}\).  For a test function \(\varphi\in C_c^\infty(\mathbb R^2)\), integration by parts on the two half-planes gives, whenever \(q+\ell\le m-1\),
\begin{align}
-\int_{\mathbb R^2}G_{q,\ell}\,\partial_r\varphi\,dr\,dz
&=
\int_{\mathbb R^2}G_{q+1,\ell}\varphi\,dr\,dz
+\int_{\mathbb R}
\left(G_{q,\ell}(0^+,z)-G_{q,\ell}(0^-,z)\right)
\varphi(0,z)\,dz
\nonumber\\
&=
\int_{\mathbb R^2}G_{q+1,\ell}\varphi\,dr\,dz.
\label{eq:parity-extension-weak-gluing}
\end{align}
\noindent The tangential derivative produces no interface term and satisfies \(\partial_zG_{q,\ell}=G_{q,\ell+1}\) weakly.  Starting from \(G_{0,0}=\mathcal E_p f\) and iterating these two identities proves that every \(G_{q,\ell}\) with \(q+\ell\le m\) is the corresponding weak derivative.  Hence \(\mathcal E_p f\in H^m(\mathbb R^2)\). \\ 

\noindent The same argument applies directly to general \(f\in H^m(\mathbb{D})\).  The Sobolev trace theorem gives the traces in equation \eqref{eq:parity-extension-trace-condition}, and the Green formula on each half-plane gives \eqref{eq:parity-extension-weak-gluing}.  The interface term again vanishes, so the piecewise formulas define all weak derivatives through order \(m\).

\hfill \\ 

\noindent\textbf{Weighted norm identity.}  Fix \(q+\ell\le m\).  By \eqref{eq:parity-extension-derivative-formula} and the even reflection of the weight,
\begin{align}
\left\|\partial_r^q\partial_z^\ell\mathcal E_p f\right\|_{L^2(\Phi^{\rm even},\mathbb R^2)}^2
&=
\int_0^\infty\!\int_{\mathbb R}
\left|\partial_r^q\partial_z^\ell f(r,z)\right|^2\Phi(r,z)\,dz\,dr
+
\int_{-\infty}^0\!\int_{\mathbb R}
\left|\partial_r^q\partial_z^\ell f(-r,z)\right|^2\Phi(-r,z)\,dz\,dr
\nonumber\\
&=
2\left\|\partial_r^q\partial_z^\ell f\right\|_{L^2(\Phi,\mathbb{D})}^2,
\label{eq:parity-extension-weighted-norm-proof}
\end{align}
\noindent where the last equality uses the change of variables \(\rho=-r\) in the negative-half-plane integral.  This proves the identity \eqref{eq:parity-extension-weighted-norm}.  Taking \(\Phi\equiv1\) and summing over all \(q+\ell\le m\) gives \eqref{eq:parity-extension-Hm-norm}.
\end{grayproof}

\hfill \\

\noindent In the estimates, this proposition is applied before every whole-plane pointwise or product bound.  After the estimate is obtained on $\mathbb R^2$, the conversion back to the half-plane uses the explicit identities
\begin{equation}
\|\mathrm E_{\rm even}f\|_{H^m(\mathbb R^2)}=\sqrt2\|f\|_{H^m(\mathbb{D})},
\qquad
\|D_r^qD_z^\ell\mathrm E_{\rm even}f\|_{L^2(\Phi^{\rm even},\mathbb R^2)}
=\sqrt2\|D_r^qD_z^\ell f\|_{L^2(\Phi,\mathbb{D})}.
\end{equation}

\hfill \\ 

\clearpage 

\begin{proposition}[Whole-plane Sobolev and interpolation constants]\label{prop:whole-plane-standard-constants}
With the Bessel norm fixed in \eqref{eq:Bessel-norm-prelim}, for $s>1$,
\begin{equation}\label{eq:sharp-Hs-Linfty-R2}
\|F\|_{L^\infty(\mathbb R^2)}
\le S_{s,2}\|F\|_{H^s_B(\mathbb R^2)},
\qquad
S_{s,2}=\frac{1}{2\sqrt\pi}\left(\frac{\Gamma(s-1)}{\Gamma(s)}\right)^{1/2}
=\frac{1}{2\sqrt{\pi(s-1)}}.
\end{equation}
This is \Cref{thm:mp-sharp-endpoint} with dimension $d=2$, exponent $r=\infty$, and Sobolev order $n=s$.  \\

\noindent Moreover, for integers $0<j<k$,
\begin{equation}\label{eq:whole-plane-Hilbert-interpolation}
\|D^jF\|_{L^2(\mathbb R^2)}
\le \|D^kF\|_{L^2(\mathbb R^2)}^{j/k}\|F\|_{L^2(\mathbb R^2)}^{1-j/k},
\end{equation}
and hence, for every $\eta>0$,
\begin{equation}\label{eq:whole-plane-Young-interpolation}
\|D^jF\|_{L^2}
\le \eta\|D^kF\|_{L^2}+C_{j,k}^{\rm int}(\eta)\|F\|_{L^2},
\qquad
C_{j,k}^{\rm int}(\eta)=\frac{k-j}{k}\left(\frac{j}{k\eta}\right)^{j/(k-j)}.
\end{equation}
\end{proposition}

\hfill 

\begin{grayproof}
\textbf{Sharp Bessel-potential Sobolev bound.}  The sharp constant is the two-dimensional endpoint constant from \Cref{thm:mp-sharp-endpoint}.  The computation below verifies the same constant with the Fourier convention~\eqref{eq:fourier-convention}.  From the inverse Fourier formula,
\begin{equation}
F(x)=\frac{1}{(2\pi)^2}\int_{\mathbb R^2}e^{ix\cdot\xi}\widehat F(\xi)\,d\xi.
\end{equation}

\hfill 

\noindent Insert the factors $(1+|\xi|^2)^{-s/2}$ and $(1+|\xi|^2)^{s/2}$ and apply the Cauchy--Schwarz inequality:
\begin{align}
|F(x)|
&\le \frac{1}{(2\pi)^2}\int_{\mathbb R^2}(1+|\xi|^2)^{-s/2}(1+|\xi|^2)^{s/2}|\widehat F(\xi)|\,d\xi\nonumber\\
&\le \frac{1}{(2\pi)^2}\left(\int_{\mathbb R^2}(1+|\xi|^2)^{-s}\,d\xi\right)^{1/2}
\left(\int_{\mathbb R^2}(1+|\xi|^2)^s|\widehat F(\xi)|^2\,d\xi\right)^{1/2}.
\end{align}

\hfill 

\noindent By the definition of the Bessel norm, the second factor is $(2\pi)\|F\|_{H^s_B}$.  Therefore
\begin{equation}
|F(x)|
\le \frac{1}{2\pi}\left(\int_{\mathbb R^2}(1+|\xi|^2)^{-s}\,d\xi\right)^{1/2}
\|F\|_{H^s_B}.
\end{equation}
The integral is finite exactly for $s>1$.  In polar coordinates,
\begin{align}
\int_{\mathbb R^2}(1+|\xi|^2)^{-s}\,d\xi
&=2\pi\int_0^\infty(1+\rho^2)^{-s}\rho\,d\rho
=\pi\int_1^\infty u^{-s}\,du
=\frac{\pi}{s-1}.
\end{align}

\hfill 

\noindent Substituting this value gives
\begin{equation}
|F(x)|\le \frac{1}{2\sqrt{\pi(s-1)}}\|F\|_{H^s_B}.
\end{equation}
Taking the supremum in $x$ proves \eqref{eq:sharp-Hs-Linfty-R2}.

\hfill \\

\noindent\textbf{Interpolation between Sobolev levels.}  Let $0<j<k$ and set $\theta=j/k$.  By Plancherel, up to the same harmless Fourier-normalization factor on both sides,
\begin{align}
\|D^jF\|_{L^2}^2
&=\int_{\mathbb R^2}|\xi|^{2j}|\widehat F(\xi)|^2\,d\xi =\int_{\mathbb R^2}\bigl(|\xi|^{2k}|\widehat F(\xi)|^2\bigr)^\theta
\bigl(|\widehat F(\xi)|^2\bigr)^{1-\theta}\,d\xi.
\end{align}

\hfill 

\noindent Hölder's inequality with exponents $1/\theta$ and $1/(1-\theta)$ gives
\begin{equation}
\|D^jF\|_{L^2}^2
\le
\left(\int |\xi|^{2k}|\widehat F|^2\right)^\theta
\left(\int |\widehat F|^2\right)^{1-\theta}.
\end{equation}
Taking square roots proves \eqref{eq:whole-plane-Hilbert-interpolation}.  \\

\noindent To obtain the additive form, we use Young's inequality with a local Young parameter \(\chi>0\):
\begin{equation}
a^\theta b^{1-\theta}
\le \theta\chi a+(1-\theta)\chi^{-\theta/(1-\theta)}b,
\qquad a,b\ge0.
\end{equation}
With $a=\|D^kF\|_{L^2}$, $b=\|F\|_{L^2}$, and $\theta=j/k$, we choose the local Young parameter by
\begin{equation}
\qquad
\chi=\frac{\eta}{\theta}.
\end{equation}
Thus the first Young term is
\begin{equation}
\qquad
\theta\chi a=\eta\|D^kF\|_{L^2}.
\end{equation}
The active second Young term is
\begin{equation}
(1-\theta)\chi^{-\theta/(1-\theta)}b
=\frac{k-j}{k}\left(\frac{j}{k\eta}\right)^{j/(k-j)}\|F\|_{L^2}
=C_{j,k}^{\rm int}(\eta)\|F\|_{L^2}.
\end{equation}
This gives \eqref{eq:whole-plane-Young-interpolation}.
\end{grayproof}

\hfill 

\noindent The pointwise and interpolation estimates are used in two distinct places in the energy argument.  First, after parity reflection, \eqref{eq:sharp-Hs-Linfty-R2} with \(s=2\) gives \(L^\infty\) control of any lower-order factor once two additional derivatives of that factor are controlled in \(L^2\).  This is the mechanism used to place coefficients such as lower-order derivatives of \(\bvar{\delta v_i}\), velocity factors, and non-top-order product factors in \(L^\infty\). Second, when a commutator produces an intermediate derivative order \(0<j<k\), \eqref{eq:whole-plane-Young-interpolation} rewrites that intermediate norm as
\begin{equation}
\|D^jF\|_{L^2}
\le
\eta\|D^kF\|_{L^2}
+
C_{j,k}^{\rm int}(\eta)\|F\|_{L^2}.
\end{equation}
The first term is a small multiple of the top-order norm and is absorbed into the damping or into the chosen top-order side of the estimate.  The second term is a low-order contribution and is bounded by the low-order energy.  Here \(D^jF\) denotes the aggregate order-\(j\) derivative norm, equivalently the sum over all multi-indices \(|\alpha|=j\), with only dimension-dependent constants changed.  \\

\noindent In \eqref{eq:whole-plane-Hilbert-interpolation}--\eqref{eq:whole-plane-Young-interpolation}, 
\(\|D^jF\|_{L^2}\) denotes the order-\(j\) derivative norm defined in Fourier variables by the multiplier \(|\xi|^j\), whereas the energy uses sums of the individual order-\(j\) coordinate-derivative norms. These quantities are related explicitly by Plancherel's theorem. Indeed, note that
\begin{equation}
|\xi|^{2j}
=
(\xi_1^2+\xi_2^2)^j
=
\sum_{q=0}^j
\binom{j}{q}
\xi_1^{2q}\xi_2^{2(j-q)}.
\end{equation}
Therefore, Plancherel's theorem gives
\begin{equation}
\|D^jF\|_{L^2}^2
=
\sum_{q=0}^j
\binom{j}{q}
\|\partial_r^q\partial_z^{j-q}F\|_{L^2}^2.
\end{equation}
Hence, by the weighted Cauchy--Schwarz inequality,
\begin{equation}
\sum_{q=0}^j
\|\partial_r^q\partial_z^{j-q}F\|_{L^2}
\le
\left(
\sum_{q=0}^j
\binom{j}{q}^{-1}
\right)^{1/2}
\|D^jF\|_{L^2},
\end{equation}
while, at order \(k\),
\begin{equation}
\|D^kF\|_{L^2}
\le
\sqrt{\binom{k}{\lfloor k/2\rfloor}}
\sum_{q=0}^k
\|\partial_r^q\partial_z^{k-q}F\|_{L^2}.
\end{equation}
Thus \eqref{eq:whole-plane-Young-interpolation} also yields the corresponding estimate for the coordinate-derivative sums: the factors above are incorporated into the interpolation constant, and the Young parameter is chosen so that the coefficient of the top-order coordinate sum is the prescribed \(\eta\).

\hfill \\ 

\hfill

\subsection{Leibniz product estimate for differentiated nonlinearities}\label{subsec:leibniz-product-estimate}

\noindent After Leibniz' rule is applied, each differentiated nonlinearity is a finite sum of products.  The estimate below uses two weighted \(L^2\) slots: one for the tested factor and one for the perturbative factor carrying the most derivatives.  Every remaining perturbative factor is lower order and is included, together with the required weight conversion, in one \(L^\infty\) multiplier.  Products containing no perturbative spatial factor are handled separately by a certified weighted norm of the fixed profile.

\hfill

\begin{proposition}[Leibniz product bound used in the energy estimates]\label{prop:whole-plane-product}
\noindent Let \(k\ge4\) and set
\begin{equation}\label{eq:explicit-Kinfty}
K_\infty:=3\sqrt2\,S_{2,2}=\frac{3}{\sqrt{2\pi}}.
\end{equation}
Let \(|\alpha|\le k\), and consider one Leibniz summand
\begin{equation}\label{eq:one-leibniz-summand-for-product}
C_\sigma\prod_{m=0}^M D^{\sigma_m}F_m,
\qquad
\sigma_0+\cdots+\sigma_M=\alpha,
\qquad
C_\sigma=\frac{\alpha!}{\sigma_0!\cdots\sigma_M!}.
\end{equation}
Let \(\mathcal P\subseteq\{0,\ldots,M\}\) be the nonempty set of indices for the perturbative spatial factors, and assume that these factors satisfy the parity and regularity hypotheses of Proposition~\ref{prop:whole-plane-reflection}.  The remaining factors are fixed smooth profiles or spatially constant coefficients with independent pointwise bounds. \\

\noindent Choose \(m_*\in\mathcal P\) so that
\begin{equation}
|\sigma_{m_*}|=\max_{m\in\mathcal P}|\sigma_m|,
\end{equation}
and let \(\wgray{\Phi_*}\) be the certified weight for \(D^{\sigma_{m_*}}F_{m_*}\). \\

\noindent For an output weight \(\wgray{\Phi_i}\), define
\begin{equation}
\mathcal M^{(i)}_{\sigma,m_*}:=
\Bigg\|
\Bigg( \ \prod_{\substack{0\le m\le M\\m\ne m_*}}D^{\sigma_m}F_m \ \Bigg)
\left(\frac{\wgray{\Phi_i}}{\wgray{\Phi_*}}\right)^{1/2}
\Bigg\|_{L^\infty(\mathbb D)}.
\label{eq:explicit-combined-product-multiplier}
\end{equation}

\noindent Then every test factor \(H\) measured in \(\wgray{\Phi_i}\) satisfies
\begin{align}
&\left|
\int_\mathbb{D} C_\sigma
\left(\prod_{m=0}^M D^{\sigma_m}F_m\right)
H\,\wgray{\Phi_i}\,dr\,dz
\right| \, \le \, 
C_\sigma\,\mathcal M^{(i)}_{\sigma,m_*}
\|D^{\sigma_{m_*}}F_{m_*}\|_{\wgray{\Phi_*}}
\|H\|_{\wgray{\Phi_i}}.
\label{eq:explicit-multifactor-certificate}
\end{align}
 
  \end{proposition}

\vspace{4mm}

\begin{grayproof}
Fix \(m_*\in\mathcal P\) with maximal derivative order.  This choice guarantees that every other perturbative factor has sufficiently few derivatives to be estimated pointwise.  Indeed, for \(m\in\mathcal P\setminus\{m_*\}\),
\begin{equation}
2|\sigma_m|
\le |\sigma_m|+|\sigma_{m_*}|
\le \sum_{q=0}^M|\sigma_q|
=|\alpha|
\le k.
\label{eq:product-maximal-factor-order-count}
\end{equation}
Hence \(|\sigma_m|\le\lfloor k/2\rfloor\).  Since \(k\ge4\), this also gives \(|\sigma_m|+2\le k\), so the two additional derivatives required by the two-dimensional \(H^2\to L^\infty\) estimate are available at energy level \(k\). \\

\noindent For each unselected perturbative factor, let \(G_m\) be the parity reflection of \(D^{\sigma_m}F_m\) supplied by Proposition~\ref{prop:whole-plane-reflection}.  Proposition~\ref{prop:whole-plane-standard-constants}, the comparison between the Bessel-potential and derivative Sobolev norms, and the reflection identity give
\begin{align}
\|D^{\sigma_m}F_m\|_{L^\infty(\mathbb{D})}
&\le S_{2,2}\|G_m\|_{H_B^2(\mathbb R^2)}
\le 3S_{2,2}\|G_m\|_{H^2(\mathbb R^2)}
=3\sqrt2\,S_{2,2}
\bigg( \ \sum_{|\gamma|\le2}
\|D^{\sigma_m+\gamma}F_m\|_{L^2(\mathbb{D})}^2 \ \bigg)^{1/2}
\nonumber\\
&=K_\infty
\bigg( \ \sum_{|\gamma|\le2}
\|D^{\sigma_m+\gamma}F_m\|_{L^2(\mathbb{D})}^2 \ \bigg)^{1/2}.
\label{eq:explicit-Linfty-factor}
\end{align}
Thus every unselected perturbative factor is controlled in \(L^\infty\).  The fixed profiles and spatially constant coefficients are controlled by their separately certified pointwise bounds. \\

\noindent The weighted allocation is now explicit:
\begin{align}
&\left(\prod_{m=0}^M D^{\sigma_m}F_m\right)H\wgray{\Phi_i} =
\Bigg[
\bigg( \ \prod_{m\ne m_*}D^{\sigma_m}F_m \bigg)
\left(\frac{\wgray{\Phi_i}}{\wgray{\Phi_*}}\right)^{1/2}
\Bigg]
\left[D^{\sigma_{m_*}}F_{m_*}\sqrt{\wgray{\Phi_*}}\right]
\left[H\sqrt{\wgray{\Phi_i}}\right].
\label{eq:product-three-factors-defined}
\end{align}

\hfill 

\noindent Applying H\"older's inequality in \(L^\infty\times L^2\times L^2\) therefore yields
\begin{align}
\bigg|
\int_\mathbb{D} C_\sigma
\bigg(\prod_{m=0}^M D^{\sigma_m}F_m\bigg)
H\,\wgray{\Phi_i}\,dr\,dz
\bigg|
& \le
C_\sigma
\bigg\|
\bigg(\prod_{m\ne m_*}D^{\sigma_m}F_m \bigg)
\left(\frac{\wgray{\Phi_i}}{\wgray{\Phi_*}}\right)^{1/2}
\bigg\|_{L^\infty(\mathbb D)} \!
\|D^{\sigma_{m_*}}F_{m_*}\|_{\wgray{\Phi_*}}
\|H\|_{\wgray{\Phi_i}}
\nonumber\\
\qquad & =
C_\sigma\,\mathcal M^{(i)}_{\sigma,m_*}
\|D^{\sigma_{m_*}}F_{m_*}\|_{\wgray{\Phi_*}}
\|H\|_{\wgray{\Phi_i}}.
\label{eq:product-holder-all-steps}
\end{align}

  \noindent This is \eqref{eq:explicit-multifactor-certificate}.  The full ratio \((\wgray{\Phi_i}/\wgray{\Phi_*})^{1/2}\) remains with the pointwise factors, which is exactly why no separate supremum estimate for the weight ratio is needed.
\end{grayproof}

\noindent Summing \eqref{eq:explicit-multifactor-certificate} over the finitely many Leibniz splits gives the bound for the full differentiated product.  If the selected factor has an intermediate derivative order, we use \eqref{eq:ns-lower-from-top} before adding the term to the final constant.  Streamfunction factors are first bounded in terms of \(\bvar{\delta\omega}\) by \Cref{thm:stream-elliptic-used}.

\clearpage 

\subsection{Streamfunction elliptic estimates}\label{subsec:streamfunction-elliptic-estimates}

\noindent This subsection constructs the streamfunction estimates needed later in the energy argument and then uses the corresponding point-value bounds to control the modulation parameters. 

\hfill

\noindent The operator
\begin{equation}
\partial_r^2+\frac{3}r\partial_r+\partial_z^2
\end{equation}

\noindent is the five-dimensional Laplacian on functions that are radial in four Euclidean variables and depend on the fifth variable \(z\).  Thus the streamfunction equation can be lifted to a Poisson equation on \(\mathbb R^5\), represented by the five-dimensional Newtonian kernel, and then reduced back to the two variables~\((r,z)\).  All estimates below are written in the reduced flat measure \(dr\,dz\), which is the measure used in the energy \(\mathfrak E_k\).

\hfill \\

\begin{proposition}[Five-dimensional kernel representation]\label{prop:five-dimensional-kernel}
\noindent Let \(\Omega_5\in C_c^\infty(\mathbb R^5)\) be radial in the first four variables, so that
\begin{equation}
\Omega_5(y,z)=\omega(|y|,z),
\qquad y\in\mathbb R^4.
\label{eq:five-dimensional-radial-source}
\end{equation}
\noindent Suppose that \(\psi\) solves
\begin{equation}
-\left(\partial_r^2+\frac{3}r\partial_r+\partial_z^2\right)\psi=\omega
\qquad (r>0),
\label{eq:reduced-streamfunction-poisson}
\end{equation}
\noindent and that the radial lift \(\Psi(y,z):=\psi(|y|,z)\) is \(C^2\) across \(y=0\) and tends to zero at infinity. \\

\noindent Then
\begin{equation}
\Delta_{y,z}\Psi
=\partial_r^2\psi+\frac{3}r\partial_r\psi+\partial_z^2\psi,
\label{eq:radial-five-dimensional-laplacian}
\end{equation}
\noindent with the identity understood distributionally at the axis.  Hence
\begin{equation}
-\Delta_{\mathbb R^5}\Psi=\Omega_5.
\label{eq:five-dimensional-poisson}
\end{equation}
\noindent Since
\begin{equation}
\Gamma_5(X)=\frac{1}{8\pi^2}|X|^{-3}
\label{eq:Gamma-five}
\end{equation}
\noindent is the fundamental solution of \(-\Delta_{\mathbb R^5}\), the decaying solution is
\begin{equation}
\Psi=\Gamma_5*\Omega_5.
\label{eq:five-dimensional-convolution}
\end{equation}
\end{proposition}

\begin{grayproof}
\noindent\textbf{Five-dimensional lifting.}  Let \(r=|y|\) with \(y\in\mathbb R^4\).  For \(r>0\), the radial Laplacian formula in four dimensions gives
\begin{equation}
\Delta_y\Psi
=\partial_r^2\psi+\frac{4-1}{r}\partial_r\psi
=\partial_r^2\psi+\frac{3}r\partial_r\psi.
\label{eq:four-dimensional-radial-laplacian}
\end{equation}
\noindent Adding \(\partial_z^2\Psi=\partial_z^2\psi\) proves \eqref{eq:radial-five-dimensional-laplacian} for \(r>0\).  \\

\noindent The assumed \(C^2\) regularity of the radial lift across \(y=0\) implies \(\partial_r\psi(0,z)=0\).  Consequently, $(3/r)\partial_r\psi$ extends regularly across $r=0$ in the weak sense. In particular, no additional singular term concentrated on the axis $y=0$ appears when the equation is interpreted distributionally. Hence the reduced equation \eqref{eq:reduced-streamfunction-poisson}, initially written away from the axis, extends across $y=0$ and gives \eqref{eq:five-dimensional-poisson} on all of $\mathbb R^5$.

\hfill \\

\noindent\textbf{Normalization of the Newtonian kernel.}  In dimension \(n\ge3\), the fundamental solution of \(-\Delta\) is
\begin{equation}
\Gamma_n(X)
=\frac{1}{n(n-2)\alpha(n)}|X|^{2-n},
\label{eq:fundamental-solution-general-n}
\end{equation}
\noindent where \(\alpha(n)\) is the volume of the unit ball in \(\mathbb R^n\) (see \citep{evans2010partial}).\\

\noindent This normalization satisfies \(-\Delta\Gamma_n=\delta_0\).  For \(n=5\),
\begin{equation}
5(5-2)\alpha(5)
=15\frac{\pi^{5/2}}{\Gamma(7/2)}
=8\pi^2,
\label{eq:Gamma-five-normalization-check}
\end{equation}
\noindent which gives the kernel \(\Gamma_5(X)=(8\pi^2)^{-1}|X|^{-3}\) in \eqref{eq:Gamma-five}. 

\hfill \\ 

\hfill 

\noindent\textbf{Convolution representation.}  Define
\begin{equation}
U(X):=(\Gamma_5*\Omega_5)(X)
=\int_{\mathbb R^5}\Gamma_5(X-Y)\Omega_5(Y)\,dY.
\label{eq:five-dimensional-potential-definition}
\end{equation}
\noindent Since \(\Omega_5\in C_c^\infty(\mathbb R^5)\), this convolution is well defined and decays at infinity.\\

\noindent For every \(\varphi\in C_c^\infty(\mathbb R^5)\), Fubini's theorem and \(-\Delta\Gamma_5=\delta_0\) give
\begin{align}
\int_{\mathbb R^5}U(X)(-\Delta\varphi)(X)\,dX
&=
\int_{\mathbb R^5}\Omega_5(Y)
\left(
\int_{\mathbb R^5}\Gamma_5(X-Y)(-\Delta\varphi)(X)\,dX
\right)dY
\end{align}
so
\begin{align}
\int_{\mathbb R^5}U(X)(-\Delta\varphi)(X)\,dX
=
\int_{\mathbb R^5}\Omega_5(Y)\varphi(Y)\,dY.
\label{eq:five-dimensional-potential-weak-identity}
\end{align}
\noindent Hence \(-\Delta U=\Omega_5\) in the distributional sense.  The lifted function \(\Psi\) satisfies the same equation by \eqref{eq:five-dimensional-poisson}.  Therefore, $\Psi-U$ satisfies
\begin{equation}
-\Delta_5(\Psi-U)=0
\end{equation}
on all of $\mathbb R^5$, so $\Psi-U$ is harmonic. Since both $\Psi$ and $U$ vanish at infinity, their difference also tends to zero at infinity. By Liouville's theorem, a harmonic function on $\mathbb R^5$ that tends to zero at infinity must be identically zero. Hence $\Psi=U$, which proves \eqref{eq:five-dimensional-convolution}.
\end{grayproof}

\hfill \\

\noindent We now identify the finite family of streamfunction expressions that must be controlled before stating the elliptic estimate. \\ 

\noindent The finite family of streamfunction expressions needed in the later linear and nonlinear terms is
\begin{equation}
\label{eq:ns-stream-family}
\mathcal S_\psi:=\left\{
\partial_{zz}, \,
\partial_{zr}, \,
 r\partial_{zz}, \,
 r\partial_{zr}, \,
 3\partial_r+r\partial_{rr}, \,
 r\partial_z, \,
 2+r\partial_r
\right\}.
\end{equation}
If \(T\in\mathcal S_\psi\), then \(T[\varphi]\) denotes the displayed differential expression applied to \(\varphi\).  For example,
\begin{equation}
    (2+r\partial_r)[\varphi]=2\varphi+r\partial_r\varphi.
\end{equation}

\noindent The finite index set \(\mathcal J_{\rm ell}\) consists of the triples \((T,i,\alpha)\), with \(T\in\mathcal S_\psi\), \(i\in\{\omega,r,z\}\), and \(|\alpha|\le k\).  The component index \(i\) specifies the output weight in which the term is estimated.  Once \(k\) is fixed, \(\mathcal J_{\rm ell}\) is finite, so every constant below is part of a finite certificate.

\hfill

\begin{proposition}[Streamfunction elliptic estimates]\label{thm:stream-elliptic-used}
For the elliptic estimates, assume \(k\ge2\).  Suppose that the point functionals and the finitely many weighted operators appearing below have been certified by the construction given after the statement.  Then the selected origin constants satisfy
\begin{align}
|\bvar{\delta\psi}(0,0)|
&\le P_{\psi,0}^{0}\mathfrak E_0+P_{\psi,0}^{k}\mathfrak H_k
\le P_{\psi,0}\mathfrak E_k,\label{eq:ell-point-theorem-psi}\\
|\partial_z\bvar{\delta\psi}(0,0)|
&\le P_{\psi_z,0}^{0}\mathfrak E_0+P_{\psi_z,0}^{k}\mathfrak H_k
\le P_{\psi_z,0}\mathfrak E_k.\label{eq:ell-point-theorem-psiz}
\end{align}
where
\begin{equation}\label{eq:ell-point-aggregate}
P_{\psi,0}:=P_{\psi,0}^{0}+\mu_k^{-1/2}P_{\psi,0}^{k},
\quad \ 
P_{\psi_z,0}:=P_{\psi_z,0}^{0}+\mu_k^{-1/2}P_{\psi_z,0}^{k}.
\end{equation}

\hfill 

\noindent The differentiated unweighted estimate is
\begin{equation}
\label{eq:ns-Cell-homogeneous}
\sum_{m=2}^{k+2}\sum_{|\alpha|=m}
\|D^\alpha\bvar{\delta\psi}\|_{L^2(\mathbb{D})}
\le C_{\rm ell}^{\rm diff}
\sum_{m=0}^{k}\sum_{|\alpha|=m}
\|D^\alpha\bvar{\delta\omega}\|_{L^2(\mathbb{D})}.
\end{equation}

\noindent The two low-order weighted Hessian estimates are
\begin{equation}
\label{eq:ns-Cwell}
\|\partial_{zz}\bvar{\delta\psi}\|_{\wgray{\Phi_z}}
\le C^{\rm w\text{-}ell}_{z,zz}\mathfrak E_k,
\qquad
\|\partial_{zr}\bvar{\delta\psi}\|_{\wgray{\Phi_r}}
\le C^{\rm w\text{-}ell}_{r,zr}\mathfrak E_k.
\end{equation}
Finally, for every \((T,i,\alpha)\in\mathcal J_{\rm ell}\), with \(m=|\alpha|\),
\begin{equation}
\label{eq:ns-Cwell-high}
\|D^\alpha T[\bvar{\delta\psi}]\|_{\wgray{\Phi_{i,m}}}
\le C^{\rm w\text{-}ell}_{T,i,m,\alpha}\mathfrak E_k.
\end{equation}
\end{proposition}

\hfill \\

\noindent The next subsections define and construct the constants needed in \Cref{thm:stream-elliptic-used} and give the certification procedure for each one.  

\hfill \\

\begin{remark*}
Proposition~\ref{prop:five-dimensional-kernel} gives the Newtonian representation for smooth, compactly supported, axis-regular sources.  We extend this representation to an arbitrary source \(\bvar{\delta\omega}\) in the perturbation class by approximation.  Using parity-preserving mollification followed by cutoff, and leveraging the regularity and decay assumptions of the perturbations, we can choose smooth axis-regular compactly supported functions \(\omega_n\) that converge to \(\bvar{\delta\omega}\) in every weighted source norm appearing in \Cref{thm:stream-elliptic-used}. For each \(n\), let \(\psi_n\) be the reduced five-dimensional Newtonian potential generated by~\(\omega_n\).  Since the Newtonian map is linear, each certified estimate in \Cref{thm:stream-elliptic-used}, applied to \(\omega_n-\omega_m\), shows that the corresponding quantities associated with \(\psi_n-\psi_m\) are Cauchy in their stated norms, which therefore determines a single streamfunction \(\bvar{\delta\psi}\).  Passing to the limit in the estimates yields exactly the same bounds and certified constants for \(\bvar{\delta\psi}\).  Passing to the limit in the equations \(-\mathcal E\psi_n=\omega_n\) after testing against smooth compactly supported functions, we get $-\mathcal E\bvar{\delta\psi}=\bvar{\delta\omega}$ in the sense of distributions on $\mathbb{D}$. The parity, axis regularity, and decay required of the perturbation streamfunction pass from \(\psi_n\) to \(\bvar{\delta\psi}\).  Finally, the limit is unique: the difference of two such solutions lifts to a decaying harmonic function on \(\mathbb R^5\), and therefore vanishes.  Thus the decaying Newtonian inverse is well defined on the full perturbation class, and every estimate in \Cref{thm:stream-elliptic-used} extends from compactly supported sources with no change in its constant.
\end{remark*}

\clearpage

\subsubsection{Certification setup and singular-integral conventions}\label{subsec:obtaining-elliptic-constants}

\noindent We now define and construct the constants appearing in \Cref{thm:stream-elliptic-used}. The goal is to turn every streamfunction quantity into a finite list of weighted estimates for \(\bvar{\delta\omega}\).  The perturbation streamfunction is the axis-regular decaying solution of
\begin{equation}\label{eq:ell-signed-equation}
-\mathcal{E}\bvar{\delta\psi}=\bvar{\delta\omega},
\qquad \quad 
\mathcal{E}:=\partial_r^2+\frac{3}r\partial_r+\partial_z^2.
\end{equation}

\hfill 

\noindent The guiding principle is simple: we rewrite each occurrence of \(\bvar{\delta\psi}\) as either a scalar functional of \(\bvar{\delta\omega}\) or an integral operator acting on a derivative of \(\bvar{\delta\omega}\).  Then we bound that functional or operator in the relevant weighted norm and convert the remaining source norm to \(\mathfrak E_0\), \(\mathfrak H_k\), or \(\mathfrak E_k\). 

\hfill \\

\paragraph{Singular-integral convention.}
\noindent Let \(U\subset\mathbb R^d\) be an open set, and let \(L(x,y)\) be smooth away from the diagonal \(x=y\).  It is precisely on this set that kernels obtained by differentiating the Newtonian potential may fail to be locally integrable.  The principal-value operator associated with \(L\) is
\begin{equation}\label{eq:ell-pv-general-definition}
(\operatorname{p.v.}\mathsf T_Lf)(x)
:=\lim_{\eta\downarrow0}
\int_{\{y\in U:\ |x-y|>\eta\}}L(x,y)f(y)\,dy.
\end{equation}
The limit is taken in the sense used by the estimate of interest: pointwise for pointwise identities, distributionally for kernel identities, and in the relevant weighted \(L^2\) norm for operator estimates.  When the kernel is locally integrable, this is just the ordinary integral operator.  When the kernel is singular at $x=y$, the integral is understood in the principal value sense by excluding the ball $|x-y|\le\eta$ and then letting $\eta\to0$.\\

\noindent Whenever a reduced singular integral is derived from a five-dimensional principal-value integral, the principal value is defined in the original five-dimensional variables: one first removes a ball of radius $\eta$ around the singular point in $\mathbb R^5$, then performs the angular integration, and finally lets $\eta\to0$. Thus the reduced operator inherits its principal-value prescription from the five-dimensional integral. This choice determines the finite local term that appears in the second-derivative formulas below, and the same convention is used for all reduced principal-value operators.  \\

\noindent The second convention concerns derivatives of kernels that are singular on the diagonal.  The superscript \(\mathrm{off}\) denotes the derivative taken only off the diagonal.  Thus we first regard the kernel as an ordinary smooth function on the set \(x\ne y\), differentiate there, and only afterwards interpret the resulting singular kernel as a principal-value operator.  For example,
\begin{equation}\label{eq:ell-pv-notation-intro}
\bigl(\operatorname{p.v.}\,\partial_{X_iX_j}^{\mathrm{off}}\Gamma_5*F\bigr)(X)
:=\lim_{\eta\downarrow0}
\int_{|X-Y|>\eta}\partial_{X_iX_j}^{\mathrm{off}}\Gamma_5(X-Y)F(Y)\,dY.
\end{equation}

\noindent  Here \(\partial_{X_iX_j}^{\mathrm{off}}\Gamma_5\) is the ordinary Hessian of \(\Gamma_5(X)=(8\pi^2)^{-1}|X|^{-3}\) computed for $X\neq0$, or equivalently for $X\neq Y$ when $\Gamma_5(X-Y)$ appears in a convolution. Because $\Gamma_5$ is singular at the origin, this pointwise Hessian does not represent the full second derivative in the sense of distributions. When the singularity at $X=Y$ is taken into account, an additional local term appears. Thus the distributional identity consists of two parts: a principal-value singular integral coming from the off-diagonal Hessian, and a Dirac mass supported at the diagonal. More precisely, the local term is $-\delta_{ij}\delta_0/5$. We keep these two contributions separate because they enter the estimates differently: the principal-value term is treated as a singular integral, while the Dirac mass contributes directly as a multiple of the source.

\hfill 

\subsubsection{Constants to be certified}\label{sssec:ell-outputs}

\noindent \Cref{thm:stream-elliptic-used} requires different groups of constants.  They control the streamfunction at the origin, the two low-order weighted Hessians, the high-order weighted terms, and the differentiated streamfunctions in unweighted \(L^2\).  The certificate stores
\begin{equation}\label{eq:ell-output-list}
P_{\psi,0}^{0},\quad P_{\psi,0}^{k},\quad P_{\psi_z,0}^{0},\quad P_{\psi_z,0}^{k},\quad C_{\rm ell}^{\rm diff},\quad C^{\rm w\text{-}ell}_{z,zz},\quad C^{\rm w\text{-}ell}_{r,zr},\quad C^{\rm w\text{-}ell}_{T,i,m,\alpha}
\ \ \text{for }(T,i,\alpha)\in\mathcal J_{\rm ell}.
\end{equation}
It also stores the constants \(\mathcal K_i^0(a,T)\) and \(\mathcal K_i^k(a,T)\) for profile-multiplied streamfunction terms. \\

\paragraph{Point outputs.}
\noindent The point constants certify
\begin{align}
|\bvar{\delta\psi}(0,0)|
&\le P_{\psi,0}^{0}\mathfrak E_0+P_{\psi,0}^{k}\mathfrak H_k,
\label{eq:ell-point-output-psi}\\
|\partial_z\bvar{\delta\psi}(0,0)|
&\le P_{\psi_z,0}^{0}\mathfrak E_0+P_{\psi_z,0}^{k}\mathfrak H_k.
\label{eq:ell-point-output-psiz}\end{align}

\hfill 

\paragraph{Weighted outputs.}
\noindent The weighted operator constants certify
\begin{align}
\|\partial_{zz}\bvar{\delta\psi}\|_{\wgray{\Phi_z}}
& \le C^{\rm w\text{-}ell}_{z,zz}\mathfrak E_k,
\qquad  \quad 
\|\partial_{zr}\bvar{\delta\psi}\|_{\wgray{\Phi_r}}
\le C^{\rm w\text{-}ell}_{r,zr}\mathfrak E_k,
\label{eq:ell-low-output-zz}\\
\|D^\alpha T[\bvar{\delta\psi}]\|_{\wgray{\Phi_{i,m}}}
& \le C^{\rm w\text{-}ell}_{T,i,m,\alpha}\mathfrak E_k, \quad \text{for } (T,i,\alpha)  \in\mathcal J_{\rm ell},\quad m=|\alpha|. 
\label{eq:ell-high-output}
\end{align}

\hfill 

\paragraph{Differentiated unweighted output.}
\noindent The unweighted differentiated constant certifies
\begin{equation}\label{eq:ell-output-Cell-homogeneous}
\sum_{m=2}^{k+2}\sum_{|\alpha|=m}
\|D^\alpha\bvar{\delta\psi}\|_{L^2(\mathbb{D})}
\le C_{\rm ell}^{\rm diff}
\sum_{m=0}^{k}\sum_{|\alpha|=m}
\|D^\alpha\bvar{\delta\omega}\|_{L^2(\mathbb{D})}.
\end{equation}

\hfill  \\

\subsubsection{Reduced kernel and complete second-derivative operators}\label{sssec:ell-kernel-representation}

\hfill 

\noindent The operator \begin{equation}
    \mathcal{E}=\partial_r^2+(3/r)\partial_r+\partial_z^2
\end{equation} 
is the ordinary Laplacian in five dimensions restricted to functions that are radial in the first four variables.  This observation is useful because the five-dimensional Poisson equation has the explicit kernel \(\Gamma_5\).  We therefore lift \((r,z)\) to \((y,z)\in\mathbb R^4\times\mathbb R\), apply the standard five-dimensional representation, and then average over the angular variables in \(y\).  The result is a two-variable kernel \(K_0(r,z,s,\zeta)\) acting on \(\bvar{\delta\omega}(s,\zeta)\).  This subsection fixes this reduced kernel and explains how its derivatives are interpreted. \\ 

\paragraph{Reduction from five to two variables.}
\noindent Fix \(e_1\in\mathbb R^4\), set \(X=(re_1,z)\), and write the source point as \(Y=(s\sigma,\zeta)\), where \(s\ge0\), \(\zeta\in\mathbb R\), and \(\sigma\) lies on the unit sphere \(S^3=\{\sigma\in\mathbb R^4:\ |\sigma|=1\}\).  \\

\noindent The lifted perturbations are
\begin{equation}\label{eq:ell-lifted-perturbations}
\widetilde{\bvar{\delta\psi}}(y,z):=\bvar{\delta\psi}(|y|,z),
\qquad
\widetilde{\bvar{\delta\omega}}(y,z):=\bvar{\delta\omega}(|y|,z).
\end{equation}
They satisfy
\begin{equation}\label{eq:ell-lifted-poisson-used}
-\Delta_{y,z}\widetilde{\bvar{\delta\psi}}
=\widetilde{\bvar{\delta\omega}},
\qquad
\widetilde{\bvar{\delta\psi}}=\Gamma_5*\widetilde{\bvar{\delta\omega}}.
\end{equation}
In the variables \(Y=(s\sigma,\zeta)\),
\begin{equation}\label{eq:ell-fiveD-measure}
dY=s^3\,ds\,d\sigma\,d\zeta,
\end{equation}
and
\begin{equation}\label{eq:ell-fiveD-distance}
|X-Y|^2=r^2+s^2-2rs(e_1\cdot\sigma)+(z-\zeta)^2.
\end{equation}

\hfill 

\noindent Since the integrand depends on \(\sigma\) only through \(e_1\cdot\sigma\), we write \(e_1\cdot\sigma=\cos\theta\).  The spherical-coordinate formula on \(S^3\) gives
\begin{equation}\label{eq:ell-S3-angular-reduction}
\int_{S^3}F(e_1\cdot\sigma)\,d\sigma
=4\pi\int_0^\pi F(\cos\theta)\sin^2\theta\,d\theta .
\end{equation}
Thus
\begin{align}
\bvar{\delta\psi}(r,z)
&=\int_0^\infty\int_{\mathbb R}
K_0(r,z,s,\zeta)\bvar{\delta\omega}(s,\zeta)\,d\zeta\,ds,
\label{eq:ell-reduced-kernel-action}
\end{align}
where
\begin{equation}\label{eq:ell-detailed-one-dimensional-kernel}
K_0(r,z,s,\zeta)
:=\frac{s^3}{2\pi}\int_0^\pi
\frac{\sin^2\theta\,d\theta}{\bigl(r^2+s^2-2rs\cos\theta+(z-\zeta)^2\bigr)^{3/2}}.
\end{equation}
The factor \(s^3\) is part of \(K_0\), so reduced operators below act with respect to the flat measure \(ds\,d\zeta\).

\hfill

\paragraph{Weights in the reduced measure.}
\noindent For a reduced weight \(\phi\), define the lifted weight
\begin{equation}\label{eq:ell-lifted-weight-convention}
\widehat\phi(y,z):=|y|^{-3}\phi(|y|,z).
\end{equation}
Then
\begin{equation}\label{eq:ell-flat-weight-lift}
\int_\mathbb{D} |f(r,z)|^2\phi(r,z)\,dr\,dz
=\frac{1}{2\pi^2}\int_{\mathbb R^5}
|\widetilde f(y,z)|^2\widehat\phi(y,z)\,dy\,dz.
\end{equation}
This identity explains why the reduced norms use the flat measure: the five-dimensional volume factor has been absorbed into the lifted weight and the reduced kernel.

\hfill

\paragraph{Regular kernel derivatives.}
\noindent The reduced kernel is first differentiated away from the main diagonal \((r,z)=(s,\zeta)\).  When the resulting kernel is singular on the diagonal, the associated operator is interpreted in the principal-value sense.  Set
\begin{equation}\label{eq:ell-q-definition}
q(r,z,s,\zeta,\theta):=r^2+s^2-2rs\cos\theta+(z-\zeta)^2.
\end{equation}
The needed derivatives of \(q\) are
\begin{equation}\label{eq:ell-q-derivatives}
\partial_rq=2(r-s\cos\theta),
\qquad
\partial_zq=2(z-\zeta),
\end{equation}
\begin{equation}\label{eq:ell-q-second-derivatives}
\partial_{rr}q=2,
\qquad
\partial_{zz}q=2,
\qquad
\partial_{rz}q=0,
\end{equation}
and all higher derivatives of \(q\) vanish. \\

\noindent Define \(P_0=1\).  For \(x\in\{r,z\}\), define recursively
\begin{equation}\label{eq:ell-kernel-recursion}
P_{\gamma+e_x}
=q\,\partial_xP_\gamma
-\left(\frac{3}{2}+|\gamma|\right)P_\gamma\partial_xq.
\end{equation}
Indeed, if \(D_{r,z}^\gamma q^{-3/2}=P_\gamma q^{-3/2-|\gamma|}\), differentiating once more gives the recursion.  Therefore the off-diagonal derivative of the reduced kernel is
\begin{equation}\label{eq:ell-diff-kernel-explicit}
D_{r,z}^{\gamma,\mathrm{off}} K_0(r,z,s,\zeta)
=\frac{s^3}{2\pi}\int_0^\pi
\sin^2\theta\,P_\gamma q^{-3/2-|\gamma|}\,d\theta,
\qquad (r,z)\ne(s,\zeta).
\end{equation}
Thus \(D_{r,z}^{\gamma,\mathrm{off}}K_0\) is an ordinary derivative off the diagonal.  The principal-value qualification is added only when this off-diagonal kernel is used as a singular integral operator. \\

\paragraph{Complete Hessian operators.}
\noindent For second derivatives, the off-diagonal kernel alone does not represent the full distributional derivative. Away from the singularity at \(X=0\), the fundamental solution \(\Gamma_5\) is smooth and harmonic, so $
\Delta \Gamma_5 = 0$ away from $X=0$. Thus, the entire distributional mass in \(-\Delta\Gamma_5\) is concentrated at the origin. When the Laplacian is decomposed into the five diagonal second derivatives,$
\Delta\Gamma_5=\sum_{i=1}^{5}\partial_{X_iX_i}\Gamma_5$,
rotational symmetry forces this local mass to be distributed equally among the five diagonal entries of the Hessian. Consequently,
\begin{equation}\label{eq:ell-distributional-hessian-Gamma}
\partial_{X_iX_j}\Gamma_5
=
\operatorname{p.v.}\!\left(
\partial_{X_iX_j}^{\mathrm{off}}\Gamma_5
\right)
-
\frac{\delta_{ij}}{5}\delta_0.
\end{equation}

\noindent The first term is the nonlocal part of the Hessian. It is obtained by differentiating \(\Gamma_5\) away from the origin and interpreting the resulting singular kernel in the principal-value sense. The second term is local and is supported entirely at \(X=0\). The factor \(\delta_{ij}\) ensures that this Dirac mass appears only in the diagonal entries. When \(i=j\), the diagonal second derivative contains the local contribution \(-\delta_0/5\). When \(i\neq j\), the mixed second derivative contains no Dirac mass. \\

\noindent We now reduce this identity to the \((r,z)\) variables.  We lift \(\bvar{\delta\omega}(r,z)\) to \(\widetilde{\bvar{\delta\omega}}(y,z)=\bvar{\delta\omega}(|y|,z)\), convolve in five dimensions with \(\Gamma_5\), and integrate over \(\sigma\in S^3\) as in \eqref{eq:ell-S3-angular-reduction}.  This angular integration produces the reduced kernel \(K_0(r,z,s,\zeta)\).  Applying the same reduction to the off-diagonal Hessian, with a small neighborhood of the singularity $X=0$ removed in the full five-dimensional space before performing the angular integration, defines
\begin{equation}\label{eq:ell-low-kernels-defined}
L_{zz}^{\mathrm{off}}:=D_{r,z}^{2e_z,\mathrm{off}}K_0,
\qquad
L_{zr}^{\mathrm{off}}:=D_{r,z}^{e_z+e_r,\mathrm{off}}K_0.
\end{equation}
The local Dirac mass reduces even more simply: when the convolution is evaluated at \((r,z)\), the term \(\delta_0(X-Y)\) just returns the source value \(\bvar{\delta\omega}(r,z)\).  Therefore the two reduced second-derivative identities used later are
\begin{align}
\partial_{zz}\bvar{\delta\psi}
&=-\frac{1}{5}\bvar{\delta\omega}
+\operatorname{p.v.}\mathsf T_{L_{zz}^{\mathrm{off}}}\bvar{\delta\omega},
\label{eq:ell-zz-pv-local}\\
\partial_{zr}\bvar{\delta\psi}
&=\operatorname{p.v.}\mathsf T_{L_{zr}^{\mathrm{off}}}\bvar{\delta\omega}.
\label{eq:ell-zr-pv-local}
\end{align}

\noindent The first formula has a local term because \(\partial_{zz}\) is the second derivative in \(z\), and thus corresponds to \(i=j=5\) in \eqref{eq:ell-distributional-hessian-Gamma}. The mixed derivative \(\partial_{zr}\) corresponds to a mixed Cartesian Hessian entry, with indices \(5\) and \(1\), so \(\delta_{ij}=0\) and there is no local Dirac term. For the norm estimate, we use the full vector \(\nabla_y\partial_z\widetilde{\bvar{\delta\psi}}\). For a radial lift, the Euclidean norm of this vector is exactly \(|\partial_{zr}\bvar{\delta\psi}|\). \\

\noindent We reserve \(L_{zz}\) and \(L_{zr}\) for the complete maps
\begin{equation}\label{eq:ell-complete-operator-definitions}
\mathsf T_{L_{zz}}:=-\frac{1}{5}I+\operatorname{p.v.}\mathsf T_{L_{zz}^{\mathrm{off}}},
\qquad
\mathsf T_{L_{zr}}:=\operatorname{p.v.}\mathsf T_{L_{zr}^{\mathrm{off}}}.
\end{equation}
Thus \(L_{zz}^{\mathrm{off}}\) and \(L_{zr}^{\mathrm{off}}\) denote only the singular-integral kernels, while \(L_{zz}\) and \(L_{zr}\) denote the full operators that are inserted into estimates. In particular, \(L_{zz}\) represents the complete \(zz\)-operator. Its definition already includes the local multiplication term \(-\frac{1}{5}f\), which comes from the Dirac mass in the diagonal Hessian. Therefore, this term should not be added separately a second time.

\hfill

\noindent Further \(z\)-derivatives are obtained by differentiating the source in these identities.  For every \(\ell\ge0\),
\begin{equation}\label{eq:ell-zz-further-z-derivatives}
\partial_z^{\ell+2}\bvar{\delta\psi}
=-\frac{1}{5}\partial_z^\ell\bvar{\delta\omega}
+\operatorname{p.v.}\mathsf T_{L_{zz}^{\mathrm{off}}}
  (\partial_z^\ell\bvar{\delta\omega}),
\end{equation}
and
\begin{equation}\label{eq:ell-zr-further-z-derivatives}
\partial_r\partial_z^{\ell+1}\bvar{\delta\psi}
=\operatorname{p.v.}\mathsf T_{L_{zr}^{\mathrm{off}}}
  (\partial_z^\ell\bvar{\delta\omega}).
\end{equation}

\hfill

\begin{proposition}[Axis-weighted Hessian bounds]\label{prop:ell-axis-pv}

\noindent For \(L\in\{zz,zr\}\), define
\begin{equation}
\widetilde{\mathsf T}_{zz}F
:=\operatorname{p.v.}\bigl((\partial_{zz}^{\mathrm{off}}\Gamma_5)\ast F\bigr),
\qquad
\widetilde{\mathsf T}_{zr}F
:=\operatorname{p.v.}\bigl((\nabla_y\partial_z\Gamma_5)\ast F\bigr).
\end{equation}
Let \(\mathcal K_L\) denote the corresponding kernel, and define
\begin{equation}
\kappa_L
:=
\begin{cases}
4/5, & L=zz,\\[2mm]
1/2, & L=zr,
\end{cases}
\qquad
B_L
:=
\begin{cases}
\dfrac{3}{2\pi^2}, & L=zz,\\[2mm]
\dfrac{15}{16\pi^2}, & L=zr.
\end{cases}
\label{eq:ell-axis-hessian-data}
\end{equation}
Choose \(B_L\) so that
\begin{equation}
|\mathcal K_{zz}(X)|
\le B_{zz}|X|^{-5},
\qquad
\|\mathcal K_{zr}(X)\|_{\mathbb R^4}
\le B_{zr}|X|^{-5},
\end{equation}
and set
\begin{align}
K_L^{\rm ax}
:={}&\kappa_L(9+2^{3/2})
+2\pi^2B_L\sqrt{\frac{2^{21}}{297}}
\frac{2^{-7}}{1-2^{-7/2}}
+\frac{256\pi^2}{\sqrt{135}}B_L
\sqrt{(2^{3/2}-1)(2^{5/2}-1)}
\frac{2^{-1}}{1-2^{-1/2}}.
\label{eq:ell-axis-pv-constant}
\end{align}
Then
\begin{equation}
\|\widetilde{\mathsf T}_LF\|_{L^2(|y|^{-3}dX)}
\le K_L^{\rm ax}\|F\|_{L^2(|y|^{-3}dX)},
\label{eq:ell-axis-pv-bound}
\end{equation}
where for \(L=zr\) the norm on the left is the \(L^2\)-norm of the Euclidean norm in \(\mathbb R^4\).\\

\noindent If \(F\) is radial in \(y\), then
\begin{equation}
\|\widetilde{\mathsf T}_{zr}F(y,z)\|_{\mathbb R^4}
=
|\partial_{zr}(\Gamma_5\ast F)(|y|,z)|.
\end{equation}

\end{proposition}

\hfill

\begin{grayproof}
\noindent We write \(w(y)=|y|^{-3}\). Fix \(L\in\{zz,zr\}\) and use \(\|\cdot\|_w\) for the corresponding weighted \(L^2\)-norm, with the Euclidean norm in \(\mathbb R^4\) when \(L=zr\). \\

\noindent First we take \(F\) with finitely many nonzero shell pieces, so the unweighted operator is defined.
We decompose \(\mathbb R^5\) into radial shells
\begin{equation}
\mathcal A_j:=\{X=(y,z):2^j\le |y|<2^{j+1}\},
\qquad F_j:=1_{\mathcal A_j}F.
\label{eq:ell-axis-shells-proof-add}
\end{equation}
Because the shells are disjoint,
\begin{equation}
\|F\|_w^2=\sum_j\|F_j\|_w^2,
\qquad
\|\widetilde{\mathsf T}_LF\|_w^2
=\sum_m\left\|1_{\mathcal A_m}\widetilde{\mathsf T}_LF\right\|_w^2.
\label{eq:ell-axis-shell-orthogonality-add}
\end{equation}

\hfill \\

\noindent\textbf{Comparable shells.}
Fix \(j\) and write the output shell as \(\mathcal A_{j+n}\). On this shell,
\(|y|^{-3/2}\le2^{-3(j+n)/2}\), while on the input shell
\(\|F_j\|_{L^2}\le2^{3(j+1)/2}\|F_j\|_w\). \\

\noindent Therefore, the unweighted multiplier bound gives
\begin{align}
\left\|1_{\mathcal A_{j+n}}\widetilde{\mathsf T}_LF_j\right\|_w
&\le 2^{-3(j+n)/2}
\left\|\widetilde{\mathsf T}_LF_j\right\|_{L^2}
\le \kappa_L2^{\frac{3}{2}(1-n)}\|F_j\|_w,
\qquad |n|\le1.
\label{eq:ell-axis-near-shell-bound-add}
\end{align}

\hfill

\noindent The three coefficients for \(n=-1,0,1\) therefore sum to
\(\kappa_L(8+2^{3/2}+1)=\kappa_L(9+2^{3/2})\).

\hfill \\

\noindent\textbf{Separated shells.}
Conjugate the operator from \(L^2(w\,dX)\) to unweighted \(L^2\). The conjugated kernel is
\begin{equation}
\mathcal H_L(X,Y)
:=|y|^{-3/2}\mathcal K_L(X-Y)|y'|^{3/2},
\qquad X=(y,z),\quad Y=(y',z').
\label{eq:ell-axis-conjugated-kernel-add}
\end{equation}

\hfill

\noindent For input supported on \(\mathcal A_j\) and output restricted to \(\mathcal A_{j+n}\), define
\begin{align}
M_{j,n}^{\rm out}
:=\sup_{X\in\mathcal A_{j+n}}
\int_{\mathcal A_j}\|\mathcal H_L(X,Y)\|\,dY, \qquad
M_{j,n}^{\rm in}
:=\sup_{Y\in\mathcal A_j}
\int_{\mathcal A_{j+n}}\|\mathcal H_L(X,Y)\|\,dX.
\label{eq:ell-axis-kernel-integrals-add}
\end{align}
Here \(\|\cdot\|\) denotes absolute value for \(L=zz\) and the Euclidean norm in \(\mathbb R^4\) for \(L=zr\).\\

\noindent Thus \(M_{j,n}^{\rm out}\) bounds the integral of the kernel norm over the input shell uniformly in the output point, while \(M_{j,n}^{\rm in}\) gives the analogous bound after integrating over the output shell. Let
\begin{equation}
G_j(Y):=|y'|^{-3/2}F_j(Y),
\qquad \|G_j\|_{L^2}=\|F_j\|_w.
\label{eq:ell-axis-conjugated-input}
\end{equation}
For each fixed \(X\), the Cauchy--Schwarz inequality gives
\begin{equation}
\left|\int_{\mathcal A_j}\mathcal H_L(X,Y)G_j(Y)\,dY\right|^2
\le
M_{j,n}^{\rm out}
\int_{\mathcal A_j}\|\mathcal H_L(X,Y)\|\,|G_j(Y)|^2\,dY.
\label{eq:ell-axis-cauchy-schwarz}
\end{equation}
Integrating in \(X\) and then reversing the order of integration gives
\begin{align}
\left\|1_{\mathcal A_{j+n}}\widetilde{\mathsf T}_LF_j\right\|_w^2
&\le
M_{j,n}^{\rm out}
\int_{\mathcal A_j}|G_j(Y)|^2
\left(\int_{\mathcal A_{j+n}}\|\mathcal H_L(X,Y)\|\,dX\right)dY
\label{eq:ell-axis-integrate-first}\\
&\le
M_{j,n}^{\rm out}M_{j,n}^{\rm in}\|F_j\|_w^2.
\label{eq:ell-axis-integrate-second}
\end{align}
Taking square roots gives
\begin{equation}
\left\|1_{\mathcal A_{j+n}}\widetilde{\mathsf T}_LF_j\right\|_w
\le \sqrt{M_{j,n}^{\rm out}M_{j,n}^{\rm in}}\,\|F_j\|_w.
\label{eq:ell-axis-integral-operator-step-add}
\end{equation}

\hfill

\noindent For \(|n|\ge2\), the radial shells are separated. Using
\begin{equation}
\|\mathcal K_L(X-Y)\|
\le B_L|X-Y|^{-5},
\end{equation}
integrating first in \(z-z'\) through
\begin{equation}
\int_{\mathbb R}(a^2+t^2)^{-5/2}\,dt=\frac{4}{3a^4},
\label{eq:ell-axis-z-integral-add}
\end{equation}
and then using four-dimensional polar coordinates gives the following explicit bounds. Write
\begin{equation}
r=|y|,\qquad \rho=|y'|,\qquad a=2^j,\qquad R=2^{j+n}.
\label{eq:ell-axis-radial-variables}
\end{equation}
Using \(|y-y'|\ge |r-\rho|\) and \(|S^3|=2\pi^2\), the integrations in \(z-z'\) and the angular variables contribute the factor
\begin{equation}
A_L:=\frac{4}{3}|S^3|B_L=\frac{8\pi^2}{3}B_L.
\label{eq:ell-axis-angular-factor}
\end{equation}
The definitions of \(M_{j,n}^{\rm out}\) and \(M_{j,n}^{\rm in}\) then give
\begin{align}
M_{j,n}^{\rm out}
&\le \sup_{R\le r<2R}A_Lr^{-3/2}
   \int_a^{2a}\frac{\rho^{9/2}}{|r-\rho|^4}\,d\rho,
\label{eq:ell-axis-out-radial}\\
M_{j,n}^{\rm in}
&\le \sup_{a\le\rho<2a}A_L\rho^{3/2}
   \int_R^{2R}\frac{r^{3/2}}{|r-\rho|^4}\,dr.
\label{eq:ell-axis-in-radial}
\end{align}
If \(n\ge2\), then \(\rho\le r/2\), so \(|r-\rho|^{-4}\le16r^{-4}\). Direct integration gives
\begin{align}
M_{j,n}^{\rm out}&\le A_L\frac{2^{21/2}}{11}2^{-11n/2},
\label{eq:ell-axis-out-positive}\\
M_{j,n}^{\rm in}&\le A_L\frac{2^{13/2}}{3}2^{-3n/2}.
\label{eq:ell-axis-in-positive}
\end{align}
If \(n\le-2\), then \(r\le\rho/2\), so \(|r-\rho|^{-4}\le16\rho^{-4}\). Direct integration gives
\begin{align}
M_{j,n}^{\rm out}&\le A_L\frac{32}{3}(2^{3/2}-1)2^{-3n/2},
\label{eq:ell-axis-out-negative}\\
M_{j,n}^{\rm in}&\le A_L\frac{32}{5}(2^{5/2}-1)2^{5n/2}.
\label{eq:ell-axis-in-negative}
\end{align}
Multiplying these bounds yields
\begin{align}
M_{j,n}^{\rm out}M_{j,n}^{\rm in}
&\le (C_L^+)^2 2^{-7n},
&&n\ge2,
\label{eq:ell-axis-kernel-product-positive}\\
M_{j,n}^{\rm out}M_{j,n}^{\rm in}
&\le (C_L^-)^2 2^{n},
&&n\le-2.
\label{eq:ell-axis-kernel-product-negative}
\end{align}
where
\begin{equation}
C_L^+:=2\pi^2B_L\sqrt{\frac{2^{21}}{297}},
\qquad
C_L^-:=\frac{256\pi^2}{\sqrt{135}}B_L
\sqrt{(2^{3/2}-1)(2^{5/2}-1)}.
\label{eq:ell-axis-far-shell-constants-add}
\end{equation}
These are exactly the product coefficients above, since
\begin{align}
C_L^+&=\frac{A_L2^{17/2}}{\sqrt{33}},
\label{eq:ell-axis-positive-coefficient}\\
C_L^-&=\frac{32A_L}{\sqrt{15}}
\sqrt{(2^{3/2}-1)(2^{5/2}-1)}.
\label{eq:ell-axis-negative-coefficient}
\end{align}
Substituting \(A_L=8\pi^2B_L/3\) gives the displayed forms. Here \(297=9\cdot33\) and \(135=9\cdot15\).

\hfill

\noindent Combining \eqref{eq:ell-axis-integral-operator-step-add}, \eqref{eq:ell-axis-kernel-product-positive}, and \eqref{eq:ell-axis-kernel-product-negative} gives
\begin{equation}
\left\|1_{\mathcal A_{j+n}}\widetilde{\mathsf T}_LF_j\right\|_w
\le
\begin{cases}
C_L^+2^{-7n/2}\|F_j\|_w,&n\ge2,\\
C_L^-2^{n/2}\|F_j\|_w,&n\le-2.
\end{cases}
\label{eq:ell-axis-far-shell-bounds-add}
\end{equation}

\hfill \\

\noindent\textbf{Summation over shells.}
Set \(a_j=\|F_j\|_w\) and
\(b_m=\|1_{\mathcal A_m}\widetilde{\mathsf T}_LF\|_w\). The preceding estimates imply
\begin{equation}
b_m\le\sum_{n\in\mathbb Z}k_na_{m-n},
\label{eq:ell-axis-discrete-convolution-add}
\end{equation}
where
\begin{equation}
k_n:=
\begin{cases}
C_L^-2^{n/2},&n\le-2,\\
8\kappa_L,&n=-1,\\
2^{3/2}\kappa_L,&n=0,\\
\kappa_L,&n=1,\\
C_L^+2^{-7n/2},&n\ge2.
\end{cases}
\label{eq:ell-axis-discrete-kernel-add}
\end{equation}

\hfill

\noindent Young's inequality on \(\ell^2(\mathbb Z)\) and \eqref{eq:ell-axis-shell-orthogonality-add} give
\begin{equation}
\|\widetilde{\mathsf T}_LF\|_w
\le\|k\|_{\ell^1}\|F\|_w.
\label{eq:ell-axis-young-add}
\end{equation}

\hfill

\noindent The two tails are geometric:
\begin{align}
\sum_{n\ge2}C_L^+2^{-7n/2}
&=C_L^+\frac{2^{-7}}{1-2^{-7/2}},\\
\sum_{n\le-2}C_L^-2^{n/2}
&=C_L^-\frac{2^{-1}}{1-2^{-1/2}}.
\end{align}
Together with the three comparable-shell coefficients, this gives
\(\|k\|_{\ell^1}=K_L^{\rm ax}\).
Finite-shell functions are dense in \(L^2(w\,dX)\), so this estimate extends \(\widetilde{\mathsf T}_L\) uniquely by continuity and proves \eqref{eq:ell-axis-pv-bound}.

\hfill \\

\noindent\textbf{Hessian data.}
The five-dimensional fundamental solution is
\(\Gamma_5(X)=(8\pi^2)^{-1}|X|^{-3}\).\\

\noindent Its off-diagonal Hessian kernel is
\begin{equation}
\partial_{X_iX_j}^{\rm off}\Gamma_5(X)
=\frac{3}{8\pi^2}
\frac{5X_iX_j-\delta_{ij}|X|^2}{|X|^7}.
\label{eq:ell-axis-hessian-kernel-add}
\end{equation}

\hfill

\noindent For a diagonal entry, \(|5X_i^2-|X|^2|\le4|X|^2\), which gives
\begin{equation}
B_{zz}=\frac{3}{2\pi^2}.
\end{equation}

\hfill

\noindent For the mixed vector kernel,
\begin{equation}
\left\|\frac{15}{8\pi^2}\frac{X_zX_y}{|X|^7}\right\|_{\mathbb R^4}
\le \frac{15}{16\pi^2}|X|^{-5},
\end{equation}
since \(|X_z|\,|X_y|\le |X|^2/2\). \\

\noindent Therefore \begin{equation} B_{zr}=\frac{15}{16\pi^2}.\end{equation}

\hfill

\noindent Using the distributional identity
\begin{equation}
\partial_{X_iX_j}\Gamma_5
=
\operatorname{p.v.}\!\left(\partial_{X_iX_j}^{\mathrm{off}}\Gamma_5\right)
-\frac{\delta_{ij}}{5}\delta_0,
\end{equation}
the Fourier multiplier of the diagonal principal-value part is
\begin{equation}
-\frac{\xi_i^2}{|\xi|^2}+\frac{1}{5},
\end{equation}
whose absolute value is at most \(4/5\).\\

\noindent The mixed term is \(-\xi_z\xi_y/|\xi|^2\), with norm
\(|\xi_z|\,|\xi_y|/|\xi|^2\le1/2\).

\hfill

\noindent Plancherel's theorem therefore gives
\(\kappa_{zz}=4/5\) and \(\kappa_{zr}=1/2\), completing the proof of \eqref{eq:ell-axis-hessian-data}.
\end{grayproof}

\hfill

\subsubsection{Point constants}\label{sssec:ell-point-constants}

\noindent We now certify the two streamfunction values at the origin that enter the modulation equations. \\ 

\noindent We compute the origin kernel, take its weighted dual norm against the $\omega$ weight, and store the resulting low-order constant.  The required output is
\begin{equation}
P_{\psi,0}^{0},\quad P_{\psi,0}^{k},\quad P_{\psi_z,0}^{0},\quad P_{\psi_z,0}^{k},
\end{equation}
with
\begin{equation}
|\bvar{\delta\psi}(0,0)|
\le P_{\psi,0}^{0}\mathfrak E_0+P_{\psi,0}^{k}\mathfrak H_k,
\qquad
|\partial_z\bvar{\delta\psi}(0,0)|
\le P_{\psi_z,0}^{0}\mathfrak E_0+P_{\psi_z,0}^{k}\mathfrak H_k.
\end{equation}

\hfill 

\noindent For the direct origin kernels, the high-order coefficients are zero, since this uses only the weighted low-order norm of \(\bvar{\delta\omega}\). \\

\hfill 

\paragraph{Origin kernels.}
\noindent We evaluate the reduced kernel at \((r,z)=(0,0)\).  Since the angular integral is explicit at the origin,
\begin{align}
k_\psi(s,\zeta)
&:=K_0(0,0,s,\zeta)
=\frac{s^3}{4(s^2+\zeta^2)^{3/2}},
\label{eq:ell-origin-kernel-psi}\\
k_{\psi_z}(s,\zeta)
&:=\partial_zK_0(0,0,s,\zeta)
=\frac{3\zeta s^3}{4(s^2+\zeta^2)^{5/2}}.
\label{eq:ell-origin-kernel-psiz}
\end{align}

\noindent Therefore
\begin{align}
\bvar{\delta\psi}(0,0)
&=\int_0^\infty\int_{\mathbb R}
k_\psi(s,\zeta)\bvar{\delta\omega}(s,\zeta)\,d\zeta\,ds,
\label{eq:ell-point-duality-psi}\\
\partial_z\bvar{\delta\psi}(0,0)
&=\int_0^\infty\int_{\mathbb R}
k_{\psi_z}(s,\zeta)\bvar{\delta\omega}(s,\zeta)\,d\zeta\,ds,
\label{eq:ell-point-duality-psiz} \end{align}

\hfill \\

\paragraph{Direct dual estimate.}
\noindent For any scalar origin kernel \(h\), we insert \(1=\wgray{\Phi_\omega}^{-1/2}\wgray{\Phi_\omega}^{1/2}\) in the pairing and apply the Cauchy--Schwarz inequality.  This gives
\begin{equation}\label{eq:ell-point-dual-template-bound}
\left|\int_0^\infty\int_{\mathbb R}h\,\bvar{\delta\omega}\,d\zeta\,ds\right|
\le
\left(\int_0^\infty\int_{\mathbb R}
\frac{|h(s,\zeta)|^2}{\wgray{\Phi_\omega(s,\zeta)}}\,d\zeta\,ds\right)^{1/2}
\|\bvar{\delta\omega}\|_{\wgray{\Phi_\omega}}.
\end{equation}

\hfill 

\noindent With \(h=k_\psi\), this gives
\begin{equation}\label{eq:ell-Ppsi-final}
P_{\psi,0}^{\rm dir}
:=\left(\int_0^\infty\int_{\mathbb R}
\frac{s^6}{16(s^2+\zeta^2)^3\wgray{\Phi_\omega(s,\zeta)}}
\,d\zeta\,ds\right)^{1/2}.
\end{equation}
With \(h=k_{\psi_z}\), it gives
\begin{equation}\label{eq:ell-Ppsiz-final}
P_{\psi_z,0}^{\rm dir}
:=\left(\int_0^\infty\int_{\mathbb R}
\frac{9\zeta^2s^6}{16(s^2+\zeta^2)^5\wgray{\Phi_\omega(s,\zeta)}}
\,d\zeta\,ds\right)^{1/2}.
\end{equation}

\hfill \\ 

\noindent The stored constants are therefore
\begin{equation}\label{eq:ell-Ppsi-selected-direct}
P_{\psi,0}^{0}:=P_{\psi,0}^{\rm dir},
\qquad
P_{\psi,0}^{k}:=0,
\end{equation}
and
\begin{equation}\label{eq:ell-Ppsiz-selected-direct}
P_{\psi_z,0}^{0}:=P_{\psi_z,0}^{\rm dir},
\qquad
P_{\psi_z,0}^{k}:=0.
\end{equation}

\hfill \\

\subsubsection{Weighted operator constants and block bounds}\label{sssec:ell-weighted-operator-constants}

\hfill 

\paragraph{Operator norm.}
\noindent For a reduced kernel \(\Upsilon\), define
\begin{equation}\label{eq:ell-kernel-operator-T}
(\mathsf T_{\Upsilon}f)(r,z)
:=\int_0^\infty\int_{\mathbb R}
\Upsilon(r,z,s,\zeta)f(s,\zeta)\,d\zeta\,ds.
\end{equation}
For regular kernels, and for principal-value kernels away from the singularity, this formula applies directly to the kernel. For the complete Hessian operators \(L_{zz}\) and \(L_{zr}\), however, \(\mathsf T_{\Upsilon}\) denotes the full operators defined in \eqref{eq:ell-complete-operator-definitions}, including any local term arising from the distributional Hessian.  Thus
\begin{equation}
\mathsf T_{L_{zz}}f=-\frac15f+\operatorname{p.v.}\mathsf T_{L_{zz}^{\mathrm{off}}}f,
\qquad
\mathsf T_{L_{zr}}f=\operatorname{p.v.}\mathsf T_{L_{zr}^{\mathrm{off}}}f.
\end{equation}

\hfill 

  \noindent If a multiplier \(a=a(r,z)\) is attached to the output, then \(a\Upsilon\) denotes the operator mapping $f$ to \( a\mathsf T_{\Upsilon}f\).  The same complete-operator convention is used, so \(aL_{zz}\) includes the local term \(-af/5\).  \\
  
  \noindent For positive input and output weights, set
\begin{equation}\label{eq:ell-target-weighted-form}
\mathsf C [\Upsilon,\phi_{\rm in},\phi_{\rm out}]
:=\sup_{f\ne0}
\frac{\|\mathsf T_{\Upsilon}f\|_{\phi_{\rm out}}}{\|f\|_{\phi_{\rm in}}}.
\end{equation}

\hfill 

\noindent This is the certified weighted operator norm used in later elliptic terms.  Once known, the bound
\begin{equation}\label{eq:ell-target-weighted-form-use}
\|\mathsf T_{\Upsilon}f\|_{\phi_{\rm out}}
\le
\mathsf C [\Upsilon,\phi_{\rm in},\phi_{\rm out}]\|f\|_{\phi_{\rm in}}
\end{equation}
can be inserted directly into the energy estimates.  

\hfill  \\ 

{
\begin{proposition}[Finite block-matrix certificate]\label{prop:ell-block-certificate}
Let \(\mathbb{D}=\bigcup_{a\in\mathcal A}\mathbb{D}_a\) be a finite disjoint measurable partition. Suppose that numbers \(c_{ab}\ge0\) have been certified for every ordered pair \((a,b)\) and satisfy
\begin{equation}\label{eq:ell-block-bound}
\|1_{\mathbb{D}_a}\mathsf T_{\Upsilon}(1_{\mathbb{D}_b}f)\|_{\phi_{\rm out}}
\le c_{ab}\|1_{\mathbb{D}_b}f\|_{\phi_{\rm in}}
\qquad\text{for every }f.
\end{equation}
Let \(C=(c_{ab})_{a,b\in\mathcal A}\).  Then

\begin{equation}\label{eq:ell-block-matrix-bound}
\mathsf C[\Upsilon,\phi_{\rm in},\phi_{\rm out}]
\le \|C\|_{\ell^2(\mathcal A)\to\ell^2(\mathcal A)}
=\sqrt{\lambda_{\max}(C^\top C)}.
\end{equation}

\end{proposition}

\begin{grayproof}
{
\noindent For each input block, set
\begin{equation}
f_b:=1_{\mathbb{D}_b}f,
\qquad
x_b:=\|f_b\|_{\phi_{\rm in}}.
\label{eq:ell-block-proof-input-vector}
\end{equation}
For each output block, set
\begin{equation}
y_a:=\|1_{\mathbb{D}_a}\mathsf T_{\Upsilon}f\|_{\phi_{\rm out}}.
\label{eq:ell-block-proof-output-vector}
\end{equation}

\hfill 

\noindent By linearity, the triangle inequality, and \eqref{eq:ell-block-bound},
\begin{equation}
y_a
\le \sum_{b\in\mathcal A}
\|1_{\mathbb{D}_a}\mathsf T_{\Upsilon}f_b\|_{\phi_{\rm out}}
\le \sum_{b\in\mathcal A}c_{ab}x_b
=(Cx)_a.
\label{eq:ell-block-proof-row}
\end{equation}
Thus \(0\le y\le Cx\) componentwise. \\

\noindent Since \(C\) and \(x\) have nonnegative entries,
\begin{equation}
\|y\|_{\ell^2}
\le \|Cx\|_{\ell^2}
\le \|C\|_{\ell^2\to\ell^2}\|x\|_{\ell^2}.
\label{eq:ell-block-proof-vector-bound}
\end{equation}
The partition is disjoint, so
\begin{equation}
\|y\|_{\ell^2}^2=\|\mathsf T_{\Upsilon}f\|_{\phi_{\rm out}}^2,
\qquad
\|x\|_{\ell^2}^2=\|f\|_{\phi_{\rm in}}^2.
\label{eq:ell-block-proof-orthogonal-sums}
\end{equation}

\hfill 

\noindent Combining \eqref{eq:ell-block-proof-vector-bound} and \eqref{eq:ell-block-proof-orthogonal-sums} proves the operator bound. The identity
\begin{equation}
    \|C\|_{\ell^2\to\ell^2}=\sqrt{\lambda_{\max}(C^\top C)}
\end{equation}
is the finite-dimensional singular-value formula.
}
\end{grayproof}
}

\noindent The proposition turns local estimates into a global operator bound.  Each block-pair constant \(c_{ab}\) measures how much input from block \(b\) can contribute to output on block \(a\).  The finite matrix \(C=(c_{ab})\) then combines all block interactions in the same way that the integral operator combines all input regions.  Bounding the \(\ell^2\) norm of \(C\) gives a rigorous global bound for \(\mathsf T_\Upsilon\) without estimating the full singular operator in one step.

\hfill \\

\paragraph{Block estimates.}

\noindent The entry \(c_{ab}\) controls the operator restricted to inputs on \(\mathbb{D}_b\) and outputs on \(\mathbb{D}_a\).  For each block pair we use one of the following estimates.  Once all entries are certified, Proposition~\ref{prop:ell-block-certificate} gives the global constant \(\mathsf C [\Upsilon,\phi_{\rm in},\phi_{\rm out}]\).
\begin{enumerate}[label=\textup{(W\arabic*)},leftmargin=*,itemsep=2.5em] 
\vspace{2mm}
\item \emph{Separated off-diagonal blocks.}  When \(\mathbb{D}_a, \, \mathbb{D}_b\) stay a positive distance from the diagonal, the kernel is nonsingular on the block pair.  For input supported on \(\mathbb{D}_b\), we insert
\begin{equation}
1=\phi_{\rm in}(s,\zeta)^{-1/2} \, \phi_{\rm in}(s,\zeta)^{1/2}
\end{equation}
in the integral defining \(1_{\mathbb{D}_a}\mathsf T_\Upsilon(1_{\mathbb{D}_b}f)\). \\

\noindent The Cauchy--Schwarz inequality in \((s,\zeta)\) gives
\begin{align}
|1_{\mathbb{D}_a}\mathsf T_\Upsilon(1_{\mathbb{D}_b}f)(r,z)|^2
&\le
\left(\int_{\mathbb{D}_b}|\Upsilon(r,z,s,\zeta)|^2\phi_{\rm in}(s,\zeta)^{-1}\,ds\,d\zeta\right)
\|1_{\mathbb{D}_b}f\|_{\phi_{\rm in}}^2 .
\end{align}

\hfill 

\noindent Multiplying by \(\phi_{\rm out}(r,z)\), integrating over \(\mathbb{D}_a\), and taking square roots gives the estimate
\begin{equation}\label{eq:ell-HS-box}
c_{ab}
\le
\left(\int_{\mathbb{D}_a}\int_{\mathbb{D}_b}
|\Upsilon(r,z,s,\zeta)|^2
\frac{\phi_{\rm out}(r,z)}{\phi_{\rm in}(s,\zeta)}
\,ds\,d\zeta\,dr\,dz\right)^{1/2}.
\end{equation}  

\vspace{4mm}

\item \emph{Blocks meeting the diagonal.} A block pair meets the diagonal if the closures of \(\mathbb{D}_a\) and \(\mathbb{D}_b\) contain points with \((r,z)=(s,\zeta)\).  On such a pair the kernel may be singular, so \eqref{eq:ell-HS-box} is not used.

\noindent On a block meeting the diagonal, we decompose the singular part into the finite sum
\begin{equation}
\sum_{\nu=1}^{N_\Upsilon}a_\nu\operatorname{p.v.}\mathsf T_{L_\nu}.
\label{eq:ell-pv-singular-sum-add}
\end{equation}

Fix one summand. Proposition~\ref{prop:ell-axis-pv} controls the principal-value operator. For this particular coefficient and block pair, we certify the complete coefficient-weight conversion
\begin{equation}
A_{\nu,ab}:=
\operatorname*{ess\,sup}_{x\in\mathbb D_a,\,y\in\mathbb D_b}
|a_\nu(x)|
\sqrt{\frac{\phi_{\rm out}(x)}{\phi_{\rm in}(y)}}<\infty.
\label{eq:ell-pv-block-weight-ratio-add}
\end{equation}
Then
\begin{align}
&\left\|1_{\mathbb D_a}a_\nu\operatorname{p.v.}\mathsf T_{L_\nu}
(1_{\mathbb D_b}f)\right\|_{\phi_{\rm out}}
\le K_{L_\nu}^{\rm ax}A_{\nu,ab}\|1_{\mathbb D_b}f\|_{\phi_{\rm in}}.
\label{eq:ell-pv-one-summand-block-bound-add}
\end{align}
Summing over the finitely many singular summands gives
\begin{equation}\label{eq:ell-pv-compact-block-bound}
c_{ab}^{\rm pv}
\le \sum_{\nu=1}^{N_\Upsilon}K_{L_\nu}^{\rm ax}A_{\nu,ab}.
\end{equation}

A regular remainder is estimated separately using \eqref{eq:ell-HS-box}. If a block is unbounded, or if the weight ratio cannot be bounded uniformly over the entire input-output block pair, we divide both blocks into dyadic annuli. We apply \eqref{eq:ell-pv-one-summand-block-bound-add} to each pair of annuli. Then we bound the interaction between every input annulus and every output annulus, and sum these bounds.

\noindent If the complete operator is \(L_{zz}\), the block also contains the local multiplication term \(-I/5\) when the input and output blocks overlap.  Its weighted block norm is bounded by
\begin{equation}\label{eq:ell-local-block-contribution}
\frac{1}{5}\operatorname*{ess\,sup}_{\mathbb{D}_a\cap\mathbb{D}_b}
\left(\frac{\phi_{\rm out}}{\phi_{\rm in}}\right)^{1/2},
\end{equation}
with default value zero when \(\mathbb{D}_a\cap\mathbb{D}_b=\emptyset\).  Thus a diagonal block for \(L_{zz}\) uses \(c_{ab}=c_{ab}^{\rm pv}\) plus the local contribution in \eqref{eq:ell-local-block-contribution}.  For \(L_{zr}\), there is no local contribution.

\vspace{2mm}

\item \emph{Separated tails.}  For tail blocks away from the diagonal, the exact kernel is replaced by an explicit majorant.  From \eqref{eq:ell-diff-kernel-explicit},
\begin{equation}
|D_{r,z}^\gamma K_0(r,z,s,\zeta)|
\le \frac{s^3}{2\pi}\int_0^\pi
\sin^2\theta\,|P_\gamma(r,z,s,\zeta,\theta)|\,
q^{-3/2-|\gamma|}\,d\theta .
\end{equation}
On the tail block, the denominator satisfies
\begin{equation}\label{eq:ell-tail-q-lower}
q(r,z,s,\zeta,\theta)
=(r-s)^2+(z-\zeta)^2+2rs(1-\cos\theta)
\ge (r-s)^2+(z-\zeta)^2.
\end{equation}
Interval arithmetic also encloses the dimensionless numerator produced by the recursion, giving a constant \(A_\gamma\) such that
\begin{equation}
\frac{1}{2\pi}\int_0^\pi
\sin^2\theta\,|P_\gamma(r,z,s,\zeta,\theta)|\,
q^{-3/2-|\gamma|}
\,d\theta
\le
A_\gamma \bigl((r-s)^2+(z-\zeta)^2\bigr)^{-(3+|\gamma|)/2}.
\end{equation}
Therefore the recursion \eqref{eq:ell-kernel-recursion} gives the tail majorant
\begin{equation}\label{eq:ell-tail-decay-simple}
|D_{r,z}^\gamma K_0(r,z,s,\zeta)|
\le A_\gamma s^3
\bigl((r-s)^2+(z-\zeta)^2\bigr)^{-(3+|\gamma|)/2}.
\end{equation}
This pointwise bound is then inserted into \eqref{eq:ell-HS-box} or into the corresponding analytic tail integral.

\end{enumerate}

\hfill \\

\noindent Applying this construction to the complete Hessian maps gives
\begin{equation}\label{eq:ell-detailed-low-operator-constants}
C^{\rm w\text{-}ell}_{z,zz}
:=\mathsf C [L_{zz},\wgray{\Phi_\omega},\wgray{\Phi_z}],
\qquad
C^{\rm w\text{-}ell}_{r,zr}
:=\mathsf C [L_{zr},\wgray{\Phi_\omega},\wgray{\Phi_r}].
\end{equation}

\paragraph{Localized near-origin fallback.}
\noindent If a global Hessian bound fails only because its target weight is singular at the origin, we localize that component alone. We fix a certified localization scale $h_{\mathsf{loc}}>0$ and a pointwise Hessian bound
\begin{equation}
\sum_{|\alpha|=2}\|D^\alpha\bvar{\delta\psi}\|_{L^\infty(\mathbb D)}
\le P_{\psi,2}^{\rm pt}\mathfrak E_k.
\label{eq:ell-pointwise-hessian-revision}
\end{equation}
\noindent To certify $P_{\psi,2}^{\rm pt}$, we first use the differentiated unweighted elliptic estimate \eqref{eq:ns-Cell-homogeneous} to control the Hessian of $\bvar{\delta\psi}$ by unweighted Sobolev norms of $\bvar{\delta\omega}$. Weight coercivity converts the required unweighted $L^2$ norms to weighted ones, interpolation supplies any intermediate derivative levels, and Proposition~\ref{prop:whole-plane-standard-constants} then gives the $L^\infty$ Hessian bound after parity reflection. \\

\noindent Inside $\mathbb D\cap B_{h_{\mathsf{loc}}}(0)$, we integrate the pointwise Hessian bound against the locally integrable target weight, which produces the factor $P_{\psi,2}^{\rm pt}(\int \Phi_i)^{1/2}$ appearing below. Outside $B_{h_{\mathsf{loc}}}(0)$, the origin singularity is removed, so we use the truncated weighted operator certificate from $\wgray{\Phi_\omega}$ to the relevant target weight. \\

\noindent When the global bound is unavailable, this gives
\begin{align}
C^{\rm w\text{-}ell}_{z,zz}
&:=P_{\psi,2}^{\rm pt}
\left(\int_{\mathbb D\cap B_{h_{\mathsf{loc}}}(0)}\Phi_z\,dr\,dz\right)^{1/2}
+\mathsf C[1_{\mathbb D\setminus B_{h_{\mathsf{loc}}}(0)}L_{zz},\wgray{\Phi_\omega},\wgray{\Phi_z}],\\
C^{\rm w\text{-}ell}_{r,zr}
&:=P_{\psi,2}^{\rm pt}
\left(\int_{\mathbb D\cap B_{h_{\mathsf{loc}}}(0)}\Phi_r\,dr\,dz\right)^{1/2}
+\mathsf C[1_{\mathbb D\setminus B_{h_{\mathsf{loc}}}(0)}L_{zr},\wgray{\Phi_\omega},\wgray{\Phi_r}].
\end{align}

\hfill \\

\subsubsection{Converting source derivatives to energy norms}\label{sssec:ell-source-conversion}

\noindent The operator estimates produce norms such as \(\|D^\sigma\bvar{\delta\omega}\|_{\wgray{\Phi_{\omega,\ell}}}\), with \(|\sigma|=\ell\).  The energy is organized instead by the low-order quantity \(\mathfrak E_0\) and the top-order quantity \(\mathfrak H_k\).   \\ 

\noindent Define \(\Phi_{\omega,0}:=\Phi_\omega\).  For \(0\le\ell\le k\), set
\begin{equation}\label{eq:ell-input-conversion-factor}
\mathcal I_\ell:=\mathcal I_\ell^0+\mu_k^{-1/2}\mathcal I_\ell^k \qquad  \text{where } \ (\mathcal I_\ell^0, \, \mathcal I_\ell^k)
:=
\begin{cases}
(1,0), & \ell=0,\\[1mm]
\bigl(I_{\omega,\ell,k}(\eta_{\omega,\ell}), \,
\eta_{\omega,\ell}\sqrt{k+1}\bigr), & 0<\ell<k,\\[1mm]
(0,1), & \ell=k,
\end{cases}
.
\end{equation}

\hfill 

\noindent Then, we can show that for every \(|\sigma|=\ell\),
\begin{equation}\label{eq:ell-interpolation-energy-conversion}
\|D^\sigma\bvar{\delta\omega}\|_{\wgray{\Phi_{\omega,\ell}}}
\le \mathcal I_\ell^0\mathfrak E_0
+\mathcal I_\ell^k\mathfrak H_k
\le \mathcal I_\ell\mathfrak E_k.
\end{equation}
The endpoint cases are direct.  For \(\ell=0\), the norm is contained in the low-order energy.  For \(\ell=k\), the norm is one of the top-order norms.  For \(0<\ell<k\), the interpolation estimate \eqref{eq:ns-lower-from-top} gives the low-order coefficient \(I_{\omega,\ell,k}(\eta_{\omega,\ell})\) and the high-order coefficient \(\eta_{\omega,\ell}\).  The factor \(\sqrt{k+1}\) converts the sum of the \(k+1\) top-order two-dimensional derivatives into the square-sum in \(\mathfrak H_k\).  The last inequality in \eqref{eq:ell-interpolation-energy-conversion} uses \(\mathfrak E_0\le\mathfrak E_k\) and \(\mathfrak H_k\le\mu_k^{-1/2}\mathfrak E_k\).

\hfill \\

\subsubsection{Low-order profile-multiplied terms}\label{sssec:ell-low-order-rows}

\noindent The low-order linear estimate contains profile-multiplied streamfunction terms, of the form
\begin{equation}
\bvar{\delta\omega}\longmapsto a(r,z)T[\bvar{\delta\psi}].
\end{equation}

\noindent Here $T$ denotes the streamfunction operator used in the corresponding row:
\begin{equation}
\{{\rm Id},\partial_r,\partial_z,\partial_{rr},\partial_{rz},\partial_{zz}\}.
\end{equation} 

\hfill 

\noindent The target component and term label are recorded in
\begin{equation}\label{eq:ell-low-order-stream-row-labels}
\mathfrak S_{\mathcal L,\psi}
:=\{(\omega,4), \, (\omega,5), \, (\omega,6),
(r,4), \, (r,5), \, (r,6), \, (r,7), \, (r,8), \, 
(z,4), \, (z,5), \, (z,6), \, (z,7), \, (z,8)\}.
\end{equation}

\hfill 

\paragraph{Row reduction.}
\noindent Fix one multiplier \(a\) and one expression \(T\).  We rewrite the row as
\begin{equation}\label{eq:ell-low-row-representation}
aT[\bvar{\delta\psi}]
=\sum_{n=1}^{N_{a,T}}
\mathsf T_{L_n}D^{\sigma_n}\bvar{\delta\omega}.
\end{equation}

\hfill 

\noindent We first expand the displayed expression \(T[\bvar{\delta\psi}]\). 
\begin{itemize}
\item If the expanded expression contains \(\partial_{zz}\bvar{\delta\psi}\) or \(\partial_{zr}\bvar{\delta\psi}\), we replace it by the complete Hessian identities \eqref{eq:ell-zz-pv-local}--\eqref{eq:ell-zr-pv-local}.  In the \(zz\) case, this includes the local multiplication term \(-\bvar{\delta\omega}/5\).
\item If extra \(z\)-derivatives appear on top of those Hessians, we use \eqref{eq:ell-zz-further-z-derivatives} and \eqref{eq:ell-zr-further-z-derivatives}.  These identities put the extra \(z\)-derivatives on the source \(\bvar{\delta\omega}\), leaving the same complete Hessian operators.
\item If the derivative expression still contains radial derivatives of \(\bvar{\delta\psi}\) that are not handled by the complete \(zz\) or \(zr\) identities, we use the previously established identities that rewrite these radial derivatives in terms of operators acting on derivatives of \(\bvar{\delta\omega}\). This expresses the remaining streamfunction derivative as a finite sum of terms involving \(\bvar{\delta\omega}\).
\end{itemize}

\hfill 

\paragraph{Elliptic row constants.}
\noindent Let \(\ell_n=|\sigma_n|\).  Define
\begin{equation}\label{eq:ell-low-row-constants}
\mathcal K_i^s(a,T)
:=\sum_{n=1}^{N_{a,T}}
\mathsf C [L_n,\wgray{\Phi_{\omega,\ell_n}},\wgray{\Phi_i}]
\mathcal I_{\ell_n}^s,
\qquad s\in\{0,k\}.
\end{equation}
Applying the operator norm to \eqref{eq:ell-low-row-representation} gives
\begin{equation}
\|aT[\bvar{\delta\psi}]\|_{\wgray{\Phi_i}}
\le
\sum_{n=1}^{N_{a,T}}
\mathsf C [L_n,\wgray{\Phi_{\omega,\ell_n}},\wgray{\Phi_i}]
\|D^{\sigma_n}\bvar{\delta\omega}\|_{\wgray{\Phi_{\omega,\ell_n}}}.
\end{equation}
Using the first inequality in \eqref{eq:ell-interpolation-energy-conversion} term by term yields
\begin{equation}
\|aT[\bvar{\delta\psi}]\|_{\wgray{\Phi_i}}
\le \mathcal K_i^0(a,T)\mathfrak E_0+
\mathcal K_i^k(a,T)\mathfrak H_k.
\end{equation}
The constant \(\mathcal K_i^0(a,T)\) multiplies the low-order energy \(\mathfrak E_0\), and \(\mathcal K_i^k(a,T)\) multiplies the top-order energy~\(\mathfrak H_k\). \\

\noindent These are the constants inserted directly into the low-order linear estimate for component \(i\).  For the velocity combinations \(r\partial_z\bvar{\delta\psi}\) and \(2\bvar{\delta\psi}+r\partial_r\bvar{\delta\psi}\), the certificate first reduces the entire displayed combination into its finite list of source terms.  Only after that reduction do we apply the triangle inequality to the resulting finite sum.  This avoids estimating \(\partial_z\bvar{\delta\psi}\), \(\bvar{\delta\psi}\), and \(\partial_r\bvar{\delta\psi}\) separately in a way that would lose cancellations already present in the combined row. 

\hfill \\

\vspace{10pt}
\subsubsection{High-order streamfunction terms}\label{sssec:ell-high-order-constants}

\noindent Fix \((T,i,\alpha)\in\mathcal J_{\rm ell}\) and put \(m=|\alpha|\).  We expand \(D^\alpha T[\bvar{\delta\psi}]\), rewrite every streamfunction derivative in terms of \(\bvar{\delta\omega}\), estimate each reduced term, and add the contributions.  The target estimate is
\begin{equation}\label{eq:ell-high-final-ineq}
\|D^\alpha T[\bvar{\delta\psi}]\|_{\wgray{\Phi_{i,m}}}
\le C^{\rm w\text{-}ell}_{T,i,m,\alpha}\mathfrak E_k.
\end{equation}
\noindent After expansion, we reduce each term to one of three forms: a direct derivative of \(\bvar{\delta\omega}\), a complete \(zz\) or \(zr\) Hessian operator, or a radial-transfer operator acting on a derivative of \(\bvar{\delta\omega}\). In each case, we identify the remaining source derivative order \(\ell\), apply the corresponding weighted bound, and multiply by the conversion factor \(\mathcal I_\ell\) from \eqref{eq:ell-interpolation-energy-conversion}. Summing these contributions gives the row constant.\\

\paragraph{Expanding the row.}
\noindent We apply \(D^\alpha\) to the displayed expression \(T[\bvar{\delta\psi}]\).  If \(T=3\partial_r+r\partial_{rr}\), we first use the elliptic equation:
\begin{equation}\label{eq:ell-Spsi-pde-replacement}
(3\partial_r+r\partial_{rr})\bvar{\delta\psi}
=-r\bvar{\delta\omega}-r\partial_{zz}\bvar{\delta\psi}.
\end{equation}
Since \(D^\alpha r=rD^\alpha+\alpha_rD^{\alpha-e_r}\), this gives
\begin{align}
D^\alpha(3\partial_r+r\partial_{rr})\bvar{\delta\psi}
={}&-rD^\alpha\bvar{\delta\omega}
-\alpha_rD^{\alpha-e_r}\bvar{\delta\omega}
-rD^{\alpha+2e_z}\bvar{\delta\psi}
-\alpha_rD^{\alpha-e_r+2e_z}\bvar{\delta\psi},
\label{eq:ell-Tstar-full-expansion}
\end{align}
where terms with \(\alpha-e_r\) are omitted when \(\alpha_r=0\).  For \(T=r\partial_z\), the product rule gives
\begin{equation}
D^\alpha(r\partial_z\bvar{\delta\psi})
=rD^{\alpha+e_z}\bvar{\delta\psi}
+\alpha_rD^{\alpha-e_r+e_z}\bvar{\delta\psi}.
\end{equation}
For \(T=2+r\partial_r\), it gives
\begin{equation}
D^\alpha(2\bvar{\delta\psi}+r\partial_r\bvar{\delta\psi})
=(2+\alpha_r)D^\alpha\bvar{\delta\psi}
+rD^{\alpha+e_r}\bvar{\delta\psi}.
\end{equation}
When all derivatives in \(D^\alpha\) fall on the streamfunction factor, the factor \(r\) remains and produces the terms such as \(rD^{\alpha+e_z}\bvar{\delta\psi}\).  When one of the \(\alpha_r\) radial derivatives falls on \(r\), it produces the additional lower-order term with coefficient \(\alpha_r\).  Since higher derivatives of \(r\) vanish, there are no further product-rule terms.  The replacement \eqref{eq:ell-Spsi-pde-replacement} is used for \(3\partial_r+r\partial_{rr}\) before applying the estimates, because the elliptic equation converts this radial second-order expression into a direct source term and a \(zz\) Hessian term.  \\ 

\paragraph{Direct source terms.}
\noindent A term of the form \(bD^\sigma\bvar{\delta\omega}\), with \(|\sigma|=\ell\), is estimated by comparing the target and input weights.  Squaring the weighted norm gives
\begin{align}
\|bD^\sigma\bvar{\delta\omega}\|_{\wgray{\Phi_{i,m}}}^2
&=\int_\mathbb{D} |D^\sigma\bvar{\delta\omega}|^2 |b|^2\wgray{\Phi_{i,m}}\,dr\,dz =\int_\mathbb{D} |D^\sigma\bvar{\delta\omega}|^2\wgray{\Phi_{\omega,\ell}}
\left(|b|^2\frac{\wgray{\Phi_{i,m}}}{\wgray{\Phi_{\omega,\ell}}}\right)\,dr\,dz \nonumber \\
&\le
\left(\operatorname*{ess\,sup}_\mathbb{D} |b|^2\frac{\wgray{\Phi_{i,m}}}{\wgray{\Phi_{\omega,\ell}}}\right)
\|D^\sigma\bvar{\delta\omega}\|_{\wgray{\Phi_{\omega,\ell}}}^2.
\end{align}
Taking square roots yields
\begin{equation}\label{eq:ell-direct-source-bound}
\|bD^\sigma\bvar{\delta\omega}\|_{\wgray{\Phi_{i,m}}}
\le
\operatorname*{ess\,sup}_\mathbb{D}
|b|\left(\frac{\wgray{\Phi_{i,m}}}{\wgray{\Phi_{\omega,\ell}}}\right)^{1/2}
\|D^\sigma\bvar{\delta\omega}\|_{\wgray{\Phi_{\omega,\ell}}}.
\end{equation}

\noindent  After this quantity is bounded, \eqref{eq:ell-interpolation-energy-conversion} converts the source norm to \(\mathfrak E_k\). The contribution of the term is therefore the certified multiplier constant times \(\mathcal I_\ell\). \\

\paragraph{Complete \(zz\) and \(zr\) terms.}
\noindent Terms with \(\partial_z^{\ell+2}\bvar{\delta\psi}\) are rewritten as
\begin{equation}\label{eq:ell-high-zz-complete}
a\partial_z^{\ell+2}\bvar{\delta\psi}
=-\frac a5\partial_z^\ell\bvar{\delta\omega}
+a\operatorname{p.v.}\mathsf T_{L_{zz}^{\mathrm{off}}}
  (\partial_z^\ell\bvar{\delta\omega}).
\end{equation}
Terms with \(\partial_r\partial_z^{\ell+1}\bvar{\delta\psi}\) are rewritten as
\begin{equation}\label{eq:ell-high-zr-complete}
a\partial_r\partial_z^{\ell+1}\bvar{\delta\psi}
=a\operatorname{p.v.}\mathsf T_{L_{zr}^{\mathrm{off}}}
  (\partial_z^\ell\bvar{\delta\omega}).
\end{equation}
The first line is exactly the complete operator \(aL_{zz}\) applied to \(\partial_z^\ell\bvar{\delta\omega}\): the term \(-a\partial_z^\ell\bvar{\delta\omega}/5\) is the local part, and the remaining term is the principal-value part. \\

\noindent The second line is \(aL_{zr}\) applied to \(\partial_z^\ell\bvar{\delta\omega}\), which has no local part.\\

\noindent Therefore a complete Hessian term contributes
\begin{equation}\label{eq:ell-high-summand-constant-clear}
\mathsf C [aL,\wgray{\Phi_{\omega,\ell}},\wgray{\Phi_{i,m}}]
\mathcal I_\ell,
\end{equation}
where \(L=L_{zz}\) or \(L=L_{zr}\).    

\hfill \\ 

\paragraph{Radial-transfer terms.}
\noindent A remaining term has a radial derivative pattern that is not one of the complete \(zz\) or \(zr\) forms above.  Concretely, after the product-rule expansion, such a term looks like \(aD^\eta\bvar{\delta\psi}\) where \(D^\eta\) contains radial derivatives but cannot be written as \(\partial_z^{\ell+2}\) or \(\partial_r\partial_z^{\ell+1}\). Such a term is reduced via
\begin{equation}\label{eq:ell-radial-transfer-representation}
aD^\eta\bvar{\delta\psi}
=\sum_{q=1}^{N_{a,\eta}}
\mathsf T_{L_{a,\eta,q}^{\rm rad}}
D^{\sigma_q}\bvar{\delta\omega},
\qquad |\sigma_q|\le |\eta|-2.
\end{equation}

\noindent Each summand in \eqref{eq:ell-radial-transfer-representation} is then estimated exactly like the other operator terms.  If the summand leaves the source derivative \(D^{\sigma_q}\bvar{\delta\omega}\), its contribution is
\begin{equation}\label{eq:ell-radial-transfer-contribution}
\mathsf C [L_{a,\eta,q}^{\rm rad},
\wgray{\Phi_{\omega,|\sigma_q|}},
\wgray{\Phi_{i,m}}]
\mathcal I_{|\sigma_q|}.
\end{equation}
The operator norm controls the map from the input weight \(\wgray{\Phi_{\omega,|\sigma_q|}}\) to the output weight \(\wgray{\Phi_{i,m}}\), and \(\mathcal I_{|\sigma_q|}\) converts the resulting source norm to \(\mathfrak E_k\) by \eqref{eq:ell-interpolation-energy-conversion}. 

\hfill \\ 

\paragraph{Assembling the row constant.}
\noindent We reduce every term in the expanded row by one of the three rules above and let \(c_{\rm term}\) denote the resulting contribution, including the appropriate interpolation factor.  Define
\begin{equation}\label{eq:ell-detailed-high-final-constant}
C^{\rm w\text{-}ell}_{T,i,m,\alpha}
:=\sum_{\text{terms in the reduced row}}c_{\rm term}.
\end{equation}
The triangle inequality gives the sum over reduced terms, \(\mathsf C \) gives the weighted operator bounds, and \eqref{eq:ell-interpolation-energy-conversion} converts each source norm to \(\mathfrak E_k\).  Therefore \eqref{eq:ell-detailed-high-final-constant} proves \eqref{eq:ell-high-final-ineq} for the fixed row.\\

\hfill \\ 

\subsubsection{The differentiated unweighted constant \texorpdfstring{\(C_{\rm ell}^{\rm diff}\)}{Cell diff}}\label{sssec:ell-unweighted-constant}

\noindent \noindent The constant $C_{\rm ell}^{\rm diff}$ comes from the standard unweighted elliptic estimate, before any weighted norms are introduced. It reflects the two-derivative smoothing of the Poisson equation: derivatives of the streamfunction of order at least two are controlled by derivatives of the vorticity source of two orders lower. \\

\noindent The construction is the same as for the weighted terms, but with both weights equal to one.  For each \(\alpha\) with \(2\le |\alpha|\le k+2\), build a finite reduction
\begin{equation}\label{eq:ell-Cell-alpha-reduction}
D^\alpha\bvar{\delta\psi}
=\sum_{n=1}^{N_\alpha}
\mathsf T_{L_{\alpha,n}}D^{\sigma_{\alpha,n}}\bvar{\delta\omega},
\qquad |\sigma_{\alpha,n}|\le |\alpha|-2.
\end{equation}
The reduction follows the same hierarchy as the weighted terms: complete \(zz\) identities when two \(z\)-derivatives are available, complete \(zr\) identities when one radial and one \(z\)-derivative are available, and radial-transfer identities for the remaining radial cases. \\ 

\noindent For each summand indexed by $n$, define
\begin{equation}\label{eq:ell-Cell-alpha-term-constant}
A_{\alpha,n}:=\mathsf C[L_{\alpha,n},1,1].
\end{equation}
Thus $A_{\alpha,n}$ is the constant used to bound the contribution of the operator $L_{\alpha,n}$. If a summand is a local term involving no singular-integral operator, its contribution is bounded directly by the absolute value of its coefficient. \\

\noindent We then collect all summands involving the same derivative $D^\sigma$ of the source and define
\begin{equation}\label{eq:ell-Cell-alpha-sigma-constant}
A_{\alpha,\sigma}
:=\sum_{\{n:\sigma_{\alpha,n}=\sigma\}}A_{\alpha,n}.
\end{equation}
Hence $A_{\alpha,\sigma}$ is the total contribution of all terms in which the derivative falling on the source is $D^\sigma$.\\

\noindent Then
\begin{equation}\label{eq:ell-Cell-alpha-q-bound}
\|D^\alpha\bvar{\delta\psi}\|_{L^2(\mathbb{D})}
\le
\sum_{|\sigma|\le |\alpha|-2}
A_{\alpha,\sigma}
\|D^\sigma\bvar{\delta\omega}\|_{L^2(\mathbb{D})}.
\end{equation}

\noindent Summing \eqref{eq:ell-Cell-alpha-q-bound} over \(\alpha\) and grouping by \(\sigma\) gives the single constant
\begin{equation}\label{eq:ell-Cell-diff-final}
C_{\rm ell}^{\rm diff}
:=\max_{|\sigma|\le k}
\sum_{\substack{2\le |\alpha|\le k+2}}
A_{\alpha,\sigma}.
\end{equation}
Indeed,
\begin{align}
\sum_{m=2}^{k+2}\sum_{|\alpha|=m}
\|D^\alpha\bvar{\delta\psi}\|_{L^2(\mathbb{D})}
&\le
\sum_{|\sigma|\le k}
\left(\sum_{2\le |\alpha|\le k+2}A_{\alpha,\sigma}\right)
\|D^\sigma\bvar{\delta\omega}\|_{L^2(\mathbb{D})}
\le C_{\rm ell}^{\rm diff}
\sum_{m=0}^{k}\sum_{|\sigma|=m}
\|D^\sigma\bvar{\delta\omega}\|_{L^2(\mathbb{D})}.
\label{eq:ell-Cell-differentiated-part}
\end{align}

\hfill  \\

\hfill 

\subsubsection{Proof of the streamfunction elliptic estimates}\label{sssec:ell-proof}

\noindent We recall \Cref{thm:stream-elliptic-used}. Suppose that the point functionals and the finitely many weighted operators appearing below have been certified by the above construction. Then the selected origin constants satisfy
\begin{align}
|\bvar{\delta\psi}(0,0)|
&\le P_{\psi,0}^{0}\mathfrak E_0+P_{\psi,0}^{k}\mathfrak H_k
\le P_{\psi,0}\mathfrak E_k,\\
|\partial_z\bvar{\delta\psi}(0,0)|
&\le P_{\psi_z,0}^{0}\mathfrak E_0+P_{\psi_z,0}^{k}\mathfrak H_k
\le P_{\psi_z,0}\mathfrak E_k.
\end{align}
where
\begin{equation}
P_{\psi,0}:=P_{\psi,0}^{0}+\mu_k^{-1/2}P_{\psi,0}^{k},
\qquad
P_{\psi_z,0}:=P_{\psi_z,0}^{0}+\mu_k^{-1/2}P_{\psi_z,0}^{k}.
\end{equation}

\hfill 

\noindent The differentiated unweighted estimate is
\begin{equation}
\sum_{m=2}^{k+2}\sum_{|\alpha|=m}
\|D^\alpha\bvar{\delta\psi}\|_{L^2(\mathbb{D})}
\le C_{\rm ell}^{\rm diff}
\sum_{m=0}^{k}\sum_{|\alpha|=m}
\|D^\alpha\bvar{\delta\omega}\|_{L^2(\mathbb{D})}.
\end{equation}
The two low-order weighted Hessian estimates are
\begin{equation}
\|\partial_{zz}\bvar{\delta\psi}\|_{\wgray{\Phi_z}}
\le C^{\rm w\text{-}ell}_{z,zz}\mathfrak E_k,
\qquad
\|\partial_{zr}\bvar{\delta\psi}\|_{\wgray{\Phi_r}}
\le C^{\rm w\text{-}ell}_{r,zr}\mathfrak E_k.
\end{equation}
Finally, for every \((T,i,\alpha)\in\mathcal J_{\rm ell}\), with \(m=|\alpha|\),
\begin{equation}
\|D^\alpha T[\bvar{\delta\psi}]\|_{\wgray{\Phi_{i,m}}}
\le C^{\rm w\text{-}ell}_{T,i,m,\alpha}\mathfrak E_k.
\end{equation}

\vspace{2mm}

\begin{grayproof}
\noindent\textbf{Roadmap.}  We convert the certified constants into the estimates in the statement.  First we use the origin kernels and the Cauchy--Schwarz inequality to control the two point values.  Next we use the finite reductions to prove the differentiated unweighted estimate.  Finally we apply the complete Hessian operator norms and the row-by-row reductions to obtain the weighted estimates.\\ 

\noindent\textbf{Step 1: origin values.}  Proposition~\ref{prop:five-dimensional-kernel} gives the Newtonian representation in \(\mathbb R^5\).  \\

\noindent Evaluating the reduced kernel at the origin gives the two scalar kernels
\begin{align}
k_\psi(s,\zeta)&:=K_0(0,0,s,\zeta),\label{eq:ell-proof-origin-kernel-psi}\\
k_{\psi_z}(s,\zeta)&:=\partial_zK_0(0,0,s,\zeta).\label{eq:ell-proof-origin-kernel-psiz}
\end{align}

 \hfill 

\noindent These direct kernels are computed explicitly in \eqref{eq:ell-origin-kernel-psi}--\eqref{eq:ell-origin-kernel-psiz}.  Therefore
\begin{align}
\bvar{\delta\psi}(0,0)
&=\int_0^\infty\int_{\mathbb R}k_\psi(s,\zeta)\bvar{\delta\omega}(s,\zeta)\,d\zeta\,ds,\\
\partial_z\bvar{\delta\psi}(0,0)
&=\int_0^\infty\int_{\mathbb R}k_{\psi_z}(s,\zeta)\bvar{\delta\omega}(s,\zeta)\,d\zeta\,ds.
\end{align}

\hfill 

\noindent We insert \(1=\wgray{\Phi_\omega}^{-1/2}\wgray{\Phi_\omega}^{1/2}\) in each integral.  The Cauchy--Schwarz inequality gives
\begin{align}
|\bvar{\delta\psi}(0,0)|
&\le \left\|k_\psi\wgray{\Phi_\omega}^{-1/2}\right\|_{L^2(\mathbb{D})}
\|\bvar{\delta\omega}\|_{\wgray{\Phi_\omega}},\\
|\partial_z\bvar{\delta\psi}(0,0)|
&\le \left\|k_{\psi_z}\wgray{\Phi_\omega}^{-1/2}\right\|_{L^2(\mathbb{D})}
\|\bvar{\delta\omega}\|_{\wgray{\Phi_\omega}}.
\end{align}

 \hfill 

\noindent  The certification subsection stores these direct bounds as
\begin{equation}
P_{\psi,0}^{0}=P_{\psi,0}^{\rm dir},\qquad
P_{\psi,0}^{k}=0,
\qquad
P_{\psi_z,0}^{0}=P_{\psi_z,0}^{\rm dir},\qquad
P_{\psi_z,0}^{k}=0.
\end{equation}

\hfill 

\noindent Therefore the point estimates are
\begin{align}
|\bvar{\delta\psi}(0,0)|
&\le P_{\psi,0}^{0}\mathfrak E_0+P_{\psi,0}^{k}\mathfrak H_k,\\
|\partial_z\bvar{\delta\psi}(0,0)|
&\le P_{\psi_z,0}^{0}\mathfrak E_0+P_{\psi_z,0}^{k}\mathfrak H_k.
\end{align}

\hfill 

\noindent Finally,
\begin{equation}
\mathfrak E_0\le \mathfrak E_k,
\qquad
\mathfrak H_k\le \mu_k^{-1/2}\mathfrak E_k.
\end{equation}

\hfill 

\noindent Substituting these two inequalities gives \eqref{eq:ell-point-theorem-psi} and \eqref{eq:ell-point-theorem-psiz}.

\hfill \\ 

\noindent\textbf{Step 2: differentiated unweighted estimate.}  Fix a multi-index \(\alpha\) with \(2\le |\alpha|\le k+2\).  The reduction described in Appendix~\ref{sssec:ell-unweighted-constant} writes the derivative of the streamfunction as a finite sum of source derivatives:
\begin{equation}
D^\alpha\bvar{\delta\psi}
=\sum_{n=1}^{N_\alpha}\mathsf T_{L_{\alpha,n}}D^{\sigma_{\alpha,n}}\bvar{\delta\omega},
\qquad |\sigma_{\alpha,n}|\le |\alpha|-2.
\end{equation}
The restriction \(|\sigma_{\alpha,n}|\le |\alpha|-2\) reflects the elliptic gain of two derivatives. \\

\noindent Each summand is certified in unweighted \(L^2\) by
\begin{equation}
\|\mathsf T_{L_{\alpha,n}}D^{\sigma_{\alpha,n}}\bvar{\delta\omega}\|_{L^2(\mathbb{D})}
\le A_{\alpha,n}
\|D^{\sigma_{\alpha,n}}\bvar{\delta\omega}\|_{L^2(\mathbb{D})}.
\end{equation}

\hfill 

\noindent We apply the triangle inequality and group the terms that contain the same source derivative.  With
\begin{equation}
A_{\alpha,\sigma}:=
\sum_{\{n:\ \sigma_{\alpha,n}=\sigma\}}A_{\alpha,n},
\end{equation}
we obtain
\begin{equation}\label{eq:ell-proof-alpha-row-bound}
\|D^\alpha\bvar{\delta\psi}\|_{L^2(\mathbb{D})}
\le
\sum_{|\sigma|\le |\alpha|-2}
A_{\alpha,\sigma}
\|D^\sigma\bvar{\delta\omega}\|_{L^2(\mathbb{D})}.
\end{equation}

\hfill 

\noindent Now sum \eqref{eq:ell-proof-alpha-row-bound} over all \(\alpha\) with \(2\le |\alpha|\le k+2\).  Since \(|\sigma|\le |\alpha|-2\), all source derivatives have order at most \(k\).  Exchanging the finite sums gives
\begin{align}
\sum_{m=2}^{k+2}\sum_{|\alpha|=m}
\|D^\alpha\bvar{\delta\psi}\|_{L^2(\mathbb{D})}
&\le
\sum_{|\sigma|\le k}
\bigg(
\sum_{\substack{2\le |\alpha|\le k+2}}
A_{\alpha,\sigma}
\bigg)
\|D^\sigma\bvar{\delta\omega}\|_{L^2(\mathbb{D})}.
\end{align}

\hfill 

\noindent By definition \eqref{eq:ell-Cell-diff-final}, each coefficient in parentheses is at most \(C_{\rm ell}^{\rm diff}\).  Therefore
\begin{align}
\sum_{m=2}^{k+2}\sum_{|\alpha|=m}
\|D^\alpha\bvar{\delta\psi}\|_{L^2(\mathbb{D})}
&\le C_{\rm ell}^{\rm diff}
\sum_{|\sigma|\le k}
\|D^\sigma\bvar{\delta\omega}\|_{L^2(\mathbb{D})} = C_{\rm ell}^{\rm diff}
\sum_{m=0}^{k}\sum_{|\sigma|=m}
\|D^\sigma\bvar{\delta\omega}\|_{L^2(\mathbb{D})}.
\end{align}

\hfill \\ 

\noindent\textbf{Step 3: low-order weighted Hessians.}  We use the principal-value and off-diagonal-kernel conventions fixed in Appendix~\ref{subsec:obtaining-elliptic-constants}.  The distributional Hessian identity in five dimensions is
\begin{equation}
\partial_{X_iX_j}\Gamma_5
=\operatorname{p.v.}(\partial_{X_iX_j}^{\mathrm{off}}\Gamma_5)-\frac{\delta_{ij}}5\delta_0.
\end{equation}

\hfill 

\noindent For the \(zz\) component this gives a principal-value operator and the local term \(-\bvar{\delta\omega}/5\).  For the mixed \(zr\) component the local term vanishes.  Angular reduction gives \eqref{eq:ell-zz-pv-local} and \eqref{eq:ell-zr-pv-local}.  Proposition~\ref{prop:ell-block-certificate}, applied with input weight \(\wgray{\Phi_\omega}\) and output weights \(\wgray{\Phi_z}\), \(\wgray{\Phi_r}\), gives
\begin{equation}
\|\partial_{zz}\bvar{\delta\psi}\|_{\wgray{\Phi_z}}
\le \mathsf C [L_{zz},\wgray{\Phi_\omega},\wgray{\Phi_z}]
\|\bvar{\delta\omega}\|_{\wgray{\Phi_\omega}},
\end{equation}
\begin{equation}
\|\partial_{zr}\bvar{\delta\psi}\|_{\wgray{\Phi_r}}
\le \mathsf C [L_{zr},\wgray{\Phi_\omega},\wgray{\Phi_r}]
\|\bvar{\delta\omega}\|_{\wgray{\Phi_\omega}}.
\end{equation}

\hfill 

\noindent If the global weighted Hessian constant for a component is finite, we obtain \eqref{eq:ns-Cwell} directly. Otherwise, we localize only that component. For the $zz$ component, we use the pointwise Hessian bound on $\mathbb D\cap B_{h_{\mathsf{loc}}}(0)$ and the truncated weighted $L_{zz}$ operator bound on $\mathbb D\setminus B_{h_{\mathsf{loc}}}(0)$:
\begin{equation}
\|1_{\mathbb D\cap B_{h_{\mathsf{loc}}}(0)}\partial_{zz}\bvar{\delta\psi}\|_{\Phi_z}
\le P_{\psi,2}^{\rm pt}
\left(\int_{\mathbb D\cap B_{h_{\mathsf{loc}}}(0)}\Phi_z\,dr\,dz\right)^{1/2}\mathfrak E_k,
\end{equation}
\begin{equation}
\|1_{\mathbb D\setminus B_{h_{\mathsf{loc}}}(0)}\partial_{zz}\bvar{\delta\psi}\|_{\Phi_z}
\le \mathsf C[1_{\mathbb D\setminus B_{h_{\mathsf{loc}}}(0)}L_{zz},\Phi_\omega,\Phi_z]\mathfrak E_k.
\end{equation}

\hfill 

\noindent Adding these two estimates gives the required $zz$ bound. The same argument gives the $zr$ bound, with $L_{zr}$ and $\wgray{\Phi_r}$ in place of $L_{zz}$ and $\wgray{\Phi_z}$.

\hfill \\ 

\hfill 

\noindent\textbf{Step 4: differentiated weighted terms.} Fix \((T,i,\alpha)\in\mathcal J_{\rm ell}\) and set \(m=|\alpha|\). \\

\noindent Appendix~\ref{sssec:ell-high-order-constants} reduces \(D^\alpha T[\bvar{\delta\psi}]\) to a finite sum of
\begin{equation}
bD^\sigma\bvar{\delta\omega},
\qquad
 aL_{zz}(D^\sigma\bvar{\delta\omega})\hbox{ or }aL_{zr}(D^\sigma\bvar{\delta\omega}),
\qquad
 \mathsf T_{L^{\rm rad}}D^\sigma\bvar{\delta\omega}.
\end{equation}
\noindent For each term, we first identify the source order \(\ell=|\sigma|\). A direct term is controlled by certifying the full multiplier
\(
|b|(\Phi_{i,m}/\Phi_{\omega,\ell})^{1/2}
\)
in \(L^\infty\). A Hessian term uses the complete weighted operator norm of \(aL_{zz}\) or \(aL_{zr}\). A radial-transfer term uses the certified operator norm of the corresponding transfer identity. In every case, \eqref{eq:ell-interpolation-energy-conversion} then converts the source norm to \(\mathfrak E_k\).

\hfill 

\noindent The row constant \(C^{\rm w\text{-}ell}_{T,i,m,\alpha}\) is the sum of these finitely many contributions, which proves \eqref{eq:ns-Cwell-high}.
\end{grayproof}

\clearpage

\subsection{Modulation constants}\label{sssec:modulation-constants}

\noindent The modulation formulas are
\begin{equation}
\bvar{\delta C}=-2\bvar{\delta\psi}_0,
\qquad
\bvar{\delta c_u}=-2(\partial_z\bvar{\delta\psi})_0.
\end{equation}
The required point estimates are
\begin{equation}
|\bvar{\delta\psi}_0|\le P_{\psi,0}\mathfrak E_k,
\qquad
|(\partial_z\bvar{\delta\psi})_0|\le P_{\psi_z,0}\mathfrak E_k.
\end{equation}
Set
\begin{align}
\label{eq:ns-MC}
M_C^{(1)}&:=2P_{\psi,0},\\
\label{eq:ns-Mcu}
M_{c_u}^{(1)}&:=2P_{\psi_z,0}.
\end{align}
Then
\begin{equation}
|\bvar{\delta C}|\le M_C^{(1)}\mathfrak E_k,
\qquad
|\bvar{\delta c_u}|\le M_{c_u}^{(1)}\mathfrak E_k.
\label{eq:euler-modulation-direct-bounds-add}
\end{equation}
Thus the modulation coefficients are linear and are used directly in the linear operator. No fixed-profile modulation contribution is added to the nonlinear operator.

\clearpage

\subsection{Ray estimate and global pointwise control}

\noindent This subsection introduces two tools used later: a ray estimate for the fixed-profile linear terms containing undifferentiated \(\bvar{\delta u}\), and the global pointwise bounds used in the nonlinear estimates.  The origin values required by the modulation formulas are $\bvar{\delta\psi}_0$ and $(\partial_z\bvar{\delta\psi})_0$.  \\

\subsubsection{Ray estimate for fixed-profile terms}

\noindent The following proposition is used only for the fixed-profile linear terms containing undifferentiated \(\bvar{\delta u}\). \\ 

\begin{proposition}[Ray estimate for fixed-profile rows]\label{prop:fixed-profile-ray-estimate}
\noindent Fix \(i\in\{\omega,r,z\}\), and let \(a\) be measurable.  Define
\begin{equation}
\mathcal C_{{\rm ray},i}(a)
:=\operatorname*{ess\,sup}_{x\in\mathbb{D}\setminus\{(0,0)\}}
\frac{\rho(x)}{\min\{\wgray{\Phi_r}(x),\wgray{\Phi_z}(x)\}}
\int_0^1 t^{-4}|a(x/t)|^2\wgray{\Phi_i}(x/t)\,dt.
\label{eq:fixed-profile-ray-constant}
\end{equation}
\noindent If \(\mathcal C_{{\rm ray},i}(a)<\infty\), then every \(f\in C^1(\mathbb{D})\) with \(f(0,0)=0\) and
\(\|\partial_r f\|_{\wgray{\Phi_r}}+\|\partial_z f\|_{\wgray{\Phi_z}}<\infty\) satisfies
\begin{equation}
\|af\|_{\wgray{\Phi_i}}
\le \mathcal C_{{\rm ray},i}(a)^{1/2}
\left(\|\partial_rf\|_{\wgray{\Phi_r}}^2+
\|\partial_zf\|_{\wgray{\Phi_z}}^2\right)^{1/2}.
\label{eq:fixed-profile-ray-ineq}
\end{equation}
\end{proposition}

\begin{grayproof}
 For $x=(r,z)\in\mathbb{D}$, the line segment $tx$, $0\le t\le1$, stays in $\mathbb{D}$.  Since $f(0,0)=0$, the fundamental theorem of calculus gives
\begin{equation}
f(x)=\int_0^1\frac{d}{dt}f(tx)\,dt
=\int_0^1 x\cdot\nabla f(tx)\,dt.
\end{equation}
The Cauchy--Schwarz inequality in the $t$ variable yields
\begin{equation}
|f(x)|^2
\le |x|^2\int_0^1|\nabla f(tx)|^2\,dt
=\rho(x)\int_0^1|\nabla f(tx)|^2\,dt.
\end{equation}

\hfill

\noindent We multiply by $|a(x)|^2\Phi_i(x)$ and integrate.  Tonelli's theorem and the change of variables $y=tx$ give $dx=t^{-2}dy$ and $\rho(x)=\rho(y)/t^2$, hence
\begin{align}
\int_\mathbb{D} |a(x)f(x)|^2\Phi_i(x)\,dx
&\le \int_0^1\int_\mathbb{D} \rho(x)|\nabla f(tx)|^2|a(x)|^2\Phi_i(x)\,dx\,dt\nonumber\\
&=\int_\mathbb{D} |\nabla f(y)|^2\rho(y)
\left(\int_0^1t^{-4}|a(y/t)|^2\Phi_i(y/t)\,dt\right)dy.
\label{eq:fixed-profile-ray-change-variables}
\end{align}
By \eqref{eq:fixed-profile-ray-constant},
\begin{equation}
\int_\mathbb{D} |af|^2\Phi_i
\le \mathcal C_{{\rm ray},i}(a)
\left(\|\partial_rf\|_{\Phi_r}^2+\|\partial_zf\|_{\Phi_z}^2\right).
\label{eq:fixed-profile-ray-proof-bound}
\end{equation}
Taking square roots gives \eqref{eq:fixed-profile-ray-ineq}.
\end{grayproof}

\hfill

\noindent We will apply Proposition~\ref{prop:fixed-profile-ray-estimate} with \(f=\bvar{\delta u}\) and
\begin{equation}
    (i,a)\in
\left\{
(\omega,2\partial_z\bar u),
(r,2\partial_{rz}\bar\psi),
(z,2\partial_{zz}\bar\psi)
\right\}.
\end{equation}
These choices are exactly the three fixed-profile linear terms containing undifferentiated \(\bvar{\delta u}\).  The global pointwise bound for \(\bvar{\delta u}\) follows separately from the energy weights.

\hfill 

\subsubsection{Global pointwise estimate}

\noindent We now collect the pointwise bounds used in the nonlinear and modulation estimates.  The proposition below gives separate constants for evolved variables, origin values, velocities, and streamfunction multipliers. \\

\begin{proposition}[Factor-specific pointwise estimates from the weighted energy]\label{prop:pointwise-weighted-energy}

\noindent Assume that each low-order weight has the positive lower bound in \eqref{eq:weight_coercivity} and that the velocity weights satisfy the far-field growth condition \eqref{eq:delta-u-tail-growth-main}. Thus the weighted $L^2$ norms control the unweighted $L^2$ norms used after reflection, and the velocity weights provide the tail control needed to recover $\bvar{\delta u}$ from its gradient.

\noindent The reflected whole-plane estimates in Propositions~\ref{prop:whole-plane-reflection}--\ref{prop:whole-plane-standard-constants}, combined with the weighted interpolation estimate when a weighted intermediate norm is needed, and the kernel elliptic estimate in \Cref{thm:stream-elliptic-used}, imply

\begin{equation}
\label{eq:ns-Cpt-global}
\|D^\alpha\bvar{\delta v_i}\|_{L^\infty(\mathbb{D})}
\le C^{\rm pt}_{i,m,\alpha}\mathfrak E_k,
\qquad i\in\{\omega,r,z\},\quad m=|\alpha|\le k-2,
\end{equation}

\noindent The constants \(C^{\rm pt}_{i,m,\alpha}\) collect the reflected Sobolev and parity-reflection constants, the interpolation constants needed for intermediate derivatives, and the coercive lower bound \(C_{\rm floor}\) used to pass from weighted to unweighted \(L^2\). Streamfunction factors are estimated separately below. \\

\noindent In particular, using the origin constants introduced in Appendix~\ref{sssec:modulation-constants},

\begin{equation}
\begin{gathered}
| (\partial_z\bvar{\delta u})_0|
\le P_{u_z,0}\mathfrak E_k,
\quad
|(\partial_z\bvar{\delta\omega})_0|
\le P_{\omega_z,0}\mathfrak E_k,\\
|\bvar{\delta\psi}_0|
\le P_{\psi,0}\mathfrak E_k,
\quad
|(\partial_z\bvar{\delta\psi})_0|
\le P_{\psi_z,0}\mathfrak E_k.
\end{gathered}
\end{equation}

\hfill

\noindent The global velocity and streamfunction factors satisfy
\begin{equation}
\begin{aligned}
\|\bvar{\delta u}\|_{L^\infty}&\le P_u\mathfrak E_k,&
\|\bvar{\delta u^r}\|_{L^\infty}&\le P_{u^r}\mathfrak E_k,&
\|\bvar{\delta u^z}\|_{L^\infty}&\le P_{u^z}\mathfrak E_k,\\
\|\partial_z\bvar{\delta\psi}\|_{L^\infty}&\le P_{\psi,z}\mathfrak E_k,&
\|r\partial_{zz}\bvar{\delta\psi}\|_{L^\infty}&\le P_{\psi,rzz}\mathfrak E_k,&
\|r\partial_{rz}\bvar{\delta\psi}\|_{L^\infty}&\le P_{\psi,rzr}\mathfrak E_k,\\
\|3\partial_r\bvar{\delta\psi}+r\partial_{rr}\bvar{\delta\psi}\|_{L^\infty}&\le P_{\psi,rrr}^{\rm ax}\mathfrak E_k.
\end{aligned}
\end{equation}

  \end{proposition}

\hfill

\begin{grayproof}
\textbf{Pointwise bounds for evolved variables.} 
Fix $i\in\{\omega,r,z\}$ and $|\alpha|=m\le k-2$. \\

\noindent Let \(\mathrm E\) denote the even or odd extension prescribed by the parity of \(\bvar{\delta v_i}\), and write
\begin{equation}
\widetilde{\bvar{\delta v_i}}:=\mathrm E\bvar{\delta v_i}.
\end{equation}
Since $m+2\le k$, Proposition~\ref{prop:whole-plane-standard-constants} with $s=2$ gives
\begin{equation}
\|D^\alpha\bvar{\delta v_i}\|_{L^\infty(\mathbb{D})}
\le
\|D^\alpha\widetilde{\bvar{\delta v_i}}\|_{L^\infty(\mathbb R^2)}
\le
S_{2,2}\|D^\alpha\widetilde{\bvar{\delta v_i}}\|_{H^2(\mathbb R^2)}.
\end{equation}

\hfill 

\noindent Expanding the $H^2$ norm gives
\begin{equation}
\|D^\alpha\widetilde{\bvar{\delta v_i}}\|_{H^2(\mathbb R^2)}
\le
\sum_{|\beta|\le2}\|D^{\alpha+\beta}\widetilde{\bvar{\delta v_i}}\|_{L^2(\mathbb R^2)}.
\end{equation}
By the reflection identity, each term in this finite sum is \(\sqrt2\) times the corresponding half-plane derivative norm.  Terms with \(|\alpha+\beta|=k\) are top-order norms and are bounded by \(\mu_k^{-1/2}\mathfrak E_k\), because \(\mathfrak E_k^2=\mathfrak E_0^2+\mu_k\mathfrak H_k^2\).  Terms with \(|\alpha+\beta|=0\) are low-order energy terms.  Terms with \(0<|\alpha+\beta|<k\) are reduced by the interpolation bound \eqref{eq:ns-all-lower-controlled}, or by the whole-plane interpolation estimate when no component weight is involved.  Adding the finitely many constants from these alternatives gives \(C^{\rm pt}_{i,m,\alpha}\), which proves \eqref{eq:ns-Cpt-global}.

\hfill

\noindent Now, for \(\Phi\in\{\Phi_\omega,\Phi_r,\Phi_z\}\), equation~\eqref{eq:weight_coercivity} gives
\begin{equation}
\|f\|_{L^2(\mathbb{D})}\le C_{\rm floor}^{-1/2}\|f\|_{\Phi}.
\end{equation}

\noindent This gives the required unweighted \(L^2\) norms for the evolved variables and every positive derivative of \(\bvar{\delta u}\).

\hfill \\ 

\noindent\textbf{The undifferentiated }$\bvar{\delta u}$\textbf{ factor.}  Since the energy controls \(\partial_r\bvar{\delta u}\) and \(\partial_z\bvar{\delta u}\), we recover the missing $L^2$ control of \(\bvar{\delta u}\) from these derivatives. \\ 

\noindent Set
\begin{equation}
C_{u,2}:=\max\left\{\frac{\rho_{\mathsf{growth}}}{C_{\rm floor}},\frac{1}{\mathsf {c_{growth}}}\right\}.
\label{eq:delta-u-L2-weight-constant}
\end{equation}

\hfill 

\noindent Then \eqref{eq:weight_coercivity} and \eqref{eq:delta-u-tail-growth-main} imply that
\begin{equation}
\rho\le C_{u,2}\min\{\Phi_r,\Phi_z\}.
\label{eq:delta-u-global-rho-comparison}
\end{equation}

\hfill 

\noindent Set \(f=\bvar{\delta u}\) and \(\mathbb{D}_R:=\mathbb{D}\cap\{|(r,z)|<R\}\).  Applying the divergence theorem to \((r,z)|f|^2\) gives
\begin{equation}
2\int_{\mathbb{D}_R}|f|^2
+2\int_{\mathbb{D}_R}f(r\partial_r f+z\partial_z f)
=
R\int_{\partial\mathbb{D}_R\cap\{|(r,z)|=R\}}|f|^2\,dS.
\label{eq:delta-u-dilation-truncated}
\end{equation}

\hfill 

\noindent The axis flux vanishes because \((r,z)\cdot n=-r=0\).  Using
\(2|f(r\partial_r f+z\partial_z f)|\le |f|^2+\rho|\nabla f|^2\), we obtain
\begin{equation}
\int_{\mathbb{D}_R}|f|^2
\le
R\int_{\partial\mathbb{D}_R\cap\{|(r,z)|=R\}}|f|^2\,dS
+
\int_{\mathbb{D}_R}\rho|\nabla f|^2.
\label{eq:delta-u-dilation-truncated-bound}
\end{equation}

\hfill 

\noindent The outer flux tends to zero by the boundary-flux assumption, so letting \(R\to\infty\) yields
\begin{equation}
\|\bvar{\delta u}\|_{L^2(\mathbb{D})}^2
\le \|\rho^{1/2}\nabla\bvar{\delta u}\|_{L^2(\mathbb{D})}^2.
\label{eq:delta-u-dilation-identity}
\end{equation}

\hfill 

\noindent Combining this with \eqref{eq:delta-u-global-rho-comparison} gives
\begin{equation}
\|\bvar{\delta u}\|_{L^2}
\le C_{u,2}^{1/2}
\left(
\|\partial_r\bvar{\delta u}\|_{\Phi_r}^2
+\|\partial_z\bvar{\delta u}\|_{\Phi_z}^2
\right)^{1/2}
\le C_{u,2}^{1/2}\mathfrak E_k.
\label{eq:delta-u-global-L2}
\end{equation}

\hfill 

\noindent Every positive derivative $D^\gamma\bvar{\delta u}$ is a derivative of either $\partial_r\bvar{\delta u}$ or $\partial_z\bvar{\delta u}$. Its $L^2$ norm is therefore supplied by the top-order energy or, at intermediate orders, by interpolation. 

\hfill \\ 

\noindent Applying the reflected \(H^2(\mathbb R^2)\hookrightarrow L^\infty(\mathbb R^2)\) estimate to \(\bvar{\delta u}\) gives
\begin{equation}
\|\bvar{\delta u}\|_{L^\infty}
\le
C\sum_{|\gamma|\le2}\|D^\gamma\bvar{\delta u}\|_{L^2}
\le P_u\mathfrak E_k,
\end{equation}

\noindent Here the zeroth-order term is controlled by \eqref{eq:delta-u-global-L2}, and the derivative terms by the energy and interpolation.

\clearpage

\noindent\textbf{Origin values of evolved variables.}  Point evaluation at the origin is bounded by the same reflected Sobolev estimate used in Appendix~\ref{sssec:modulation-constants}.  \\

\noindent For example, to bound \((\partial_z\bvar{\delta u})_0\) or \((\partial_z\bvar{\delta\omega})_0\), we reflect the corresponding first derivative to \(\mathbb R^2\) and apply \(H^2(\mathbb R^2)\hookrightarrow L^\infty(\mathbb R^2)\).  This requires \(L^2\) control of derivatives of the evolved variable up to order three.  Thus \(k\ge3\) supplies these origin values through the low-order energy, the top-order energy, and interpolation.  The resulting constants are \(P_{u_z,0}\) and \(P_{\omega_z,0}\).

\hfill \\

\hfill 

\noindent\textbf{Streamfunction point values and multipliers.}  For $\bvar{\delta\psi}_0$ and $(\partial_z\bvar{\delta\psi})_0$, Appendix~\ref{sssec:ell-point-constants} uses the direct dual-kernel estimate and gives
\begin{equation}
|\bvar{\delta\psi}_0|\le P_{\psi,0}\mathfrak E_k,
\qquad |(\partial_z\bvar{\delta\psi})_0|\le P_{\psi_z,0}\mathfrak E_k.
\end{equation}

\hfill

\noindent The kernel derivative \(\partial_zK_0\) has a weak diagonal singularity. We therefore use the Euclidean distance between \((r,z)\) and \((s,\zeta)\), split at distance one, and certify the single constant
\begin{equation}
\begin{aligned}
P_{\psi,z}:=\sup_{(r,z)\in\mathbb D}\Bigg[&
C^{\rm pt}_{\omega,0,0}
\int_{|(r,z)-(s,\zeta)|\le1}
|\partial_zK_0(r,z,s,\zeta)|\,d\zeta\,ds\\
&+\left(
\int_{|(r,z)-(s,\zeta)|>1}
\frac{|\partial_zK_0(r,z,s,\zeta)|^2}{\Phi_\omega(s,\zeta)}\,d\zeta\,ds
\right)^{1/2}\Bigg].
\end{aligned}
\label{eq:global-psiz-pointwise-row}
\end{equation}
\noindent The first integral is controlled by \(\|\bvar{\delta\omega}\|_{L^\infty}\), and the second by the weighted Cauchy--Schwarz inequality.\\

\noindent Hence
\begin{equation}
\|\partial_z\bvar{\delta\psi}\|_{L^\infty}
\le P_{\psi,z}\mathfrak E_k.
\end{equation}

\hfill 

\noindent The velocity factors
\begin{equation}
\bvar{\delta u^r}=-r\partial_z\bvar{\delta\psi},
\qquad
\bvar{\delta u^z}=2\bvar{\delta\psi}+r\partial_r\bvar{\delta\psi}
\end{equation}
and the remaining streamfunction multipliers are treated as members of $\mathcal S_\psi$ in \eqref{eq:ns-stream-family}.  The required terms in $\mathcal J_{\rm ell}$ are first controlled by \eqref{eq:ns-Cwell-high}.  The reflected Sobolev estimate, together with the previously certified weighted elliptic and coefficient-conversion bounds, then gives \(P_{u^r}\), \(P_{u^z}\), \(P_{\psi,rzz}\), \(P_{\psi,rzr}\), and \(P_{\psi,rrr}^{\rm ax}\).
\end{grayproof}

\clearpage

\noindent Proposition~\ref{prop:pointwise-weighted-energy} converts energy control into the pointwise bounds needed in nonlinear products.  Evolved variables are handled by reflection, whole-plane Sobolev embedding, and interpolation.  Streamfunction variables are handled by the kernel formula together with the weighted elliptic certificate.  In a typical nonlinear product, one perturbative factor is placed in \(L^\infty\) with size \(\mathcal O(\mathfrak E_k)\), while the other two factors are the weighted \(L^2\) term from the equation and the weighted \(L^2\) test factor from the energy pairing.  This gives the cubic bounds. \\

\noindent The same proposition controls the origin values entering \(\bvar{\delta C}\) and \(\bvar{\delta c_u}\).  Both modulation coefficients are linear, so these bounds enter only the linear modulation terms and the ordinary quadratic products already present in the nonlinear equations.

\hfill \\ 

\noindent The next identity is the estimate used whenever a transport term is paired with the same quantity it differentiates.  It  moves the derivative off the factor paired with \(\partial_\ell f\) and onto the coefficient and the weight.  This prevents the energy estimate from requiring one derivative more than is present in the energy.

\hfill

\begin{proposition}[Weighted transport identity]\label{prop:weighted-transport}
\noindent Let \(G\) be real-valued and let \(f\in L^2(\Phi_i)\) be sufficiently regular for the integrations below. Assume the boundary flux vanishes and
\begin{equation}
\left\|\frac{\partial_\ell(G\Phi_i)}{\Phi_i}\right\|_{L^\infty(\mathbb D)}<\infty.
\label{eq:weighted-transport-coefficient-assumption}
\end{equation}
\noindent If \(\Phi_i\) is singular at the origin, assume in addition that \(G(0,0)=0\) and that \(\nabla G\) is bounded in a neighborhood of the origin. Then
\begin{equation}
\left|\left\langle-G\partial_\ell f,f\right\rangle_{\Phi_i}\right|
\le\frac12
\left\|\frac{\partial_\ell(G\Phi_i)}{\Phi_i}\right\|_{L^\infty(\mathbb D)}
\|f\|_{\Phi_i}^2.
\label{eq:prop-weighted-transport}
\end{equation}
\end{proposition}

\begin{grayproof}
\noindent We first suppose that \(\Phi_i\) is regular at the origin. Since \(2f\partial_\ell f=\partial_\ell(f^2)\),
\begin{align}
\left\langle-G\partial_\ell f,f\right\rangle_{\Phi_i}
&=-\frac12\int_{\mathbb D}G\,\partial_\ell(f^2)\Phi_i\,dr\,dz\\
&=\frac12\int_{\mathbb D}f^2\,\partial_\ell(G\Phi_i)\,dr\,dz.
\label{eq:weighted-transport-ibp-step}
\end{align}
\noindent The second equality is integration by parts. The boundary term is zero by the assumed boundary-flux condition, including the axis and far-field limits. Factoring out \(\Phi_i\) gives
\begin{align}
\left|\left\langle-G\partial_\ell f,f\right\rangle_{\Phi_i}\right|
&\le\frac12\int_{\mathbb D}|f|^2\Phi_i
\left|\frac{\partial_\ell(G\Phi_i)}{\Phi_i}\right|\,dr\,dz\\
&\le\frac12
\left\|\frac{\partial_\ell(G\Phi_i)}{\Phi_i}\right\|_{L^\infty(\mathbb D)}
\|f\|_{\Phi_i}^2.
\end{align}

\noindent Now suppose that \(\Phi_i\) is singular at the origin. We choose a smooth radial cutoff \(\chi_s\) with
\begin{equation}
\chi_s=0\quad\text{when }r^2+z^2\le s^2,
\qquad
\chi_s=1\quad\text{when }r^2+z^2\ge4s^2,
\qquad
|\nabla\chi_s|\le C/s.
\end{equation}
\noindent The cutoff removes the singular point, so the preceding integration by parts applies to \(\chi_sG\). It gives
\begin{align}
-\int_{\mathbb D}\chi_sG(\partial_\ell f)f\Phi_i
={}&\frac12\int_{\mathbb D}\chi_s f^2\partial_\ell(G\Phi_i)
+\frac12\int_{\mathbb D}f^2G\Phi_i\partial_\ell\chi_s.
\label{eq:weighted-transport-cutoff-step}
\end{align}
\noindent For the first term, \eqref{eq:weighted-transport-coefficient-assumption} gives the integrable bound
\begin{equation}
\left|\chi_s f^2\partial_\ell(G\Phi_i)\right|
\le
\left\|\frac{\partial_\ell(G\Phi_i)}{\Phi_i}\right\|_{L^\infty(\mathbb D)}
|f|^2\Phi_i.
\end{equation}
\noindent Hence dominated convergence gives the first term in \eqref{eq:weighted-transport-ibp-step} as \(s\to0\).\\

\noindent It remains to show that the cutoff term tends to zero. We choose a neighborhood of the origin on which \(|\nabla G|\le C_G\). Since \(G(0,0)=0\), the mean value theorem gives, for \(s\) small enough,
\begin{equation}
|G(r,z)|\le C_G\sqrt{r^2+z^2}.
\end{equation}
\noindent The derivative \(\nabla\chi_s\) is supported where \(s\le\sqrt{r^2+z^2}\le2s\). Therefore
\begin{equation}
|G\partial_\ell\chi_s|
\le C C_G
\end{equation}
\noindent on its support. Consequently,
\begin{equation}
\left|\int_{\mathbb D}f^2G\Phi_i\partial_\ell\chi_s\right|
\le C C_G
\int_{\{s\le\sqrt{r^2+z^2}\le2s\}}|f|^2\Phi_i\,dr\,dz.
\end{equation}
\noindent The last integral tends to zero because \(f\in L^2(\Phi_i)\) and the annuli shrink to the origin.\\

\noindent Letting \(s\to0\) proves \eqref{eq:weighted-transport-ibp-step}, and the same \(L^\infty\) estimate then gives \eqref{eq:prop-weighted-transport}.
\end{grayproof}

\hfill 

\noindent Proposition~\ref{prop:weighted-transport} is the weighted integration-by-parts identity used whenever a nonlinear transport coefficient multiplies $\partial_\ell f$ and the term is paired with $f$ in the weighted inner product.  It is used in the low-order nonlinear estimate and in the principal transport part of the high-order nonlinear estimate.

\hfill

\noindent The high-order linear proof also differentiates the weighted transport identity after commuting $D^\alpha$ through coefficients and weights.  For that reason we record the top-order weight constants that control first and second derivatives of the higher-order weight functions:
\begin{equation}
\label{eq:ns-top-weight-constants}
W^{(1)}_{i,k,\ell}:=\left\|\wgray{\Phi_{i,k}}^{-1}\partial_\ell\wgray{\Phi_{i,k}}\right\|_{L^\infty(\mathbb{D})},
\qquad
W^{(2)}_{i,k,\ell m}:=\left\|\wgray{\Phi_{i,k}}^{-1}\partial_{\ell m}\wgray{\Phi_{i,k}}\right\|_{L^\infty(\mathbb{D})}.
\end{equation}

\noindent These constants are the high-order analogues of the low-order transport constants.   

\clearpage

\section{Low-order linear estimate}\label{app:low-order-linear}

\noindent We now prove the low-order linear estimate \eqref{eq:ns-CL}. The transport and selected zeroth-order component couplings are kept together in a signed symmetric $3\times3$ quadratic form after the component weights are inserted. Its negativity is certified by the matrix damping bound. The coupling from vorticity to $\partial_z u$ proportional to $2\bar u$ is kept outside the matrix because independent singular weights do not guarantee that its matrix multiplier is finite, and thus estimated directly. The remaining terms are handled according to their structure: undifferentiated $\bvar{\delta u}$ terms by the ray estimate, streamfunction terms by the weighted elliptic estimates, and linear modulation terms by the modulation bounds.

\hfill \\

\subsection{Splitting the low-order linear operator}\label{app:low-order-linear-setup}

\hfill 

\noindent The low-order linear pairing is
\begin{equation}
\label{eq:low-linear-pairing-def}
\begin{aligned}
\mathcal Q_{\mathcal L}^{\rm low}(\bvar{\delta})
&:=\langle \mathcal L_\omega(\bvar{\delta}),\bvar{\delta v_\omega}\rangle_{\wgray{\Phi_\omega}}
+\langle \mathcal L_r(\bvar{\delta}),\bvar{\delta v_r}\rangle_{\wgray{\Phi_r}}
+\langle \mathcal L_z(\bvar{\delta}),\bvar{\delta v_z}\rangle_{\wgray{\Phi_z}}.
\end{aligned}
\end{equation}

\hfill

\noindent Recall that \(\mathcal Q_{\mathcal L}^{\rm low}\) contains the terms obtained when the residual amplitude offset multiplies a perturbation.  Then, \(\bvar{\delta C}\) and \(\bvar{\delta c_u}\) are exactly linear and therefore enter the linear operator in full.

\hfill 

\noindent The same transport field appears in all three components.  It is
\begin{equation}
\label{eq:low-linear-profile-transport-coeffs}
G_r:=\lambda r+\bar u^r,
\qquad
G_z:=\lambda z+\bar C+\bar u^z.
\end{equation}

\hfill

\noindent Using \(\bar u^r=-r\partial_z\bar\psi\), \(\bar u^z=2\bar\psi+r\partial_r\bar\psi\), and the normalization conditions,
\begin{equation}
G_r(0,0)=0,
\qquad
G_z(0,0)=\bar C+2\bar\psi_0=0.
\label{eq:low-linear-origin-stagnation-add-current}
\end{equation}
Thus the origin is a fixed point of the meridional transport field.

\hfill  \\

\noindent For each \(i\in\{\omega,r,z\}\), we separate the transport part and non-transport terms:
\begin{equation}
\label{eq:low-linear-operator-split}
\mathcal L_i(\bvar{\delta})
=-G_r\partial_r\bvar{\delta v_i}-G_z\partial_z\bvar{\delta v_i}
+\sum_m B_{i,m}(\bvar{\delta}).
\end{equation}

\hfill 

\noindent The complete list of non-transport terms is
\begin{align}
\label{eq:low-linear-B-omega-list}
B_{\omega,1}&=(\bar c_u-\lambda
+\cures)\bvar{\delta v_\omega},
&B_{\omega,2}&=2\partial_z\bar u\,\bvar{\delta u},\nonumber\\
B_{\omega,3}&=2\bar u\,\bvar{\delta v_z},
&B_{\omega,4}&=-2\partial_z\bar\omega\,\bvar{\delta\psi},\nonumber\\
B_{\omega,5}&=r\partial_r\bar\omega\,\partial_z\bvar{\delta\psi},
&B_{\omega,6}&=-r\partial_z\bar\omega\,\partial_r\bvar{\delta\psi},\nonumber\\
B_{\omega,7}&=-\partial_z\bar\omega\,\bvar{\delta C},
&B_{\omega,8}&=\bar\omega\,\bvar{\delta c_u},
\end{align}
\begin{align}
\label{eq:low-linear-B-r-list}
B_{r,1}&=(2\partial_z\bar\psi+\bar c_u-\lambda-\partial_r\bar u^r
+\cures)\bvar{\delta v_r},
&B_{r,2}&=-\partial_r\bar u^z\,\bvar{\delta v_z},\nonumber\\
B_{r,3}&=2\partial_{rz}\bar\psi\,\bvar{\delta u},
&B_{r,4}&=(2\bar u+r\partial_r\bar u)\partial_{rz}\bvar{\delta\psi},\nonumber\\
B_{r,5}&=(3\partial_r\bar u+r\partial_{rr}\bar u)\partial_z\bvar{\delta\psi},
&B_{r,6}&=-r\partial_z\bar u\,\partial_{rr}\bvar{\delta\psi},\nonumber\\
B_{r,7}&=-(3\partial_z\bar u+r\partial_{rz}\bar u)\partial_r\bvar{\delta\psi},
&B_{r,8}&=-2\partial_{rz}\bar u\,\bvar{\delta\psi},\nonumber\\
B_{r,9}&=-\partial_{rz}\bar u\,\bvar{\delta C},
&B_{r,10}&=\partial_r\bar u\,\bvar{\delta c_u},
\end{align}
\begin{align}
\label{eq:low-linear-B-z-list}
B_{z,1}&=-\partial_z\bar u^r\,\bvar{\delta v_r},
&B_{z,2}&=(2\partial_z\bar\psi+\bar c_u-\lambda-\partial_z\bar u^z
+\cures)\bvar{\delta v_z},\nonumber\\
B_{z,3}&=2\partial_{zz}\bar\psi\,\bvar{\delta u},
&B_{z,4}&=(2\bar u+r\partial_r\bar u)\partial_{zz}\bvar{\delta\psi},\nonumber\\
B_{z,5}&=(r\partial_{rz}\bar u)\partial_z\bvar{\delta\psi},
&B_{z,6}&=-r\partial_z\bar u\,\partial_{rz}\bvar{\delta\psi},\nonumber\\
B_{z,7}&=-r\partial_{zz}\bar u\,\partial_r\bvar{\delta\psi},
&B_{z,8}&=-2\partial_{zz}\bar u\,\bvar{\delta\psi},\nonumber\\
B_{z,9}&=-\partial_{zz}\bar u\,\bvar{\delta C},
&B_{z,10}&=\partial_z\bar u\,\bvar{\delta c_u}.
\end{align}

\hfill \\

\noindent The matrix part contains the transport terms and the five labels
\begin{equation}
\label{eq:low-linear-selected-set-fixed}
\mathcal S_{\mathcal L}
=
\{B_{\omega,1},B_{r,1},B_{r,2},B_{z,1},B_{z,2}\}.
\end{equation}
\noindent The three diagonal labels contribute to the diagonal entries, while $B_{r,2}$ and $B_{z,1}$ contribute to the symmetric $r$--$z$ entry. The coupling
\begin{equation}
B_{\omega,3}=2\bar u\,\bvar{\delta v_z}
\end{equation}
\noindent is estimated directly. Its matrix coefficient would contain $2\bar u(\Phi_\omega/\Phi_z)^{1/2}$, and the independent low-order singularity rates do not by themselves make this multiplier finite. The five-label matrix and the direct cost of $B_{\omega,3}$ are therefore certified separately. \\

\noindent The direct index sets are
\begin{align}
\mathcal I_\omega&=\{2,3,4,5,6,7,8\},
\label{eq:low-linear-direct-index-set-omega}\\
\mathcal I_r&=\{3,4,5,6,7,8,9,10\},
\label{eq:low-linear-direct-index-set-r}\\
\mathcal I_z&=\{3,4,5,6,7,8,9,10\}.
\label{eq:low-linear-direct-index-set-z}
\end{align}

\hfill

\noindent We refer to the terms estimated individually outside the signed matrix the direct terms. They fall into four fixed groups.
\begin{itemize}
    \item The terms $B_{\omega,2}$, $B_{r,3}$, and $B_{z,3}$ contain undifferentiated $\bvar{\delta u}$ and are controlled by Proposition~\ref{prop:fixed-profile-ray-estimate}.
    \item The term $B_{\omega,3}=2\bar u\,\bvar{\delta v_z}$ is controlled directly by the weighted $L^2$ norm of the fixed profile $\bar u$ and the global pointwise estimate for $\bvar{\delta v_z}$.
    \item The terms $B_{\omega,4}$--$B_{\omega,6}$, $B_{r,4}$--$B_{r,8}$, and $B_{z,4}$--$B_{z,8}$ contain $\bvar{\delta\psi}$ or its derivatives and use the weighted elliptic estimates for $\bvar{\delta\psi}=\Gamma_5*\bvar{\delta\omega}$.
    \item The terms $B_{\omega,7}$, $B_{\omega,8}$, $B_{r,9}$, $B_{r,10}$, $B_{z,9}$, and $B_{z,10}$ contain $\bvar{\delta C}$ or $\bvar{\delta c_u}$ and use the linear modulation bounds.
\end{itemize}

\hfill \\

\noindent With this fixed split,
\begin{equation}
\label{eq:low-linear-sign-remainder-split}
\mathcal L_i(\bvar{\delta})
=\mathcal L_{i,{\rm mat}}^{\rm low}(\bvar{\delta})
+\sum_{m\in\mathcal I_i}B_{i,m}(\bvar{\delta}),
\end{equation}

\noindent where \(\mathcal L_{i,{\rm mat}}^{\rm low}\) contains the transport \((G_r,G_z)\) terms and the labels from \(\mathcal S_{\mathcal L}\) that belong to component \(i\).

\hfill \\

\subsection{Low-order linear estimate}\label{app:low-order-linear-combined}

\noindent The estimate is written in the three weighted point variables
\begin{equation}
\label{eq:low-linear-xi-def}
\bvar{d\xi_\omega}:=\wgray{\Phi_\omega}^{1/2}\bvar{\delta v_\omega},
\qquad
\bvar{d\xi_r}:=\wgray{\Phi_r}^{1/2}\bvar{\delta v_r},
\qquad
\bvar{d\xi_z}:=\wgray{\Phi_z}^{1/2}\bvar{\delta v_z},
\end{equation}
and set
\begin{equation}
\label{eq:low-linear-xi-column}
\bvar{d\xi}(r,z):=
\left(\bvar{d\xi_\omega},\bvar{d\xi_r},\bvar{d\xi_z}\right)(r,z)^\top .
\end{equation}

\hfill \\

\noindent The theorem combines one global matrix bound with the remaining direct terms. \\

\begin{theorem}[Low-order linear estimate]\label{prop:low-order-linear}
\noindent After weighted integration by parts,
\begin{equation}
\sum_{i\in\{\omega,r,z\}}
\langle \mathcal L_{i,{\rm mat}}^{\rm low}(\bvar{\delta}),\bvar{\delta v_i}\rangle_{\Phi_i}
=
\int_{\mathbb D}
\bvar{d\xi}^{\top}\mathbb M_{\mathcal L}\bvar{d\xi}\,dr\,dz.
\label{eq:low-linear-sign-matrix-form}
\end{equation}
\noindent Assume
\begin{equation}
\lambda_{\max}\!\left(\mathbb M_{\mathcal L}(r,z)\right)
\le -\Lambda_{\mathcal L}^{\rm mat}<0
\qquad\text{for a.e. }(r,z)\in\mathbb D.
\label{eq:low-linear-sign-certificate}
\end{equation}

\noindent Assume also that each direct label satisfies the term bound
\begin{equation}
\label{eq:low-linear-remainder-norm-bound}
\|B_{i,m}(\bvar{\delta})\|_{\wgray{\Phi_i}}
\le \kappa^0_{i,m}\mathfrak E_0+\kappa^k_{i,m}\mathfrak H_k,
\qquad i\in\{\omega,r,z\},
\quad m\in\mathcal I_i.
\end{equation}
For every direct term with \(\kappa^k_{i,m}>0\), choose a positive balancing constant \(\eta_{i,m}\), and set
\begin{equation}
\mathcal Y_i:=\{m\in\mathcal I_i:\ \kappa^k_{i,m}>0\}.
\end{equation}
Define
\begin{equation}
\label{eq:low-linear-total-losses}
A_{\mathcal L}^{0}:=
\sum_{i\in\{\omega,r,z\}}\sum_{m\in\mathcal I_i}\kappa^0_{i,m}
+
\sum_{i\in\{\omega,r,z\}}\sum_{m\in\mathcal Y_i}\eta_{i,m},
\qquad
A_{\mathcal L}^{k}:=
\sum_{i\in\{\omega,r,z\}}\sum_{m\in\mathcal Y_i}
\frac{(\kappa^k_{i,m})^2}{4\eta_{i,m}}.
\end{equation}
Set
\begin{equation}
\label{eq:low-linear-damping-positive}
\DLamLlow:=\Lambda_{\mathcal L}^{\rm mat}-A_{\mathcal L}^{0},
\qquad
\DLamLhigh:=A_{\mathcal L}^{k}.
\end{equation}

\noindent If \(\DLamLlow>0\), then
\begin{equation}
\label{eq:low-linear-final-estimate}
\mathcal Q_{\mathcal L}^{\rm low}(\bvar{\delta})
\le -\DLamLlow \, \mathfrak E_0^2+
\DLamLhigh \, \mathfrak H_k^2.
\end{equation}
\end{theorem}

\begin{grayproof}
\noindent\textbf{Matrix part.} The matrix part consists only of the transport terms and the five labels in \(\mathcal S_{\mathcal L}\).  After integration by parts, the transport terms are also pointwise in \(\bvar{d\xi}\).  For each component, the weighted transport identity gives
\begin{equation}
\label{eq:low-linear-transport-signed}
\begin{aligned}
\left\langle -G_r\partial_r\bvar{\delta v_i}-G_z\partial_z\bvar{\delta v_i},
\bvar{\delta v_i}\right\rangle_{\wgray{\Phi_i}}
&=\frac{1}{2}\int_\mathbb{D}
\bigl(\partial_r(G_r\wgray{\Phi_i})+
\partial_z(G_z\wgray{\Phi_i})\bigr)|\bvar{\delta v_i}|^2\,dr\,dz \\
&=\int_\mathbb{D}
\frac{\partial_r(G_r\wgray{\Phi_i})+
\partial_z(G_z\wgray{\Phi_i})}{2\wgray{\Phi_i}}
|\bvar{d\xi_i}|^2\,dr\,dz.
\end{aligned}
\end{equation}

{\noindent The signed axial transport contribution is retained through $G_z\partial_z\log\Phi_i$, with $\partial_zG_z$ appearing in the divergence term.}

\hfill  \\

\noindent The boundary terms vanish by the admissibility and boundary-flux hypotheses, together with the prescribed parity at the axis, and density extends the identity to the weighted class.  A diagonal matrix label \(B_{i,m}=a\bvar{\delta v_i}\) adds \(a\) to the corresponding diagonal entry of \(\mathbb M_{\mathcal L}\).  
      
\hfill 

\noindent For a cross matrix label \(B_{i,m}=a\bvar{\delta v_j}\), with \(i\ne j\), we use the combined coefficient
\begin{equation}
\label{eq:low-linear-cross-weighted}
\int_{\mathbb D}a\bvar{\delta v_j}\bvar{\delta v_i}\Phi_i\,dr\,dz
=
\int_{\mathbb D}a\left(\frac{\Phi_i}{\wgray{\Phi_j}}\right)^{1/2}
\bvar{d\xi_i}\bvar{d\xi_j}\,dr\,dz.
\end{equation}

\hfill 

\noindent The matrix certificate therefore gives
\begin{equation}
\sum_i
\langle \mathcal L_{i,{\rm mat}}^{\rm low}(\bvar{\delta}),\bvar{\delta v_i}\rangle_{\Phi_i}
\le
-\Lambda_{\mathcal L}^{\rm mat}\mathfrak E_0^2.
\label{eq:low-linear-matrix-integrated}
\end{equation}

\hfill \\

\noindent\textbf{Direct part.}  Now fix one direct label \(B_{i,m}\) with \(m\in\mathcal I_i\).  The Cauchy--Schwarz inequality gives
\begin{equation}
\left|\langle B_{i,m}(\bvar{\delta}),\bvar{\delta v_i}\rangle_{\wgray{\Phi_i}}\right|
\le
\|B_{i,m}(\bvar{\delta})\|_{\wgray{\Phi_i}}
\|\bvar{\delta v_i}\|_{\wgray{\Phi_i}}.
\end{equation} 
Using \eqref{eq:low-linear-remainder-norm-bound} and
\(\|\bvar{\delta v_i}\|_{\wgray{\Phi_i}}\le\mathfrak E_0\), this becomes
\begin{equation}
\label{eq:low-linear-remainder-pairing}
\begin{aligned}
\left|\langle B_{i,m}(\bvar{\delta}),\bvar{\delta v_i}\rangle_{\wgray{\Phi_i}}\right|
&\le
\left(\kappa^0_{i,m}\mathfrak E_0+
\kappa^k_{i,m}\mathfrak H_k\right)\mathfrak E_0 =
\kappa^0_{i,m}\mathfrak E_0^2+
\kappa^k_{i,m}\mathfrak E_0\mathfrak H_k .
\end{aligned}
\end{equation}

\hfill 

\noindent If \(\kappa^k_{i,m}=0\), the term contributes only
\(\kappa^0_{i,m}\mathfrak E_0^2\).  If \(\kappa^k_{i,m}>0\), we choose a positive balancing constant \(\eta_{i,m}\) and apply Young's inequality:
\begin{equation}
\label{eq:low-linear-young}
\kappa^k_{i,m}\mathfrak E_0\mathfrak H_k
\le \eta_{i,m}\mathfrak E_0^2+
\frac{(\kappa^k_{i,m})^2}{4\eta_{i,m}}\mathfrak H_k^2.
\end{equation}
 
\hfill 

 \noindent  Summing all direct terms and adding \eqref{eq:low-linear-matrix-integrated} gives
\begin{equation}
\mathcal Q_{\mathcal L}^{\rm low}(\bvar{\delta})
\le
- \bigl( \Lambda_{\mathcal L}^{\rm mat}-A _{\mathcal L}^{0}\bigr) \mathfrak E_0^2
+ A_{\mathcal L}^{k}\mathfrak H_k^2.
\end{equation}

\end{grayproof}

\subsection{Direct term rules}\label{app:low-order-linear-c5}

\noindent This subsection gives the four rules used for the direct terms in \eqref{eq:low-linear-direct-index-set-omega}--\eqref{eq:low-linear-direct-index-set-z}. Each direct term produces two certified constants \((\kappa^0_{i,m},\kappa^k_{i,m})\) in \eqref{eq:low-linear-remainder-norm-bound}.

\begin{enumerate}[label=(R\arabic*),leftmargin=.38in,itemsep=3em,topsep=2em]
\item \emph{Terms proportional to $\bvar{\delta u}$.} For $B_{\omega,2}$, $B_{r,3}$, and $B_{z,3}$, we write $B_{i,m}=a\bvar{\delta u}$ with the corresponding fixed profile coefficient $a$. Proposition~\ref{prop:fixed-profile-ray-estimate} gives
\begin{equation}
\|a\bvar{\delta u}\|_{\wgray{\Phi_i}}
\le \mathcal C_{{\rm ray},i}(a)^{1/2}
\left(\|\partial_r\bvar{\delta u}\|_{\wgray{\Phi_r}}^2
+\|\partial_z\bvar{\delta u}\|_{\wgray{\Phi_z}}^2\right)^{1/2}
\le \mathcal C_{{\rm ray},i}(a)^{1/2}\mathfrak E_0.
\label{eq:low-linear-rule-u-ray-detail}
\end{equation}
\noindent Hence
\begin{equation}
\kappa^0_{i,m}=\mathcal C_{{\rm ray},i}(a)^{1/2},
\qquad
\kappa^k_{i,m}=0.
\label{eq:low-linear-rule-u-ray}
\end{equation}

\item \emph{The direct vorticity--$z$ coupling.} For
\begin{equation}
B_{\omega,3}=2\bar u\,\bvar{\delta v_z},
\end{equation}
\noindent weighted H\"older and Proposition~\ref{prop:pointwise-weighted-energy} give
\begin{equation}
\|B_{\omega,3}\|_{\Phi_\omega}
\le 2\|\bar u\|_{\Phi_\omega}\|\bvar{\delta v_z}\|_{L^\infty(\mathbb D)}
\le 2C^{\rm pt}_{z,0,0}\|\bar u\|_{\Phi_\omega}\mathfrak E_k.
\label{eq:low-linear-rule-Bomega3-direct}
\end{equation}
\noindent Using $\mathfrak E_k\le \mathfrak E_0+\sqrt{\mu_k}\,\mathfrak H_k$ gives the admissible choice
\begin{equation}
\kappa^0_{\omega,3}=2C^{\rm pt}_{z,0,0}\|\bar u\|_{\Phi_\omega},
\qquad
\kappa^k_{\omega,3}=2\sqrt{\mu_k}\,C^{\rm pt}_{z,0,0}\|\bar u\|_{\Phi_\omega}.
\label{eq:low-linear-rule-Bomega3-constants}
\end{equation}
\noindent Its low- and high-order costs are recorded among the direct losses.

\item \emph{Streamfunction terms.} The labels $B_{\omega,4}$--$B_{\omega,6}$, $B_{r,4}$--$B_{r,8}$, and $B_{z,4}$--$B_{z,8}$ have the form $B_{i,m}=a\,T\bvar{\delta\psi}$, where $T$ is exactly the streamfunction expression displayed in that label. Thus $T$ is one of
\begin{equation}
\{{\rm Id},\partial_r,\partial_z,\partial_{rr},\partial_{rz},\partial_{zz}\}.
\end{equation}
\noindent Since $\bvar{\delta\psi}=\Gamma_5*\bvar{\delta\omega}$, we certify the actual weighted operator needed by that row:
\begin{equation}
\|a\,T(\Gamma_5*\bvar{\delta\omega})\|_{\wgray{\Phi_i}}
\le \kappa^0_{i,m}\mathfrak E_0+\kappa^k_{i,m}\mathfrak H_k.
\label{eq:low-linear-rule-stream}
\end{equation}
\noindent The multiplier $a$, the streamfunction operator $T$, and the target weight $\wgray{\Phi_i}$ are therefore certified together rather than through a separate universal weight-comparison constant. Appendix~\ref{sssec:ell-low-order-rows} denotes the resulting low/high row constants by $\mathcal K_i^0(a,T)$ and $\mathcal K_i^k(a,T)$. These are the values entered as $\kappa^0_{i,m}$ and $\kappa^k_{i,m}$.

\item \emph{Linear modulation terms.} If a direct label is $a\bvar{\delta C}$ or $a\bvar{\delta c_u}$, we combine the scalar modulation bound with the weighted norm of the fixed coefficient:
\begin{equation}
\label{eq:low-linear-rule-modulation}
\begin{array}{ll}
\kappa^0_{i,m}=\|a\|_{\wgray{\Phi_i}}\kappa_{\bvar{\delta C}}^{0},
\quad
\kappa^k_{i,m}=\|a\|_{\wgray{\Phi_i}}\kappa_{\bvar{\delta C}}^{k}, & B_{i,m}=a\bvar{\delta C},\\[1.4ex]
\kappa^0_{i,m}=\|a\|_{\wgray{\Phi_i}}\kappa_{\bvar{\delta c_u}}^{0},
\quad
\kappa^k_{i,m}=\|a\|_{\wgray{\Phi_i}}\kappa_{\bvar{\delta c_u}}^{k}, & B_{i,m}=a\bvar{\delta c_u}.
\end{array}
\end{equation}
\end{enumerate}

\hfill  \\ 

\hfill

\subsection{Explicit numerical certificate}\label{app:low-order-linear-explicit-damping}

\noindent This subsection gives the finite data needed to certify the low-order linear estimate.\\

\noindent The construction first forms the symmetric matrix \(\mathbb M_{\mathcal L}\) from the transport contribution and the selected matrix labels in \eqref{eq:low-linear-selected-set-fixed}.  In the standard case, its damping constant is certified on the full half-plane by combining a finite-box enclosure with the tail closure in Appendix~\ref{sssec:low-order-linear-tail-closure}. \\

\noindent If strict matrix damping fails only on sufficiently small exception sets, Appendix~\ref{sssec:localized-matrix-certificate} gives a localized matrix certificate that lowers the available matrix damping and records the corresponding high-order cost. \\

\noindent After the matrix part is fixed by one of these two routes, the remaining direct terms in \eqref{eq:low-linear-direct-index-set-omega}--\eqref{eq:low-linear-direct-index-set-z} give the losses \(A_{\mathcal L}^{0}\) and \(A_{\mathcal L}^{k}\).  In the standard global-matrix case, the final constants are
\begin{equation}
\label{eq:low-linear-explicit-outputs}
\DLamLlow=\Lambda_{\mathcal L}^{\rm mat}-A_{\mathcal L}^{0},
\qquad
\DLamLhigh=A_{\mathcal L}^{k}.
\end{equation}
Here \(\Lambda_{\mathcal L}^{\rm mat}\) is the damping produced by the matrix on all of \(\mathbb{D}\), while \(A_{\mathcal L}^{0}\) and \(A_{\mathcal L}^{k}\) are the low-order and high-order losses produced by the direct terms.

\hfill \\

\subsubsection{Matrix construction}\label{sssec:low-order-linear-matrix-construction}

\noindent The matrix is written in the component order \((\omega,r,z)\). Only the transport terms and the terms in \(\mathcal S_{\mathcal L}\) enter this matrix. The remaining terms are kept for the direct-term estimate below. \\

\paragraph{Transport coefficients.}
We compute
\begin{equation}
G_r=\lambda r+\bar u^r,
\qquad
G_z=\lambda z+\bar C+\bar u^z,
\end{equation}
where
\begin{equation}
    \bar u^r=-r\partial_z\bar\psi,
\qquad
\bar u^z=2\bar\psi+r\partial_r\bar\psi
\end{equation}

\hfill

\hfill \\

\noindent The derivatives needed for the transport diagonal are
\begin{align}
\partial_r\bar u^r&=-\partial_z\bar\psi-r\partial_{rz}\bar\psi,
&
\partial_z\bar u^r&=-r\partial_{zz}\bar\psi,\\
\partial_r\bar u^z&=3\partial_r\bar\psi+r\partial_{rr}\bar\psi,
&
\partial_z\bar u^z&=2\partial_z\bar\psi+r\partial_{rz}\bar\psi.
\end{align}

\hfill

\noindent This gives
\begin{equation}
\partial_rG_r=\lambda+\partial_r\bar u^r,
\qquad
\partial_zG_z
=\lambda+\partial_z\bar u^z.
\end{equation}

\hfill \\

\noindent For each component \(i\in\{\omega,r,z\}\), set
\begin{equation}
D_i
:=\frac{1}{2}\left(
\partial_rG_r+\partial_zG_z
+G_r\partial_r\log\wgray{\Phi_i}
+G_z\partial_z\log\wgray{\Phi_i}
\right).
\end{equation}

\hfill \\

\noindent Initialize the matrix with
\begin{equation}
\mathbb M_{\mathcal L}=\operatorname{diag}(D_\omega,D_r,D_z).
\end{equation}

\hfill \\

\paragraph{Diagonal matrix labels.}
\noindent The constant modulation offset \(\cures\) is absorbed into the reference amplitude coefficient and therefore contributes to the zeroth-order linear coefficients. We add each displayed coefficient to the indicated diagonal entry:
\begin{equation}
\begin{array}{lll}
(\omega,1):& (\omega,\omega) &
\bar c_u+\cures-\lambda,\\[0.6ex]
(r,1):& (r,r) &
2\partial_z\bar\psi+\bar c_u+\cures-\lambda-\partial_r\bar u^r,\\[0.6ex]
(z,2):& (z,z) &
2\partial_z\bar\psi+\bar c_u+\cures-\lambda-\partial_z\bar u^z.
\end{array}
\end{equation}

\hfill \\

\paragraph{Cross matrix labels.}
\noindent Both $r$--$z$ couplings contribute to the same symmetric off-diagonal entry. We keep each coefficient together with the weight conversion required by its source and target components:
\begin{equation}
\begin{array}{lll}
(r,2):& (r,z),(z,r) &
-\partial_r\bar u^z
\left(\dfrac{\Phi_r}{\Phi_z}\right)^{1/2},\\[1.8ex]
(z,1):& (r,z),(z,r) &
-\partial_z\bar u^r
\left(\dfrac{\Phi_z}{\Phi_r}\right)^{1/2}.
\end{array}
\end{equation}
\noindent In the symmetric quadratic form, one half of each listed coefficient is placed in each off-diagonal entry. The coupling $B_{\omega,3}=2\bar u\,\bvar{\delta v_z}$ remains outside the matrix because the corresponding multiplier $2\bar u(\Phi_\omega/\Phi_z)^{1/2}$ may not be finite.\\

\paragraph{Matrix damping.}
\noindent The preceding entries define the pointwise matrix. A damping constant $\Lambda_{\mathcal L}^{\rm mat}>0$ is certified by
\begin{equation}
\label{eq:low-linear-explicit-Lambda-mat}
\lambda_{\max}\bigl(\mathbb M_{\mathcal L}(r,z)\bigr)
\le -\Lambda_{\mathcal L}^{\rm mat}
\qquad\text{for almost every }(r,z)\in\mathbb D.
\end{equation}
\noindent The finite-box enclosure is combined with the tail estimate below. Any exceptional region of positive measure is handled by Appendix~\ref{sssec:localized-matrix-certificate}.

\hfill \\

\subsubsection{Tail closure for the matrix damping bound}
\label{sssec:low-order-linear-tail-closure}

{
\noindent Define the numerical enclosure box by
\begin{equation}
\label{eq:low-linear-tail-box}
\mathbb{D}_{\rm box}:=[0,R_{\rm box}]\times[-R_{\rm box},R_{\rm box}]
\end{equation}
This is the box on which the pointwise matrix is enclosed numerically.  This subsubsection gives an analytic bound on \(\mathbb{D}\setminus\mathbb{D}_{\rm box}\).  Combining the two bounds produces one global matrix damping constant.  If localized exception sets remain inside the box, we use the certificate in the next subsubsection instead.

\hfill \\

\noindent Throughout the tail calculation,
\begin{equation}
\label{eq:low-linear-tail-normalization}
\lambda=\frac{1}{2},
\qquad
\bar c_u=-1.
\end{equation}
The profile-dependent constant \(\bar C\) remains symbolic until the numerical certificate is evaluated. \\ 

\paragraph{Exact tail matrix.}
\noindent Assume that the steady profile and every profile derivative entering \(\mathbb M_{\mathcal L}\) vanish exactly on \(\mathbb{D}\setminus\mathbb{D}_{\rm box}\).  There,
\begin{equation}
\bar u=\bar\omega=\bar\psi=0,
\qquad
\bar u^r=\bar u^z=0,
\qquad
G_r=\frac r2,
\qquad
G_z=\frac z2+\bar C.
\end{equation}
All selected cross-couplings vanish.    The residual amplitude offset remains because it is a spatial constant multiplying the perturbation. \\

\noindent Hence the tail matrix is diagonal, with entry
\begin{equation}
\label{eq:low-linear-tail-diagonal-entry}
-1+\cures
+\frac14(r\partial_r+z\partial_z)\log\wgray{\Phi_i}
+\frac12(\bar C)\partial_z\log\wgray{\Phi_i},
\qquad i\in\{\omega,r,z\}.
\end{equation}
Thus the tail eigenvalue bound reduces to one scalar inequality for each weight. 

\hfill \\ 

\paragraph{Scalar tail condition.}
\noindent It is sufficient to find \(\gamma_{\rm tail}>0\) such that
\begin{equation}
\label{eq:low-linear-tail-scalar-condition}
\sup_{\mathbb{D}\setminus\mathbb{D}_{\rm box}}
\left[
(r\partial_r+z\partial_z)\log\wgray{\Phi_i}
+2(\bar C)\partial_z\log\wgray{\Phi_i}
\right]
\le 4(1-\cures-\gamma_{\rm tail})
\end{equation}
for every \(i\in\{\omega,r,z\}\). \\

\noindent Substitution into \eqref{eq:low-linear-tail-diagonal-entry} gives
\begin{equation}
\label{eq:low-linear-tail-eigenvalue-bound}
\lambda_{\max}\bigl(\mathbb M_{\mathcal L}(r,z)\bigr)
\le-\gamma_{\rm tail}
\qquad\text{on }\mathbb{D}\setminus\mathbb{D}_{\rm box}.
\end{equation}

\hfill 

\noindent If the finite-box computation gives
\begin{equation}
\lambda_{\max}\bigl(\mathbb M_{\mathcal L}(r,z)\bigr)
\le-\Lambda_{{\mathcal L},{\rm box}}^{\rm mat}\qquad\text{for a.e. }(r,z)\in\mathbb D_{\rm box}.
\end{equation}
then any choice satisfying
\begin{equation}
\label{eq:low-linear-tail-global-Lambda-choice}
0<\Lambda_{\mathcal L}^{\rm mat}
\le\min\left\{\Lambda_{{\mathcal L},{\rm box}}^{\rm mat},\gamma_{\rm tail}\right\}
\end{equation}
is a global matrix damping constant.

\hfill \\

\paragraph{Radial tail weights.}
\noindent Suppose that the tail weight is radial in
\(\rho=r^2+z^2\):
\begin{equation}
\label{eq:low-linear-tail-radial-weight}
\wgray{\Phi_i}(r,z)=\phi_i(\rho).
\end{equation}
Then
\begin{equation}
\label{eq:low-linear-tail-radial-chain-rule}
\partial_r\log\wgray{\Phi_i}
=2r\frac{d}{d\rho}\log\phi_i,
\qquad
\partial_z\log\wgray{\Phi_i}
=2z\frac{d}{d\rho}\log\phi_i,
\end{equation}
and the left side of \eqref{eq:low-linear-tail-scalar-condition} becomes
\begin{equation}
\label{eq:low-linear-tail-radial-expression}
\left(1+\frac{2z}{\rho}(\bar C)\right)
\left(2\rho\frac{d}{d\rho}\log\phi_i(\rho)\right).
\end{equation}

\noindent  The first factor contains the complete signed axial drift correction.  Since \(\sqrt\rho\ge R_{\rm box}\) and \(|z|\le\sqrt\rho\) on the tail,
\begin{equation}
\label{eq:low-linear-tail-finite-box-factor}
1-\frac{2|\bar C|}{R_{\rm box}}
\le
1+\frac{2z}{\rho}(\bar C)
\le
1+\frac{2|\bar C|}{R_{\rm box}}.
\end{equation}

\hfill 

\noindent The convenient sufficient assumptions
\begin{equation}
\label{eq:low-linear-tail-basic-phase-condition}
\cures<1,
\qquad
R_{\rm box}>2|\bar C|
\end{equation}
have separate roles.  The first leaves a positive zeroth-order margin.  The second keeps the correction factor in \eqref{eq:low-linear-tail-radial-expression} positive throughout the tail.

\hfill \\ 

\paragraph{Admissible model rates.}
\noindent The implemented weights have a positive floor \(C_{\rm floor}>0\).  Because the correction factor tends to one, a necessary asymptotic guide is
\begin{equation}
\limsup_{\rho\to\infty}
2\rho\frac{d}{d\rho}\log\phi_i(\rho)
<4(1-\cures).
\label{eq:low-linear-tail-asymptotic-rate-condition}
\end{equation}
The finite-box certificate must also use the upper factor in \eqref{eq:low-linear-tail-finite-box-factor}. \\ 

\begin{itemize}[leftmargin=1.3em,itemsep=2.5em,topsep=0.1em]
\item \textbf{Exponential decay.}  If
\begin{equation}
\phi_i(\rho)=C_{\rm floor}+a\exp(-b\rho^p),
\qquad a,b,p>0,
\end{equation}
then \(d(\log\phi_i)/d\rho\le0\).  Under \eqref{eq:low-linear-tail-basic-phase-condition}, the left side of \eqref{eq:low-linear-tail-scalar-condition} is nonpositive.  Any
\(0<\gamma_{\rm tail}<1-\cures\)
therefore works.

\item \textbf{Algebraic decay.}  If
\begin{equation}
\phi_i(\rho)=C_{\rm floor}+a\rho^{-p},
\qquad a,p>0,
\end{equation}
then the same conclusion holds.

\item \textbf{Algebraic growth.}  If
\begin{equation}
\phi_i(\rho)=C_{\rm floor}+a\rho^p,
\qquad a,p>0,
\end{equation}
then
\begin{equation}
\label{eq:low-linear-tail-algebraic-growth-rate}
2\rho\frac{d}{d\rho}\log\phi_i(\rho)
=2p\frac{a\rho^p}{C_{\rm floor}+a\rho^p}
\le2p,
\end{equation}
and the expression converges to \(2p\).  Hence the asymptotic condition is \(p<2(1-\cures)\).  A sufficient finite-box condition is
\begin{equation}
\label{eq:low-linear-tail-algebraic-growth-threshold}
p\left(1+\frac{2|\bar C|}{R_{\rm box}}\right)
<2(1-\cures).
\end{equation}
The strict inequality leaves room to choose \(\gamma_{\rm tail}>0\).

\item \textbf{Positive exponential growth.}  If
\begin{equation}
\phi_i(\rho)=C_{\rm floor}+a\exp(b\rho^p),
\qquad a,b,p>0,
\end{equation}
then
\begin{equation}
2\rho\frac{d}{d\rho}\log\phi_i(\rho)
=2bp\rho^p\frac{a\exp(b\rho^p)}{C_{\rm floor}+a\exp(b\rho^p)},
\end{equation}
\noindent which is unbounded as $\rho\to\infty$. Hence the scalar tail certificate cannot hold.
\end{itemize}

\hfill \\

\begin{center}
\renewcommand{\arraystretch}{1.17}
\setlength{\tabcolsep}{4pt}
\textbf{Summary of radial model rates.}\vspace{0.5em}

\newcolumntype{Y}{>{\raggedright\arraybackslash}X}
\begin{tabularx}{\textwidth}{@{}lXYY@{}}
\toprule
\textbf{Rate}
& \textbf{Model \(\phi_i(\rho)\)}
& \textbf{Asymptotic condition}
& \textbf{Finite-tail condition} \\
\midrule
Exponential decay
& \(C_{\rm floor}+a\exp(-b\rho^p)\)
& \(a,b,p>0\), with \(\cures<1\)
& \eqref{eq:low-linear-tail-basic-phase-condition} \\
\addlinespace[0.15em]
Algebraic decay
& \(C_{\rm floor}+a\rho^{-p}\)
& \(a,p>0\), with \(\cures<1\)
& \eqref{eq:low-linear-tail-basic-phase-condition} \\
\addlinespace[0.15em]
Algebraic growth
& \(C_{\rm floor}+a\rho^p\)
& \(p<2(1-\cures)\)
& \eqref{eq:low-linear-tail-algebraic-growth-threshold} \\
\addlinespace[0.15em]
Exponential growth
& \(C_{\rm floor}+a\exp(b\rho^p)\)
& Not possible
& Not possible \\
\bottomrule
\end{tabularx}
\end{center}

\hfill \\ 

\noindent For the \(r,z\) weights, \eqref{eq:delta-u-tail-growth-main} also requires \(p\ge1\).  Thus an algebraically growing tail is compatible precisely when one can choose
\begin{equation}
1\le p<
\frac{2(1-\cures)}{1+2|\bar C|/R_{\rm box}}.
\label{eq:low-linear-tail-compatible-growth-range}
\end{equation}
The interval in \eqref{eq:low-linear-tail-compatible-growth-range} is nonempty if and only if
\begin{equation}
1+\frac{2|\bar C|}{R_{\rm box}}
<2(1-\cures).
\label{eq:low-linear-tail-compatible-growth-nonempty}
\end{equation}
As \(R_{\rm box}\to\infty\), this reduces to the necessary condition \(\cures<1/2\).

\hfill \\
\hfill
}

\subsubsection{Localized matrix certificate}\label{sssec:localized-matrix-certificate}\label{app:localized-matrix-certificate}

\noindent The global matrix certificate uses a uniform negative upper bound for \(\lambda_{\max}(\mathbb M_{\mathcal L})\) on the whole half-plane. This is sufficient, but not necessary for the integrated energy estimate. An isolated point at which the pointwise matrix margin fails causes no loss in the matrix integral, because a set of measure zero carries no \(L^2\) mass. The relevant issue is a neighborhood on which the uniform negative margin is unavailable. The localized certificate allows finitely many such neighborhoods, provided that the weighted energy that can concentrate in them is quantitatively controlled.

\hfill

\noindent Small area or size by itself does not provide this control, since an \(L^2\) perturbation may place most of its mass in a set of arbitrarily small measure. The localized certificate therefore uses two complementary inputs: a uniform damping bound outside the non-damping region and, on each component of that region, an estimate of the localized weighted energy in terms of the global low-order and top-order energies. The part controlled by \(\mathfrak E_0^2\) reduces the coefficient of the low-order damping term, while the part controlled by \(\mathfrak H_k^2\) appears as an additional top-order term that must be absorbed by the differentiated damping. \\ 

\paragraph{Localized non-damping regions.}
Let the localized non-damping region be a finite union of pairwise disjoint measurable pieces indexed by a finite set \(\mathcal I_{\mathcal L}^{\rm stagn}\):
\begin{equation}
\mathbb{D}_{\mathcal L}^{\rm stagn}
=
\bigcup_{\ell\in\mathcal I_{\mathcal L}^{\rm stagn}}B_\ell,
\qquad
\mathbb{D}
=
\mathbb{D}_{\mathcal L}^{\rm damp}\cup\mathbb{D}_{\mathcal L}^{\rm stagn},
\qquad
\mathbb{D}_{\mathcal L}^{\rm damp}\cap\mathbb{D}_{\mathcal L}^{\rm stagn}=\emptyset.
\label{eq:localized-matrix-region-split}
\end{equation}
The localized certificate below is used only for additional non-damping pieces away from the normalized origin. Accordingly, each $B_\ell$ is chosen bounded and small enough that \(\operatorname{dist}(B_\ell,(0,0))>0\). The origin itself is handled by the weighted matrix estimate with $G_r(0,0)=G_z(0,0)=0$.

\hfill \\ 

\paragraph{Certified inputs.}
On the damping region, we certify the pointwise matrix bound
\begin{equation}
\lambda_{\max}\bigl(\mathbb M_{\mathcal L}(r,z)\bigr)
\le -\Lambda_{\mathcal L}^{\rm damp}
\qquad
(r,z)\in \mathbb{D}_{\mathcal L}^{\rm damp},
\qquad
\Lambda_{\mathcal L}^{\rm damp}>0.
\label{eq:localized-matrix-damp-bound}
\end{equation}
On the non-damping pieces, we certify upper bounds
\begin{equation}
\lambda_{\max}\bigl(\mathbb M_{\mathcal L}(r,z)\bigr)
\le \Lambda_{\mathcal L}^{{\rm stagn},\ell}
\qquad
(r,z)\in B_\ell,
\qquad
\Lambda_{\mathcal L}^{{\rm stagn},\ell}\ge 0.
\label{eq:localized-matrix-stagn-bound}
\end{equation}
Thus the matrix is strictly damping on \(\mathbb{D}_{\mathcal L}^{\rm damp}\), while on each \(B_\ell\) it may lose damping or produce bounded growth.

\hfill \\

\noindent The additional input is the concentration estimate
\begin{equation}
\int_{B_\ell}|\bvar{d\xi}|^2\,dr\,dz
\le
m_{\mathcal L}^{0,{\rm stagn},\ell}\mathfrak E_0^2
+
m_{\mathcal L}^{k,{\rm stagn},\ell}\mathfrak H_k^2,
\qquad \ell\in\mathcal I_{\mathcal L}^{\rm stagn},
\label{eq:localized-matrix-stagn-concentration}
\end{equation}
for nonnegative constants \(m_{\mathcal L}^{0,{\rm stagn},\ell}\) and \(m_{\mathcal L}^{k,{\rm stagn},\ell}\). The first term quantifies the localized mass that is controlled at low order and therefore decreases the effective coefficient of \(\mathfrak E_0^2\). The second quantifies the remaining localized mass controlled by \(\mathfrak H_k^2\) and therefore contributes an additional top-order term.

\hfill \\

\paragraph{Effective matrix estimate.}
For each non-damping piece, set
\begin{equation}
\gamma_\ell
:=
\Lambda_{\mathcal L}^{\rm damp}+\Lambda_{\mathcal L}^{{\rm stagn},\ell}.
\label{eq:localized-matrix-piece-cost}
\end{equation}

\hfill 

\noindent Using \eqref{eq:localized-matrix-region-split} and the identity \(\int_\mathbb{D}|\bvar{d\xi}|^2=\mathfrak E_0^2\),
\begin{align}
\int_\mathbb{D} \bvar{d\xi}^{\top}\mathbb M_{\mathcal L}\bvar{d\xi}\,dr\,dz
&\le
-\Lambda_{\mathcal L}^{\rm damp}
\int_{\mathbb{D}_{\mathcal L}^{\rm damp}}|\bvar{d\xi}|^2\,dr\,dz
+
\sum_{\ell\in\mathcal I_{\mathcal L}^{\rm stagn}}\Lambda_{\mathcal L}^{{\rm stagn},\ell}
\int_{B_\ell}|\bvar{d\xi}|^2\,dr\,dz
\nonumber\\
&=
-\Lambda_{\mathcal L}^{\rm damp}\mathfrak E_0^2
+
\sum_{\ell\in\mathcal I_{\mathcal L}^{\rm stagn}}\gamma_\ell
\int_{B_\ell}|\bvar{d\xi}|^2\,dr\,dz.
\label{eq:localized-matrix-split-derivation}
\end{align}
Define
\begin{align}
\Lambda_{\mathcal L}^{\rm mat,stagn}
&:=
\Lambda_{\mathcal L}^{\rm damp}
-
\sum_{\ell\in\mathcal I_{\mathcal L}^{\rm stagn}}\gamma_\ell m_{\mathcal L}^{0,{\rm stagn},\ell},
\label{eq:localized-matrix-stagn-Lambda-choice}\\
A_{\mathcal L}^{k,\rm stagn}
&:=
\sum_{\ell\in\mathcal I_{\mathcal L}^{\rm stagn}}\gamma_\ell m_{\mathcal L}^{k,{\rm stagn},\ell}.
\label{eq:localized-matrix-stagn-Ak-choice}
\end{align}
Then \eqref{eq:localized-matrix-stagn-concentration} gives
\begin{equation}
\int_\mathbb{D} \bvar{d\xi}^{\top}\mathbb M_{\mathcal L}\bvar{d\xi}\,dr\,dz
\le
-\Lambda_{\mathcal L}^{\rm mat,stagn}\mathfrak E_0^2
+
A_{\mathcal L}^{k,\rm stagn}\mathfrak H_k^2.
\label{eq:localized-matrix-estimate}
\end{equation}
This localized estimate has two effects: it lowers the matrix damping from \(\Lambda_{\mathcal L}^{\rm damp}\) to \(\Lambda_{\mathcal L}^{\rm mat,stagn}\), and it adds the top-order term \(A_{\mathcal L}^{k,\rm stagn}\mathfrak H_k^2\).

\hfill \\ 

\paragraph{Closing conditions.}
After the direct terms are included, we use
\begin{equation}
\DLamLlow
:=
\Lambda_{\mathcal L}^{\rm mat,stagn}-A_{\mathcal L}^{0},
\qquad
\DLamLhigh
:=
A_{\mathcal L}^{k}+A_{\mathcal L}^{k,\rm stagn}.
\label{eq:localized-linear-updated-constants}
\end{equation}

\hfill 

\noindent The remaining low-order damping is positive if and only if
\begin{equation}
\sum_{\ell\in\mathcal I_{\mathcal L}^{\rm stagn}}\gamma_\ell m_{\mathcal L}^{0,{\rm stagn},\ell}
<
\Lambda_{\mathcal L}^{\rm damp}-A_{\mathcal L}^{0}.
\label{eq:localized-matrix-margin-condition}
\end{equation}
Thus the total reduction caused by the low-order parts of the concentration estimates must be smaller than the damping margin available before the localized non-damping region is introduced.

\hfill \\

\noindent The localized estimate introduces the additional top-order loss
\(A_{\mathcal L}^{k,\rm stagn}\mathfrak H_k^2\). This loss must be absorbed by the differentiated linear estimate from Appendix~\ref{app:high-order-linear}.  The differentiated damping comes from the complete high-order matrix, including the first transport commutators:
\begin{equation}
\lambda_{\max}\!\left(\mathbb M_{\mathcal L}^{(k),\rm full}(r,z)\right)
\le-\Lambda_{D^\alpha L}^{\rm full}<0
\qquad\text{for every }(r,z)\in\mathbb{D}.
\label{eq:localized-requires-full-high-certificate-add}
\end{equation}

\hfill \\
\noindent This condition is different from the low-order matrix bound.  At a point where \((G_r,G_z)=0\), the term \((G_r,G_z)\!\cdot\nabla\log\Phi_{i,k}\) vanishes, but the coefficients involving \(\nabla (G_r,G_z)\) remain in the differentiated matrix and can change the sign of the high-order estimate.

\hfill \\
\noindent With the remainder constants from Appendix~\ref{app:high-order-linear}, set
\begin{equation}
\DLamDLhigh
:=\Lambda_{D^\alpha L}^{\rm full}-A_{\rm lin}^{k,\rm rem},
\qquad
\DLamDLlow:=A_{\rm lin}^{0,\rm rem}.
\label{eq:localized-full-high-constants-add}
\end{equation}

\hfill 

\noindent If \(\DLamDLhigh>0\), we choose \(0<\mu_k\le1\) so that
\begin{equation}
A_{\mathcal L}^{k}+A_{\mathcal L}^{k,\rm stagn}
<\mu_k\DLamDLhigh,
\qquad
\mu_k\DLamDLlow
<\Lambda_{\mathcal L}^{\rm mat,stagn}-A_{\mathcal L}^{0}.
\label{eq:localized-matrix-mu-interval-full-add}
\end{equation}

\hfill \\
\noindent The first inequality absorbs all top-order losses generated by the low-order estimate.  The second keeps the remaining low-order damping positive.  If the complete high-order matrix itself fails to damp on a neighborhood, the same localization argument cannot simply be repeated at order \(k\), because \(\mathfrak H_k\) does not control concentration of its own top derivatives.  A separate top-order coercive estimate is then required.

\hfill

\paragraph{Sufficient local concentration estimates.}
To apply the localized certificate, it is enough to prove \eqref{eq:localized-matrix-stagn-concentration} on each piece \(B_\ell\). Since every \(B_\ell\) is bounded and stays a positive distance from the origin, and the weights are singular only at the origin,
\(
\sup_{B_\ell}\wgray{\Phi_i}<\infty
\)
for every \(i\in\{\omega,r,z\}\). The weighted energy can therefore be controlled by ordinary local estimates for the components \(\bvar{\delta v_i}\).

\hfill 

\noindent One option is to use pointwise estimates. If
\begin{equation}
\|\bvar{\delta v_i}\|_{L^\infty(B_\ell)}
\le
P_{i,\ell}^{0}\mathfrak E_0+P_{i,\ell}^{k}\mathfrak H_k,
\qquad i\in\{\omega,r,z\},
\label{eq:localized-matrix-interior-pointwise-input}
\end{equation}
then \(|\bvar{d\xi}|^2=\sum_i\wgray{\Phi_i}|\bvar{\delta v_i}|^2\) and the identity \((a+b)^2\le2a^2+2b^2\) give \eqref{eq:localized-matrix-stagn-concentration},
\begin{equation}
\int_{B_\ell}|\bvar{d\xi}|^2\,dr\,dz
\le
m_{\mathcal L}^{0,{\rm stagn},\ell}\mathfrak E_0^2
+
m_{\mathcal L}^{k,{\rm stagn},\ell}\mathfrak H_k^2,
\qquad \ell\in\mathcal I_{\mathcal L}^{\rm stagn},
\end{equation}
with
\begin{align}
m_{\mathcal L}^{0,{\rm stagn},\ell}
&:=
2|B_\ell|\sum_{i\in\{\omega,r,z\}}
\Bigl(\sup_{B_\ell}\wgray{\Phi_i}\Bigr)(P_{i,\ell}^{0})^2,
\label{eq:localized-matrix-interior-m0}\\
m_{\mathcal L}^{k,{\rm stagn},\ell}
&:=
2|B_\ell|\sum_{i\in\{\omega,r,z\}}
\Bigl(\sup_{B_\ell}\wgray{\Phi_i}\Bigr)(P_{i,\ell}^{k})^2.
\label{eq:localized-matrix-interior-mk}
\end{align}
Thus shrinking \(B_\ell\) improves the constants, provided that the weight suprema and the pointwise constants remain uniformly controlled.

\hfill

\noindent Alternatively, if a local Sobolev, interpolation, restriction, or Poincar\'e estimate gives
\begin{equation}
\int_{B_\ell}|\bvar{\delta v_i}|^2\,dr\,dz
\le
Q_{i,\ell}^{0}\mathfrak E_0^2
+
Q_{i,\ell}^{k}\mathfrak H_k^2,
\qquad i\in\{\omega,r,z\},
\label{eq:localized-matrix-local-L2-input}
\end{equation}
then \eqref{eq:localized-matrix-stagn-concentration} holds with
\begin{equation}
m_{\mathcal L}^{0,{\rm stagn},\ell}
:=
\sum_{i\in\{\omega,r,z\}}
\Bigl(\sup_{B_\ell}\wgray{\Phi_i}\Bigr)Q_{i,\ell}^{0},
\qquad
m_{\mathcal L}^{k,{\rm stagn},\ell}
:=
\sum_{i\in\{\omega,r,z\}}
\Bigl(\sup_{B_\ell}\wgray{\Phi_i}\Bigr)Q_{i,\ell}^{k}.
\label{eq:localized-matrix-local-L2-constants}
\end{equation}

\hfill \\
\noindent\textbf{Summary.} The localized certificate closes once three things are verified: the concentration estimate \eqref{eq:localized-matrix-stagn-concentration}, the positive low-order margin \eqref{eq:localized-matrix-margin-condition}, and the complete differentiated matrix certificate \eqref{eq:localized-requires-full-high-certificate-add}.  The same \(\mu_k\) must satisfy both inequalities in \eqref{eq:localized-matrix-mu-interval-full-add}. 

\hfill \\

\subsubsection{Direct terms}

\noindent After the matrix contribution has been fixed, either through the global constant \(\Lambda_{\mathcal L}^{\rm mat}\) or through the localized matrix certificate in Appendix~\ref{sssec:localized-matrix-certificate}, it remains to evaluate the direct terms in \eqref{eq:low-linear-direct-index-set-omega}--\eqref{eq:low-linear-direct-index-set-z}. The five labels in \(\mathcal S_{\mathcal L}\) are included in the signed damping matrix. \\ 

\paragraph{Undifferentiated velocity labels.}
For \((\omega,2)\), \((r,3)\), and \((z,3)\), we apply Proposition~\ref{prop:fixed-profile-ray-estimate} with the corresponding fixed coefficient and target energy weight. This gives

\begin{align}
\kappa^0_{\omega,2}
&=\mathcal C_{{\rm ray},\omega}(2\partial_z\bar u)^{1/2},\qquad \kappa^k_{\omega,2}=0,
\label{eq:low-linear-ray-constants}\\
\kappa^0_{r,3}
&=\mathcal C_{{\rm ray},r}(2\partial_{rz}\bar\psi)^{1/2},\qquad \kappa^k_{r,3}=0,
\label{eq:low-linear-ray-constants-r}\\
\kappa^0_{z,3}
&=\mathcal C_{{\rm ray},z}(2\partial_{zz}\bar\psi)^{1/2},\qquad \kappa^k_{z,3}=0.
\label{eq:low-linear-ray-constants-z}
\end{align}

\hfill

\paragraph{Direct vorticity--\(z\) coupling.}
\noindent The coupling
\begin{equation}
B_{\omega,3}=2\bar u\,\bvar{\delta v_z}
\end{equation}
\noindent is estimated outside the signed matrix. Putting it in the matrix would require $2\bar u(\Phi_\omega/\Phi_z)^{1/2}$ to be bounded, which is not implied by the independent singularity rates of $\wgray{\Phi_\omega}$ and $\wgray{\Phi_z}$. Instead, we keep $\bar u$ in the $\wgray{\Phi_\omega}$ norm and use the pointwise estimate for $\bvar{\delta v_z}$. Proposition~\ref{prop:pointwise-weighted-energy} gives
\begin{equation}
\|B_{\omega,3}\|_{\Phi_\omega}
\le2\|\bar u\|_{\Phi_\omega}\|\bvar{\delta v_z}\|_{L^\infty(\mathbb D)}
\le2C^{\rm pt}_{z,0,0}\|\bar u\|_{\Phi_\omega}\mathfrak E_k.
\label{eq:low-linear-Bomega3-direct-explicit}
\end{equation}
\noindent Thus one may use
\begin{equation}
\kappa^0_{\omega,3}=2C^{\rm pt}_{z,0,0}\|\bar u\|_{\Phi_\omega},
\qquad
\kappa^k_{\omega,3}=2\sqrt{\mu_k}\,C^{\rm pt}_{z,0,0}\|\bar u\|_{\Phi_\omega}.
\label{eq:low-linear-Bomega3-constants-explicit}
\end{equation}
Its contribution is therefore included in the direct losses \(A_{\mathcal L}^0\) and \(A_{\mathcal L}^k\), not in the matrix eigenvalue.

\hfill \\

\paragraph{Streamfunction labels.}
The streamfunction terms are
\begin{equation}
(\omega,4),(\omega,5),(\omega,6),
(r,4),(r,5),(r,6),(r,7),(r,8),
(z,4),(z,5),(z,6),(z,7),(z,8).
\end{equation}

\noindent Each term has the form \(a\,T(\Gamma_5*\bvar{\delta\omega})\).  \\

\noindent We certify
\begin{equation}
\|a\,T(\Gamma_5*\bvar{\delta\omega})\|_{\wgray{\Phi_i}}
\le
\mathcal K_i^0(a,T)\mathfrak E_0+
\mathcal K_i^k(a,T)\mathfrak H_k.
\end{equation}

\hfill 

\noindent For \(s=0,k\), we enter
\begin{align}
\kappa^s_{\omega,4}&=\mathcal K_\omega^s(-2\partial_z\bar\omega,{\rm Id}),
&
\kappa^s_{\omega,5}&=\mathcal K_\omega^s(r\partial_r\bar\omega,\partial_z),\nonumber\\
\kappa^s_{\omega,6}&=\mathcal K_\omega^s(-r\partial_z\bar\omega,\partial_r),
&
\kappa^s_{r,4}&=\mathcal K_r^s(2\bar u+r\partial_r\bar u,\partial_{rz}),\nonumber\\
\kappa^s_{r,5}&=\mathcal K_r^s(3\partial_r\bar u+r\partial_{rr}\bar u,\partial_z),
&
\kappa^s_{r,6}&=\mathcal K_r^s(-r\partial_z\bar u,\partial_{rr}),\nonumber\\
\kappa^s_{r,7}&=\mathcal K_r^s(-(3\partial_z\bar u+r\partial_{rz}\bar u),\partial_r),
&
\kappa^s_{r,8}&=\mathcal K_r^s(-2\partial_{rz}\bar u,{\rm Id}),\nonumber\\
\kappa^s_{z,4}&=\mathcal K_z^s(2\bar u+r\partial_r\bar u,\partial_{zz}),
&
\kappa^s_{z,5}&=\mathcal K_z^s(r\partial_{rz}\bar u,\partial_z),\nonumber\\
\kappa^s_{z,6}&=\mathcal K_z^s(-r\partial_z\bar u,\partial_{rz}),
&
\kappa^s_{z,7}&=\mathcal K_z^s(-r\partial_{zz}\bar u,\partial_r),\nonumber\\
\kappa^s_{z,8}&=\mathcal K_z^s(-2\partial_{zz}\bar u,{\rm Id}).
\end{align}

\hfill

\noindent If a multiplier is split into scalar summands, we add the certified \(\mathcal K\)-constants for the summands.

\hfill \\

\paragraph{Modulation labels.}
The modulation terms are built from the linear modulation bounds.  First we certify the origin evaluation constants, for \(s=0,k\):

\begin{align}
  |\bvar{\delta\psi}_0|
 &\le P_{\psi,0}^{0}\mathfrak E_0+P_{\psi,0}^{k}\mathfrak H_k,\\
 | (\partial_z\bvar{\delta\psi})_0 |
 &\le P_{\psi_z,0}^{0}\mathfrak E_0+P_{\psi_z,0}^{k}\mathfrak H_k.
\end{align}

\hfill

\noindent Set
\begin{align}
\kappa_{\bvar{\delta C}}^{s}
&=2P_{\psi,0}^{s},\\
\kappa_{\bvar{\delta c_u}}^{s}
&=2P_{\psi_z,0}^{s}.
\end{align}

\hfill  \\

\noindent The modulation terms are \((\omega,7)\), \((\omega,8)\), \((r,9)\), \((r,10)\), \((z,9)\), and \((z,10)\).  For \(s=0,k\), enter
\begin{align}
\kappa^s_{\omega,7}
&=\|\partial_z\bar\omega\|_{\wgray{\Phi_\omega}}\kappa_{\bvar{\delta C}}^{s},
&
\kappa^s_{r,9}
&=\|\partial_{rz}\bar u\|_{\wgray{\Phi_r}}\kappa_{\bvar{\delta C}}^{s},\nonumber\\
\kappa^s_{z,9}
&=\|\partial_{zz}\bar u\|_{\wgray{\Phi_z}}\kappa_{\bvar{\delta C}}^{s},
&
\kappa^s_{\omega,8}
&=\|\bar\omega\|_{\wgray{\Phi_\omega}}\kappa_{\bvar{\delta c_u}}^{s},\nonumber\\
\kappa^s_{r,10}
&=\|\partial_r\bar u\|_{\wgray{\Phi_r}}\kappa_{\bvar{\delta c_u}}^{s},
&
\kappa^s_{z,10}
&=\|\partial_z\bar u\|_{\wgray{\Phi_z}}\kappa_{\bvar{\delta c_u}}^{s}.
\end{align}

\clearpage

\subsubsection{Assembly}

\noindent The direct-term set used in the final sums is ordered by component and derivative count:
\begin{equation}
\label{eq:low-linear-direct-set-present}
\begin{aligned}
\mathfrak S_{\mathcal L}:={}&
\mathfrak S_{\mathcal L,\omega}\cup\mathfrak S_{\mathcal L,r}\cup\mathfrak S_{\mathcal L,z},\\
\mathfrak S_{\mathcal L,\omega}:={}&\{(\omega,2),(\omega,3),(\omega,4),(\omega,5),(\omega,6),(\omega,7),(\omega,8)\},\\
\mathfrak S_{\mathcal L,r}:={}&\{(r,3),(r,4),(r,5),(r,6),(r,7),(r,8),(r,9),(r,10)\},\\
\mathfrak S_{\mathcal L,z}:={}&\{(z,3),(z,4),(z,5),(z,6),(z,7),(z,8),(z,9),(z,10)\}.
\end{aligned}
\end{equation}
This set is exactly the complement of the five matrix labels within the non-transport list. 
\hfill  \\ 

\noindent The terms that have a top-order coefficient are
\begin{equation}
\mathfrak S_{\mathcal L}^{k}
:=
\{(i,m)\in\mathfrak S_{\mathcal L}:\kappa^k_{i,m}>0\}.
\end{equation}

\hfill 

\noindent The direct losses are
\begin{equation}
\label{eq:low-linear-explicit-A0-Ak}
\begin{aligned}
A_{\mathcal L}^{0}
&=
\sum_{(i,m)\in\mathfrak S_{\mathcal L}}\kappa^0_{i,m}
+
\sum_{(i,m)\in\mathfrak S_{\mathcal L}^{k}}\eta_{i,m},\\[0.8ex]
A_{\mathcal L}^{k}
&=
\sum_{(i,m)\in\mathfrak S_{\mathcal L}^{k}}
\frac{(\kappa^k_{i,m})^2}{4\eta_{i,m}}.
\end{aligned}
\end{equation}

\hfill \\

\noindent These formulas come from applying the Cauchy--Schwarz inequality to each direct term and then splitting the mixed term.  Namely,
\begin{equation}
\|B_{i,m}(\bvar{\delta})\|_{\wgray{\Phi_i}}\mathfrak E_0
\le
\kappa^0_{i,m}\mathfrak E_0^2
+\kappa^k_{i,m}\mathfrak E_0\mathfrak H_k,
\end{equation}
and, for \((i,m)\in\mathfrak S_{\mathcal L}^{k}\),
\begin{equation}
\kappa^k_{i,m}\mathfrak E_0\mathfrak H_k
\le
\eta_{i,m}\mathfrak E_0^2
+\frac{(\kappa^k_{i,m})^2}{4\eta_{i,m}}\mathfrak H_k^2.
\end{equation}

\hfill \\

\noindent terms outside \(\mathfrak S_{\mathcal L}^{k}\) have \(\kappa^k_{i,m}=0\), so they contribute only \(\kappa^0_{i,m}\) to \(A_{\mathcal L}^{0}\). 

\hfill \\

\noindent In the standard global-matrix case, combining \eqref{eq:low-linear-explicit-Lambda-mat} and \eqref{eq:low-linear-explicit-A0-Ak}, the certified constants are
\begin{equation}
\label{eq:low-linear-explicit-final-constants}
\DLamLlow=\Lambda_{\mathcal L}^{\rm mat}-A_{\mathcal L}^{0},
\qquad
\DLamLhigh=A_{\mathcal L}^{k}.
\end{equation}

\hfill

\noindent If the localized matrix certificate in Appendix~\ref{sssec:localized-matrix-certificate} is used, the final constants are instead the updated constants in \eqref{eq:localized-linear-updated-constants}.  This keeps the low-order linear estimate in the same form while using the localized matrix contribution from \eqref{eq:localized-matrix-estimate}. \\

\hfill 

\noindent The low-order requirement is \(\DLamLlow>0\). \\

\noindent Combining the low- and high-order linear estimates leaves the top-order coefficient \(\mu_k\DLamDLhigh-\DLamLhigh\).  It is positive exactly when
\begin{equation}
\DLamLhigh<\mu_k\DLamDLhigh.
\label{eq:low-linear-exact-high-loss-absorption}
\end{equation}

\clearpage

\section{High-Order Linear Estimate}\label{app:high-order-linear}

\noindent This appendix presents the high-order linear estimate. The main difference with the low-order estimate is the differentiated transport term: when one derivative falls on a transport coefficient, the resulting commutator is still top order and must be included in the signed damping matrix.  Terms in which two or more derivatives fall on a transport coefficient are estimated as remainders.  The nonlocal linear terms are also kept outside the pointwise matrix and are controlled by the streamfunction estimates.

\hfill

\noindent We use the same transport field \((G_r,G_z)\) defined in \eqref{eq:low-linear-profile-transport-coeffs}.  By \eqref{eq:low-linear-origin-stagnation-add-current},
\begin{equation}
G_r(0,0)=G_z(0,0)=0.
\label{eq:high-linear-origin-stagnation-add-current}
\end{equation}
As before

\begin{equation}
\partial_rG_r=\lambda+\partial_r\bar u^r,
\qquad
\partial_zG_r=\partial_z\bar u^r,
\qquad
\partial_rG_z=\partial_r\bar u^z,
\qquad
\partial_zG_z=\lambda+\partial_z\bar u^z.
\label{eq:high-linear-explicit-G-gradient-add}
\end{equation}

\hfill \\
\noindent To see the commutators directly, we write \(\alpha=(q,k-q)\), with \(0\le q\le k\).  In the Leibniz expansion of
\begin{equation}
    D^\alpha(G_r\partial_r f+G_z\partial_z f),
\end{equation} the first transport commutators are the terms in which exactly one of the \(k\) derivatives falls on \(G_r\) or \(G_z\). There are exactly four such terms: derivatives can fall in either the \(r\)- or \(z\)-direction, and on either component of \((G_r,G_z)\).  Each leaves exactly \(k\) derivatives on \(f\), so none of them may be put into a lower-order remainder.

\hfill \\

\noindent\textbf{Complete top-order matrix.} For \(0\le q\le k\), set
\begin{equation}
V_{i,q}:=\partial_r^q\partial_z^{k-q}\bvar{\delta v_i},
\qquad i\in\{\omega,r,z\}.
\label{eq:high-linear-top-derivative-array-add}
\end{equation}
\noindent When a neighboring index falls outside \(0\le q\le k\), the corresponding term below is simply absent.  Thus, at \(q=0\) there is no \(V_{i,q-1}\) term, and at \(q=k\) there is no \(V_{i,q+1}\) term.

\hfill \\
\noindent The complete top-order transport expression is
\begin{equation}
\begin{aligned}
&-G_r\partial_rV_{i,q}-G_z\partial_zV_{i,q}
-\bigl(q\,\partial_rG_r+(k-q)\,\partial_zG_z\bigr)V_{i,q} 
-q\,\partial_rG_z\,V_{i,q-1}
-(k-q)\,\partial_zG_r\,V_{i,q+1}.
\end{aligned}
\label{eq:high-linear-complete-transport-add}
\end{equation}
\noindent The first two terms transport the top derivative.  The remaining coefficient-derivative contributions are precisely the first transport commutators: \(\partial_rG_r\) and \(\partial_zG_z\) act on \(V_{i,q}\), while \(\partial_rG_z\) and \(\partial_zG_r\) couple \(V_{i,q}\) to the neighboring derivative indices.  The omitted transport terms are exactly those in which at least two derivatives fall on \((G_r,G_z)\). Such terms contain at most \(k-1\) derivatives of the perturbation and are treated in the remainder estimate below.  The local zeroth-order and cross-component linear terms are not transport terms and are listed separately in \eqref{eq:high-linear-explicit-top-local-matrix-add}--\eqref{eq:high-linear-explicit-top-local-z-add}.

\hfill \\
\noindent The remaining local top-order terms are, term by term,
\begin{align}
\omega\text{-term}:\quad
& (\bar c_u-\lambda+\cures)V_{\omega,q}+2\bar u\,V_{z,q},
\label{eq:high-linear-explicit-top-local-matrix-add}\\
r\text{-term}:\quad
& \left(2\partial_z\bar\psi+\bar c_u-\lambda-\partial_r\bar u^r+\cures\right)V_{r,q}
-\partial_r\bar u^z\,V_{z,q},
\label{eq:high-linear-explicit-top-local-r-add}\\
z\text{-term}:\quad
& -\partial_z\bar u^r\,V_{r,q}
+\left(2\partial_z\bar\psi+\bar c_u-\lambda-\partial_z\bar u^z+\cures\right)V_{z,q}.
\label{eq:high-linear-explicit-top-local-z-add}
\end{align}
\noindent These are exactly the local terms for which all \(k\) derivatives remain on the perturbation.  If a derivative lands on one of their spatially varying coefficients, the resulting perturbation factor has order at most \(k-1\) and is treated as a remainder.

\clearpage

\noindent As in the low-order construction, it is convenient to write the signed quadratic form directly as a symmetric matrix.  Order the weighted top derivatives by component and derivative index:
\begin{equation}
\mathbf V_k:=
\left(
\Phi_{\omega,k}^{1/2}V_{\omega,0},\ldots,\Phi_{\omega,k}^{1/2}V_{\omega,k},
\Phi_{r,k}^{1/2}V_{r,0},\ldots,\Phi_{r,k}^{1/2}V_{r,k},
\Phi_{z,k}^{1/2}V_{z,0},\ldots,\Phi_{z,k}^{1/2}V_{z,k}
\right)^{\!\top}.
\label{eq:high-linear-weighted-vector-add}
\end{equation}

\hfill 

\noindent The matrix \(\mathbb M_{\mathcal L}^{(k),\rm full}\) is indexed by the pairs \((i,q)\), with \(i\in\{\omega,r,z\}\) and \(0\le q\le k\):
\begin{equation}
\mathbb M_{\mathcal L}^{(k),\rm full}
=
\left[
\left(\mathbb M_{\mathcal L}^{(k),\rm full}\right)_{(i,q),(j,p)}
\right]_{\substack{i,j\in\{\omega,r,z\}\\0\le q,p\le k}},
\qquad
\mathbb M_{\mathcal L}^{(k),\rm full}
=\left(\mathbb M_{\mathcal L}^{(k),\rm full}\right)^{\!\top}.
\label{eq:high-linear-explicit-block-matrix-add}
\end{equation}

\hfill

\paragraph{Diagonal entries.} The diagonal entries are obtained exactly as in the low-order matrix: weighted integration by parts contributes the transport divergence and weight-gradient terms, the first commutators contribute the terms proportional to \(q\) and \(k-q\), and equations \eqref{eq:high-linear-explicit-top-local-matrix-add}--\eqref{eq:high-linear-explicit-top-local-z-add} contribute the local diagonal coefficients.  Explicitly, for \(0\le q\le k\),
\begin{align}
\left(\mathbb M_{\mathcal L}^{(k),\rm full}\right)_{(\omega,q),(\omega,q)}
={}&
\frac{1}{2}\left(
\partial_rG_r+\partial_zG_z
+G_r\partial_r\log\Phi_{\omega,k}
+G_z\partial_z\log\Phi_{\omega,k}
\right)\nonumber\\
& \quad +\bar c_u-\lambda+\cures
-q\,\partial_rG_r-(k-q)\,\partial_zG_z,
\label{eq:high-linear-explicit-diagonal-omega-add}\\
\left(\mathbb M_{\mathcal L}^{(k),\rm full}\right)_{(r,q),(r,q)}
={}&
\frac{1}{2}\left(
\partial_rG_r+\partial_zG_z
+G_r\partial_r\log\Phi_{r,k}
+G_z\partial_z\log\Phi_{r,k}
\right)\nonumber\\
& \quad  +2\partial_z\bar\psi+\bar c_u-\lambda-\partial_r\bar u^r+\cures
-q\,\partial_rG_r-(k-q)\,\partial_zG_z,
\label{eq:high-linear-explicit-diagonal-r-add}\\
\left(\mathbb M_{\mathcal L}^{(k),\rm full}\right)_{(z,q),(z,q)}
={}&
\frac{1}{2}\left(
\partial_rG_r+\partial_zG_z
+G_r\partial_r\log\Phi_{z,k}
+G_z\partial_z\log\Phi_{z,k}
\right)\nonumber\\
& \quad  +2\partial_z\bar\psi+\bar c_u-\lambda-\partial_z\bar u^z+\cures
-q\,\partial_rG_r-(k-q)\,\partial_zG_z.
\label{eq:high-linear-explicit-diagonal-z-add}
\end{align}

\hfill

\paragraph{Same-component off-diagonal entries.} For each component \(i\in\{\omega,r,z\}\), the only same-component off-diagonal entries are the neighboring derivative couplings.  For \(0\le q<k\),
\begin{equation}
\left(\mathbb M_{\mathcal L}^{(k),\rm full}\right)_{(i,q),(i,q+1)}
=
\left(\mathbb M_{\mathcal L}^{(k),\rm full}\right)_{(i,q+1),(i,q)}
=
-\frac{1}{2}\left[(k-q)\,\partial_z\bar u^r+(q+1)\,\partial_r\bar u^z\right].
\label{eq:high-linear-explicit-neighbor-entry-add}
\end{equation}

\hfill

\paragraph{Cross-component entries.} The cross-component entries are obtained as in the low-order matrix by keeping each coefficient together with its required weight conversion and then symmetrizing.  For \(0\le q\le k\),
\begin{align}
\left(\mathbb M_{\mathcal L}^{(k),\rm full}\right)_{(\omega,q),(z,q)}
&=
\left(\mathbb M_{\mathcal L}^{(k),\rm full}\right)_{(z,q),(\omega,q)}
=
\bar u\left(\frac{\Phi_{\omega,k}}{\Phi_{z,k}}\right)^{1/2},
\label{eq:high-linear-explicit-cross-blocks-add}\\
\left(\mathbb M_{\mathcal L}^{(k),\rm full}\right)_{(r,q),(z,q)}
&=
\left(\mathbb M_{\mathcal L}^{(k),\rm full}\right)_{(z,q),(r,q)}=
-\frac{1}{2}\left[
\partial_r\bar u^z\left(\frac{\Phi_{r,k}}{\Phi_{z,k}}\right)^{1/2}
+\partial_z\bar u^r\left(\frac{\Phi_{z,k}}{\Phi_{r,k}}\right)^{1/2}
\right].
\label{eq:high-linear-explicit-cross-rz-add}
\end{align}

\hfill 

\noindent For the $\omega$--$z$ entry, the quantity to certify is the complete multiplier $\bar u \sqrt{ \Phi_{\omega,k}/\Phi_{z,k}}$. Near the origin, $\wgray{\Phi_{i,k}}=1+\rho^k\wgray{\Phi_i}$ and $\rho^k\wgray{\Phi_i}\to0$, so both high-order weights tend to $1$ and this multiplier has no singular weight-conversion obstruction there. On the rest of the finite box and on the tail, the complete product is enclosed directly. \\

\noindent There are no \(\omega\)-\(r\) entries.  There are also no cross-component entries with different derivative indices \(q\ne p\).  Together with symmetry, \eqref{eq:high-linear-explicit-diagonal-omega-add}--\eqref{eq:high-linear-explicit-cross-rz-add} therefore specify every entry of \(\mathbb M_{\mathcal L}^{(k),\rm full}\).

\hfill \\

\paragraph{Damping certificate.} The required high-order damping certificate is
\begin{equation}
\operatorname*{ess\,sup}_{(r,z)\in\mathbb D} \lambda_{\max}\!\left(\mathbb M_{\mathcal L}^{(k),\rm full}(r,z)\right) \le-\Lambda_{D^\alpha L}^{\rm full}<0.
\label{eq:high-linear-full-matrix-certificate-add}
\end{equation}

\hfill

\noindent All differentiated linear terms not represented by this matrix are required to satisfy
\begin{equation}
\left|
\sum_{i\in\{\omega,r,z\}}\sum_{|\alpha|=k}
\left\langle D^\alpha\mathcal L_i,D^\alpha\bvar{\delta v_i}\right\rangle_{\Phi_{i,k}}
-
\int_\mathbb{D} \mathbf V_k^{\top}\mathbb M_{\mathcal L}^{(k),\rm full}\mathbf V_k\,dr\,dz
\right|
\le
A_{\rm lin}^{k,\rm rem}\mathfrak H_k^2
+A_{\rm lin}^{0,\rm rem}\mathfrak E_0^2.
\label{eq:high-linear-full-remainder-bound-add}
\end{equation}

\hfill \\
\begin{theorem}[High-order linear estimate]\label{prop:high-order-linear}
\noindent Assume \eqref{eq:high-linear-full-matrix-certificate-add}. Assume \eqref{eq:high-linear-full-remainder-bound-add} and
\begin{equation}
A_{\rm lin}^{k,\rm rem}<\Lambda_{D^\alpha L}^{\rm full}.
\label{eq:high-linear-full-euler-absorption-add}
\end{equation}
\noindent Then there are constants \(\DLamDLhigh>0\) and \(\DLamDLlow\ge0\) such that
\begin{equation}
\sum_{i\in\{\omega,r,z\}}\sum_{|\alpha|=k}
\left\langle D^\alpha\mathcal L_i,D^\alpha\bvar{\delta v_i}\right\rangle_{\Phi_{i,k}}
\le
\DLamDLlow \, \mathfrak E_0^2-\DLamDLhigh \, \mathfrak H_k^2.
\label{eq:high-order-linear-final-full-add}
\end{equation}

  \end{theorem}

\hfill 

\begin{grayproof}
\noindent The proof proceeds in 5 steps. \\

\noindent \textbf{1. Separate the top-order transport terms.}
\noindent For \(\alpha=(q,k-q)\), the exact Leibniz expansion is
\begin{equation}
\begin{aligned}
D^\alpha\!\left(-G_r\partial_r\bvar{\delta v_i}-G_z\partial_z\bvar{\delta v_i}\right)
={}&
-G_r\partial_rV_{i,q}-G_z\partial_zV_{i,q}
-\bigl(q\,\partial_rG_r+(k-q)\,\partial_zG_z\bigr)V_{i,q}\\
&-q\,\partial_rG_z\,V_{i,q-1}
-(k-q)\,\partial_zG_r\,V_{i,q+1}\\
&-
\sum_{\substack{0<\beta\le\alpha\\ |\beta|\ge2}}
\binom{\alpha}{\beta}
\left[
(D^\beta G_r)\partial_rD^{\alpha-\beta}\bvar{\delta v_i}
+(D^\beta G_z)\partial_zD^{\alpha-\beta}\bvar{\delta v_i}
\right].
\end{aligned}
\label{eq:high-linear-exact-transport-leibniz-add}
\end{equation}
\noindent The first two lines are exactly the retained transport expression \eqref{eq:high-linear-complete-transport-add}.  In the last line the perturbation has derivative order
\begin{equation}
k+1-|\beta|\le k-1,
\label{eq:high-linear-remainder-order-add}
\end{equation}
\noindent so those terms belong to the remainder. \\

\noindent The same distinction applies to local coefficients: the terms in \eqref{eq:high-linear-explicit-top-local-matrix-add}--\eqref{eq:high-linear-explicit-top-local-z-add} keep all \(k\) derivatives on the perturbation, while any term with a derivative on a spatially varying local coefficient has at most \(k-1\) perturbation derivatives.

\hfill \\

\noindent\textbf{2. Compute the signed top-order form.}
\noindent For each \(i,q\), weighted integration by parts gives
\begin{equation}
\left\langle
-(G_r,G_z)\!\cdot\nabla V_{i,q},V_{i,q}
\right\rangle_{\Phi_{i,k}}
=
\frac{1}{2}\int_\mathbb{D}
\left(
\partial_rG_r+\partial_zG_z
+G_r\partial_r\log\Phi_{i,k}
+G_z\partial_z\log\Phi_{i,k}
\right)
\left|\Phi_{i,k}^{1/2}V_{i,q}\right|^2\,dr\,dz.
\label{eq:high-linear-top-transport-ibp-add}
\end{equation}

\hfill 

\noindent The two first commutators that stay on \(V_{i,q}\) contribute
\begin{equation}
-q\,\partial_rG_r-(k-q)\,\partial_zG_z.
\label{eq:high-linear-diagonal-commutator-add}
\end{equation}
\noindent The other two first commutators couple neighboring derivative indices and give \eqref{eq:high-linear-explicit-neighbor-entry-add} after symmetrization.  The local cross-component terms give \eqref{eq:high-linear-explicit-cross-blocks-add}--\eqref{eq:high-linear-explicit-cross-rz-add}. \\

\noindent Therefore the entire retained top-order pairing is
\begin{equation}
\int_\mathbb{D}
\mathbf V_k^{\top}\mathbb M_{\mathcal L}^{(k),\rm full}\mathbf V_k\,dr\,dz.
\label{eq:high-linear-full-matrix-identity-add}
\end{equation}

\hfill 

\noindent Since \(\int_\mathbb{D}|\mathbf V_k|^2\,dr\,dz=\mathfrak H_k^2\), the certificate \eqref{eq:high-linear-full-matrix-certificate-add} gives
\begin{equation}
\int_\mathbb{D}
\mathbf V_k^{\top}\mathbb M_{\mathcal L}^{(k),\rm full}\mathbf V_k\,dr\,dz
\le
-\Lambda_{D^\alpha L}^{\rm full}\mathfrak H_k^2.
\label{eq:high-linear-full-signed-bound-add}
\end{equation}

\hfill \\
\noindent At a point where \((G_r,G_z)=0\), the weight-gradient terms in the diagonal entries vanish, but the first-commutator terms involving \(\nabla (G_r,G_z)\) remain.  

\hfill \\

\hfill 

\noindent\textbf{3. Bound the remaining terms.}
\noindent The transport remainders have the form
\begin{equation}
\binom{\alpha}{\beta}
(D^\beta (G_r,G_z))\!\cdot\nabla D^{\alpha-\beta}\bvar{\delta v_i},
\qquad 2\le|\beta|\le k,
\label{eq:high-linear-lower-transport-remainder-add}
\end{equation}
\noindent while a typical local coefficient remainder has the form
\begin{equation}
\binom{\alpha}{\beta}(D^\beta a)D^{\alpha-\beta}\bvar{\delta v_j},
\qquad 1\le|\beta|\le k.
\label{eq:high-linear-lower-zero-order-remainder-add}
\end{equation}
\noindent In both cases the perturbation factor has order at most \(k-1\). \\

\noindent We put the differentiated steady coefficient in \(L^\infty\). \\   

\noindent For \(0<|\gamma|<k\), interpolation \eqref{eq:ns-lower-from-top} and the Cauchy--Schwarz inequality over the \(k+1\) derivatives of total order \(k\) give
\begin{equation}
\|D^\gamma\bvar{\delta v_j}\|_{\Phi_{j,|\gamma|}}
\le
\eta\sqrt{k+1}\,\mathfrak H_k
+I_{j,|\gamma|,k}(\eta)\mathfrak E_0,
\qquad 0<|\gamma|<k.
\label{eq:linear-intermediate-reduction}
\end{equation}

\hfill 

\noindent The order-zero endpoint is controlled directly by \(\mathfrak E_0\). For a source component \(j\) and target component \(i\), we keep the coefficient together with the required weight conversion:
\begin{equation}
\left|\left\langle aD^\gamma f,D^k g\right\rangle_{\Phi_{i,k}}\right|
\le
\left\|a\left(\frac{\Phi_{i,k}}{\Phi_{j,|\gamma|}}\right)^{1/2}\right\|_{L^\infty(\mathbb D)}
 \|D^\gamma f\|_{\Phi_{j,|\gamma|}}
\|D^k g\|_{\Phi_{i,k}}.
\label{eq:high-linear-weight-bridge-used-add}
\end{equation} 

\hfill

 \noindent For a transport remainder in \eqref{eq:high-linear-lower-transport-remainder-add}, the perturbation factor has derivative level $k+1-|\beta|$. Thus the certificate uses the complete multiplier
\begin{equation}
\left\|
(D^\beta G_\ell)
\left(\frac{\Phi_{i,k}}{\Phi_{i,k+1-|\beta|}}\right)^{1/2}
\right\|_{L^\infty(\mathbb D)},
\qquad 2\le|\beta|\le k,\quad \ell\in\{r,z\}.
\label{eq:high-linear-transport-level-bridge-add}
\end{equation} 

\hfill 

\noindent Streamfunction terms are first reduced by \Cref{thm:stream-elliptic-used}.  After these reductions, each remainder is bounded by a combination of \(\mathfrak H_k^2\) and \(\mathfrak E_0\mathfrak H_k\).  We split the mixed term by
\begin{equation}
C\mathfrak E_0\mathfrak H_k
\le
   \eta_{\rm lin} \mathfrak H_k^2+   \frac{C^2}{4\eta_{\rm lin}} \mathfrak E_0^2,
\qquad    \eta_{\rm lin} >0.
\label{eq:linear-young-used}
\end{equation} 

\hfill 

\noindent Summing the certified term bounds yields \eqref{eq:high-linear-full-remainder-bound-add}.  The essential point is that no absolute-value estimate is applied until after the complete signed matrix has been extracted.

\hfill \\

\hfill 

\noindent\textbf{4. Conclude the estimate.}
\noindent Combining \eqref{eq:high-linear-full-signed-bound-add} and \eqref{eq:high-linear-full-remainder-bound-add} gives
\begin{equation}
\sum_{i,\ |\alpha|=k}
\left\langle D^\alpha\mathcal L_i,D^\alpha\bvar{\delta v_i}\right\rangle_{\Phi_{i,k}}
\le
-\left(\Lambda_{D^\alpha L}^{\rm full}-A_{\rm lin}^{k,\rm rem}\right)\mathfrak H_k^2
+A_{\rm lin}^{0,\rm rem}\mathfrak E_0^2.
\label{eq:high-linear-euler-combined-add}
\end{equation}
\noindent Hence
\begin{equation}
\DLamDLhigh
:=\Lambda_{D^\alpha L}^{\rm full}-A_{\rm lin}^{k,\rm rem}>0,
\qquad
\DLamDLlow:=A_{\rm lin}^{0,\rm rem},
\label{eq:high-linear-full-euler-constants-add}
\end{equation}
\noindent which proves \eqref{eq:high-order-linear-final-full-add}.

  \end{grayproof}

\clearpage

\section{Low-Order Nonlinear Estimate}\label{app:low-order-nonlinear}

\noindent The low-order nonlinear estimate uses only the variables already present in the low-order energy.  Each nonlinear term is treated by one of two templates.  Transport terms are integrated by parts in the weighted inner product.  Non-transport product terms are estimated directly by placing one factor in \(L^\infty\) and the other factors in the weighted \(L^2\) norms controlled by the energy. \\

\noindent The modulation coefficients
\(
\bvar{\delta C}=-2\bvar{\delta\psi}_0
\)
and
\(
\bvar{\delta c_u}=-2(\partial_z\bvar{\delta\psi})_0
\)
are exactly linear in the perturbation. Their fixed-profile contributions are therefore included in the linear operator, with no additional fixed-profile modulation terms in the nonlinear operator. 

\hfill \\

\noindent If \(a=a(r,z,t)\) is a scalar transport coefficient, Proposition~\ref{prop:weighted-transport} gives the estimate
\begin{equation}
\label{eq:ns-transport-template}
|\langle a\partial_\ell f,f\rangle_{\wgray{\Phi_i}}|
\le \frac{1}{2}\left\|\partial_\ell a+a\frac{\partial_\ell\wgray{\Phi_i}}{\wgray{\Phi_i}}\right\|_{L^\infty}
\|f\|_{\wgray{\Phi_i}}^2.
\end{equation}

\hfill

\noindent For a non-transport product term, we use weighted Hölder in the form
\begin{equation}
\label{eq:ns-source-template}
\int_\mathbb{D}|ABC|\wgray{\Phi_i}
\le
\left\|A\left(\frac{\Phi_i}{\wgray{\Phi_j}}\right)^{1/2}\right\|_{L^\infty(\mathbb D)}
\|B\|_{\wgray{\Phi_j}}\|C\|_{\wgray{\Phi_i}}.
\end{equation}

\hfill

\noindent The pointwise factor is controlled by Proposition~\ref{prop:pointwise-weighted-energy}, streamfunction factors are first controlled by \Cref{thm:stream-elliptic-used}, and modulation factors are controlled by \eqref{eq:euler-modulation-direct-bounds-add}. When \(i\ne j\), the certificate bounds the complete multiplier
\(
A(\Phi_i/\wgray{\Phi_j})^{1/2}
\)
for the actual summand.

\hfill 

\noindent Since
\begin{equation}
\bvar{\delta C}=-2\bvar{\delta\psi}_0,
\qquad
\bvar{\delta c_u}=-2(\partial_z\bvar{\delta\psi})_0,
\end{equation}
both modulation coefficients are exactly linear and are included in the linear operator. Hence there are no additional fixed-profile modulation source terms in the nonlinear operator.

\hfill

\begin{theorem}[Low-order nonlinear estimate]\label{prop:low-order-nonlinear}
The nonlinear terms satisfy

\begin{equation}
\left|\langle\mathcal N_\omega,\bvar{\delta\omega}\rangle_{\wgray{\Phi_\omega}}
+\langle\mathcal N_r,\partial_r\bvar{\delta u}\rangle_{\wgray{\Phi_r}}
+\langle\mathcal N_z,\partial_z\bvar{\delta u}\rangle_{\wgray{\Phi_z}}\right|
\le \DLamNthree \, \mathfrak E_k^3.
\label{eq:low-order-nonlinear-final}
\end{equation}

\hfill 

\end{theorem}

\begin{grayproof}
\noindent\textbf{Roadmap.} We separate the nonlinear pairings into transport terms and all other products. A transport term $a\partial_\ell\bvar{\delta v_i}$ paired with $\bvar{\delta v_i}$ is handled by \eqref{eq:ns-transport-template}. Every other product is handled by \eqref{eq:ns-source-template}. In either case, two factors remain in weighted $L^2$, while the remaining multiplier is bounded by $\mathcal O(\mathfrak E_k)$.

\hfill

 \noindent\textbf{Transport summands.}  Let \(\mathcal T_0\) be the finite set of transport summands in \(\mathcal N_\omega\), \(\mathcal N_r\), and \(\mathcal N_z\).  An element \(T\in\mathcal T_0\) records one additive summand together with the paired component \(i\), the direction \(\ell\), and the coefficient \(a_T\).  Applying \eqref{eq:ns-transport-template} gives
\begin{equation}
|\langle a_T\partial_\ell\bvar{\delta v_i},\bvar{\delta v_i}\rangle_{\wgray{\Phi_i}}|
\le A_T\|\bvar{\delta v_i}\|_{\wgray{\Phi_i}}^2,
\qquad
A_T:=\frac{1}{2}\left\|\partial_\ell a_T+a_T\frac{\partial_\ell\wgray{\Phi_i}}{\wgray{\Phi_i}}\right\|_{L^\infty}.
\label{eq:low-transport-summand}
\end{equation}

\hfill 

\noindent For axial transport, we keep the centered coefficient
\begin{equation}
a_T=-(\bvar{\delta C}+\bvar{\delta u^z})
\label{eq:low-transport-centered-axial-add}
\end{equation}
\noindent intact. \\

\noindent Since $\bvar{\delta C}=-2\bvar{\delta\psi}_0$ and $\bvar{\delta u^z}(0,0)=2\bvar{\delta\psi}_0$, one has $a_T(0,0)=0$. This cancellation is used before the factor $\partial_z\Phi_i/\Phi_i$ is estimated in \eqref{eq:low-transport-summand}. \\

\noindent The coefficients \(a_T\) are made from \(\bvar{\delta C}\), \(\bvar{\delta u^r}\), \(\bvar{\delta u^z}\), and the streamfunction expressions defining these velocities. Since \(\bvar{\delta C}=-2\bvar{\delta\psi}_0\) is exactly linear, the modulation estimates, pointwise estimates, and streamfunction estimates give a constant \(A_T^{(1)}\) such that
\begin{equation}
A_T\le A_T^{(1)}\mathfrak E_k.
\label{eq:low-transport-coefficient-split}
\end{equation}

\hfill 

\noindent Since \(\|\bvar{\delta v_i}\|_{\wgray{\Phi_i}}\le \mathfrak E_k\), the summand \(T\) contributes
\begin{equation}
A_T^{(1)}\mathfrak E_k^3.
\label{eq:low-transport-summand-contribution}
\end{equation}

\hfill

\noindent\textbf{Non-transport product summands.} Let $\mathcal S_0$ be the finite set of the remaining product summands. For $S\in\mathcal S_0$, we write $B_S$ for the perturbative factor measured in $\wgray{\Phi_j}$, set $C_S=\bvar{\delta v_i}$ for the tested factor, and we place all remaining factors in $A_S$. Then \eqref{eq:ns-source-template} gives
\begin{equation}
\int_\mathbb{D}|A_SB_SC_S|\wgray{\Phi_i}
\le
\left\|A_S\left(\frac{\wgray{\Phi_i}}{\wgray{\Phi_j}}\right)^{1/2}\right\|_{L^\infty(\mathbb D)}
\|B_S\|_{\wgray{\Phi_j}}\|C_S\|_{\wgray{\Phi_i}}.
\label{eq:low-source-summand}
\end{equation}

\hfill

\noindent We apply Proposition~\ref{prop:pointwise-weighted-energy}, \Cref{thm:stream-elliptic-used}, and the modulation bounds to the actual factors in $A_S$. Keeping any required weight conversion inside the same multiplier gives the certified bound
\begin{equation}
\left\|A_S\left(\frac{\Phi_i}{\wgray{\Phi_j}}\right)^{1/2}\right\|_{L^\infty(\mathbb D)}
\le C_S^{(1)}\mathfrak E_k
\label{eq:low-source-pointwise-split}
\end{equation}
for every summand. The two weighted factors are controlled by $\mathfrak E_k$, with any fixed constant absorbed into the summand coefficient. 

\hfill

\noindent Thus each non-transport product summand contributes
\begin{equation}
B_S^{(1)}\mathfrak E_k^3.
\label{eq:low-source-summand-contribution}
\end{equation}

\hfill \\

\noindent\textbf{Finite summation.} Since the nonlinearities contain finitely many summands, the sets $\mathcal T_0$ and $\mathcal S_0$ are finite. Define
\begin{equation}
\DLamNthree:=\sum_{T\in\mathcal T_0}A_T^{(1)}
+\sum_{S\in\mathcal S_0}B_S^{(1)}.
\label{eq:low-nonlinear-constants-defined}
\end{equation}

\hfill 

\noindent Summing \eqref{eq:low-transport-summand-contribution} and \eqref{eq:low-source-summand-contribution} over all summands gives \eqref{eq:low-order-nonlinear-final}. 
 \end{grayproof}

\clearpage

\section{High-Order Nonlinear Estimate}\label{app:high-order-nonlinear}

\noindent The top-order nonlinear estimate has one potential derivative loss.  Differentiating a transport term \(a\partial_\ell\bvar{\delta v_i}\) produces \(a\partial_\ell D^\alpha\bvar{\delta v_i}\), which contains one derivative beyond the energy.  Weighted integration by parts removes that derivative.  Every remaining term is then estimated by the same simple rule: we keep two factors in weighted \(L^2\), place all lower-order factors in a combined \(L^\infty\) multiplier, and keep the weight conversion inside that multiplier. \\ 

\begin{theorem}[High-order nonlinear estimate]\label{prop:high-order-nonlinear}

Assume \(k\ge4\) and that the standing certified product, interpolation, elliptic, modulation, and profile bounds hold.  Then

\begin{equation}
\mu_k\left|
\sum_{|\alpha|=k}\langle D^\alpha\mathcal N_\omega,D^\alpha\bvar{\delta\omega}\rangle_{\wgray{\Phi_{\omega,k}}}
+\sum_{|\alpha|=k}\langle D^\alpha\mathcal N_r,D^\alpha\partial_r\bvar{\delta u}\rangle_{\wgray{\Phi_{r,k}}}
+\sum_{|\alpha|=k}\langle D^\alpha\mathcal N_z,D^\alpha\partial_z\bvar{\delta u}\rangle_{\wgray{\Phi_{z,k}}}
\right|
\le \DLamDNthree \, \mathfrak E_k^3.
\label{eq:high-order-nonlinear-final}
\end{equation}

\end{theorem}

\hfill 

\begin{grayproof}
\noindent Fix a component \(i\in\{\omega,r,z\}\) and a multi-index \(|\alpha|=k\).  The tested factor is always \(D^\alpha\bvar{\delta v_i}\), and the top-order energy gives
\begin{equation}
\sqrt{\mu_k}\,
\|D^\alpha\bvar{\delta v_i}\|_{\wgray{\Phi_{i,k}}}
\le \mathfrak E_k.
\label{eq:high-nonlinear-one-top-factor-energy}
\end{equation}

\hfill 

\noindent The weighted product estimate used below is written directly in terms of the energy weights.  If a Leibniz summand contains a perturbative factor \(D^\gamma\bvar{\delta v_j}\) of order \(s=|\gamma|\), and \(A\) denotes the product of the remaining factors in that summand, then
\begin{align}
&\left|
\int_\mathbb{D} A\,D^\gamma\bvar{\delta v_j}\,
D^\alpha\bvar{\delta v_i}\,\wgray{\Phi_{i,k}}\,dr\,dz
\right| \le
\left\|A
\left(\frac{\wgray{\Phi_{i,k}}}{\wgray{\Phi_{j,s}}}\right)^{1/2}
\right\|_{L^\infty}
\|D^\gamma\bvar{\delta v_j}\|_{\wgray{\Phi_{j,s}}}
\|D^\alpha\bvar{\delta v_i}\|_{\wgray{\Phi_{i,k}}}.
\label{eq:high-nonlinear-two-slot-holder}
\end{align}

\hfill 

\noindent The ratio of weights stays inside the pointwise multiplier.  Streamfunction factors are first converted to vorticity estimates by \Cref{thm:stream-elliptic-used}.  For \(0\le s\le k\), the low-order energy, the top-order energy, and the certified interpolation estimate imply
\begin{equation}
\sqrt{\mu_k}\,
\|D^\gamma\bvar{\delta v_j}\|_{\wgray{\Phi_{j,s}}}
\le C_{j,s}\mathfrak E_k,
\qquad |\gamma|=s,
\label{eq:high-nonlinear-selected-factor-energy-conversion}
\end{equation}

\noindent with the corresponding elliptic constant included when the original factor contains \(\bvar{\delta\psi}\). \\

\noindent Every nonlinear summand is quadratic because $\bvar{\delta C}$ and $\bvar{\delta c_u}$ are exactly linear and their fixed-profile contributions are included in the linear operator. Thus the remaining multiplier is bounded by \(C^{(1)}\mathfrak E_k\). Together with the two weighted \(L^2\) factors in \eqref{eq:high-nonlinear-two-slot-holder}, this contributes a cubic term.

\hfill \\ 

\noindent\textbf{Principal transport terms.}

\hfill 

\noindent For an axial transport term, we take $a=-(\bvar{\delta C}+\bvar{\delta u^z})$ throughout the principal term and its commutators. Its origin cancellation is retained before the logarithmic derivative of the weight is estimated.\\

\noindent We consider a nonlinear transport product \(a\partial_\ell\bvar{\delta v_i}\).  Leibniz' rule gives
\begin{equation}
D^\alpha(a\partial_\ell\bvar{\delta v_i})
=a\partial_\ell D^\alpha\bvar{\delta v_i}
+\sum_{0<\beta\le\alpha}\binom{\alpha}{\beta}
(D^\beta a)D^{\alpha-\beta}\partial_\ell\bvar{\delta v_i}.
\label{eq:high-nonlinear-transport-leibniz}
\end{equation}
Only the first term appears to contain one derivative beyond the energy.  Set \(f=D^\alpha\bvar{\delta v_i}\) and apply Proposition~\ref{prop:weighted-transport} with the weight \(\wgray{\Phi_{i,k}}\):
\begin{equation}
\left\langle a\partial_\ell f,f\right\rangle_{\wgray{\Phi_{i,k}}}
=-\frac12\int_\mathbb{D}
\left(
\partial_\ell a
+a\frac{\partial_\ell\wgray{\Phi_{i,k}}}{\wgray{\Phi_{i,k}}}
\right)|f|^2\wgray{\Phi_{i,k}}\,dr\,dz.
\label{eq:high-nonlinear-principal-transport-identity}
\end{equation}

\hfill

\noindent The pointwise estimates, the modulation bounds, and \Cref{thm:stream-elliptic-used} give
\begin{equation}
\frac12\left\|
\partial_\ell a
+a\frac{\partial_\ell\wgray{\Phi_{i,k}}}{\wgray{\Phi_{i,k}}}
\right\|_{L^\infty}
\le A_{i,\alpha,\ell}^{(1)}\mathfrak E_k.
\label{eq:high-nonlinear-principal-coeff-split}
\end{equation}

\hfill 

\noindent Combining this with \(\mu_k\|f\|_{\wgray{\Phi_{i,k}}}^2\le\mathfrak E_k^2\) gives the required cubic bound.  This weighted integration by parts removes the only possible derivative loss.

\hfill \\

\noindent\textbf{Transport commutators.}
A commutator term is
\begin{equation}
\binom{\alpha}{\beta}(D^\beta a)
D^{\alpha-\beta}\partial_\ell\bvar{\delta v_i},
\qquad 0<\beta\le\alpha.
\label{eq:high-nonlinear-commutator-model}
\end{equation}
The second factor has derivative order
\begin{equation}
s=k-|\beta|+1.
\label{eq:high-nonlinear-commutator-complementary-order}
\end{equation}

\hfill 

\noindent If \(|\beta|\le k-2\), then \(D^\beta a\) is a lower-order factor and is placed in \(L^\infty\), while the order-\(s\) factor remains in the weight \(\wgray{\Phi_{i,s}}\).  The certified pointwise estimate, including the required weight conversion, is

\begin{equation}
\left\|D^\beta a
\left(\frac{\wgray{\Phi_{i,k}}}{\wgray{\Phi_{i,s}}}\right)^{1/2}
\right\|_{L^\infty}
\le A_{i,\alpha,\beta}^{(1)}\mathfrak E_k.
\label{eq:high-nonlinear-forward-multiplier}
\end{equation}
The order-\(s\) factor is controlled by \eqref{eq:high-nonlinear-selected-factor-energy-conversion}, so \eqref{eq:high-nonlinear-two-slot-holder} gives a cubic contribution.

\hfill 

\noindent If \(|\beta|\ge k-1\), then \(s\le2\).  In this range the complementary factor is placed in \(L^\infty\), and \(D^\beta a\) is kept in its usual weighted \(L^2\) norm.  Let \(\wgray{\Phi_{j,s_a}}\) be the certified weight for \(D^\beta a\), after applying the streamfunction estimate first when necessary.  Here \(s_a\) is the derivative level in that certified estimate.  Then
\begin{align}
\|D^\beta a\|_{\wgray{\Phi_{j,s_a}}}
&\le B_{i,j,\alpha,\beta}^{(1)}\mathfrak E_k,
\label{eq:high-nonlinear-reverse-weighted-bound}\\
\left\|
D^{\alpha-\beta}\partial_\ell\bvar{\delta v_i}
\left(\frac{\wgray{\Phi_{i,k}}}{\wgray{\Phi_{j,s_a}}}\right)^{1/2}
\right\|_{L^\infty}
&\le B_{i,j,\alpha,\beta}^{\rm pt}\mathfrak E_k.
\label{eq:high-nonlinear-reverse-pointwise-bound}
\end{align}

\hfill 

\noindent The pointwise estimate uses at most two derivatives of the complementary factor, and \(k\ge4\) leaves the two additional derivatives needed for the \(H^2\to L^\infty\) bound.  Applying the weighted H\"older's inequality with the weights displayed in \eqref{eq:high-nonlinear-reverse-weighted-bound}--\eqref{eq:high-nonlinear-reverse-pointwise-bound} again gives

\begin{equation}
\mu_k\binom{\alpha}{\beta}
\left|\left\langle
(D^\beta a)D^{\alpha-\beta}\partial_\ell\bvar{\delta v_i},
D^\alpha\bvar{\delta v_i}
\right\rangle_{\wgray{\Phi_{i,k}}}\right|
\le C_{i,\alpha,\beta}^{(3)}\mathfrak E_k^3.
\label{eq:high-nonlinear-commutator-contribution}
\end{equation}

\hfill 

\noindent\textbf{Non-transport products.}
We apply Leibniz' rule to one nonlinear source product and fix a resulting summand.  We choose a perturbative spatial factor of maximal derivative order, say \(D^{\sigma_*}\bvar{\delta v_j}\), and set \(s_*=|\sigma_*|\).  \\

\noindent Proposition~\ref{prop:whole-plane-product} gives the direct estimate
\begin{align}
&\left|
\left\langle
C_\sigma\prod_mD^{\sigma_m}F_m,
D^\alpha\bvar{\delta v_i}
\right\rangle_{\wgray{\Phi_{i,k}}}
\right|
\le C_\sigma
\left\|
\left(\prod_{m\ne *}D^{\sigma_m}F_m\right)
\left(\frac{\wgray{\Phi_{i,k}}}{\wgray{\Phi_{j,s_*}}}\right)^{1/2}
\right\|_{L^\infty}
\|D^{\sigma_*}\bvar{\delta v_j}\|_{\wgray{\Phi_{j,s_*}}}
\|D^\alpha\bvar{\delta v_i}\|_{\wgray{\Phi_{i,k}}}.
\label{eq:high-nonlinear-source-product-direct}
\end{align}

\hfill 

\noindent All unselected perturbative factors are lower order and are therefore controlled pointwise. \\

\noindent Fixed profiles use their certified pointwise bounds, and streamfunction factors are first reduced using \Cref{thm:stream-elliptic-used}. \\

\noindent Since the original summand is quadratic, the combined pointwise multiplier satisfies
\begin{equation}
\left\|
\left(\prod_{m\ne *}D^{\sigma_m}F_m\right)
\left(\frac{\wgray{\Phi_{i,k}}}{\wgray{\Phi_{j,s_*}}}\right)^{1/2}
\right\|_{L^\infty}
\le E_{i,\alpha,\sigma}^{(1)}\mathfrak E_k.
\label{eq:high-nonlinear-source-combined-multiplier}
\end{equation}
Equations \eqref{eq:high-nonlinear-one-top-factor-energy}, \eqref{eq:high-nonlinear-selected-factor-energy-conversion}, and \eqref{eq:high-nonlinear-source-product-direct} then give a contribution bounded by \(C^{(3)}\mathfrak E_k^3\).\\

\noindent There is no Leibniz summand with no perturbative spatial factor coming from a quadratic modulation coefficient, because the fixed-profile contributions of $\bvar{\delta C}$ and $\bvar{\delta c_u}$ are linear and already included in the linear operator. Thus no fixed-profile quadratic modulation term is present in the high-order estimate.

\hfill \\

\noindent\textbf{Summation.}

\hfill 

\noindent For fixed \(k\), the displayed nonlinearities have finitely many Leibniz summands.  For each summand \(J\), the preceding argument provides a certified nonnegative constant \(B_J^{(3)}\) such that
\begin{equation}
|\mathcal C_J|
\le B_J^{(3)}\mathfrak E_k^3.
\label{eq:high-nonlinear-summand-bound}
\end{equation}

\noindent Define
\begin{equation}
\DLamDNthree:=\sum_J B_J^{(3)}.
\label{eq:high-nonlinear-constants-defined}
\end{equation}

\noindent Summing over components, multi-indices, transport terms, commutators, and source products proves

\begin{equation}
\mu_k\left|
\sum_{|\alpha|=k}\langle D^\alpha\mathcal N_\omega,D^\alpha\bvar{\delta\omega}\rangle_{\wgray{\Phi_{\omega,k}}}
+\sum_{|\alpha|=k}\langle D^\alpha\mathcal N_r,D^\alpha\partial_r\bvar{\delta u}\rangle_{\wgray{\Phi_{r,k}}}
+\sum_{|\alpha|=k}\langle D^\alpha\mathcal N_z,D^\alpha\partial_z\bvar{\delta u}\rangle_{\wgray{\Phi_{z,k}}}
\right|
\le \DLamDNthree \, \mathfrak E_k^3.
\end{equation}

\end{grayproof}

\clearpage

\section{Proofs of the Stability Theorem}\label{appx: nonlinear stability proof}

\hfill 

\noindent We now prove the single-radius stability Theorem~\ref{thm:main-stability}.\\

\begin{theorem}[Single-radius stability]\label{thm:main-stability_appx}
Assume that the following perturbation equations hold for admissible perturbations:
\begin{align}
\partial_t\bvar{\delta\omega}
&=\mathcal L_\omega(\bvar{\delta})+
\mathcal N_\omega(\bvar{\delta})+\mathcal R_\omega, \label{eq: Stability theorem perturbation 1_appx}\\
\partial_t\partial_r\bvar{\delta u}
&=\partial_r\mathcal L_u(\bvar{\delta})+
\mathcal N_{r}(\bvar{\delta})+\mathcal R_r,\label{eq: Stability theorem perturbation 2_appx}\\
\partial_t\partial_z\bvar{\delta u}
&=\partial_z\mathcal L_u(\bvar{\delta})+
\mathcal N_{z}(\bvar{\delta})+\mathcal R_z. \label{eq: Stability theorem perturbation 3_appx}
\end{align}
Assume that the following estimates have been certified with finite constants:

\begin{equation}
\begin{aligned}
&\langle \mathcal L_\omega,\bvar{\delta\omega}\rangle_{\wgray{\Phi_\omega}}
+\langle \partial_r\mathcal L_u,\partial_r\bvar{\delta u}\rangle_{\wgray{\Phi_r}}
+\langle \partial_z\mathcal L_u,\partial_z\bvar{\delta u}\rangle_{\wgray{\Phi_z}}
\le -\DLamLlow \, \mathfrak E_0^2+\DLamLhigh \, \mathfrak H_k^2.
\end{aligned}
\label{eq:ns-thm-low-linear_appx}
\end{equation}

\vspace{2mm}

\begin{equation}
\begin{aligned}
&\sum_{i\in\{\omega,r,z\}}\sum_{|\alpha|=k}
\langle D^\alpha\mathcal L_i,D^\alpha\bvar{\delta v_i}\rangle_{\wgray{\Phi_{i,k}}}
\le \DLamDLlow \, \mathfrak E_0^2-\DLamDLhigh \, \mathfrak H_k^2.
\end{aligned}
\label{eq:ns-thm-high-linear_appx}
\end{equation}

\vspace{2mm}

\begin{equation}
\begin{aligned}
&\left|\langle\mathcal N_\omega,\bvar{\delta\omega}\rangle_{\wgray{\Phi_\omega}}
+\langle\mathcal N_{r},\partial_r\bvar{\delta u}\rangle_{\wgray{\Phi_r}}
+\langle\mathcal N_{z},\partial_z\bvar{\delta u}\rangle_{\wgray{\Phi_z}}
\right|
\le \DLamNthree \, \mathfrak E_k^3.
\end{aligned}
\label{eq:ns-thm-low-nonlinear_appx}
\end{equation}

\vspace{2mm}

\begin{equation}
\begin{aligned}
&\mu_k\left|
\sum_{|\alpha|=k}\langle D^\alpha\mathcal N_\omega,D^\alpha\bvar{\delta\omega}\rangle_{\wgray{\Phi_{\omega,k}}}
+\sum_{|\alpha|=k}\langle D^\alpha\mathcal N_r,D^\alpha\partial_r\bvar{\delta u}\rangle_{\wgray{\Phi_{r,k}}}
+\sum_{|\alpha|=k}\langle D^\alpha\mathcal N_z,D^\alpha\partial_z\bvar{\delta u}\rangle_{\wgray{\Phi_{z,k}}}
\right|
\le \DLamDNthree \, \mathfrak E_k^3.
\end{aligned}
\label{eq:ns-thm-high-nonlinear_appx}
\end{equation}

\vspace{2mm}

\begin{equation}
\begin{aligned}
&\left|
\sum_{i\in\{\omega,r,z\}}\langle \mathcal R_i,\bvar{\delta v_i}\rangle_{\wgray{\Phi_i}}
+\mu_k\sum_{i\in\{\omega,r,z\}}\sum_{|\alpha|=k}
\langle D^\alpha\mathcal R_i,D^\alpha\bvar{\delta v_i}\rangle_{\wgray{\Phi_{i,k}}}
\right|
\le \DLamR \, \mathfrak E_k.
\end{aligned}
\label{eq:ns-thm-residual_appx}
\end{equation}

\hfill \\

\noindent Choose $\mu_k \in (0,1]$ and
\begin{equation}
\label{eq:ns-mu-choice_appx}
\DLamStab(\mu_k):=
\min\left\{
\DLamLlow-\mu_k\DLamDLlow,
\ \DLamDLhigh-\frac{\DLamLhigh}{\mu_k}
\right\}>0.
\end{equation}
After this choice of $\mu_k$ is fixed, write $\DLamStab$ for the positive number $\DLamStab(\mu_k)$. \\

\noindent Define the combined nonlinear constant by
\begin{equation}
\label{eq:ns-N2N3_appx}
\DLamThree:=\DLamNthree+\DLamDNthree.
\end{equation}

\noindent Assume that $\delta_*>0$ satisfies

\begin{equation}
\label{eq:ns-smallness_appx}
\DLamThree\delta_*+\frac{\DLamR}{\delta_*}<\DLamStab.
\end{equation}

\hfill 

\noindent Then every solution with $\mathfrak E_k(0)<\delta_*$ remains in the ball $\mathfrak E_k(t)<\delta_*$ for as long as the solution exists and the certified estimates apply.
\end{theorem}

\newpage

\begin{grayproof} \Lean \\

\noindent We pair the three perturbation equations with $(\bvar{\delta\omega},\partial_r\bvar{\delta u},\partial_z\bvar{\delta u})$ in the low-order weights and add the $k$th-order pairings with prefactor $\mu_k$. \\

\noindent The two linear estimates in the theorem, with the high-order estimate multiplied by $\mu_k$, give
\begin{align}
\label{eq:ns-close-linear-add_appx}
&-\DLamLlow \, \mathfrak E_0^2+\DLamLhigh \, \mathfrak H_k^2
+\mu_k\bigl(-\DLamDLhigh \, \mathfrak H_k^2+\DLamDLlow \, \mathfrak E_0^2\bigr) \\ &  \qquad \qquad \qquad =-(\DLamLlow-\mu_k\DLamDLlow)\mathfrak E_0^2-(\mu_k\DLamDLhigh-\DLamLhigh)\mathfrak H_k^2. \nonumber
\end{align}

\hfill

\noindent By \eqref{eq:ns-mu-choice_appx},
\begin{equation}
\DLamLlow-\mu_k\DLamDLlow\ge \DLamStab,
\qquad
\mu_k\DLamDLhigh-\DLamLhigh\ge \mu_k\DLamStab.
\end{equation}

\hfill 

\noindent Therefore
\begin{equation}
\label{eq:ns-close-linear-gamma_appx}
-(\DLamLlow-\mu_k\DLamDLlow)\mathfrak E_0^2-(\mu_k\DLamDLhigh-\DLamLhigh)\mathfrak H_k^2
\le -\DLamStab\mathfrak E_0^2-\mu_k\DLamStab\mathfrak H_k^2
= -\DLamStab\mathfrak E_k^2,
\end{equation}
since $\mathfrak E_k^2=\mathfrak E_0^2+\mu_k\mathfrak H_k^2$.

\hfill  \\

\noindent The low-order nonlinear estimate \eqref{eq:ns-thm-low-nonlinear_appx}, the high-order nonlinear estimate \eqref{eq:ns-thm-high-nonlinear_appx}, and \eqref{eq:ns-N2N3_appx} give

\begin{equation}
\label{eq:ns-close-nonlinear_appx}
\left|
\sum_{i\in\{\omega,r,z\}}
\langle\mathcal N_i,\bvar{\delta v_i}\rangle_{\wgray{\Phi_i}}
+\mu_k\sum_{i\in\{\omega,r,z\}}\sum_{|\alpha|=k}
\langle D^\alpha\mathcal N_i,D^\alpha\bvar{\delta v_i}\rangle_{\wgray{\Phi_{i,k}}}
\right| \le \DLamThree\mathfrak E_k^3.
\end{equation}

\hfill \\ 

\noindent Combining \eqref{eq:ns-close-linear-gamma_appx}, \eqref{eq:ns-close-nonlinear_appx}, and the residual bound \eqref{eq:ns-thm-residual_appx} yields

\begin{equation}
\label{eq:ns-energy-squared-ineq_appx}
\frac12\frac{d}{dt}\mathfrak E_k^2
\le -\DLamStab\mathfrak E_k^2+\DLamThree\mathfrak E_k^3+
\DLamR \mathfrak E_k.
\end{equation}

\hfill 

\noindent It remains to show that the energy ball is forward invariant.  We suppose that \(\mathfrak E_k(0)<\delta_*\) and argue by contradiction. Define the first boundary-contact time by
\begin{equation}
t_*:=\inf\{t>0:\mathfrak E_k(t)=\delta_*\}
\label{eq:ns-closing-first-contact-time}
\end{equation}
 \\ 

\noindent Set \(F(t):=\mathfrak E_k(t)^2\). Then \(F(t)<\delta_*^2\) for \(t<t_*\) and \(F(t_*)=\delta_*^2\), so the left derivative satisfies
\begin{equation}
F'(t_*)=
\lim_{h\downarrow0}\frac{F(t_*)-F(t_*-h)}{h}\ge0.
\label{eq:ns-closing-first-contact-left-derivative}
\end{equation}

\hfill 

\noindent On the other hand, evaluating \eqref{eq:ns-energy-squared-ineq_appx} at \(t_*\) gives

\begin{align}
\frac12F'(t_*)
&\le -\DLamStab\delta_*^2
+\DLamThree\delta_*^3
+\DLamR \delta_*
=-\delta_*^2\left(
\DLamStab-\DLamThree\delta_*
-\frac{\DLamR}{\delta_*}\right)<0
\label{eq:euler-boundary-inward-derivative}
\end{align}
by \eqref{eq:ns-smallness_appx}, which contradicts \eqref{eq:ns-closing-first-contact-left-derivative}.  \\

\noindent Therefore the ball
\(\mathfrak E_k<\delta_*\) is forward invariant, which proves the stated
nonlinear stability. 
\end{grayproof}

\clearpage

\end{document}